\documentclass[12pt,reqno]{amsart}
\usepackage{amsmath,amssymb,amsbsy,amsfonts,amsthm,latexsym,amsopn,amstext,
            amsxtra,euscript,amscd,verbatim,epsf,colordvi,graphics}
\usepackage{Young}
\usepackage{tabmac}
\usepackage{color}
\usepackage{pgf,tikz}
\usetikzlibrary{decorations.pathreplacing}
\usetikzlibrary{arrows}

\usepackage{subfigure}

\usepackage{hyperref}

\newtheorem{theorem}{Theorem}[section]
\newtheorem{corollary}[theorem]{Corollary}

\newtheorem{proposition}[theorem]{Proposition}

\theoremstyle{definition}
\newtheorem{definition}{Definition}[section]

\newtheorem{remark}{Remark}[section]
\newtheorem{example}{Example}[section]

\numberwithin{equation}{section}

\numberwithin{table}{section}

\numberwithin{figure}{section}
\newtheorem{alg}[theorem]{Algorithm}

\newcommand{\olga}{\color{red}}

\def \N{\mathbb{N}}

\def\cA{{\mathcal A}}
\def\cB{{\mathcal B}}

\def\cX{{\mathcal X}}
\def\cY{{\mathcal Y}}

\def\ts{{\hskip .1in}}

\def\rf#1{\left\lceil#1\right\rceil}

\def\({\left(}
\def\){\right)}
\def\[{\left[}
\def\]{\right]}
\def\<{\langle}
\def\>{\rangle}

\title[Symmetries of Littlewood--Richardson coefficients]{ A uniform action of the dihedral group  $\mathbb Z_2\times { D}_3$ \\on Littlewood--Richardson coefficients }

\author{Olga Azenhas}
\address{University of Coimbra, CMUC, Department of Mathematics, Portugal}
\email{oazenhas@mat.uc.pt}

\author{Alessandro Conflitti}
\email{alessandro.conflitti@gmail.com}

\author{Ricardo Mamede}
\address{University of Coimbra, CMUC, Department of Mathematics, Portugal}
\email{mamede@mat.uc.pt}

\thanks{ This work was partially supported by the Centre for Mathematics, University of Coimbra (funded by the Portuguese Government through FCT/MCTES, DOI 10.54499/UIDB/00324/2020). This work has been partially conducted when the second author was a member of CMUC, Centre for Mathematics, University of Coimbra.}



\begin{document}

\begin{abstract} We show that the dihedral group $\mathbb Z_2\times D_3$ of order twelve acts faithfully on the set $\mathcal{LR}$, either consisting of Littlewood--Richardson (LR) tableaux, or their companion tableaux,  or Knutson--Tao (KT) hives or   Knutson--Tao--Woodward (KTW) puzzles,
via involutions which  simultaneously  conjugate   or shuffle a Littlewood-Richardson (LR) triple  of partitions.
 The  action of $\mathbb Z_2\times { D}_3$  carries a linear time    index two subgroup $H\simeq D_3$ action,  where  an involution which goes from $H$ into the other
coset of $H$ is  difficult in the sense that it is not manifest neither exhibited by simple means.   Pak and Vallejo have earlier made this observation with respect to the subgroup  of index two in the symmetric group $\mathfrak{ S}_3$ consisting of cyclic permutations which  $H$ extends.
   The other half   LR symmetries, not in the range of the  $H$-action,  are hidden  and consist of commutativity and conjugation symmetries.  Their  exhibition is reduced to  the action of a remaining generator of  $\mathbb Z_2\times D_3$, which belongs to the other coset of $H$, and  enables to  reduce in linear time all known LR commuters and transposers to each other, and to the  Luzstig-- Sch\"utzenberger involution. A KT hive is specified  by superimposing the companion tableau pair of an LR tableau, and  its $\mathbb Z_2\times D_3$--symmetries   are exhibited  via the corresponding LR companion tableau pair.  The action  of $\mathbb Z_2\times D_3$ on KTW puzzles, naturally in bijection with Purbhoo mosaics,  is consistent with the migration map on  mosaics which translates to \textit{jeu de taquin} slides or tableau-switching on LR tableaux. Their $H$--symmetries are  reduced to simple procedures on a KTW  puzzle via label swapping together with simple reflections of an equilateral triangle, that is, puzzle dualities, and  rotations on an equilateral triangle. Finally, the $\mathfrak{S}_3$--symmetries under this action,  distributed in the two $H$-cosets,  are consistent with the Thomas-Yong carton rule based on the infusion involution, a specific governance of \emph{jeu de taquin} slides in the tableau switching.
\end{abstract}

\maketitle

\section{Introduction}

This paper aims to fulfill the study by Pak and Vallejo with
 LR conjugation involutions (LR transposers), which have not been considered in \cite{PV2}, and thus to give a complete and uniform picture of the
LR symmetries under the  action of the dihedral group $\mathbb Z_2\times { D}_3$. Namely, one shows that all LR transposers known up to date coincide. Furthermore, the LR transposers and LR commuters are linear time reducible to each other, in particular, to Lusztig-Sch\"utzenberger involution. (This involution is realized by the Sch\"utzenberger
\textit{evacuation} \cite{st,ful}
  on tableaux of normal shape, and by the Benkart-Sottile-Stroomer \textit{reversal} \cite{bss} on
 tableaux of skew-shape.) This  amounts to show  the linear cost computational complexity of the involutions exhibiting the $H$--symmetries which in a Knutson-Tao-Woodward (KTW) puzzle \cite{knutson} (puzzle for short) consist of rotations and simple reflections  on an equilateral triangle, { where} the latter are together with   label swapping of the puzzle. This is consistent with the coincidence of LR commuters known up to date \cite{az17}.  {  We refer to  the recent survey by Pak \cite{IPak} for a  comprehensive account of combinatorial interpretations of LR coefficients}.

The Berenstein--Zelevinsky (BZ) triangles \cite{bz} give an interpretation of  LR coefficients
   which manifest all  $\mathfrak{S}_3$--symmetries except the commutativity, that is, the swapping of two entries in the LR triple  (see \cite[Remarks (a)]{bz}), and  the conjugation symmetry is also hidden in BZ triangles.
Pak and Vallejo have defined in~\cite{PV1} bijections, given by explicit {\em linear maps}, between LR tableaux, Knutson--Tao (KT)
hives ~\cite{knutson0}, and BZ triangles, which combined with the symmetries of BZ triangles give all the $\mathfrak{S}_3$-symmetries
except  the commutativity. The conjugation symmetry is not considered in their work.
As pointed out in~\cite{PV1},~\cite[Section 7.6]{PV2}, regarding to the   $\mathfrak{S}_3$--symmetries of LR triples,  those defined by the index two subgroup $R$   in $ \mathfrak{S}_3$, consisting of  cyclic permutations, can be given by easily computed involutions in every combinatorial model.  The    analysis of the  symmetries of LR triples in BZ triangles, under the aforesaid action of $\mathbb Z_2\times \mathfrak{S}_3$ \cite[Remarks (a)]{bz}, suggests that a larger subgroup $ H\supseteq R$ of  symmetries of LR triples, with index two, in $\mathbb Z_2\times \mathfrak{S}_3$, can be given by easily computed involutions in every combinatorial model.
The LR coefficients are preserved in linear time by the action of $ H$. The $H$--symmetries  are  manifest in a KTW puzzle. They are exhibited via rotations and  reflections of the dihedral  group $D_3\simeq \mathfrak{S}_3$ on an equilateral triangle together with  the label swapping of the puzzle. These involutions incur to simultaneously transpose and commute entries of the LR triple of partitions.  We thus consider the isomorphic group $\mathbb{Z}_2\times H\simeq \mathbb{Z}_2\times D_3$ - action on  the set $\mathcal{LR}$ either consisting of Littlewood-Richardson (LR) tableaux, or their companion tableaux \cite{hk},  (see also \cite{akt,tka}),  or  hives,  or  puzzles.

  The symmetries outside of $H$ are linear time reducible either to the commutativity or to the transposition symmetries, and are given by the action of a remaining generator of $\mathbb Z_2\times { D}_3$ in the other coset of $H$.
More precisely, the other half   LR symmetries, realized by the action of the elements in the other coset of  $H$, are hidden  and consist of commutativity and conjugation symmetries.  They are given by the action of a remaining generator,  realized by the \emph{reversal} involution \cite{bss} or the Sch\"utzenberger \textit{evacuation} involution, which enables to  reduce in linear time the bijections exhibiting commutativity and conjugation  symmetries to each other.
    Since all known LR commuters are involutions and coincide, this incurs the coincidence and the involutory nature of  all known LR transposers. The  LR commuters and the LR transposers are linear time reducible to the Lusztig-Sch\"utzenberger  involution. 

    Explicit constructions of the corresponding  symmetry involutions  on  companion tableaux  of LR skew tableaux are provided as well. 
    To pass  from symmetries of LR  (skew) tableaux to symmetries of  companion tableaux we use  the action  of the longest element of a symmetric group on a crystal by sending  an LR tableau to its reversal, and Lascoux's double crystal graph structure on biwords by relating the so called left and right strings of the double crystal graph \cite{double}. This analysis also incurs in an explicit relationship between  two interlocking Gelfand--Tsetlin (GT) patterns in a hive, that is, the left and the right companion of an LR tableau \cite{hk}.
    Finally, one takes the  linear time index two subgroup $H$ action on Knutson-Tao-Woodward puzzles and Purbhoo mosaics \cite{mosaic} and
explains what operations on LR (skew) tableaux they translate to and back.

\subsection{ LR coefficients as structure coefficients and symmetries}

Schur functions $s_\lambda$  where $\lambda$ runs over all Young shapes (partitions)
form  a linear $\mathbb{Z}$-basis for the ring $\Lambda$ of symmetric functions in  countably many variables  $x_1,x_2,\dots$. The Littlewood--Richardson (LR)  coefficients are the structure coefficients in the Schur function product with respect to this basis.
The product $s_{\mu}s_{\nu}$ in $\Lambda$ is therefore a non--negative integral
linear combination of Schur functions $s_{\lambda}$, \begin{equation}\label{schur}s_{\mu}s_{\nu}=\sum_{\lambda}
c_{\mu\;\nu}^{\lambda}s_{\lambda},\end{equation}
 \noindent where the  structure constants $c_{\mu\;\nu}^{\lambda}$ in this Schur basis expansion,    depending only on the three partitions
$\mu,\nu$ and $\lambda$, are called  Littlewood--Richardson (LR) coefficients
~\cite{lr,mac,sa,st}. Let  $\Lambda_d:= \mathbb{Z}[x_1,\dots,x_d]^{\mathfrak{S}_d}$ be the ring of symmetric polynomials in the variables $x_1,\dots, x_d$. In representation theory the Schur polynomials in $\Lambda_d$ occur as
characters of the general linear group $GL_d(\mathbb{C})$, and they correspond to characters of the symmetric group via the Frobenius map. Schur functions also have an intersection-theoretic
interpretation as representatives of Schubert classes in the cohomology ring of a Grassmannian.
 More precisely, the product of two Schubert cohomology classes $ S_\mu. S_\nu$ on a Grassmannian $Gr(d, \mathbb{C}^n )$
is also known to be a positive combination of other Schubert classes, and  the coefficients can be computed using
the Littlewood-Richardson rule \cite{lr} or the more symmetric puzzle rule from  Knutson-Tao-Woodward \cite{knutson},
$$ S_\mu. S_\nu = \sum_{
\textit{puzzles $P$}} S_\lambda,$$
$\lambda$ the southern  side of $P$ , read left to right
where the sum is taken over puzzles with NW side labeled $\mu$, NE side labeled $\nu$, both from left
to right. (See Section \ref{sec:puzzle} for details on puzzles.)
Let  $ H^*(Gr(d, \mathbb{C}^n ))$ the cohomology of the Grassmannian $Gr(d, \mathbb{C}^n )$. The map
$\Lambda_d\rightarrow H^*(Gr(d, \mathbb{C}^n ))$
which sends $s_\mu$ to $S_\mu$ if $\mu$ fits the rectangle $d\times (n-d)$ and to $0$ otherwise, is a ring epimorphism.

Let $0< d< n$ be fixed integers, and fix Young shapes   $\mu,\nu, \lambda$ contained in the ambient non empty rectangle rectangle  $D:=d\times (n-d)$ according to the French convention. Let $c_{\mu\;\nu\,\lambda^\vee}:=c_{\mu\;\nu}^{\lambda}$ where $\lambda^\vee$    is the Young shape defined by the  set complement of $\lambda$ in that rectangle. (See subsection \ref{sec:pre}.)
Henceforth, $c_{\mu\;\nu\,\lambda^\vee}$ is the coefficient of $s_{d\times(n-d)}$ in the Schur expansion of $s_\mu s_\nu s_{\lambda^\vee}$ in $\Lambda$, and clearly $c_{\mu,\nu,\lambda}$ is invariant under the shuffling of the partition triple  $(\mu,\nu,\lambda)$:
\begin{equation}\label{schur2}c_{\mu,\nu,\lambda}= c_{\nu\;\mu\;\lambda}=c_{\mu\;\lambda\;\nu}=c_{\lambda\;\nu\;\mu}\end{equation}
\begin{equation}\label{schur2a}c_{\mu,\nu,\lambda}=c_{\lambda\;\mu\;\nu}=c_{\nu\;\lambda\;\mu}.\end{equation}

Let $\lambda^t$ denote the conjugate or
transpose of the partition $\lambda$ with ambient rectangle $(n-d)\times d$. While the $\mathfrak{S}_3$-symmetries
  \eqref{schur2}, \eqref{schur2a} are obvious from the Schur expansion of $s_\mu s_\nu s_{\lambda^\vee}$, it is not  the case from the Schur  expansion of $s_{\mu^t}s_{\nu^t}s_{{\lambda^\vee}^t}$ where $c_{\mu^t\;\nu^t\,{\lambda^\vee}^t}$ is the coefficient of $s_{(n-d)\times d}$ in that expansion, that the LR conjugation
symmetry,
\begin{equation}\label{schurconjug}c_{\mu\;\nu\;\lambda}=c_{\mu^t\;\nu^t\;\lambda^t}\end{equation} holds.  This last symmetry  is shown {\em via} the involutive $\mathbb{Z}$-automorphism $\omega$ of the $\mathbb{Z}$--algebra $\Lambda$ of symmetric functions \cite{mac,st}.

On the other hand,    LR coefficients enumerate several combinatorial objects depending on three partitions $\mu$, $\nu$ and $\lambda$,
and the combinatorics of their symmetries is rather uniform in the sense that in all  combinatorial means the commutativity \eqref{schur2}
and the conjugation \eqref{schurconjug} symmetries are {\em hidden}, see~\cite{az1,akt,bz,DK3,knutson0,PV1,PV2,tka}. Interestingly,  commuting and transposing simultaneously gives, $c_{\mu\;\nu}^{\lambda}=c_{\nu^t\;\mu^t}^{\lambda^t}$,  a symmetry  revealed by {\em simple} means.  
In  (KTW) puzzles \cite{knutson}, it means the {\em puzzle duality}, that is, one gets this symmetry {\em via} the vertical reflection  of a puzzle with label swapping.
 This  {\em simple} involution,  denoted $\spadesuit$, subsection \ref{sec:H}, is in turn   translated to LR tableaux through  {\em simple} operations in Definition \ref{D:spade},  as well as to   LR companion tableaux (for the definition, see  subsection \ref{subsec:lr}) in Algorithm \ref{alg:spade}.
The vertical reflection with labelling swapping followed with a clockwise rotation of a  KTW puzzle by $\pi/3$ and $2\pi/3$ radians  gives the two remaining {\em puzzle dualities} $c_{\mu\;\nu\;\lambda}=c_{\lambda^t\;\nu^t\;\mu^t}$ and $c_{\mu\,\nu\,\lambda}=c_{\mu^t\,\lambda^t\,\nu^t}$. They are again  translated to LR tableaux or companion LR tableaux through simple involutions, denoted $\blacklozenge$  and $\clubsuit$  respectively. (For  definitions, see algorithms   \ref{alg}, \ref{alg:blacklozenge},  subsection \ref{sec:H}, and Definition \ref{D:club}.)  The three symmetries consisting of puzzle dualities and the three symmetries consisting of  puzzle rotations
$c_{\mu\;\nu\;\lambda}=c_{\lambda\;\mu\;\nu}=c_{\nu\;\lambda\;\mu}$ \eqref{schur2a} are exhibited by the faithful  action of a two index subgroup $H\simeq {D}_3$  \eqref{hgroup} of the dihedral group $\mathbb Z_2\times { D}_3$  on KTW puzzles.


\subsection{The set $\mathcal{LR}$ and a representation of $\mathbb Z_2\times { D}_3$ in $\mathfrak{S}_{\mathcal{LR}}$}
 Let $0< d< n$ be fixed integers. Throughout, a partition is identified with its Young diagram which fits inside, according to the French convention,  a nonempty rectangle $D:=d\times (n-d)$. Let $\binom{[n]}{d}$ denote the set of binary words consisting of $d$ ones and $n-d$ zeroes.
 Our partitions in the ambient rectangle space $D$ are
identified with  the $01$--words in $\binom{[n]}{d}$ as follows: the positions of the zeroes and ones in a
$01$--word are respectively the positions of the horizontal and
vertical steps along the boundary of the corresponding Young
diagram, starting in the right lower corner of the rectangle and ending up at the upper left corner. In particular, the empty partition $\emptyset$ is identified with $0^{n-d}1^d$, and $D$ with $1^{d}0^{n-d}$.(See picture \eqref{rectangle} in  Section \ref{sec:pre}).

Given $\mu,\nu,\lambda$ be partitions contained in $D$,  $\rm LR_{\mu,\nu}^\lambda=\rm LR(\mu,\nu,\lambda^\vee)$
is the set of LR
tableaux of shape $\lambda/\mu$ and content $\nu$ and   its cardinality is $c_{\mu\;\nu}^{\lambda}=c_{\mu\;\nu\;\lambda^\vee}$. See subsection \ref{subsec:standard} for the definition. In particular, $c_{\mu,\nu,\lambda^\vee}>0$ only if $\mu,\nu\subset \lambda$ and the sum of all parts of $\mu,\nu$ and $\lambda^\vee$ is $|\mu|+|\nu|+|\lambda^\vee|=n$.
Let $\mathcal{LR}_{d,n}$ be the set of all LR tableaux or KTW puzzles of size $n$, that is,
\begin{equation}\label{lrset}\displaystyle\mathcal{LR}_{d,n}=\bigsqcup_{(\mu,\nu,\lambda)} \rm LR_{\mu,\nu}^{\lambda^\vee},
\end{equation}  where $(\mu,\nu,\lambda)\in \binom{[n]}{d}^3\cup \binom{[n]}{n-d}^3$.
   (For some  partition-triples $(\mu,\nu,\lambda)$, we may have $\rm LR_{\mu,\nu}^{\lambda^\vee}=\emptyset$ in which case $c_{\mu,\nu,\lambda}=0$.) For simplicity, we  write $\mathcal{LR}:=\mathcal{LR}_{d,n}$.

Given  $(\mu,\nu,\lambda)\in \binom{[n]}{d}^3\cup \binom{[n]}{n-d}^3$,
the LR coefficients $c_{\mu\;\nu\,\lambda}$ are invariant under the following action of the dihedral group $\mathbb Z_2\times {D}_3$: the non--identity element of $\mathbb Z_2$ transposes simultaneously $\mu$, $\nu$ and $\lambda$, and the reflections in ${D}_3$ swap two entries in the triple $(\mu, \nu,
\lambda)$.
Denoting by $\tau$ the non--identity element of the cyclic group $\mathbb Z_2$, and  by $\varsigma_{1}$,  $\varsigma_{2}$ two swaps  or reflections  of the dihedral group ${D}_3$ of order six, consider
 $\mathbb Z_2\times {D}_3$, the dihedral group of order twelve, as the free group generated by the involutions $\tau$, $\varsigma_1$, $\varsigma_2$  subject to the relations inherited from $\mathbb{ Z}_2$ and ${D}_3$ as a Coxeter group, and such that $\tau$ commutes with both $\varsigma_{1}$ and $\varsigma_{2}$,
 \begin{equation}\label{dihedraltwelve}\mathbb Z_2\times {D}_3=\langle \tau,\varsigma_{1}, \varsigma_{2}\mid \tau^2=\varsigma_{1}^2=\varsigma_{2}^2=(\varsigma_{1} \varsigma_{2})^3=1=(\tau\varsigma_{1})^2=(\tau\varsigma_{2})^2\rangle.\end{equation}

  Let $H$ be the two index subgroup of
$\mathbb Z_2\times {D}_3$, defined by
\begin{eqnarray}\label{hgroup}
H:=&\langle\tau\varsigma_1,\tau\varsigma_2\rangle&=\{1,\tau\varsigma_1,\tau\varsigma_2,
\varsigma_1\varsigma_2,\varsigma_2\varsigma_1,\tau\varsigma_1\varsigma_2\varsigma_1\}\\
=&
\langle\tau\varsigma_i,\tau\varsigma\rangle& =\{1,\tau\varsigma_i,\tau\varsigma,
\varsigma_i\varsigma,\varsigma\varsigma_i,\tau\varsigma_i\varsigma\varsigma_i\}, \;i=1,2,\nonumber
\end{eqnarray}
where $\varsigma=\varsigma_1\varsigma_2\varsigma_1= \varsigma_2\varsigma_1\varsigma_2$.
The subgroup $H\simeq D_3$  contains the two index cyclic group $R=\langle \varsigma_1\varsigma_2 \rangle=\{1,\varsigma_1\varsigma_2,\varsigma_2\varsigma_1\}$ of $D_3$.   Since $H$ is a  subgroup  of index 2
in $\mathbb Z_2\times \mathfrak{S}_3$, $H$ is  normal,  the
quotient group $\left(\mathbb Z_2\times \mathfrak{S}_3\right)/H$ is cyclic and every element different from the
identity is a generator, {\em i.e.}  $\zeta H=H\zeta\neq H$, for $\zeta=\varsigma_1,$ $\varsigma_2,\tau$, $\varsigma_1\varsigma_2\varsigma_1
$, $\tau\varsigma_2\varsigma_1$, $\tau\varsigma_1\varsigma_2\notin H$. Therefore, as a set, $\mathbb Z_2\times {D}_3=H\sqcup \zeta H$, for any $\zeta\notin H$, and
$H$ affords also the following  presentations of $\mathbb Z_2\times {D}_3\simeq \mathbb Z_2\times H$ useful for our purposes:
\begin{equation}\label{dihedral1} \langle\tau,\tau\varsigma_1,\tau\varsigma_2:\tau^2=(\tau\varsigma_1)^2=(\tau\varsigma_2)^2=
(\tau\varsigma_1\tau\varsigma_2)^3=(\tau\varsigma_1\tau)^2=(\tau\varsigma_2\tau)^2=1\rangle.\end{equation}
  \begin{equation} \label{eq:groupres}\langle\tau,\tau\varsigma_i,\tau\varsigma:\tau^2=(\tau\varsigma_i)^2=(\tau\varsigma)^2=
(\tau\varsigma_i\tau\varsigma)^3=(\tau\varsigma_i\tau)^2=(\tau\varsigma\tau)^2=1\rangle,\; i=1,2.\end{equation}
\begin{equation}\label{eq:groupres3} \langle\varsigma_i,\tau\varsigma_1,\tau\varsigma_2:\varsigma_i^2=(\tau\varsigma_2)^2=(\tau\varsigma_1)^2=
(\tau\varsigma_1\tau\varsigma_2)^3=(\tau\varsigma_i\varsigma_i)^2=(\tau\varsigma_j\tau\varsigma_i\varsigma_i)^2=1\rangle,\;1 \le j\neq i\le 2.
\end{equation}
\begin{equation}\label{eq:groupres4} \langle\varsigma,\tau\varsigma_1,\tau\varsigma_2:\varsigma^2=(\tau\varsigma_2)^2=(\tau\varsigma_1)^2=
(\tau\varsigma_1\tau\varsigma_2)^3=1
\rangle,
\end{equation}
\begin{equation}\label{eq:groupres2} \langle\varsigma_i,\tau\varsigma_i,\tau\varsigma:\varsigma_i^2=(\tau\varsigma_i)^2=(\tau\varsigma)^2=
(\tau\varsigma_i\tau\varsigma)^3=(\varsigma_i\tau\varsigma_i)^2=(\tau\varsigma\tau)^2=1\rangle,\;i=1,2.
\end{equation}

   We show that the group $\mathbb Z_2\times {D}_3$ acts faithfully on the set $\mathcal{LR}$ trough involutions which    conjugate or shuffle the entries of a partition-triple $(\mu,\nu,\lambda)$.
  Using the presentation \eqref{eq:groupres},   the following is an injective  group homomorphism where $\mathfrak{S}_{\mathcal{LR}}$ is the group of all bijections from $\mathcal{LR}$ to itself or permutations of the set $\mathcal{LR}$
\begin{eqnarray}
\varpi:\mathbb Z_2\times {D}_3&\longrightarrow& \mathfrak{S}_{\mathcal{LR}}\\
\tau&\mapsto&\varrho\nonumber\\
\tau\varsigma_1\varsigma_2\varsigma_1&\mapsto&\blacklozenge\nonumber\\
\tau\varsigma_1&\mapsto&\spadesuit\nonumber
\end{eqnarray}
realized by the involutions $\blacklozenge $, $\spadesuit$ and $\varrho$ on $\mathcal{LR}$. (See section \ref{sec:commutertransposer} for definitions.) More precisely, $\blacklozenge $ is an involution on $ {\rm LR}(\mu,\nu,\lambda)\sqcup {\rm LR}(\lambda^t,\nu^t,\mu^t)$ \eqref{orthogonal} and on the union of KTW puzzles of boundary $(\mu,\nu,\lambda)$ or $(\lambda^t,\nu^t,\mu^t)$, subsection \ref{sec:H}, which agrees with Tao's bijection between KTW puzzles and LR tableaux, see Example \ref{ex:taoblacklozenge}; $\spadesuit$ denotes the vertical puzzle duality in subsection \ref{sec:H}, and the involution on $ {\rm LR}(\mu,\nu,\lambda)\sqcup {\rm LR}(\nu^t,\mu^t,\lambda^t)$ given by the procedure in Definition \ref{D:spade}; and $\varrho=\bullet\blacklozenge \,\eta$ in  Theorem \ref{plr*}, with $\bullet$  the rotation map \eqref{rotation}, and $\eta$ the Lusztig-Sch\"utzenberger involution, is an involution in $ {\rm LR}(\mu,\nu,\lambda)\sqcup {\rm LR}(\mu^t,\nu^t,\lambda^t)$.

 Observe that for any nonempty partition $\gamma$, $\rm LR(\emptyset,\gamma,\gamma^\vee)$ consists of the sole straight LR tableau of shape and weight $\gamma$, $Y(\gamma)$, also called Yamanouchi tableau of shape $\gamma$. For any $g\neq h\in \mathbb Z_2\times {D}_3$, $\varpi(g)(T)\neq \varpi(h)(T)$ for some $T\in LR(\epsilon, \delta, \alpha)$ where $(\epsilon, \delta, \alpha)$ is a shuffle or a shuffle and a transposition of  the entries of $(\emptyset,\gamma,\gamma^\vee)$ providing $\gamma\neq \gamma^\vee$.
The monomorphism $\varpi$ shows that $\mathbb Z_2\times {D}_3$ is a group of symmetries of LR coefficients regarding to the conjugation and the shuffling of the  boundary partition-triple  of an element in $\mathcal{LR}$, and thereby
\begin{equation}\mathbb Z_2\times {D}_3\simeq \langle\spadesuit,\blacklozenge,\varrho\rangle:= \langle\spadesuit,\blacklozenge,{\varrho}:{\varrho}^2=\spadesuit^2=\blacklozenge^2=
(\spadesuit\blacklozenge)^3=(\spadesuit{\varrho})^2=(\blacklozenge{ \varrho})^2=1\rangle.
\end{equation}

\subsection{The $H$-action and the action of a remaining generator of $\mathbb Z_2\times {D}_3$ on $\mathcal{LR}$} The involutions exhibiting the $H$-symmetries of LR triples, define a faithfully  group action of $H$ on KTW puzzles and on LR tableaux.
As for the computational complexity  of our involutions, we  study the  invariance of LR coefficients, under the action of the  two-index subgroup $H$ on  the set  $\mathcal{LR}$, where
\begin{equation}H=\langle\tau\varsigma_1,\tau\varsigma_2\rangle
=\langle\tau\varsigma_1,\tau\varsigma_1\varsigma_2\varsigma_1\rangle
\simeq \langle\spadesuit,\blacklozenge\rangle=
 \{\spadesuit,\blacklozenge:\spadesuit^2=\blacklozenge^2={\bf 1}=
(\spadesuit\blacklozenge)^3\}\simeq {{D}_3}.
\end{equation}

 The $H$-invariance of LR coefficients is proved through the exhibition of simple involutions $\spadesuit,\blacklozenge,\clubsuit$ on KTW puzzles, called {\em puzzle dualities},  and simultaneously on LR tableaux. Tao's bijection shows how they do translate to each other.  Example \ref{ex:taoblacklozenge} clearly illustrates this translation  for the involution $\blacklozenge$. Puzzle dualities on KTW puzzles are the diagonal reflections (linear maps) together with 0 and 1 label swapping. Puzzle dualities $\spadesuit$ and $\clubsuit$ on a LR tableau $T$ is obtained  by an hybrid  switching:   a pair consisting of a row strict Yamanouchi tableau (the transpose of a Yamanouchi tableau) and $T$  that $T$ extends  \cite[Section 2, p.22]{bss}; and a pair consisting of $T$ and  the transpose of a  Yamanouchi tableau in the anti-normal form that extends $T$.
In Appendix \ref{a:mosaic}, one considers the index two subgroup $H$ action on  puzzles and Purbhoo mosaics and
explains what operations on LR (skew) tableaux they translate to.

\medskip

 The  action of the elements in the coset $\zeta H$ exhibit  the {\em hidden} symmetries consisting of the \textit{LR transposer} (an involution exhibiting $c_{\mu\,\nu\,\lambda}=c_{\mu^t\,\nu^t\,\lambda^t}$) and the \textit{LR commuters} (an involution exhibiting $c_{\mu\;\nu\;\lambda}=c_{\epsilon\,\delta\,\theta}$ with $(\epsilon,\delta,\theta)$ obtained by commuting two  entries in $ (\mu,\nu,\lambda)$ )  respectively. (See Section \ref{sec:check} for  details.)
 \medskip

Recalling that BZ triangles, hives and LR tableaux are related through linear bijections ~\cite{PV1}, { we prove that the  involutions exhibiting the $H$-symmetries or equalities}
$$c_{\mu\;\nu\;\lambda}=c_{\mu^t\;\lambda^t\,\nu^t},\;c_{\mu\;\nu\;\lambda}=c_{\nu^t\;\mu^t\;\lambda^t},\;
c_{\mu\;\nu\;\lambda}=c_{\lambda\;\mu\;\nu},\; c_{\mu\;\nu\;\lambda}=c_{\nu\;\lambda\;\mu},\;
c_{\mu\;\nu\;\lambda}=c_{\lambda^t\;\nu^t\;\mu^t},
$$
are intrinsically easy to exhibit in every model under consideration.

While the $H$-symmetries are easy to exhibit,
 the symmetries under the action of $\varsigma_1 H=\varsigma_2 H=\varsigma_1\varsigma_2\varsigma_1H=$ $\tau H=\tau\varsigma_2\varsigma_1H=\tau\varsigma_1\varsigma_2H$, {\em i.e.}  the commutativity and the conjugation symmetries, giving the equalities
\begin{eqnarray} c_{\mu\;\nu\;\lambda}=c_{\nu\;\mu\;\lambda},
\;c_{\mu\;\nu\;\lambda}=c_{\mu\;\lambda\;\nu},\;c_{\mu\;\nu\;\lambda}=c_{\lambda\;\nu\;\mu},\label{outH}\\
c_{\mu\,\nu\,\lambda}=c_{\mu^t\;\nu^t\;\lambda^t},\;c_{\mu\;\nu\;\lambda}=c_{\lambda^t\;\mu^t\;\nu^t},\;
c_{\mu\,\nu\,\lambda}=c_{\nu^t\;\lambda^t\,\mu^t},\label{outH1}\end{eqnarray}
are difficult to exhibit.   The symmetries outside of $H$  are reduced to the action of the elements in   the coset $\varsigma H$ where $\varsigma$ is an LR commuter or an LR transposer and exhibit any of the six symmetries in \eqref{outH}.  The LR commuters and transposers in $\varsigma H$ are therefore linear  time equivalent to each other and can be reduced to the
Luzstig--Sch\"utzenberger involution realized by the Sch\"utzenberger evacuation or by the Benkart-Sottille-Stroomer reversal. The computational complexity analysis is uniform in the aforementioned combinatorial models.


\subsection{Organization}  The rest of this paper
is structured into seven sections as follows.
In the next section, divided into seven subsections, we review the needed tableau calculus for computational complexity,  linear reduction and linear equivalence of bijections. In Section \ref{crystal}, divided into four subsections, we make the connection between the tableau calculus in  previous section  and  tableau crystals.
In particular, we notice the pair of companion Gelfand-Tsetlin patterns or companion tableaux associated with a Littlewood-Richardson (LR) tableau, the linear bijection to pass from (skew) LR tableaux to their left and right companions of normal shape and recall its crystal characterization as  highest weight  and lowest elements of a crystal respectively. LR companion pais specify hives.
In Section \ref{s:maps}, divided into four subsections, we consider linear reduction and linear equivalence of bijections, and introduce the linear time bijections rotation and orthogonal transpose (composition of rotation with transposition) on LR tableaux. On this regard, algorithms are \ref{alg} and \ref{alg:blacklozenge} are defined. In Section \ref{sec:commutertransposer}, also divided into four subsections, we consider linear reduction and linear equivalence of bijections. In particular we study LR transposers  coincidence and  linear equivalence to an LR commuter which in turn  are linear time equivalent to the Lusztig--Sch\"utzenberger involution.
In Section \ref{sec:puzzle}  we study the $\mathbb Z_2\times \mathfrak{S}_3$--symmetries and the subgroup $\mathcal{H}$  of KT puzzle  dualities and rotations. In Section \ref{sec:hives},  the action of $\mathbb{Z}_2\times\mathfrak{S}_3$ on LR tableaux is translated to  the LR companion pairs in Proposition \ref{th:rho1} and Corollary \ref{cor:rho1}, and on hives is described in Theorem \ref{actionLRcompanion}. Finally in Appendix \ref{a:mosaic} we encompass our analysis with migration in Puhrboo mosaics \cite{mosaic}.

\medskip

\emph{  This work is based on the DMUC preprint \cite{acmcmuc24}. Part of it has  appeared as  extended abstract \cite{acm}  in the FPSAC 2009  proceedings or  in the preprint \cite{acmfpsac}. This paper was intended to be the full version.}

\section{Preliminaries on tableau combinatorics }

\subsection{Young (skew) diagrams and  linear transformations}
\label{sec:pre} Throughout we fix a non empty rectangle $D$ of size $d\times (n-d)$, $n>d>0$ as an ambient space.
A \emph {partition} (or straight shape, normal shape or normal form) $\lambda$ is a finite weakly decreasing sequence of
non--negative integers
$\lambda=(\lambda_1,\lambda_2,\ldots,\lambda_d),$ with at most $d$ parts (positive entries),
$\lambda_1\geq\lambda_2\geq\cdots\geq\lambda_d\ge 0$. We assume $\lambda_1\le n-d$. The
number of  parts is the \emph{length} $\ell(\lambda)\le d$,  and the \emph{ weight} is
$|\lambda|=\lambda_1+\lambda_2+\cdots+\lambda_d\le n$. We use lower case Greek letters such as $\lambda$
 to represent partitions. Often we drop the commas and parentheses when writing a partition. For instance, 2210 is the partition $(2,2,1,0)$.  
 We think of $\mathbb{Z}\times \mathbb{Z}$ as consisting of boxes, and number the rows and columns of $\mathbb{Z}\times \mathbb{Z}$ so that rows number increase bottom to top and columns number left to right. Compass directions are defined as usual according to the canonical basis of $\mathbb{R}^2$. The Young
diagram  of $\lambda$ (normal shape) is the collection of boxes $
\{(i,j)\in\mathbb{Z}^{2}|\; 1\leq i\leq d, 1\leq j\leq
\lambda_i\}$ in French convention.
We do not make distinction between a partition
$\lambda$ and its Young diagram whose diagram 
fits inside, according
to the French convention, the lower left corner of the rectangle $D$
anchored at the origin of $\mathbb{Z}\times \mathbb{Z}$.
In particular, we regard this rectangle $D$ as a
Young diagram  with $d$ parts  of   size $n-d$, and hence all our Young diagrams are subsets of the
Young diagram   $D$. 

 We let  $\binom{[n]}{d}$ denote the set of $01$--words of length $n$,  with $d$ $1$'s and $n-d$ $0$'s. There is a natural action of the symmetric group  $\mathfrak{S}_n$ on this set. In particular, the longest permutation or \emph{reverse permutation} ${\rm rev}$ acts on $\binom{[n]}{d}$ by reversing the words. Our partitions in $D$ are
identified with  the $01$--words in $\binom{[n]}{d}$ as follows: the positions of the zeroes and ones in a
$01$--word are respectively the positions of the horizontal and
vertical steps along the boundary of the corresponding Young
diagram, starting in the right lower corner of the rectangle and ending up at the upper left corner. In particular, the empty partition $\emptyset$ is identified with $0^{n-d}1^d$, and $D$ with $1^{d}0^{n-d}$.
In
 the example below, with $d=4$, $n=10$, the partition
$\lambda=4210$, depicted in green, is identified with
the $01$--word
${\color{red}00}{\color{blue}1}{\color{red}
00}{\color{blue}1}{\color{red}0}{\color{blue}1}
{\color{red}0} {\color{blue}1}\in \binom{[10]}{4}$,
\begin{align}\label{rectangle}
\includegraphics[scale=0.42,trim = 0cm 0.1cm 0cm
0cm,clip]{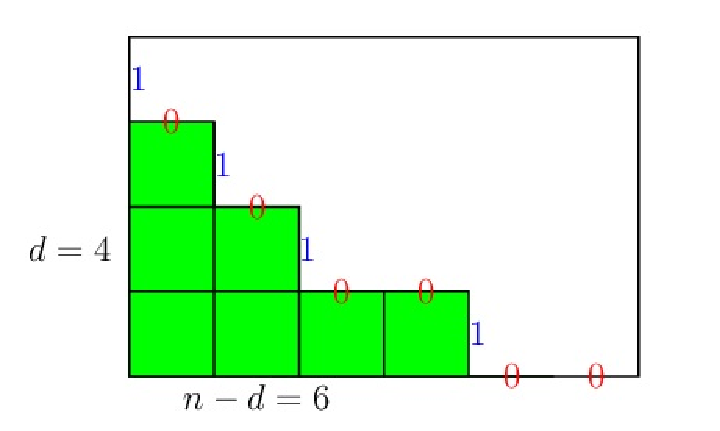}
\end{align}
 Reverting each word in $\binom{[n]}{d}$ gives the complement of each partition in $D$. The {\em complement} of $\lambda=(\lambda_1,\dots,\lambda_d)$ is the partition     $\lambda^\vee=(n-d-\lambda_d,\dots,n-d-\lambda_1)$. Equivalently, rotate by $\pi $ radians, the set complement of $\lambda$ in $D$, and against the lower left corner to obtain the  normal shape of $\lambda^\vee$.
 Since reversing is an involution on $\binom{[n]}{d}$, $(\lambda^\vee)^\vee=\lambda$. In particular, $\emptyset^\vee=D$ and $D^\vee=\emptyset$. Considering $\lambda$ in \eqref{rectangle},  the left hand side of \eqref{01word}, gives the complement of $\lambda$ in \eqref{rectangle}, $\lambda^\vee=6542=$ $
{\color{blue}1}{\color{red}0}{\color{blue}1}{\color{red}0}{\color{blue}1}{\color{red}00}{\color{blue}1}{\color{red}00}$
the reverse of the $01$--word of $\lambda$ in \eqref{rectangle}. 

Indeed $\binom{[n]}{d}$ and $\binom{[n]}{n-d}$ are in bijection. The transposition operation on partitions realizes such a bijection. The {\em transpose} (or \emph{conjugate}) of the partition $\lambda$  is the partition $\lambda^t$ obtained by reflecting $\lambda$ about the line $y=x$.
The {\em transpose} (or \emph{conjugate}) of a partition $\lambda$ in $D$,  $\lambda^t$, is given  by
the $01$--word with $n-d$ $1$'s and $d$ zeroes,  obtained from $\lambda$ by reversing its $01$ word and
swapping zeroes and  ones. In particular,  $D^t$ denotes the $(n-d)\times d$ rectangle, and  $\lambda^t$ is contained $D^t$.
Since transposing is an involution on the set $\binom{[n]}{d}\cup\binom{[n]}{n-d}$, $(\lambda^t)^t=\lambda$.
Considering $\lambda$ in \eqref{rectangle},   the middle picture in \eqref{01word}, illustrates $\lambda^t=321100={\color{red}0}{\color{blue}1}{\color{red}0}{\color{blue}1}{\color{red}0}{\color{blue}11}
{\color{red}0} {\color{blue}11}$.
The {\em complement transpose} (or {\em transpose complement}) of $\lambda$,
$\lambda^{\vee\,t}=\lambda^{t\,\vee}$, is  the complement of the transpose (or the transpose of the complement) of  $\lambda$, identified with the $01$--word
obtained from $\lambda$ by swapping zeroes and  ones. The right hand side of  \eqref{01word} illustrates $\lambda^{t\,\vee}={\color{blue}1}{\color{blue}1}{\color{red}0}{\color{blue}1}{\color{blue}1}{\color{red}0}{\color{blue}1}{\color{red}0}
{\color{blue}1}{\color{red}0}$ depicted in green.
\begin{eqnarray}\label{01word}\parbox{5cm}{\begin{tikzpicture}[scale=.35]
\node  at (-1.5,2) {$\lambda^{\vee}=$};
\draw (0,0) -- (6,0) -- (6,4) -- (0,4) -- (0,0);
\filldraw[fill=green, draw=black]
(0,0) rectangle (1,1) (1,0) rectangle (2,1) (2,0) rectangle (3,1) (3,0) rectangle (4,1) (4,0) rectangle (5,1) (5,0) rectangle (6,1)
(0,1) rectangle (1,2) (1,1) rectangle (2,2) (2,1) rectangle (3,2) (3,1) rectangle (4,2) (4,1) rectangle (5,2)
(0,2) rectangle (1,3) (1,2) rectangle (2,3) (2,2) rectangle (3,3) (3,2) rectangle (4,3)
(0,3) rectangle (1,4) (1,3) rectangle (2,4);
\end{tikzpicture}}\parbox{4cm}{\begin{tikzpicture}[scale=.4]
\node  at (-1.5,2) {$\lambda^t=$};
\draw (0,0) -- (4,0) -- (4,6) -- (0,6) -- (0,0);
\filldraw[fill=green, draw=black]
(0,0) rectangle (1,1) (1,0) rectangle (2,1) (2,0) rectangle (3,1)
(0,1) rectangle (1,2) (1,1) rectangle (2,2)
(0,2) rectangle (1,3)
(0,3) rectangle (1,4);
\end{tikzpicture}}\parbox{4cm}{\begin{tikzpicture}[scale=.4]
\node  at (-2,2) {$(\lambda^{\vee})^t=$};
\draw (0,0) -- (4,0) -- (4,6) -- (0,6) -- (0,0);
\filldraw[fill=green, draw=black]
(0,0) rectangle (1,1) (1,0) rectangle (2,1) (2,0) rectangle (3,1) (3,0) rectangle (4,1)
(0,1) rectangle (1,2) (1,1) rectangle (2,2) (2,1) rectangle (3,2) (3,1) rectangle (4,2)
(0,2) rectangle (1,3) (1,2) rectangle (2,3) (2,2) rectangle (3,3)
(0,3) rectangle (1,4) (1,3) rectangle (2,4) (2,3) rectangle (3,4)
(0,4) rectangle (1,5) (1,4) rectangle (2,5)
(0,5) rectangle (1,6);
\end{tikzpicture}.} \end{eqnarray}

A partition $\mu$ is said to be contained in a partition $\lambda$
if the Young diagram of $\mu$ is contained in the Young diagram of
$\lambda$. In this case, one defines the {\em skew shape} (or skew partition) $\lambda/\mu$ to
be the set $\{(i,j)\in\mathbb{Z}^{2}|\ts (i,j)\in \lambda,
(i,j)\notin\mu\}$ of boxes in the Young diagram of $\lambda$ that remains
after one removes those boxes corresponding to $\mu$. It is  convenient to identify $\lambda/\mu$ with the nonnegative vector $\lambda-\mu$. When $\mu$ is the null partition $\emptyset$, the
skew--diagram $\lambda/\mu$ equals the   normal shape $\lambda$. The
number of boxes in $\lambda/\mu$ is $|\lambda/\mu|=|\lambda|-|\mu|$.
 The {\em antinormal shape} (or \emph{antinormal form}) of $\lambda$ is the skew
shape $D/\lambda^\vee$. Equivalently, $\lambda$ is  $\pi$ radians-rotated  and placed against the upper right corner of $D$. We also think of it  as  the reverse of $\lambda$   fitting the upper right corner of $D$.

The {\em transpose} (conjugate) of $\lambda/\mu$ is defined to be  $(\lambda/\mu)^t:=\lambda^t/\mu^t$, or as  the image of $\lambda/\mu$ under the linear transformation
$(i,j)\mapsto (j,i)$.

 The  {\em rotate} of $\lambda/\mu$, $(\lambda/\mu)^\bullet$,  is  the image of $\lambda/\mu$ under the linear transformation
$(i,j)\mapsto (d-i+1, n-d-j+1)$.
 Equivalently
$(\lambda/\mu)^{\bullet}=\mu^\vee/\lambda^\vee$. In particular, $(D/\mu)^\bullet=\mu^\vee$, and $\lambda^\bullet:=(\lambda/\emptyset)^{\bullet}=D/\lambda^\vee$  is  the anti-normal shape of $\lambda$ and we think of $\lambda^\bullet$ as the reverse of $\lambda$.
The {\em orthogonal transpose} or the {\em rotate transpose} 
   is the composition of the transposition and the
rotation maps $\bullet t=t \bullet$.
The {\em rotate
transpose shape} $(\lambda/\mu)^{\bullet t}=$ $(\mu^\vee)^t/(\lambda^\vee)^t$ $=(\mu^t)^\vee/(\lambda^t)^\vee$ $=(\lambda/\mu)^{t\bullet }$ is then the image of
$\lambda/\mu$  under the linear transformation $(i,j)\mapsto (n-d-j+1,d-i+1)$.
 In particular, $(D/\mu)^{\bullet  t}=\mu^{\vee t}$ and $\lambda/\emptyset)^{\bullet t}=D^t/(\lambda^t)^\vee$
is the anti-normal shape of $\lambda^t$. For instance, if
$\mu=21$ and $\lambda=4210$  as above, we have
\begin{equation}\label{partoperation}\parbox{3.8cm}{\begin{tikzpicture}[scale=.35]
\node  at (-1.7,2) {$\lambda/\mu=$};
\draw (0,0) -- (6,0) -- (6,4) -- (0,4) -- (0,0);
\filldraw[fill=green, draw=black]
 (2,0) rectangle (3,1) (3,0) rectangle (4,1)
 (1,1) rectangle (2,2)
(0,2) rectangle (1,3);
\end{tikzpicture}}\parbox{4.5cm}{\begin{tikzpicture}[scale=.35]
\node  at (-2.4,2) {$(\lambda/\mu)^\bullet=$};
\draw (0,0) -- (6,0) -- (6,4) -- (0,4) -- (0,0);
\filldraw[fill=green, draw=black]
(2,3) rectangle (3,4) (3,3) rectangle (4,4)
 (4,2) rectangle (5,3)
(5,1) rectangle (6,2);
\end{tikzpicture}}\parbox{4cm}{\begin{tikzpicture}[scale=.35]
\node  at (-2.4,2) {$(\lambda/\mu)^t=$};
\draw (0,0) -- (4,0) -- (4,6) -- (0,6) -- (0,0);
\filldraw[fill=green, draw=black]
(2,0) rectangle (3,1)
 (1,1) rectangle (2,2)
(0,2) rectangle (1,3)
(0,3) rectangle (1,4);
\end{tikzpicture}}\parbox{4cm}{\begin{tikzpicture}[scale=.35]
\node  at (-2.6,2) {$(\lambda/\mu)^{t\bullet}=$};
\draw (0,0) -- (4,0) -- (4,6) -- (0,6) -- (0,0);
\filldraw[fill=green, draw=black]
(1,5) rectangle (2,6)
 (2,4) rectangle (3,5)
(3,3) rectangle (4,4)
(3,2) rectangle (4,3);
\end{tikzpicture}.} \end{equation}

\subsection{ Tableaux and Littlewood-Richardson tableaux  }\label{subsec:standard}

 A (\emph{semistandard}) {\it Young tableau} (SSYT) $ T$ of shape $\lambda/\mu$ on the alphabet $[n]$, is  a filling
of the boxes of the  skew diagram $\lambda/\mu$ with positive
integers  such that the entries are strictly  increasing in each column from
bottom to top, and  weakly increasing in each row from left to right. When $\mu$
is the empty partition we say that $T$ has normal shape $\lambda$.
 The (\emph{row reading}) \emph { word} $w(T)$ of a Young tableau $T$ is the
sequence of positive integers obtained by reading the entries of $ T$  right--to--left, the rows of $ T$,
from bottom to top. The \emph {column word} $w_{col}(T)$ is the word
obtained by reading the entries of $T$, from right to left along each column, starting in
the rightmost column and moving upwards. The nonnegative vector $m=(m_1,\ldots,m_n)$ is said to be the \emph{content} (or \emph{weight})  of $T$ if it is the
content (weight)  of its word, that is,  $m_k$ is the number of $k$'s in $T$.
 Denote
by $YT(\lambda/\mu,m)$ the set of Young tableaux of shape $\lambda/\mu$ and content $m$. If $\mu=\emptyset$ then we write $YT(\lambda,m)$ to indicate the Young
tableaux of normal shape $\lambda$ and weight $m$.

Let $\delta_r$ denote the \emph{staircase shape} partition $\delta_r=(r,r-1,\dots,1)$. A word $w$ of length $r$ can be identified with the \emph{diagonally-shaped} tableau with word $w$, a
semistandard tableau of shape $\delta_r/\delta_{r-1}$ with row reading word $w$.

 For any skew diagram $\lambda/\mu$ a {\it Littlewood--Richardson} {(LR)} tableau is a Young tableau of shape $\lambda/\mu$ such that any prefix of its word contains at least as many letters $i$ as letters $i+1$, for
all $i$. A word such that every prefix satisfies this
property is called a {\it lattice permutation, ballot} or   \emph{Yamanouchi}. Its content is always a partition. The column word of an LR--tableau is also a
Yamanouchi word of the same content. When $\mu=0$ we get the the LR tableau of straight shape $\nu$ or the
\emph{ Yamanouchi tableau} $Y(\nu)$,  the unique tableau of shape and weight
$\nu$, that is, the tableau of shape $\nu$ where each row
$i$ is  filled with $\nu_i$ $i$'s from right to left along each column, starting in
the rightmost column and moving upwards. See Example \ref{yamanouchi}.

For a given finite alphabet, say $[r]$,  a word is said to be \emph{opposite} or \emph{anti-}\emph{Yamanouchi} (\emph{anti-ballot}) if every suffix contains at least as many letters $i$ as letters $i-1$, for all $i\le r$. Its content it is always the reverse of  a partition. A Young tableau of shape $\lambda/\mu$ whose word is anti-ballot is called  an \emph{opposite} (\emph{anti-}) \emph{LR tableau}. When $\mu=0$, we get the opposite LR tableau of straight shape $\nu$ or the
\emph{ opposite Yamanouchi tableau} $Y(\nu^\bullet)$,  the unique tableau of shape $\nu$ and reverse weight
$\nu^\bullet$.  See Example \ref{oppositeyamanouchi}.

Given the rectangle $D$ with $\mu,\nu,\lambda\subseteq D$, the \emph{ boundary data} of a LR tableau of shape $\lambda/\mu$ and weight $\nu$ is $(\mu,\nu,\lambda^\vee)$, and $\rm{LR}^\lambda_{\mu,\nu}=\rm{LR}(\mu,\nu,\lambda^{\vee})$ denotes the set of LR tableaux with that boundary data. Let $c_{\mu,\nu,\lambda^\vee}$ denote the cardinal of $\rm {LR}(\mu,\nu,\lambda^{\vee})$. We define  $\rm LR(\mu,\nu^\bullet,\lambda)$ to be the set
 of {\em opposite LR tableaux} of shape $\lambda^\vee/\mu$ and weight
$\nu^\bullet$, and  { $c_{\mu,\nu^\bullet,\lambda}$} denotes its cardinality.
  In Example \ref{yamanouchi}, for $n=9$ and  $d=3$, $T$  is an LR tableau in $\rm LR(210,532,320)$ with
Yamanouchi word $w(T)=1111221332$.
\begin{example}\label{yamanouchi} Let $n=9$ and  $d=3$. $
T=\begin{tikzpicture}[scale=.5]
\filldraw[fill=green, draw=green]
(0,0) rectangle (1,1) (1,0) rectangle (2,1) (0,1) rectangle (1,2);
\draw (0,0) -- (6,0) -- (6,3) -- (0,3) -- (0,0);
\filldraw[fill=yellow, draw=yellow]
(3,2) rectangle (4,3) (4,2) rectangle (5,3) (5,2) rectangle (6,3)
(4,1) rectangle (5,2) (5,1) rectangle (6,2);
\draw (0,0) -- (6,0) -- (6,3) -- (0,3) -- (0,0);
\draw
(2,0) rectangle (3,1) (3,0) rectangle (4,1) (4,0) rectangle (5,1) (5,0) rectangle (6,1)
(1,1) rectangle (2,2) (2,1) rectangle (3,2) (3,1) rectangle (4,2)
(0,2) rectangle (1,3) (1,2) rectangle (2,3) (2,2) rectangle (3,3);
\draw (0,0) -- (6,0) -- (6,3) -- (0,3) -- (0,0);
\node at (2.5,0.5) {1};\node at (3.5,0.5) {1};\node at (4.5,0.5) {1};\node at (5.5,0.5) {1};
\node at (1.5,1.5) {1};\node at (2.5,1.5) {2};\node at (3.5,1.5) {2};
\node at (0.5,2.5) {2};\node at (1.5,2.5) {3};\node at (2.5,2.5) {3};
\end{tikzpicture}$  is an LR tableau  with boundary data $(\mu,\nu,\lambda^{\vee})$, where
$\mu=210$, $\nu=532$ and $\lambda^{\vee}=320$. Its word is $w(T)=1111221332$ and its column word is
$w_{col}(T)=1112123132$ both of content the partition $\nu=532$. The Yamanouchi tableau $Y(\nu)=\begin{Young} 3&3&&&&\cr 2&2&2&&&\cr 1&1&1&1&1&\cr
\end{Young}.$
\end{example}

\begin{example}\label{oppositeyamanouchi} An \emph{ opposite} LR tableau $ S=\begin{Young} 3&3&3&&&\cr &2&2&2&&\cr &&1&1&3&3\cr
\end{Young}$ on the alphabet $[3]$, with the same boundary data as in Example \ref{yamanouchi}, and opposite ballot word $w(S)=3311222333$ of content  $\nu^\bullet=235$, the reverse partition $\nu$. The opposite Yamanouchi tableau
$Y(\nu^\bullet)=\begin{Young} 3&3&&&&\cr 2&2&3&&&\cr 1&1&2&3&3&\cr
\end{Young}$.
\end{example}

The {\em standard order} of the boxes on a semistandard Young tableau
 is given by the numerical ordering of the labels with priority, in the case of equality,
given by rule northwest=smaller, southeast=larger.
A Young tableau with $s$ boxes is {\it standard} if it is filled
with the numbers 1 through $s$  without repetitions. The \emph{standardization} of a semistandard tableau $T$ of content $m$, denoted by $\widehat T$, is the enumeration of the
 labeled boxes according to the standard order of the boxes in $T$.
The standardization $\widehat w$ of a word $w$ is
defined accordingly, from right to left. For instance, the
standardization of the tableau $ T$ above in Example \ref{yamanouchi} is $$
\widehat{T}=\begin{tikzpicture}[scale=.5]
\filldraw[fill=green, draw=green]
(0,0) rectangle (1,1) (1,0) rectangle (2,1) (0,1) rectangle (1,2);
\draw (0,0) -- (6,0) -- (6,3) -- (0,3) -- (0,0);
\filldraw[fill=yellow, draw=yellow]
(3,2) rectangle (4,3) (4,2) rectangle (5,3) (5,2) rectangle (6,3)
(4,1) rectangle (5,2) (5,1) rectangle (6,2);
\draw (0,0) -- (6,0) -- (6,3) -- (0,3) -- (0,0);
\draw
(2,0) rectangle (3,1) (3,0) rectangle (4,1) (4,0) rectangle (5,1) (5,0) rectangle (6,1)
(1,1) rectangle (2,2) (2,1) rectangle (3,2) (3,1) rectangle (4,2)
(0,2) rectangle (1,3) (1,2) rectangle (2,3) (2,2) rectangle (3,3);
\draw (0,0) -- (6,0) -- (6,3) -- (0,3) -- (0,0);
\node at (2.5,0.5) {2};\node at (3.5,0.5) {3};\node at (4.5,0.5) {4};\node at (5.5,0.5) {5};
\node at (1.5,1.5) {1};\node at (2.5,1.5) {7};\node at (3.5,1.5) {8};
\node at (0.5,2.5) {6};\node at (1.5,2.5) {9};\node at (2.5,2.5) {10};
\end{tikzpicture},$$  and
$\widehat{w(T)}:=w(\widehat T)=5432871(10)96$. If $w=w_1 w_2\dots
w_s$ is a word and $\alpha$ is a permutation in the symmetric group
$\mathfrak{S}_{s}$, define $\alpha w= w_{\alpha(1)}\dots
w_{\alpha(s)}$. In the case $T$ is standard we have
$w_{col}(\widehat T)= rev\, w({\widehat T}^t)$, with $rev$ the
longest permutation in $\mathfrak{S}_{s}$. The transposition of  a standard tableau $T$ is still a standard tableau {  written $T^t$}. If $T$ is  semistandard then $T^t$ means a tableau strictly increasing eastward and weakly increasing northward.

\subsection{The recording matrix of a tableau}\label{sec:recordingmatrix}
Given the tableau $T\in YT(\lambda/\mu,m)$, where $m=(m_1,\dots,m_n)$ and $\ell(\lambda)\le n$,
let $M$ $=$ $(M_{ij})$ be the $n\times n$ matrix
with non--negative integer entries such that
$M_{ij}$ is the number of $j's$ in the $i$th row of $T$, called
the {\em recording matrix} of $T$ \cite{lee2,PV2}.
The recording matrix of a Young tableau of normal shape is an upper
triangular matrix. Observe that in the case of a non straight shape the recording matrix determines
the skew tableau only up to a parallel shift on the skew shape $\lambda/\mu$, see~\cite{lee2,PV2}.

\subsection{The linear involution rotation map $\bullet$ on SSYT's and  LR tableaux}\label{rotate}
Given the alphabet $[d]$, and an integer $i$ in $[d]$,  let
$i^\bullet:=d-i+1$ the {\em complement} of $i$ with respect to $[d]$.
Given the word $w=w_1w_2\cdots w_{s}$, over the alphabet
$[d]$, of weight $m=(m_1,\ldots,m_d)$,
$w^\bullet:=w_{s}^\bullet\cdots w_2^\bullet w_1^\bullet$ is the {\em dual or reverse
complement word}
 of $w$ and  ${\rm rev}\,m=(m_d,\ldots,m_1)$, the reverse of $m$, its weight.
  Indeed $w^{\bullet\bullet}=w$.
  We next extend the map $\bullet$ on words to skew tableaux recalling that a word is identified with  a diagonally shaped tableau.

  Given a Young tableau $\rm T$ of shape $\lambda/\mu$ and weight
$m$ and an ambient rectangle,  the {\em rotate} or \emph{dual} of $T$, $\bullet( T)$, is defined to be the Young tableau of
shape $(\lambda/\mu)^\bullet$ and reverse weight ${\rm rev}\, m$, obtained from $\rm T$   by rotating $\pi$ radians the shape $\lambda/\mu$
  while complementing each entry, that is, replacing each entry $i$ with  $i^\bullet$. The word of ${\rm T^\bullet}$ satisfies $w(\bullet(T))=w( T)^\bullet$, and
${\bullet\bullet( T)}=T$. The {\em rotation map} is involutive and commutes with standardization $\bullet({\widehat T })=\widehat{\bullet(T)}$.
 It is also a {linear}  map since
 $M$ is the
recording matrix of $T$ if and only if  $M^\bullet:=P_{\rm rev}MP_{\rm rev}$ is  the recording matrix of
$\bullet(T)$.
Given the finite alphabet $[d]$, the rotation map $\bullet:T\rightarrow \bullet(T),\;w(T)\mapsto w(T)^\bullet $
 is a {  linear involution
 on the set of semistandard Young tableaux} over $[d]$, which swaps the inner border with the outer border and reverses the weight.

 If  $w$ is a Yamanouchi word of weight
$\nu=(\nu_1,\dots,\nu_d)$, write $\nu^t=(\nu_1^t,\dots,\nu_{\nu_1}^t)$ and
observe that $w$ is a shuffle of the words $12\cdots\nu^t_i$ for
$i=1,\dots,\nu_1$. Similarly, if $w$ is an opposite Yamanouchi word of weight $\nu^\bullet$, $w$ is a shuffle of the words $d-\nu_i^t+1\cdots
d-1  d$, for  $i=1,\dots,\nu_1$.
Thus, fixing such a shuffle for $w$, we may obtain $w^\bullet$,  in the Yamanouchi case,  by first  replacing, for each $i=1,\dots,\nu_1$, the word $12\cdots\nu^t_i$  in the shuffle of $w$ with   $d d-1\cdots d-\nu_i^t+1$ and then reversing the resulting word; and, in the opposite Yamanouchi case, by first  replacing, for each $i=1,\dots,\nu_1$, the word $d-\nu_i^t+1\cdots d-1 d$  in the shuffle of $w$ with   $\nu^t_i\cdots 2 1$  and then reversing the resulting word.
(The result does not depend on the chosen shuffle.) The dual Yamanouchi word $w^\bullet$, in the former case is a shuffle of the words $d-\nu_i^t+1\cdots d-1  d$, for  $i=1,\dots,\nu_1$, that is, an opposite Yamanouchi word, and a Yamanouchi word in the latter.
Thus, a word is {\em opposite Yamanouchi } if and only its  \emph{dual word is Yamanouchi}.
Therefore, for a skew diagram $\lambda/\mu$, an {\it opposite Littlewood--Richardson} {(LR)} tableau is a Young tableau of shape $\lambda/\mu$ such that its dual is an LR tableau or its word is the dual of a Yamanouchi word. When $\mu=0$, one obtains the rotate Yamanouchi tableau $\bullet(Y(\nu))$.

 The rotation $\bullet$ of $\rm LR(\mu,\nu^\bullet,\lambda)$,  the set of opposite LR tableaux of shape $\lambda/\mu$ and content $\nu^\bullet$, gives $\rm LR(\lambda,\nu,\mu)$ the set of LR tableaux of shape $(\lambda^\vee/\mu)^\bullet=\mu^\vee/\lambda$
and content $\nu$, and {\em vice-versa}.
\begin{proposition}\label{rotation}
The rotation map
 \begin{equation}\label{eq:rotation}\bullet:\rm LR(\mu,\nu,\lambda)\cup \rm LR(\lambda,\nu^\bullet,\mu)\rightarrow\rm
\rm LR(\lambda,\nu^\bullet,\mu)\cup LR(\mu,\nu,\lambda),\;T\mapsto \bullet( T),
\end{equation}
is a linear  involution on $\rm LR(\mu,\nu,\lambda)\cup\rm\rm LR(\lambda,\nu^\bullet,\mu)$ that transforms the LR tableau $T$ into  its dual $\bullet(T)$, the opposite LR tableau  of  shape $(\lambda^\vee/\mu)^\bullet$ and content $\nu^\bullet$, and {\em vice versa}. It exhibits the symmetry $c_{\mu,\nu,\lambda}=c_{\lambda,\nu^\bullet,\mu}$.
\end{proposition}
\begin{example}\label{ex:rotateyam} Let $d=3$ and $n=7$. Given $\nu=421$, the Yamanouchi $Y(\nu)$ and  the rotate Yamanouchi tableau $\bullet({Y(\nu)})$ are displayed below
$$Y(\nu)=\begin{Young} 3&&&\cr 2&2&&\cr 1&1&1&1\cr
\end{Young},\;\;w=1111223,\; \;\bullet(Y(\nu)) =\begin{Young} 3&3&3&3\cr &&2&2\cr &&&1\cr
\end{Young},\,w^\bullet=1223333.$$
\end{example}
For $d=4$ and $n=8$,
 $\nu=4210$,  $\bullet(Y(\nu)) =\begin{Young} 4&4&4&4\cr &&3&3\cr &&&2\cr &&&\cr
\end{Young},\,w^\bullet=2334444.$

\subsection{Jeu de taquin and Burge correspondence }\label{subsec:taquinburge}
 \subsubsection{Jeu de taquin}Whenever partitions $\nu\subseteq\mu\subseteq\lambda$, we say that the shape
$\lambda/\mu$ extends the shape $\mu/\nu$. An {\it inside corner} of
$\lambda/\mu$ is a box in  the diagram $\mu$ such that the boxes above and to the right (if any) are not in $\mu$.
When a box extends $\lambda/\mu$, this box is called an {\it outside corner}.
Let $\rm T$ be a (semi standard) Young tableau of shape $\lambda/\mu$, and let $b$ be an inside corner for $
T$. A {\it contracting slide}, see~\cite{bss,schu} of $T$ into the
box $b$ is performed by moving the empty box at $b$ through $\rm T$,
successively interchanging it with the neighbouring integers to the
north and east according to the following rules:
$(i)$ if the empty box has only one neighbour,  interchange it with that
neighbour;
    $(ii)$ if it has two unequal neighbours, interchange it with the
    smaller one; and
    $(iii)$ if it has two equal neighbours, interchange it with that one to
    the north.
The empty box moves in this fashion until it  becomes an outside
corner. This contracting slide can be reversed by performing an
analogous procedure over the outside corner, called an {\it
expanding slide}. This procedure is  known as Sch\"utzenberger {\it jeu de taquin}. Performing  contracting slides over successive inside
corners  in $\mu$  reduces $ T$ to a tableau $T^{\rm n}:=\texttt{rect}({T})$ of
normal shape,  called the {\it rectification} or {\em the normal form} of $ T$. The rectification of $T$ is independent of the particular sequence of inside corners used, \cite{thomas}, and so
${\rm rect}({T})$ is well defined. When $T\in \rm LR(\mu,\nu,\lambda)$,
${\rm rect}T=Y(\nu)$, and if $T\in \rm LR(\mu,\nu^\bullet,\lambda)$, ${\rm rect} T=Y(\nu^\bullet)$.

Similarly, inside the rectangle $D$, there exists exactly one tableau of anti-normal shape { $T^{\rm a}:=\texttt{ arect}T$} produced by the \emph{reverse jeu de taquin} by performing expanding slides over each successive  outside corner  in $D/\lambda$ ~\cite{bss}, called the {\em anti-normal form} (or the {\em contre-tableau} or \emph{anti-rectification}) of $T$. In particular, if $\mu=0$, the anti-rectification of $T$ produces $\lambda^\vee$. Applying reverse
{\em jeu de taquin}  slides to the canonical LR tableau or Yamanouchi tableau $Y(\nu)$ of shape $\nu$ inside the $D$ rectangle, we obtain  its anti-normal form ${\rm arect}Y(\nu)$. As ${\rm arect}Y(\nu)$ fits the
upper right corner of  $D$, ${\rm arect}Y(\nu)$ is the LR tableau  of anti-normal shape $D/\nu^\vee$  and content $\nu$.
For instance,
\begin{equation}\label{anti} Y(\lambda)=\begin{Young} 3&&&&\cr 2&2&&&\cr 1&1&1&1&\cr
\end{Young}\longleftrightarrow {\rm arect}Y(\lambda) =\begin{Young} &1&1&2&3\cr &&&1&2\cr &&&&1\cr
\end{Young}.\end{equation}
The (anti) rectification of a word means the (anti) rectification of the diagonally shaped tableau with that word.
\subsubsection{Burge correspondence a variation of Robinson-Schensted-Knuth correspondence} \label{burge} We consider a \emph{variation of the RSK--correspondence} on a  two-line array known as the {\em Burge correspondence}, see~\cite{burge}, and \cite[Appendix A.4.1]{ful}, where the ordering on the two-line array $W=\left(\begin{smallmatrix}a_1&a_2&\dots&a_N\cr
b_1&b_2&\dots&b_N\end{smallmatrix}\right)$ of positive integers is such that $a_i<a_{i+1}$, or $a_i=a_{i+1}$ and $b_i\ge b_{i+1}$, called \emph{Burge array}. The procedure to transform the biword $W$ into the semistandard tableau pair $(P(W),Q(W))$ of the same shape is the column bump $(\cdots(b_1\leftarrow b_2)\cdots)\leftarrow b_N$ and place in $a_1,\dots,a_N$ respectively. That is, $Q(W)$ is defined to be the semistandard tableau of shape $\lambda$ such that if $P(W)$ has shape $\lambda$ and $b_i$ is inserted in $(\cdots(b_1\leftarrow b_2)\cdots)\leftarrow b_{i-1}$ to create a node in $\lambda$ then we fill the node with $a_i$.  As usual, when the first line is the permutation $12\cdots N$ (the standardization of the first row of $W$), we identify $W$ with the second line word $w:=b_1b_2\cdots b_N$. In this case, $P(w)=P(W)$ and $Q(w)= \widehat{ Q(W)}$. The insertion tableau in the Burge  and in the  RSK correspondence coincide as follows:  the column bumping $(\cdots(b_1\leftarrow b_2)\cdots)\leftarrow b_N$  is equal to the row bumping $(\cdots(b_N\leftarrow b_{N-1})\cdots)\leftarrow b_1$  giving $P(w)$. Instead of biwords we may consider matrices of nonnegative integers. The biword $W$ may also be described by the $m\times t$ matrix whose $(i,j)$ entry is the number of times $\binom{i}{j}$, $i\in [m]$ and $j\in[t]$, occurs in the array. The Burge correspondence is then a correspondence between matrices $A$ with nonnegative entries and ordered pairs $(P,Q)$ of tableaux of the same shape.

 From now on when we refer to the RSK-correspondence we mean the Burge correspondence. Thanks to RSK--correspondence a word $w$ is uniquely determined by a tableau pair $(P(w),$ $ Q(w))$
 of the same normal shape, where ${P}(w)$ is the \emph{insertion tableau}
obtained by \emph{column insertion} of the letters of $w$ from left to right, and $Q(w)$ standard, called the \emph{$Q$--symbol} or \emph{recording tableau} of $w$. Reciprocally, every tableau pair $(P,Q)$ of the same shape with $Q$ standard determines a unique word on the alphabet of $P$ and with same weight.
Given the alphabet $[d]$, the RSK correspondence gives a bijection between $\textsf{words(k)}$ the set of words of length $k\ge 0$ and pairs of tableaux \begin{equation}\label{rsk}\textsf{RSK}: \textsf{words}(k)\longrightarrow \bigsqcup_{|\lambda|=k} \rm SSYT(\lambda)\times \rm SYT(\lambda), \;\textsf{RSK}(w)=(P(w),Q(w)).
\end{equation}
Insertion can be
translated into the language of Knuth elementary transformations on a word \cite{ful}.
Two words $w$ and $v$ are said \emph{Knuth equivalent} if one can be transformed into another by a sequence of Knuth moves. Two words are Knuth equivalent if and only if they have the
same insertion tableau. Each Knuth class is in bijection with the
set of all standard tableaux with  normal shape given by the unique tableau of normal shape in
that Knuth class.

 Two tableaux $ T$ and $ R$ of arbitrary shape are {\it Knuth equivalent},
written $ T\equiv R$, if and only if $P(w(T))=P(w(R))$. Since row and column words of $T$ are Knuth equivalent, one also has  $P(w(T))=P(w_{col}(T))$ \cite{ful}.
 Sch\"utzenberger sliding  in a skew tableau $T$  preserves the Knuth class of its word. Thereby,  $P(w(T))=T^{\rm n}$, and $T\equiv R$ if and only if $T^{\rm n}= R^{\rm n}$, i.e. one of them can be
transformed into the other
 with a sequence of contracting and expanding {\em jeu de taquin}
slides. In particular, inside $D$, there exists exactly one tableau $T^{\rm a}$ of anti-normal shape Knuth equivalent to $T$.  Hence $T\equiv R$ if and only if $T^{\rm n}=R^{\rm n}$, equivalently, $T\equiv R$ if and only if $T^{\rm a}=R^{\rm a}$. Recall that $ T\equiv R$ if and only if  $\widehat T\equiv \widehat R$, and  $P(w(\widehat T))=\widehat {P(w(T)})$, $Q(w(\widehat T))=Q(w(T))$.

\subsubsection{The recording tableau in Burge correspondence and LR  tableaux} \label{burgecomp}
Under Burge correspondence, there is a bijection between Burge arrays $\left(\begin{smallmatrix}{\bf y}\\
w
\end{smallmatrix}\right)$, where $w$ is a Yamanouchi word of weight $\nu$, and tableau pairs $(Y( \nu),G)$ where $G$ is of shape $\nu$ and weight $\rm wt(\bf y)$.
Let $w=w_1
w_2\cdots w_s$ be  a Yamanouchi word of content $\nu$ such that $\textsf{RSK}(\left(\begin{smallmatrix}{\bf y}\\
w
\end{smallmatrix}\right))=(Y(\nu),G)$ for some Burge array $\left(\begin{smallmatrix}{\bf y}\\
w
\end{smallmatrix}\right)$,  and put the number $k$ in
the $w_k$th row of the diagram $\nu$.
The  labels of the $i$th row
are the $k$'s such that $w_k=i$, thus its length is  $\nu_i$ and its
shape is $\nu$. We denote this standard tableau of shape $\nu$ by $U(w)$. Hence, $\textsf{RSK}(w)=(Y(\nu),U(w))$ where $U(w)=\widehat G$.

\begin{definition}Given the partition $\nu$,  \emph{$\mathbb{Y}(\nu)$ denotes the set of Yamanouchi words of weight $\nu$}.
\end{definition}
Any tableau pair $(Y(\nu),P)$ with $P$ standard of shape $\nu$ produces a Yamanouchi word of weight $\nu$, and thus $P=U(w)$.
Then the
 map $w\mapsto U(w)$ defines a bijection   between Yamanouchi words of content $\nu$ and standard Young tableaux with  shape $\nu$. Hence $|\mathbb{Y}(\nu)|=|  \textsf{SYT}|$.
In Example \ref{yamanouchi},  $w=1111221332$, a  Yamanouchi word of content $\nu=532$,
 gives \begin{equation}\label{U} U(w)=
 \begin{Young} 8&9&&&\cr 5&6&10&&\cr 1&2&3&4&7\cr
\end{Young},\end{equation}
 where the entries of the $i$th row are the positions  of the $i$'s in the word  of $T$ (according to the LR numbering).

Let  $T\in \rm LR(\mu,\nu,\lambda^\vee)$. 
 We may associate to $T$ two biwords (or matrices) $W^{\lambda/\mu}$ and $W^\nu$  consisting of the same biletters   but ordered differently.  Consider the words ${\bf y}=1^{\lambda_1-\mu_1}$ $2^{\lambda_2-\mu_2}$ $
\dots{\ell(\lambda)}^{\lambda_{\ell(\lambda)}-\mu_{\ell(\lambda)}}$ of weight $\lambda/\mu$ and  ${\bf x}=1^{\nu_1}2^{\nu_2}\dots
{\ell(\nu)}^{\nu_{\ell(\nu)}}$ of weight $\nu$, both of length $|\nu|=|\lambda|-|\mu|$, and put
\begin{eqnarray} W^{\lambda/\mu}:=\left(\begin{smallmatrix}{\bf y}\\
w(T)
\end{smallmatrix}\right),&&
W_{\nu}:=\left(\begin{smallmatrix}{\bf{g}}\\
{\bf x}
\end{smallmatrix}\right)\label{biword}
\end{eqnarray}
where $W_{\nu}$  is a reordering of the biletters of $W^{\lambda/\mu}$
such that ${\bf g}=g_1g_2\dots g_{|\nu|}$ satisfies $g_i\ge g_{i+1}$ whenever $x_i=x_{i+1}$. The first, as a matrix, is a lower triangular matrix, and the second, as a matrix, is the transpose of the former thus an upper triangular matrix. See Example \ref{ex:burge}.
The symmetry of Burge correspondence \cite{ful,Lot,double} gives the following result.
\begin{proposition}\label{prop:burge} Let  $T\in \rm LR(\mu,\nu,\lambda^\vee)$ with the Burge arrays \eqref{biword}.
Then under RSK correspondence one has:

  $(a)$ $W^{\lambda/\mu}\mapsto(Y(\nu),G)$ and $W_{\nu}\mapsto(G,Y(\nu))$ where $Q(W_\nu)=Y(\nu)=P(w(T))$ and $Q(W^{\lambda/\mu})=G=P({\bf g})$.

$(b)$ $w(G)={\bf g}$, that is, $W_{\nu}=\left(\begin{smallmatrix}{w(G)}\\
{\bf x}
\end{smallmatrix}\right)$.

 $(c)$ {$\widehat G=\widehat{Q(W^{\lambda/\mu})}=Q(w(T))=U(w(T))$}, that is, the
$Q$-symbol of a Yamanouchi word $w$, with respect to Burge correspondence, is $U(w)$. That is,
$G$ is a semistandard Young tableau of shape $\nu$ and content $\lambda/\mu$, such that  each row $i$  tell us in which rows of $T$ the $i$'s are filled in.
\end{proposition}

 \begin{example}\label{ex:burge} Let $T$ be the LR tableau in Example \ref{yamanouchi}.
 The Burge correspondence gives: $W^{\lambda/\mu}=\left(\begin{smallmatrix}1^4&2^3&3^3\\
&w(T)&
\end{smallmatrix}\right)=\left(\begin{smallmatrix}1&1&1&1&2&2&2&3&3&3\\
1&1&1&1&2&2&1&3&3&2
\end{smallmatrix}\right)=\left(\begin{smallmatrix}4&0&0\\
1&2&0\\
0&1&2
\end{smallmatrix}\right)\rightarrow (Y(\nu), G)$,
and $W_{\nu}=\left(\begin{smallmatrix}&w(G)&\\
1^5&2^3&3^2
\end{smallmatrix}\right)=\left(\begin{smallmatrix}2&1&1&1&1&3&2&2&3&3\\
1&1&1&1&1&2&2&2&3&3
\end{smallmatrix}\right)=\left(\begin{smallmatrix}4&1&0\\
0&2&1\\
0&0&2
\end{smallmatrix}\right)\rightarrow(G, Y(\nu))$
where    $Y(\nu)=P(w(T))=Q(W^\nu)$ where \begin{equation}\label{GG} G=\begin{Young} 3&3&&&&\cr 2&2&3&&&\cr 1&1&1&1&2&\cr
\end{Young}\end{equation} of normal shape $\nu=532$ and
weight $\lambda/\mu=643-210=433$. The standardization of $G$ gives $\widehat G=U(w(T))$ in \eqref{U}.
\end{example}

\subsection{Evacuation, reverse complementation, rotation and RSK}

Given a tableau $\rm T$ of normal shape, the {\em Sch\"utzenberger evacuation} of $T$, ${\rm evac T}$, is  a tableau with the shape  and reverse weight of $T$ that can be characterized in different ways:
 the normal form of the rotation of
 $ T$, ${T^\bullet}^{\rm n}$;  the insertion tableau of the word
$w(T^\bullet)=w({\rm T})^\bullet$ \cite{ful};  or  the rotation of  the anti--normal form $T^{\rm a}$. Thus $\rm evac T={T^{\rm a}
}^\bullet=T^{\bullet \rm n}=P(w(T)^\bullet)$ and { $T^{\rm a}:=\texttt{arect} (T)=\bullet \texttt{evac}(T)$}. 
Indeed $\rm evac\, \rm evac T=T$. Given the Yamanouchi tableau $Y(\nu)$, its opposite satisfies $Y(\nu^\bullet)={Y(\nu)^{\rm a}}^\bullet={Y(\nu)^\bullet}^{\rm n}= \rm evac Y(\nu)$. For instance, using Example \ref{ex:rotateyam}
$$Y(\nu)=\begin{Young} 3&&&\cr 2&2&&\cr 1&1&1&1\cr
\end{Young}\;\;w=1111223,\; \;Y(\nu)^\bullet =\begin{Young} 3&3&3&3\cr &&2&2\cr &&&1\cr
\end{Young},\,w^\bullet=1223333,$$
$$ Y(\nu)^{\rm a} =\begin{Young} 1&1&2&3\cr &&1&2\cr &&&1\cr
\end{Young},\;\;
 Y(\nu^\bullet)=\begin{Young} 3&&&\cr 2&3&&\cr 1&2&3&3\cr
\end{Young}={Y(\nu)^{\rm a}}^\bullet={Y(\nu)^\bullet}^{\rm n}=\rm evac Y(\nu)=P(w^\bullet).
$$
\begin{proposition}\label{duality}{\em [Duality of Burge correspondence \cite[Appendix A 4.1]{ful}.]} Let $w$ and $u$ be two words. Under Burge correspondence

$(a)$ The word $w$ corresponds to the tableau-pair $(P,Q)$ if and only if $w^\bullet$ corresponds to $(\rm evac \,P, \rm evac \,Q)$.

 $(b)$ For any tableau $T$, $w(T)$ corresponds to the tableau-pair $(T^{\rm n}, Q)$ and $w(T)^\bullet=w(T^\bullet)$  to the pair $({T^\bullet}^{\rm n}={\rm evac T^{\rm n}},\rm evac\,Q)$.

  $(c)$   $w\equiv u$ if and only if $w^\bullet\equiv u^\bullet$, and $Q(u)=Q(v)$ if and only if $Q(u^\bullet)=Q(v^\bullet)$.
  Similarly, $\rm rev\,w$ corresponds to $(P^t,Q^{Et})$.
  \end{proposition}
  \subsection{Tableau switching and reversal involution 
}\label{subsec:switching}
 In this subsection we recall Haiman's result \cite[Theorem 2.13]{h2}: a skew tableau is uniquely determined by the skew shape,  dual Knuth class and  Knuth class (rectification or anti normal form).

Two tableaux $ T$ and $ R$ of the same shape are said to be {\it dual Knuth
equivalent}, written ${T}\overset{d}{\equiv}{R}$, if for some sequence
of contracting slides or/and expanding slides that can be applied to
one of them, can also be applied to the other, and the sequence of
shape changes is the same for both, see~\cite{h2}. In fact if two tableaux on the same shape have the same shape changes for some sequence of {\em  jeu de taquin} slides they  also have the same shape changes for any other. Hence, dual Knuth equivalent tableaux have the same (skew) shape as well as the same shape of their normal forms and  the same anti--normal shape of their anti--normal forms.
Moreover, two tableaux of the same normal shape or anti--normal shape are dual Knuth equivalent \cite[Proposition 2.14]{h2}.
 \emph{Dual Knuth equivalence on tableaux of the same shape} may
also be characterized by the $Q$ symbols or recording tableaux of their words in the RSK correspondence:
 ${T}\overset{d}{\equiv}{R}$ if and only if ${Q}(w({T}))={Q}(w({R}))$. In addition, either row  or column words may be used, ${Q}(w({T}))={Q}(w({R}))$ if and only if ${Q}(w_{col}({T}))={Q}(w_{col}({R}))$.  \emph{Dual Knuth equivalence on words} of the \emph{same length} is defined by identifying  two words of the same length with the same $Q$-symbol. Alternatively, we may identify a word of length $r$ with  the diagonal shape tableau $\delta_r/\delta_{r-1}$ with that word, where $\delta_i=(i,i-1,\dots,2,1)$, for $i=r-1,r$, and apply the definition of dual Knuth equivalence on tableaux of the same shape.

Let
$ S$ and $ T$ be tableaux such that $\rm T$ extends $\rm S$, that is, the outer border of $S$ is the inner border of $T$, and
consider the set union $\rm S\cup T$.
The {\em tableau switching}, see~\cite{bss}, may be presented
as a procedure based on {\em jeu de taquin} elementary moves, for moving two tableaux past one another,  transforming $\rm S\cup \rm T$
into $\rm A\cup \rm B$, where $B$ is a tableau Knuth equivalent to $
T$ which extends $\rm A$, and $A$ is a tableau Knuth equivalent to $
S$. We write $\rm S\cup \rm T\overset {\bf s}\longrightarrow \rm A\cup
\rm B$. In particular, if  $ S$ is of normal shape,  $\rm A=\rm
T^{\rm n}$, and $ S= B^{\rm n}$. Switching of $S$ with $T$ may be
described as follows:  $\widehat T$ is a set of instructions telling
where expanding slides can be applied to $S$. (Similarly, $\widehat S$ is a set of instructions telling
where contracting slides can be applied to $T$.) Moreover, switching commutes with standardization. Switching and
dual Knuth equivalence are related as in the theorem below. It combines tableau switching \cite{bss} with Haiman dual equivalence \cite[Corollaries 2.8, 2.9]{h2}.
\begin{theorem} \label{t1}
 \cite[Theorem 4.3]{bss},  \cite[Corollaries 2.8, 2.9]{h2}. Let $T$ and $U$ be tableaux with the same shape and dual equivalent and let $W$ be a tableau
which extends $T$ (or $T$ extends).
 If $ T\cup
      W\overset{s}\longrightarrow  Z\cup  X  $ and $ U\cup  W\overset{s}\longrightarrow  Z\cup
     Y$, then $ X\overset{d}\equiv  Y$. If $ W\cup T
      \overset{s}\longrightarrow  Z\cup  X  $ and $ W\cup  U\overset{s}\longrightarrow  Z\cup
     Y$, then $ X\overset{d}\equiv  Y$.
\end{theorem}
\begin{theorem}\label{t2}
    \cite[Theorem 2.13]{h2}. Let $\mathcal{D}$ be a dual Knuth equivalence class and $\mathcal{K}$
    be a Knuth equivalence class, both corresponding to the same normal shape (that is, the elements of $\mathcal{D}$ rectify to the normal shape of the unique tableau of normal shape in $\mathcal{K}$).
    Then, there is a unique tableau in $\mathcal{D}\cap\mathcal{K}$ which is the unique tableau in $\mathcal{D}$ that rectifies to the normal shape of the unique tableau of normal shape in $\mathcal{K}$. Tableau switching ${\bf s}$ allows to compute ${\mathcal D}\cap {\mathcal K}$.
\end{theorem}

\begin{alg} \cite{bss} \label{alg:KD} Computation of ${\mathcal D}\cap {\mathcal K}$ with ${\mathcal D}$ and $ {\mathcal K}$ corresponding to the same normal shape. Let
$Q\in {\mathcal D}$ and let $V\in {\mathcal K}$ be the unique tableau
with normal shape in this Knuth class ($V$ and $Q^{\rm n}$ have the same normal shape), and $W$ any tableau of normal shape  that $Q$ extends:

Step 1. Compute \begin{equation}\begin{matrix}\label{linearcoststep}
W\cup Q &&W\cup X\\
{\scriptstyle \bf s}\!\downarrow&&\uparrow\!{\scriptstyle \bf s}\\
Q^{\rm n}\cup Z&\rightarrow&V\cup Z.
\end{matrix}\end{equation}

Step 2.  $X\overset{d}{\equiv} Q$,
$X\overset{}{\equiv}V$,
and $\mathcal{D}\cap\mathcal{K}=\{X\}$ where $X$ is the only tableau in $\mathcal{D}$ that rectifies to $V$.
 \end{alg}

 In particular, if $\mathcal{K}$ is the Yamanouchi Knuth class given by the normal shape corresponding to $\mathcal{D}$, $X$ is the only LR tableau in $\mathcal{D}$  whose content is the normal shape corresponding to $\mathcal{D}$.

 \remark \label{re:KD} In Algorithm \ref{alg:KD}, $\mathcal{D}$ and $\mathcal{K}$ also correspond to the same anti-normal shape. If in \eqref{linearcoststep} we  consider $V^{\rm a}\in \mathcal{K}$,   the  anti-normal form of $V$, and $W$ any tableau with anti-normal shape that extends $Q\in \mathcal{D}$, then     $Q\cup W  \overset{\bf s} \rightarrow  Z\cup Q^{\rm a}\rightarrow Z\cup V^{\rm a} \overset{\bf s} \rightarrow X\cup W$ ($Q^{\rm a}$ and  $V^{\rm a}$ have the same anti-normal shape)  to obtain $\mathcal{D}\cap\mathcal{K}=\{X\}$. Note that $X\overset{d}{\equiv} Q^{\rm a}\overset{d}{\equiv} Q$ and
$X\overset{}{\equiv}V^{\rm a}{\equiv}V$. { Note also that $Q^{\rm a}=\bullet  \texttt{evac}Q^{\rm n}$.}

\begin{corollary} LR
tableaux  form a complete transversal for the set of dual Knuth
equivalence classes. Moreover, the LR coefficient $c_{\mu\nu}^\lambda$
counts the number of dual Knuth equivalence classes of tableaux of shape $\lambda/\mu$ whose
rectification has shape  $\nu$.
\end{corollary}

\subsubsection{The reversal involution} \begin{definition} \cite{bss} Given a tableau ${T}$ of any shape, the {\em reversal
} of $T$,  ${T}^e$, is  defined to be the unique tableau  Knuth equivalent to ${\rm T}^\bullet$ (${\rm T}^\bullet\equiv T^{\bullet{\rm n}}=\rm evac T^{\rm n}$), and
dual Knuth equivalent to ${\rm T}$. In other words, $T^e$ is the unique tableau dual equivalent to $T$ that rectifies to the evacuation of $T^{\rm n}$, that is, $T^{e\,\rm n}=\rm evac T^{\rm n}$. If $T$ has normal shape,
$\rm evac T=T^e$.
\end{definition}

 Algorithm \ref{alg:KD}  calculates
$T^e=[\rm evac T^{\rm\, n}]_K\cap[T]_{dK}$ (by abuse of notation we omit the brackets in $\{T^e\}$),
 where ~~~$[$\; $ ]_K$ denotes Knuth
class and $[$\; $ ]_{dK}$ dual Knuth class. That is, we choose any $W$ of straight shape that  $T$ extends, to form  $W\cup T$, we rectify $T$, using $W$ as  a set of \emph{jeu taquin} instructions  to obtain $T^{\rm\, n}\cup Z$, replace $T^{\rm\, n}$ with $\rm evac T$ to obtain $\rm evac T\cup Z$. Then  by reverse \emph{jeu de taquin} instructed by $Z$, we obtain $W\cup T^e$. Alternatively, if $\textsf{RSK(w(T))}=(T^{\rm n},Q))$ and since $\textsf{RSK(w($T^e$))}=(\rm evac T^{\rm n},Q)$
 then $T^e$
 can be calculated as $\textsf{RSK}^{-1}(\rm evac T^{\rm n},Q)=w(T^e)$.

The mapping $T\rightarrow T^e$ is called {\em reversal} and is an involution on the set of SSYT which preserves the shape and reverses the weight.
  Observe that $T^{e\,e}=[T^{e\,\rm n\,E}]_K\cap[T^e]_{dK}=[T^{\bullet\rm n\,E}]_K\cap[T]_{dK}=[T^{\rm n\,E E}]_K\cap[T]_{dK}=T$.

\begin{corollary}\label{cor:reversal} Let ${T}\in{\rm LR}(\mu,\nu,\lambda^\vee)$. The reversal of $T$, ${T}^e=[Y(\nu^\bullet)]_K\cap[T]_{dK}$ is the only opposite LR tableau in ${ \rm LR}(\mu,\nu^\bullet,\lambda^\vee)$ dual Knuth equivalent to $T$. Opposite  LR
tableaux  form another complete transversal for the set of dual Knuth
equivalence classes. The LR coefficient $c_{\mu\nu}^\lambda=c_{\mu\nu^\bullet}^\lambda$
also counts the number of dual Knuth equivalence classes of tableaux of shape $\lambda/\mu$ whose
anti--normal form has shape  $\nu^\bullet$.
\end{corollary}

\section{Crystals of tableaux,  Luzstig-Sch\"utzenberger involution and hives}\label{crystal}
Kashiwara and Nakashima \cite{NK 94} has shown that semistandard tableaux can be arranged into crystals. We recall  briefly and refer to \cite{kwon,bumschi}.

A $\mathfrak{gl}_d$-crystal is a finite set $B$ along with maps
$\rm wt :B\rightarrow \mathbb{Z}^r$
$e_i, f_i :B\rightarrow B \cup \{0\}$
for $i = 1,\dots,d$, obeying the following axioms for any $b, b'\in B$,

(i) if $e_i(b) \neq 0$ then $\rm wt(e_i(b)) = wt(b) + \alpha_i$,

(ii) if $f_i(b) \neq 0 $ then $\rm wt(f_i(b)) = wt(b)-\alpha_i$,

(iii) $b' = e_i(b)$ if and only if $b = f_i(b')$, and

(iv) if $b, b' \in B$ such that $e_i(b) = f_i(b') = 0$ and $f^k_i (b) = b'$ for some $k\ge  0$,
then $\rm wt(b') = s_i \rm wt(b)$,

\noindent where $\alpha_i=\epsilon_i-\epsilon_{i+1}$ are the $A_{d-1}$ simple roots,  and $s_i$ is the simple transposition of $\mathfrak{S}_d$, $i=1,\dots,d-1$.
The crystals that we deal with also allow to define length functions $\varepsilon_i,\varphi_i:B\rightarrow \mathbb{Z}$
$i=1,\dots,d-1$,
$$\varepsilon_i(b)=\max\{a:e_i^a(b)\neq 0\},\quad \varphi_i(b)=\max\{a:f_i^a(b)\neq 0\}.$$

Let $B_d=\{1,\dots,d\}$ be the standard $\mathfrak{gl}_d$-crystal consisting of the words of  a sole letter on the alphabet $[d]$ whose coloured crystal graph is $$1\overset{1}\longrightarrow 2\overset{2}\longrightarrow \cdots \overset{d-1}\longrightarrow d-1\overset{d-1}\longrightarrow d.$$ The Kashiwara   raising operators $f_i$ and lowering operators $e_i$ are defined for $i\in I=[d-1]$ as follows: $f_i(i)=i+1$,  $f_{d-1}(d)=0$, and $e_i(i+1)=i$, $e_1(1)=0$, otherwise, the letters are unchanged.
 The weight $\rm wt(b)=\epsilon_b$, for $b=1,\dots, d$,  the canonical basis of $\mathbb{R}^d$. The highest (lowest) weight element of $B_d$ is the word $1$ ($d$), and the highest (lowest)  weight is $\epsilon_1$ ($\epsilon_d$).

The tensor product of crystals allows us to define the crystal
$\mathcal{W}_d=\bigsqcup\limits_{k> 0} {B_d}^{\otimes k}\sqcup\{\emptyset\}$ of all finite words on $[ d]$ where $\emptyset$  is the empty word and the vertex $w_1\otimes \dots \otimes w_k\in {B_d}^{\otimes k}$ is identified with the word $w=w_1\cdots w_k$ of length $k$ on $[ d]$. We describe the crystal  of ${B_d}^{\otimes k}$ as the crystal structure on the set $\textsf{word(k)}$ of all words of length $k$ on $[d]$. Following  the tensor product rule, the action of the Kashiwara \emph{raising} and \emph{lowering operators} $e_i$ and $f_i$ on $w$, for $i\in[d-1]$, is  given by the $i$-\emph{signature rule} \cite{NK 94,kwon} which is induced from those operators on $B_d$. We substitute each $i$ by $+$  and each $i+1$ by $-$, and erase the letter in any other case. Then successively erase any pair $+-$ until all the remaining letters form a word  $\textsf{sign(w)}_i=-^a +^b$. We define $\varphi_i(w):=b$ and $\varepsilon_i(w):=a$.  If $a=0$, $e_i(w)=0$, and if $a>0$, $e_i$  changes to $i$ the letter $i+1$ of $w$ corresponding to the rightmost unbracketed $-$ (i.e., not erased), whereas if $b=0$, $f_i(w)=0$, and if $b>0$, $f_i$  changes to $i+1$ the letter  $i$ corresponding to the leftmost unbracketed $+$.

The  crystal ${B_d}^{\otimes k}$, as a graph, is the union of connected components. The connected components of
${B_d}^{\otimes k}$ are the coplactic classes or dual Knuth classes in the RSK correspondence that identify words with the same recording tableau in $SYT(\lambda)$ for some $\lambda$. For each standard tableau $Q\in SYT(\lambda)$ there is an embedding of $SSYT(\lambda)$ in ${B_d}^{\otimes k}$,
$$\textsf{read}_Q=\textsf{RSK}^{-1}(.,Q):SSYT(\lambda)\longrightarrow {B_d}^{\otimes k}.$$
Furthermore, given $w\in{B_d}^{\otimes k}$ there exists some $\lambda$ and $R\in SYT(\lambda)$ such that $\textsf{read}_S=w$.

The crystal $B_d(\lambda)$ is defined as the crystal structure on
the set $\textsf{SSYT}(\lambda)$ on the alphabet $[d]$ induced by the map $\textsf{read}_Q$  for any $Q\in   \textsf{SYT}(\lambda)$ which does not depend on the choice of $Q$. We have a crystal isomorphism afforded by RSK correspondence $${B_d}^{\otimes k}\longrightarrow \bigsqcup_{|\lambda|=k} B(\lambda)\times \textsf{SYT}(\lambda),\;w\mapsto (P(w),Q)$$
Choose a word $w$ on $[d]$ such that the shape of $P(w)$ is $\lambda$. If we replace every word of its coplactic class with its insertion tableau we obtain the crystal of tableaux ${B}_d(\lambda)$ that has all semistandard tableaux of shape $\lambda$ on the alphabet $[d]$,
\begin{equation}\label{crystaldecomp}{B_d}^{\otimes k}\approx \bigsqcup_{|\lambda|=k} B(\lambda)^{|\textsf{SYT}(\lambda)|}\approx \bigsqcup_{|\lambda|=k} B(\lambda)^{|\mathbb{Y}(\lambda)|},\end{equation}
where $\mathbb{Y}(\lambda)$ denotes the set of Yamanouchi words of weight $\lambda$.
 Each connected component of ${B_d}^{\otimes k}$ has a unique highest weight word which is a Yamanouch word and a unique lowest weight word which is the reversal of that Yamanouchi word. The highest weight element of $B(\lambda)$ is the Yamanouchi tableau $Y(\lambda)$,
and the lowest weight element $Y(\lambda^\bullet)=\rm evac Y(\lambda)$. Two connected components are isomorphic if and only if they have the same highest (lowest) weight \cite{Kash 95}.
	Two words on $[d]$   belong to the same connected component of $\mathcal{W}_d$  if and only if they are dual equivalent. This means that both words are obtained  from the same highest weight word, through a  sequence  of crystal operators $f_i$, or one is obtained from  another by some sequence of crystal operators $f_i$ and $e_j$, $i, j \in [d-1]$.
	Also, two words $w_1,w_2$ on $[d]$ are Knuth equivalent if and only if they occur in the same place in two isomorphic connected components of $\mathcal{W}_d$, that is, they are obtained from two highest words with the same weight through a same sequence of crystal operators. Crystal operators preserve dual Knuth classes and commute with any admissible sequence of  \emph{jeu de taquin} moves.

\subsection{Crystal of a skew-tableau}

For $\mu\subseteq\lambda\subseteq D$, 
let  $B(\lambda/\mu)$ be the set of all semi-standard tableaux of shape $\lambda/\mu$ on the alphabet $[d]$. The column reading of each tableau in $B(\lambda/\mu)$ embeds it in $B^{\otimes |\lambda|-|\mu|}$ the crystal of words of length $|\lambda|-|\mu|$ on the alphabet $[d]$. From \eqref{crystaldecomp} it decomposes

$$B(\lambda/\mu)\cong \bigsqcup_{\begin{smallmatrix}\nu,\;\ell(\nu)\le d\\
T\in \rm LR_{\mu,\nu}^{\lambda}\end{smallmatrix}} B(T)\cong \bigsqcup_{\nu,\;\ell(\nu)\le d} B(\nu)^{c_{\mu,\nu,\lambda}},$$
 where $B(T)$  is the { crystal} connected component of $B(\lambda/\mu)$  containing  the LR tableau $T\in  \rm LR(\mu,\nu,\lambda^\vee)$ for some partition $\nu\subseteq \lambda$. Each crystal $B(T)$ with $T\in \rm LR(\mu,\nu,\lambda^\vee)$ is $\mathfrak{gl}_d$-crystal isomorphic to $B(\nu)$. It  consists of all tableaux with the same shape as $T$ on the alphabet $[d]$ whose normal forms define the set $B(\nu)$, with  highest weight element $Y(\nu)\equiv T$. That is, for each $T\in LR(\mu,\nu,\lambda^\vee)$, $B(T)$ is  a dual Knuth class with highest weight element $T$. Let $T^{\rm low}$ be the lowest weight element of $B(T)$, the unique opposite LR tableau in $B(T)$  Knuth equivalent to $Y(\nu^\bullet)$.

\subsection{The rotated crystal graph and the reversal}
Consider the set of skew-tableaux \linebreak $B(\lambda/\mu)^{\texttt{rotate}}$ of shape $\mu^\vee/\lambda^\vee$ as the image of $B(\lambda/\mu)$ under the map $\texttt{rotate}$, $ U\mapsto \texttt{rotate}(U)$, where $\texttt{rotate}(U)$ is obtained from $U\in B(\lambda/\mu)$ under rotation   by $\pi$ radians while dualizing its word. The map $\texttt{rotate}:B(\lambda/\mu)\rightarrow B(\lambda/\mu)^{\texttt{rotate}}$ is a set bijection preserving the connected components but not a crystal isomorphism. The set $B(\lambda/\mu)^{\texttt{rotate}}=B(\mu^\vee/\lambda^\vee)$ has a crystal structure by flipping upside down each connected component of the crystal $B(\lambda/\mu)$, reverting the arrows and applying the operation $\texttt{rotate}$ to the vertices. If $T$ is the highest weight of a connected component of $B(\lambda/\mu)$ then $\texttt{rotate}(T^{\rm low})\equiv Y(\nu)$ and $\texttt{rotate}(T)\equiv \texttt{rotate}(T)^{ \rm n}=Y(\texttt{rotate}(\nu))$ are respectively the highest and the lowest weights elements of a same connected  component of $B(\lambda/\mu)^{\texttt{rotate}}$.
The crystals $B(\lambda/\mu)\approxeq B(\lambda/\mu)^\bullet$ are isomorphic because they have the same multiset of highest weights   but the isomorphism  is not canonical. Reversal is a set involution on each connected component of $B(\lambda/\mu)$,
$$e:B(\lambda/\mu)\rightarrow B(\lambda/\mu),\; T\mapsto {e}(T)\equiv \texttt{rotate}(T),$$
 is the unique element of  the connected component $B(T)$ of $B(\lambda/\mu)$ containing $T$, Knuth equivalent to $\texttt{rotate}(T)$.
In particular, evacuation is a set involution on $B(\lambda)$
  $${\rm evac}:B(\lambda)\rightarrow B(\lambda),\; T\mapsto {\rm evac}(T)=\texttt{rotate}(T)^{\rm n},$$
that is, is the unique element of $B(\lambda)$ Knuth equivalent to $T^\bullet$. Reversal   interchanges the lowest and highest weight elements in each connected component.

\begin{example}\label{ex:ls} Lusztig-Sch\"utzenberger involution: evacuation $E$ on  crystals of tableaux with normal shape:
$\qquad\scriptsize\textsf{evacuation} :\begin{smallmatrix}2&\\ 1&1\end{smallmatrix}\overset{\scriptsize\textsf{rotate}}\rightarrow\begin{smallmatrix}3&3\\ &2\end{smallmatrix}\overset{\scriptsize\textsf{rectification}}\rightarrow \begin{smallmatrix}2&\\ 3&3\end{smallmatrix}$

\medskip
\begin{tikzpicture}
 \node  at (-0.3,1){{\scriptsize$f_1=$}};\draw[-] [draw=blue, thick] (0,1) -- (0.2,1);
 \node  at (0.7,1){{\scriptsize$f_2=$}};\draw[-] [draw=red, thick] (1,1) -- (1.2,1);
 \end{tikzpicture}

\medskip
$\begin{tikzpicture}
\node at (0,20){\begin{tikzpicture}
  [scale=.85,auto=left,every node/.style={circle,fill=blue!20}]
  \node (n0) at (4,10) {$\begin{smallmatrix}2&\\ 1&1\end{smallmatrix} $};
  \node (n1l) at (2,8.8)  {$\begin{smallmatrix}3&\\ 1&1\end{smallmatrix}$};
  \node (n1r) at (6,8.8)  {$\begin{smallmatrix}{\color{brown} 2}&\\ \bf {\color{brown}1}&{\color{brown}2}\end{smallmatrix}$};
  \node (n2l) at (4.7,7)  {$\begin{smallmatrix}3&{}\\ {1}&2\end{smallmatrix}$};
  \node (n2r) at (3.2,7) {$\begin{smallmatrix}2&{}\\ {1}&3\end{smallmatrix}$};
  \node (n3l) at (1.7,5.9)  {$\begin{smallmatrix}3&\\ {1}&3\end{smallmatrix}$};
  \node (n3r) at (6.5,5.8)  {$\begin{smallmatrix}{\color{brown}3}&{}\\ {\color{brown}2}&{\color{brown}2}\end{smallmatrix}$};
\node (n4l) at (4,4.2)  {$\begin{smallmatrix}3&{}\\ {2}&3\end{smallmatrix}$};

    \foreach \from/\to in {n0/n1r,n1l/n2l,n2l/n3r,n3l/n4l}
    \draw[->] [draw=blue,  thick](\from) -- (\to);
    \foreach \from/\to in {n0/n1r}
     \draw[->] [draw=blue,  thick](\from) -- (\to);
     \foreach \from/\to in {n0/n1l,n1r/n2r,n2r/n3l}
     \draw[->] [draw=red,  thick](\from) -- (\to);

    \foreach \from/\to in {n0/n1l,n1r/n2r,n2r/n3l, n3r/n4l}
 \draw [->] [draw=red,  thick](\from) -- (\to);
\end{tikzpicture}};
\end{tikzpicture}
\begin{tikzpicture}
\node at (0,20){\begin{tikzpicture}
  [scale=.85,auto=left,every node/.style={circle,fill=blue!20}]
  \node (n0) at (4,10) {$\begin{smallmatrix}3&3\\ &2\end{smallmatrix} $};
  \node (n1l) at (2,8.8)  {$\begin{smallmatrix}3&3\\ &1\end{smallmatrix}$};
  \node (n1r) at (6,8.8)  {$\begin{smallmatrix}{\color{brown}2}&{\color{brown}3}\\ &{\color{brown}2}\end{smallmatrix}$};
  \node (n2l) at (4.7,7)  {$\begin{smallmatrix}2&{3}\\ {}&1\end{smallmatrix}$};
  \node (n2r) at (3.2,7) {$\begin{smallmatrix}1&{3}\\ {}&2\end{smallmatrix}$};
  \node (n3l) at (1.7,5.9)  {$\begin{smallmatrix}1&{3}\\ {}&1\end{smallmatrix}$};
  \node (n3r) at (6.5,5.8)  {$\begin{smallmatrix}2&{2}\\ {}&1\end{smallmatrix}$};
\node (n4l) at (4,4.2)  {$\begin{smallmatrix}1&{2}\\ {}&1\end{smallmatrix}$};

    \foreach \from/\to in {n0/n1r,n1l/n2l,n2l/n3r,n3l/n4l}
    \draw[<-] [draw=red,  thick](\from) -- (\to);
    \foreach \from/\to in {n0/n1r}
     \draw[<-] [draw=red,  thick](\from) -- (\to);
     \foreach \from/\to in {n0/n1l,n1r/n2r,n2r/n3l}
     \draw[<-] [draw=blue,  thick](\from) -- (\to);

    \foreach \from/\to in {n0/n1l,n1r/n2r,n2r/n3l, n3r/n4l}
 \draw [<-] [draw=blue,  thick](\from) -- (\to);
\end{tikzpicture}};
\end{tikzpicture}
\begin{tikzpicture}
\node at (0,20){\begin{tikzpicture}
  [scale=.85,auto=left,every node/.style={circle,fill=blue!20}]
  \node (n0) at (4,10) {$\begin{smallmatrix}3&\\ 2&3\end{smallmatrix} $};
  \node (n1l) at (2,8.8)  {$\begin{smallmatrix}3&\\ 1&3\end{smallmatrix}$};
  \node (n1r) at (6,8.8)  {$\begin{smallmatrix}{\color{brown}3}&\\ {\color{brown}2}&{\color{brown}2}\end{smallmatrix}$};
  \node (n2l) at (4.7,7)  {$\begin{smallmatrix}2&\\ {}1&3\end{smallmatrix}$};
  \node (n2r) at (3.2,7) {$\begin{smallmatrix}3&{}\\ {1}&2\end{smallmatrix}$};
  \node (n3l) at (1.7,5.9)  {$\begin{smallmatrix}3&{}\\ {1}&1\end{smallmatrix}$};
  \node (n3r) at (6.5,5.8)  {$\begin{smallmatrix}2&{}\\ {1}&2\end{smallmatrix}$};
\node (n4l) at (4,4.2)  {$\begin{smallmatrix}2&{}\\ {1}&1\end{smallmatrix}$};

    \foreach \from/\to in {n0/n1r,n1l/n2l,n2l/n3r,n3l/n4l}
    \draw[<-] [draw=red,  thick](\from) -- (\to);
    \foreach \from/\to in {n0/n1r}
     \draw[<-] [draw=red,  thick](\from) -- (\to);
     \foreach \from/\to in {n0/n1l,n1r/n2r,n2r/n3l}
     \draw[<-] [draw=blue,  thick](\from) -- (\to);

    \foreach \from/\to in {n0/n1l,n1r/n2r,n2r/n3l, n3r/n4l}
 \draw [<-] [draw=blue,  thick](\from) -- (\to);
\end{tikzpicture}};
\end{tikzpicture}$

$$\scriptsize B((2,1,0),3)\underset{\textsf{rotate}}\rightarrow B((2,2)/(1),3)\underset{\textsf{rectification}}\rightarrow \scriptsize B((2,1,0),3)$$

$$ 
\texttt{evac}( \begin{smallmatrix}{\color{brown} 1}&{\color{brown}2}\\ \bf {\color{brown}2}&\end{smallmatrix})
=\begin{smallmatrix}{\color{brown}2}&{\color{brown}2}\\ {\color{brown}3}&\end{smallmatrix}
$$

\end{example}
\subsection{The action of the symmetric group on a crystal and partial Sch\"utzenberger involutions}\label{groupaction}Recall the signature rule on a crystal of words.
Given  the word $w=w_1\cdots w_k\in {B_d}^{\otimes k}$ with $i$-signature $-^a +^b=x_{j_1}\cdots x_{j_a}\cdots x_{j_{a+b}}$, one defines the \textit{crystal (Kashiwara)  reflection operator} $\sigma_i$ on $w$ by putting $\sigma_i(i^a \;{i+1}^b)=i^b \;{i+1}^a=x'_{j_1}\cdots x'_{j_a}\cdots x'_{j_{a+b}}$ and $\sigma_i(w)=y_1\cdots y_k$ where $y_j=x'_j$ if $j\in \{j_1,\dots,j_{a+b}\}$ and $y_j=x_j$ otherwise. The operators $\sigma_i$, $i\in[d-1]$ are involutions and one defines an action of the symmetric group $\mathfrak{S}_d$ on $\mathcal{W}_d$ by acting on its  connected components isomorphic to $B(\lambda)$, for some $\lambda$. They commute with any meaningful sequence of {\em jeu de taquin} moves.
The subgraph obtained from ${B_d}(\lambda)$ by erasing all edges of colour $\neq i$ is a disjoint union of $i$-strings of various lengths
$$\bullet\overset{i}\dashrightarrow\bullet\overset{i}\dashrightarrow\cdots\bullet\overset{i}\dashrightarrow\bullet$$
The operator $\sigma_i$ is the involution on $B(\lambda)$ which reverses each $i$-string, that is, $\sigma_i(w)$ is the vertex on the $i$ string of $w$ such that $\varepsilon_i(\sigma_i(w))=\varphi_i(w)$. It coincides with the action of the partial Sch\"utzenberger involution on the alphabet $\{i,i+1\}$ on the $i$-string \cite{att}. The $i$-string is itself a crystal graph with highest weight the top and lowest weight the bottom of the $i$-string, thus $\sigma_i$ interchanges the highest with the lowest weight.
Each simple transposition $s_i\in\mathfrak{S}_d$ acts in a $i$-string so that $\rm wt(\sigma_i w)=s_i \rm wt(w)$. Therefore,
the action of the longest Weyl group element (in  type $A$, the reverse permutation) on a connected component of $\mathcal{W}_d$ isomorphic to $B(\lambda)$ agrees with the  action of the Sch\"utzenberger's involution (or evacuation) on Yamanouchi or opposite Yamanouchi  tableaux  and with the reversal on Yamanouchi or opposite Yamanouchi words. The highest weight element is the Yamanouchi tableau $Y(\nu)$,
and the lowest weight element $\sigma_0Y(\nu)=Y(\nu^\bullet)$, with $\sigma_0\in\mathfrak{S}_d$.
Crystal operators and  crystal reflection  operators acting on words (for
definitions see~\cite{Lot,plaxique,kwon,bumschi}) preserve Knuth equivalence and the $Q$--symbol, and, henceforth, also the dual Knuth class, when acting on the word of a tableau.
 Let $w$ be a Yamanouchi word of weight
$\nu$, with
$\ell(\nu)\le d$,
and let $\sigma_i$ denote the
crystal reflection  operator (for definitions see \cite[Section 5.5]{Lot} or  \cite[Definition  2.35]{bumschi})
 acting on the subword
 over the alphabet $\{i,i+1\}$, for  $1\le i<d$. The  operators $\sigma_i$ satisfy the Coxeter relations { of the symmetric group}.
If  $\omega_0:=s_{i_N}\cdots s_{i_1}$ is the
$\mathfrak{S}_d$ long element, put $\sigma_0:=\sigma_{i_N}\cdots\sigma_{i_1}$.
  Then {$\sigma_0 w$ is the opposite Yamanouchi word of weight $\nu^\bullet$}, $\sigma_0 w\equiv w^\bullet\equiv Y(\nu^\bullet)$, and $Q(\sigma_0w)=Q(w)=U(w)$ (whereas $Q(w^\bullet)=\rm evac Q(w)$).

 For $\mu\subseteq\lambda\subseteq D$, 
let  $B(\lambda/\mu)$ be the set of all semi-standard tableaux of shape $\lambda/\mu$ on the alphabet $[d]$.

$$B(\lambda/\mu)\simeq \bigoplus_{\begin{smallmatrix}\nu \subseteq \lambda\\
T\in \rm LR_{\mu,\nu}^{\lambda}\end{smallmatrix}} B(T),$$
 where $B(T)$  is the {crystal} connected component of $B(\lambda/\mu)$  containing  the LR tableau $T\in  \rm LR(\mu,\nu,\lambda^\vee)$ for some partition $\nu\subseteq \lambda$. Each crystal $B(T)$ with $T\in \rm LR(\mu,\nu,\lambda^\vee)$ is $\mathfrak{gl}_d$-crystal isomorphic to $B(\nu)$.
  That is, $B(T)$ is  a dual Knuth class. Since $\sigma_0T$ is dual Knuth equivalent to $T$, the  lowest weight element of $B(T)$ is the reversal
  LR tableau $eT=\sigma_0T\equiv Y(\nu^\bullet)$ in $\rm LR(\mu,\nu^\bullet,\lambda)$ and $T^{e\,e}=\sigma_0\sigma_0 T=T$.

 For LR tableaux we may provide another procedure  used in ~\cite{az} to calculate the reversal. This procedure is illuminated by  the action of the longest permutation of the symmetric group $\mathfrak{S}_d$ on the crystal $B_d(T)$ of a skew tableau $T$ of shape $\lambda/\mu$ and content $\nu$
 ~\cite{kwon, bumschi} over the alphabet $\{1,\dots,d\}$, with $\ell(\nu)\le d$.

\subsection{ Left and right LR companion tableaux,  crystals and hives}
\label{subsec:lr}


The \emph{recording matrix $M$ of an LR tableau} $T\in  \rm LR(\mu,\nu,\lambda^\vee)$ is the
$d\times d$ lower triangular matrix  identified with the LHS of \eqref{biword}.  Its transposition $M^t$ is the upper triangular matrix identified with the RHS of  \eqref{biword}. $M^t$  is the recording  matrix of the semistandard Young tableau $G$ of shape $\nu$ and weight $\lambda/\mu$, the recording tableau of $\textsf{RSK}(M)=(Y(\nu),G)$ in Proposition \ref{prop:burge}.  The  semistandard Young tableau $G$ of shape $\nu$ and content $\lambda/\mu$ is such that  each row $i$  tell us in which rows of $T$ the $i$'s are filled in, Proposition \ref{prop:burge}, $(c)$, and is
 called the \emph{right Gelfand-Tsetlin} (GT) pattern, or the \emph{right LR  companion tableau} of $T$. It satisfies $\widehat G=U(w(T))$.

A GT pattern, $G$, is a triangular array of
non-negative integers $G=(\nu_j^{(i)})_{1\leq j\leq i\leq d}$   displayed
as below (for more details and references therein we refer to \cite{akt,tka}):
\begin{equation}
\begin{array}{cccccccccccccccccc}
\nu_1^{(d)}&&\nu_2^{(d)}&&\cdots&&\nu_{d-1}^{(d)}&&\nu_n^{(d)}&\cr
&\nu_1^{(d-1)}&&\nu_2^{(d-1)}&&\cdots&&\nu_{d-1}^{(d-1)}\cr
&&\cdots&&\cdots&&\cdots\cr
&&&\nu_1^{(2)}&&\nu_2^{(2)}\cr
&&&&\nu_1^{(1)}\cr
\end{array}
\end{equation}
where entries satisfy the betweenness conditions $\nu_j^{(i+1)}\geq \nu_j^{(i)}\geq \nu_{j+1}^{(i+1)}$
for $1\leq j\leq i<d$. The $i$th row, enumerated from bottom to top, necessarily constitutes a partition $\nu^{(i)}$
of length $\leq i$. Such a GT pattern is said to be of \emph{shape} (or \emph{type})  $\nu^{(d)}$ and of weight $\gamma=(\gamma_1,\gamma_2,\ldots,\gamma_n)$
where $\gamma_i=|\nu^{(i)}|-|\nu^{(i-1)}|$ for $i=1,2,\ldots,d$, with $|\nu^{(0)}|=0$. There is a natural bijection between GT patterns  and semistandard tableaux of the same shape and weight: a semistandard tableau of shape $\nu$ is a nested sequence of partitions $\nu^{(1)}\subseteq \nu^{(2)}\subseteq\cdots\subseteq \nu^{(d)}$ where $\nu^{(i)}$ specify the shape of that part of the semistandard tableau consisting of entries $\leq i$,
that is to say having $\nu_j^{(i)}-\nu_j^{(i-1)}$ entries $i$ in row $j$ for $1\leq j\leq i\leq d$.

 \begin{definition}\label{companiondef} Denote by $\rm LR_{\nu,\lambda/\mu}$ the set of right LR companion tableaux of $\rm LR(\mu,\nu,\lambda^\vee)$. The elements of $\rm LR_{\nu,\lambda/\mu}$ are called  right LR companion tableaux of shape $\nu$ and content $\lambda/\mu$.
 \end{definition}

 Thanks to Proposition \ref{prop:burge}, the  linear transformation $M\mapsto                                                                                                                                                                                                                                                  M^t$, where $M$ is the recording matrix of an LR tableau, defines  a linear  bijection  $\iota$  between LR tableaux of boundary data $(\mu,\nu,\lambda)$ and their right companions. Let \begin{equation}\label{rightcompanionmap}\iota: {\rm LR(\mu,\nu,\lambda^\vee)}\rightarrow {\rm LR_{\nu,\lambda/\mu}},\;T\mapsto \iota(T)=G,\end{equation} where $G$ is the semistandard tableau of shape $\nu$ and content $\lambda/\mu$ with recording matrix $M^t$ given that $T$ has recording matrix $M$.

  The \emph{left LR companion tableau} (or \emph{left Gelfand-Tsetlin pattern})  of $T\in  \rm LR(\mu,\nu,\lambda^\vee)$ is
 defined  to be the semistandard tableau $L$ of shape $\mu$ and content $\rm rev(\lambda/\nu)$ which records the sequence of partitions
giving the shapes occupied by the entries $<r$ in rows $r,r+1,\dots,d$ of $T$, for $r=1,2,\ldots,d$ (for more details we refer to \cite{akt,tka} and references therein).
Given the recording matrix $M=$ $(M_{ij})_{1\le i,j\le d}$ of $T$, the semistandard tableau $L$ is given by the nested sequence of partitions
\begin{eqnarray}\label{left}(\mu_d+\sum_{j=1}^{d-1}m_{dj})\subseteq\cdots\subseteq(\mu_r+\sum_{j=1}^{r-1}m_{rj},\dots,\mu_{d-1}+\sum_{j=1}^{r-1}m_{d-1,j},
\mu_{d}+\sum_{j=1}^{r-1}m_{d,j})\subseteq \cdots\subseteq\mu,\nonumber\\\qquad\end{eqnarray}
where $\mu=(\mu_1,\dots,\mu_d)$ and $1\le r\le d$.
This is equivalent to the linear transformation given in  \cite[Proposition 12]{PV2} to obtain the recording matrix  of the LR left companion $L$ of $T$ given its recording matrix.
\begin{definition}
Denote by $\rm LR^-_{\nu,\lambda/\mu}$ the \emph{set of left LR companion tableaux} of $\rm LR(\mu,\nu,\lambda^\vee)$. \end{definition}
Clearly the sets $\rm LR^-_{\mu,\lambda/\nu}$, $\rm LR_{\nu,\lambda/\mu}$ and $\rm LR(\mu,\nu,\lambda^\vee)$ are mutually in linear bijection.

  Given a semistandard tableau $G$ of shape $\nu$ and content $\gamma=\lambda-\mu$  the transpose of its recording matrix determines
 an LR skew tableau only up to a parallel shift of the skew shape. Given a semistandard tableau $G$ of shape $\nu$ and content $\lambda/\mu$ how do we check whether $G$ is in $\rm LR_{\nu,\lambda/\mu}$? The same question for $L$ of shape $\mu$ and content ${\rm rev}(\lambda-\nu)$ to be in $\rm LR^-_{\mu,\lambda/\nu}$.

  The   elements of $\rm LR_{\nu,\lambda/\mu}$  are those tableaux (GT patterns) $G=$  $\nu^1\supseteq\nu^2\supseteq\cdots\supseteq\nu^d=(\nu_1,\dots,\nu_d)$ and content $\lambda/\mu$   \cite[Theorem 1]{GZpolyed} (we refer to \cite[Section 2.4]{akt} for more details) such that
  \begin{equation}\label{gtright}
  \sum_{k=1}^i(\nu_k^{(j)}-\nu_k^{(j-1)})-\sum_{k=1}^{i-1}(\nu_{k}^{j-1}-\nu_k^{j-2})\le \mu_{j-1}-\mu_j,\,1\le i<j\le d.
  \end{equation}

   Alternatively, let $B(\mu)$ and $B(\nu)$ be the crystals of all semistandard tableaux of shape $\mu$  and $\nu$, on the alphabet $[d]$, respectively. Let us consider $B(\mu)\otimes B(\nu)$ and their highest weight elements. Given $G\in B(\nu)$ of weight $\lambda/\mu$,    $Y(\mu)\otimes G$ is the highest weight element of weight $\lambda$ of  a connected component, isomorphic to $B(\lambda)$, of $B(\mu)\otimes B(\nu)$ if and only if $\varepsilon_j(G)\le \mu_{j-1}-\mu_j$, for all $1< j\le d$, where $\varepsilon_j(G)=\max\{k\in\mathbb{Z}_{\ge 0}:e_i^k(G)\neq 0\}$ and $e_i$ is a raising operator \cite[Appendix]{Nak05}, \cite{kwon}. The only if part of the former statement is equivalent to
   \eqref{gtright} holds.  Therefore $G\in \rm LR_{\nu,\lambda/\mu}$ if and only if $Y(\mu)\otimes G$ is a highest weight element of $B(\lambda$  if and only if $G$ is  a vertex of $B(\nu)$ of weight $\lambda/\mu$ such that $\varepsilon_j(G)\le \mu_{j-1}-\mu_j$, for all $1< j\le d$ if and only if  \eqref{gtright} holds.
   Equivalently, $\rm LR_{\nu,\lambda/\mu}$ is the set of semistandard tableaux $G$ of shape $\nu$ and content $\lambda/\mu$  satisfying the equation $Y(\mu).G=Y(\lambda)$ where "$.$'' refers to  the column insertion of $G$ in $Y(\mu)$ \cite{thomas1,kwon}.

 The   elements of $\rm LR^-_{\mu,\lambda/\nu}$  are those tableaux (GT patterns) $L=$  $\mu^1\supseteq\mu^2\supseteq\cdots\supseteq\mu^d=(\mu_1,\dots,\mu_d)$ and
  content $\rm rev(\lambda/\nu)$  such that  \cite{BZphy}
   (see also \cite[Section 2.4]{akt} for details)
   \begin{equation}\label{gtleft}
  \sum_{k=j}^d(\mu_{k-i}^{(d-i)}-\mu_{k-i+1}^{(d-i+1)})-\sum_{k=j+1}^{d}(\mu_{k-i-1}^{(d-i-1)}-\mu_{k-i}^{(d-i)})\le \nu_{i}-\nu_{i+1},\,1\le i<j\le  d.
  \end{equation}

   Similarly, under considerations of the lowest weight elements of $B(\mu)\otimes B(\nu)$, given $L\in B(\mu)$ of weight $\rm rev(\lambda/\nu)$,
    $L\otimes Y({\rm rev}\nu)$ is the lowest weight element of weight $\rm rev\lambda$ of a connected component, isomorphic to $B(\lambda)$, of  $B(\mu)\otimes B(\nu)$ if and only if $\varphi_{d-i}(L)\le \nu_i-\nu_{i+1}$, for $1\le i<d$,
    where $\varphi_j(L)=\max\{k\in\mathbb{Z}_{\ge 0}:f_j^k(L)\neq 0\}$ and $f_j$ is a lowering operator  \cite{kwon}. The only if part of the former statement is equivalent to \eqref{gtleft}. Therefore $L\in \rm LR^-_{\nu,\lambda/\mu}$ if and only if $L\otimes Y({\rm rev}\nu)$ is a lowest weight element of $B(\lambda$  if and only if $L$ is  a vertex of $B(\mu)$ of weight $\rm rev\lambda/\nu$ such that $\varphi_{d-i}(L)\le \nu_i-\nu_{i+1}$, for all $1< i< d$ if and only if  \eqref{gtleft} holds.
    Equivalently,  $L$ of shape $\mu$ and content $\rm rev(\lambda/\nu)$  satisfies the equation $L\otimes Y({\rm rev}\nu)=Y(\rm rev\lambda)$ where "$.$'' refers to  the column insertion of $Y(\rm rev \nu)$ in $L$ \cite{kwon}.

  \subsubsection{LR companion pairs and hives} Thanks to \cite{hk}, a pair $(L_\mu,G_\nu)$ of semistandard tableaux of shapes $\mu$ and $\nu$ and weights $rev(\lambda/\nu)$ and $\lambda/\nu$ respectively, is said to be a \emph{LR companion pair} of $ \rm LR(\mu,\nu,\lambda^\vee)$ if and only if $L_\mu\otimes Y({\rm rev}\nu)$ and $Y(\mu)\otimes G_\nu$ are the lowest respectively the highest weight elements of a connected component, isomorphic to $B(\lambda)$, of $B(\mu)\otimes B(\nu)$ (see \cite[subsection 12.1]{akt} for more details). That is   $Y(\mu)\otimes G$ and $L\otimes Y(\rm rev \nu)$ have the same recording tableau $Q$ (an LR tableau of shape $\lambda/\mu$ and content $\nu$) \cite{thomas1,kwon}. Alternatively, $(L,G)$ satisfy  certain linear equalities which can be expressed by the triangle condition in a hive \cite[Section 2.4]{akt}.

   In \cite[Proposition 12]{PV2}, for  $ T\in \rm LR(\mu,\nu,\lambda^\vee)$, it is given  a linear  map \cite[Proposition 12]{PV2}  to transform the recording matrix of $T$ into the left companion of $T$. Thereby given a right companion  the corresponding left companion can be obtained from that by a linear transformation.
   Later in Section \ref{sec:hives} we give another relation using the action of $\mathfrak{S}_3\times\mathbb{Z}_2$ on $\mathcal{LR}$.
\begin{example}\label{ex:companion}
Consider $T\in \rm LR(210,532,320)$ in Example \ref{yamanouchi}. The recording matrix of
$T$ is $M=\left(\begin{smallmatrix}
4&0&0\\
1&2&0\\
0&1&2
\end{smallmatrix}\right)$. Its transposition  $M^t=\left(\begin{smallmatrix}
4&1&0\\
0&2&1\\
0&0&2
\end{smallmatrix}\right)$ encodes the right companion tableau
 \begin{equation}\label{G} G=\begin{Young} 3&3&&&&\cr 2&2&3&&&\cr 1&1&1&1&2&\cr
\end{Young}\end{equation} of normal shape $\nu=532$ and
weight $\lambda/\mu=643-210=433$. The standardization of $G$ gives $\widehat G=U(w(T))$ in \eqref{U}.
The left companion of $T$ is given by the sequence of shapes $(1)\subseteq (20)\subseteq (210)$, equivalently,
\begin{equation}\label{L} L=\begin{Young} &&&&&\cr 3&&&&&\cr 1&2&&&&\cr
\end{Young}\end{equation} of normal shape $\mu=210$ and
weight $\rm rev\lambda/\nu=\rm rev(643-532)=111$.
\end{example}


\subsubsection{Companion tableau of an opposite LR tableau}
The \emph{recording matrix $M$} of an \emph{opposite LR tableau} $T\in  \rm LR(\mu,\nu^\bullet,\lambda^\vee)$ is a
$d\times d$ upper triangular matrix $M$. Its transposition $M^t$, a lower triangular matrix, is the recording matrix of the antitableau $H$ of shape $\nu^\bullet$ and weight $\lambda/\mu$ such that the row $i$ entries of $H$ tell us which rows of $T$ are filled with $i$'s. The antitableau $H$ is the companion of $T$. In Example \ref{oppositeyamanouchi} the recording matrix of $S$ is
$M(S)=\left(\begin{smallmatrix}2&0&2\cr
0&3&0\cr
0&0&3\cr
\end{smallmatrix}\right)$ and $M(S)^t=\left(\begin{smallmatrix}2&0&0\cr
0&3&0\cr
2&0&3\cr
\end{smallmatrix}\right)$ with  companion tableau $H=\begin{Young}1&1&3&3&3\cr
&&2&2&2\cr
&&&1&1\cr
\end{Young}
$

\section{ The linear rotation and orthogonal transpose maps on LR tableaux and companions}
\label{s:maps}
\subsection{Linear reduction and linear equivalence of bijections}
\label{S:linearcost}
We follow closely~\cite{PV2}
for this section. Using ideas and techniques of Theoretical Computer
Science, see~\cite{AHU,clrs}, each bijection can be seen as an
algorithm having one type of combinatorial objects as \emph{input},
and another as \emph{output}. We define a \emph{correspondence} as
a one--to--one map established by a bijection; therefore, obviously
several different defined bijections can produce the same
correspondence. In this way one can think of a correspondence as a
function which is computed by the algorithm, viz. the bijection. The
computational complexity is, roughly, the number of steps in the
bijection. Two bijections are \emph{identical} if and only if they
define the same correspondence. Obviously one task can be performed
by several different algorithms, each one having its own
computational complexity, see~\cite{AHU,clrs}. For example we recall that
there are several ways to multiply large integers, from naive
algorithms, e.g. the Russian peasant algorithm, to that ones using
FFT (Fast Fourier Transform), e.g. Sch\"{o}nhage--Strassen
algorithm; see e.g.~\cite{gg} for a comprehensive and update
reference. Formally, a function $f$ reduces linearly to $g$, if it
is possible to compute $f$ in time linear in the time it takes to
compute $g$; $f$ and $g$ are linearly equivalent if $f$ reduces
linearly to $g$ and vice versa. This defines an equivalence relation
on functions, which can be translated into a linear equivalence on
bijections.

Let $D = (d_{1}, \ldots , d_{n})$ be an array of integers,
and let $m = m(D) := \max_{i} d_{i}$. The \emph{bit--size} of $D$,
denoted by $\langle D \rangle$, is the
amount of space required to store $D$; for simplicity from now on we
assume that $\langle D \rangle = n \rf{\log_{2} m + 1}$. We view a
bijection $\delta : \cA \longrightarrow \cB$ as an algorithm which
inputs $A \in \cA$ and outputs $B = \delta\(A\) \in \cB$. We need to
present Young tableaux as arrays of integers so that we can store
them and compute their bit--size. Suppose $A \in YT(\mu,
m,\lambda)$: a way to encode $A$ is through its recording matrix
$\(c_{i,j}\)$, which is defined by $c_{i,j} =a_{i,j} - a_{i,j-1}$;
in other words, $c_{i,j}$ is the number of $j$'s in the $i$--th row
of $A$; this is the way Young tableaux will be presented in the
input and output of the algorithms. Finally, we say that a map
$\gamma : \cA \longrightarrow \cB$ is \emph{size--neutral} if the
ratio $\frac{\langle \gamma(A)\rangle}{\langle A \rangle}$ is
bounded for all $A \in \cA$. Throughout the paper we consider only
size--neutral maps, so we can investigate the linear equivalence of
maps comparing them by the number of times other maps are used,
without be bothered by the timing. In fact, if we drop the condition
of being size--neutral, it can happen that a map increases the
bit--size of combinatorial objects, when it transforms the input
into the output, and this affects the timing of its subsequent
applications. Let $\cA$ and $\cB$ be two possibly infinite sets of
finite integer arrays, and let $\delta : \cA \longrightarrow \cB$ be
an explicit map between them. We say that $\delta$ has {\em linear cost}
if $\delta$ computes $\delta\(A\) \in \cB$ in linear time
$O\(\langle A \rangle \)$ for all $A \in \cA$. There are many ways
to construct new bijections out of existing ones: we call such
algorithms \emph{circuits} and we define below several of them that
we need.

\begin{description}

\item{$\circ$}
Suppose $\delta_{1} : \cA_{1} \longrightarrow \cX_{1}$, $\gamma :
\cX_{1} \longrightarrow \cX_{2}$ and $\delta_{2} : \cX_{2}
\longrightarrow \cB$, such that $\delta_{1}$ and $\delta_{2}$ have
linear cost, and consider $\chi=\delta_{2} \circ \gamma \circ
\delta_{1} : \cA \longrightarrow \cB$. We call this circuit
\emph{trivial} and denote it by $I\(\delta_{1}, \gamma,
\delta_{2}\)$.

\item{$\circ$}
Suppose $\gamma_{1} : \cA \longrightarrow \cX$ and $\gamma_{2} : \cX
\longrightarrow \cB$, and let $\chi=\gamma_{2} \circ \gamma_{1} :
\cA \longrightarrow \cB$. We call this circuit \emph{sequential} and
denote it by $S\(\gamma_{1}, \gamma_{2}\)$.

\item{$\circ$}
Suppose $\delta_{1} : \cA \longrightarrow \cX_{1} \times \cX_{2}$,
$\gamma_{1} : \cX_{1} \longrightarrow \cY_{1}$, $\gamma_{2} :
\cX_{2} \longrightarrow \cY_{2}$, and $\delta_{2} : \cY_{1} \times
\cY_{2} \longrightarrow \cB$, such that $\delta_{1}$ and
$\delta_{1}$ have linear cost. Consider $\chi=\delta_{2} \circ
\(\gamma_{1} \times \gamma_{2}\) \circ \delta_{1} : \cA
\longrightarrow \cB$: we call this circuit \emph{parallel} and
denote it by $P\(\delta_{1},\gamma_{1},\gamma_{2},\delta_{2}\)$.

\end{description}

For a fixed bijection $\alpha$, we say that $\beth$ is an
$\alpha$\emph{--based ps--circuit} if one of the following holds:

\begin{description}

\item{$\bullet$}
$\beth=\delta$, where $\delta$ is a bijection having linear cost.

\item{$\bullet$}
$\beth=I\(\delta_{1}, \alpha, \delta_{2}\)$, where
$\delta_{1},\delta_{2}$ are bijections having linear cost.

\item{$\bullet$}
$\beth=P\(\delta_{1},\gamma_{1},\gamma_{2},\delta_{2}\)$, where
$\gamma_{1},\gamma_{2}$ are $\alpha$--based ps--circuits and
$\delta_{1},\delta_{2}$ are bijections having linear cost.

\item{$\bullet$}
$\beth=S\(\gamma_{1}, \gamma_{2}\)$, where $\gamma_{1}, \gamma_{2}$
are $\alpha$--based ps--circuits.

\end{description}

In other words, $\beth$ is an $\alpha$--based ps--circuit if there
is a parallel--sequential algorithm which uses only a finite number
of linear cost maps and a finite number of application of map
$\alpha$. The $\alpha$--cost of $\beth$ is the number of times the
map $\alpha$ is used; we denote it by $s\(\beth\)$.

Let $\gamma : \cA \longrightarrow \cB$ be a map produced by the
$\alpha$--based ps--circuit $\beth$. We say that $\beth$ computes
$\gamma$ at cost $s\(\beth\)$ of $\alpha$. A map $\beta$ is
\emph{linearly reducible} to $\alpha$, write $\beta \hookrightarrow
\alpha$, if there exist a finite $\alpha$--based ps--circuit $\beth$
which computes $\beta$. In this case we say that $\beta$ can be
computed in at most $s\(\beth\)$ cost of $\alpha$. We say that maps
$\alpha$ and $\beta$ are linearly equivalent, write $\alpha \sim
\beta$, if $\alpha$ is linearly reducible to $\beta$, and $\beta$ is
linearly reducible to $\alpha$. We recall, gluing together, results
proved in Section 4.2 of~\cite{PV2}.

\begin{proposition}
Suppose $\alpha_{1} \hookrightarrow \alpha_{2}$ and $\alpha_{2}
\hookrightarrow \alpha_{3}$, then $\alpha_{1} \hookrightarrow
\alpha_{3}$. Moreover, if $\alpha_{1}$ can be computed in at most
$s_{1}$ cost of $\alpha_{2}$, and $\alpha_{2}$ can be computed in at
most $s_{2}$ cost of $\alpha_{3}$, then $\alpha_{1}$ can be computed
in at most $s_{1}s_{2}$ cost of $\alpha_{3}$. Suppose $\alpha_{1}
\sim \alpha_{2}$ and $\alpha_{2} \sim \alpha_{3}$, then $\alpha_{1}
\sim \alpha_{3}$ Suppose $\alpha_{1} \hookrightarrow \alpha_{2}
\hookrightarrow \ldots \hookrightarrow \alpha_{n} \hookrightarrow
\alpha_{1}$, then $\alpha_{1} \sim \alpha_{2} \sim \ldots \sim
\alpha_{n} \sim \alpha_{1}$.
\end{proposition}

\subsection{The linear involution rotation on LR    and  companion tableaux}\label{comprotated}

 Given a Yamanouchi word $w=w_1\cdots w_s$ of weight $\nu$,
define the standard tableau
 $U(w^\bullet)$ of  anti-normal shape $\nu^\bullet$ such that the label $k$ is in
row $i$ if and only if $w^\bullet_{s-i+1}=k$ where $w^\bullet=w_s^\bullet\cdots w_1^\bullet$. Thus $U(w^\bullet)=U(w)^\bullet$

 and  $w\mapsto U(w)^\bullet$ defines a bijection  between Yamanouchi words and standard tableaux  of  anti-normal shape given by the reverse content of the  Yamanouchi word.

If $T\in \rm LR(\mu,\nu,\lambda)$ and  $M$ is the recording matrix of $T$, $M^\bullet:=P_{rev}MP_{rev}$, the $\pi$ radians rotation of $M$, is the recording matrix of $\bullet T\in LR(\lambda,\nu^\bullet,\mu)$. Since $M^{\bullet\, t}=M^{t\bullet}$,  if $G$ is the companion of  $T$ then $M^{\bullet\, t}$ is the recording matrix of $\bullet G$ which we define to be   the companion of $\bullet T$.
Therefore  $U(w(T)^\bullet)={U(w(T))\,}^\bullet$ $=\bullet {\widehat G\,}=\widehat {\bullet G}$.
{  Let $\rm LR_{\nu^\bullet,(\lambda/\mu)^\bullet}$ be the set of companion tableaux of the opposite LR tableaux of shape $\lambda/\mu$ and content $\nu^\bullet$. Then $\rm LR_{\nu^\bullet,(\lambda/\mu)^\bullet}:=\bullet \rm LR_{\nu,\lambda/\mu}$.}
The map $\bullet$ is an involution on the set $\rm LR(\mu,\nu,\lambda)\sqcup \rm LR(\lambda,\nu^\bullet,\mu)$ and on the set of  LR companion tableaux
 $\rm LR_{\nu,\lambda/\mu}\sqcup \rm LR_{\nu^\bullet,(\lambda/\mu)^\bullet}$.

\begin{proposition} If $G$ is the right companion tableau of $T$, the {\em companion tableau} of $\bullet T$ is the semistandard tableau $\bullet G$ of anti-normal shape $\nu^\bullet$, the rotated of $G$,  whose  $i$-th row  tell us in which rows of $\bullet T$ the $i$'s are filled in.
\end{proposition}

\begin{example} \label{gt} Let $T\in LR(\mu,\nu,\lambda^\vee)$, in  Example \ref{yamanouchi}, with recording matrix
$M=\left(\begin{smallmatrix}
4&0&0\\
1&2&0\\
0&1&2
\end{smallmatrix}\right)$, and companion tableau
$\label{G1} G=\begin{Young} 3&3&&&&\cr 2&2&3&&&\cr 1&1&1&1&2&\cr
\end{Young}$ of shape $\nu=532$ and weight $\lambda/\mu=433$ encoded by $M^t$.
Then $\bullet T=\begin{Young} 3&3&3&3&&\cr &&2&2&3&\cr &&&1&1&2\cr
\end{Young}$ is in $LR(\lambda^\vee,\nu^\bullet,\mu)$ with dual Yamanouchi  word $w^\bullet=2113223333$, and has   recording matrix  $M^\bullet:=P_{rev}MP_{rev}=\left(\begin{smallmatrix}
2&1&0\\
0&2&1\\
0&0&4
\end{smallmatrix}\right)$, with  $P_{rev}$ the $rev$ permutation matrix. The matrix ${M^\bullet}^t=M^{t\bullet}=\left(\begin{smallmatrix}
2&0&0\\
1&2&0\\
0&1&4
\end{smallmatrix}\right)$ encodes $\bullet G=\begin{Young} &2&3&3&3&3\cr &&&1&2&2\cr &&&&1&1\cr
\end{Young}$ of anti normal shape $\nu^\bullet=235$ and weight $(\lambda/\mu)^\bullet=\mu^\vee/\lambda^\vee=334=(\lambda-\mu)_{rev}$. The antitableau  $\bullet G$  gives explicit information  of $\bullet T$ namely the row $i$ entries of $\bullet G$ tell us which rows of $\bullet T$ are filled with $i$'s. One  has

$\widehat {\bullet G}=\bullet{\widehat G\,}={U(w)\,}^\bullet=U(w^\bullet)=\begin{Young} &4&7&8&9&10\cr &&&1&5&6\cr &&&&2&3\cr
\end{Young}$.
\end{example}

\subsection{The 
linear cost orthogonal transpose involution $\blacklozenge$}\label{blacklozenge}
There is another simple bijection, denoted  $\blacklozenge$,
between LR tableaux  of conjugate weights and  rotate conjugate shapes \cite{az,az1,za}.
Given a \emph{Yamanouchi} word $w$ on the alphabet $[d]$ of weight
$\nu=(\nu_1,\dots,\nu_d)$, let $\nu^t=(\nu_1^t,\dots,\nu_{\nu_1}^t)$.
We define $w^\blacklozenge$ to be the
Yamanouchi word, on the alphabet $[\nu_1]$, of weight $\nu^t$ obtained by replacing in $w$ the subword of length $\nu_i$,
 consisting  of all letters $i$,
     with the subword $12\cdots \nu_i$, for each $i=1,\dots,d$.   When $w$ is \emph{opposite Yamanouchi} word, $w^\blacklozenge$ is defined to be the
\emph{opposite Yamanouchi} word, on the alphabet $[\nu_1]$, of weight $\nu^{t\bullet}$ obtained by replacing in $w$ for each
     $i$, for each $i=1,\dots,d$,  the subword of length $\nu_i$,
 consisting  of all letters $d-i+1$
     with the subword $ (\nu_1-\nu_i+1)\cdots(\nu_1-1)\,\nu_1$.

     {If $w$ is a \emph{Yamanouchi}  word, the word $w^{\blacklozenge \bullet}$ is calculated by first replacing in $w$ each string $i^{\nu_i}$ with $1 2\cdots \nu_i$ in $w$, for $i=1,\dots,d$, to obtain $w^\blacklozenge$. Then, after reversing $w^\blacklozenge$, for $i=1,\dots,d$, our string $1 2\cdots \nu_i$ is transformed into $\nu_i\cdots 21$ which we replace with $(\nu_1-\nu_i+1)\cdots(\nu_1-1)\,\nu_1$.
     On the other hand, the word $w^{\bullet\blacklozenge }$ is calculated by first replacing in $w$ each string $i^{\nu_i}$ with $(d-i+1)^{\nu_i}$, for $i=1,\dots,d$, and then reversing the word to obtain $w^\bullet$. Then, for $i=1,\dots,d$, our string $(d-i+1)^{\nu_i}$ in $w^\bullet$ is transformed into $(\nu_1-\nu_i+1)\cdots(\nu_1-1)\,\nu_1$. Thereby, $w^{\bullet\blacklozenge}=w^{\blacklozenge\bullet}$ is an \emph{opposite Yamanouchi} word of weight $\nu^{t\bullet}=\nu^{\bullet t}$  \eqref{partoperation}.}
     Henceforth, if $w$ is \emph{Yamanouchi}, the  word $w^{\blacklozenge\bullet}$ can be obtained in just  one single
     step by replacing in $w$, for each
     $i$,  the subword of length $\nu_i$,
 consisting  of all letters $i$
     with the subword $\nu_1\,(\nu_1-1)\cdots (\nu_1-\nu_i+1)$,  and then reversing the resulting word. See Example \ref{ex:bullblack}.
     If $w$ is \emph{ opposite Yamanouchi } word, we also have $w^{\bullet\blacklozenge}=w^{\blacklozenge\bullet}$ by reducing to the previous case because every \emph{opposite Yamanouchi} word is the dual of some \emph{Yamanouchi} word and $\bullet$ is an involution. The Yamanouchi word $w^{\blacklozenge\bullet}$  is obtained by replacing each string $i^{\nu_i}$ with $\nu_i\cdots 21$ and then reverse the resulting word.

     The map  $w\mapsto w^\blacklozenge$ defines a bijection between  (opposite) Yamanouchi words  of  conjugate content.
Clearly,  $U(w^\blacklozenge)=U(w)^t$  has shape $\nu^t$, and
      $U({w^\bullet}^\blacklozenge)=U(w^{\blacklozenge\bullet})=U(w)^{\bullet t}=U(w)^{t\bullet }$ has shape $\nu^\bullet$.

The operation $\blacklozenge$ on Yamanouchi or opposite Yamanouchi words is now extended to LR or opposite LR tableaux in the sense that $\blacklozenge$ can be seen as defined on diagonally-shaped LR or opposite LR tableaux.
The orthogonal transpose  map $\blacklozenge$ is defined on LR and opposite LR tableaux as follows. Given ${\rm T}\in {\rm LR}(\mu,\nu,\lambda)$ (respectively ${\rm
LR}(\mu,\nu^\bullet,\lambda)$) with (opposite) Yamanouchi word $w$ , the { {\em orthogonal  transpose}} of $T$, $\blacklozenge{\rm
T}$, is   the (opposite) LR tableau of shape
$(\lambda/\mu)^{t \bullet}=(\mu^\vee)^t/(\lambda^\vee)^t$ and (opposite) Yamanouchi column word $w^{\blacklozenge}$  of weight $\nu^t$ ($\nu^{\bullet t}$). It is obtained from $\rm
T$ by replacing the word $w$  with $w^{\blacklozenge}$,  and then transpose and rotate
the shape $\lambda/\mu$ by $\pi$ radians,

\begin{equation}\label{orthogonal}\begin{array}{ccccc}
{\scriptsize{\blacklozenge}}:{\rm LR}(\mu,\nu,\lambda) \cup  {\rm
LR}(\mu,\nu^\bullet,\lambda))&\longrightarrow& {\rm
LR}(\lambda^{t},\nu^t,\mu^{t})\cup{\rm
LR}(\lambda^{t},\nu^{\bullet t},\mu^{t}),\\
T&\mapsto&  T^\blacklozenge,\;\,w_{col}(T^\blacklozenge)=w(T)^\blacklozenge.
\end{array}\end{equation} The map $\blacklozenge$ is an involution  on  ${\rm LR}(\mu,\nu,\lambda) \cup   {\rm
LR}(\lambda^{t},\nu^t,\mu^{t})$ and on ${\rm
LR}(\mu,\nu^\bullet,\lambda)\cup {\rm
LR}(\lambda^{t},\nu^{\bullet t},\mu^{t})$ that transposes  inner shape, outer shape and weight, and simultaneously swaps the inner and the outer shapes.  In  subsection \ref{S:comp}, one proves that $\blacklozenge$ is a linear time involution. It exhibits the symmetries $c_{\mu,\nu,\lambda}=c_ {\lambda^t,\nu^t,\mu^t}$, $c_{\mu,\nu^\bullet,\lambda}=c_ {\lambda^t,\nu^{\bullet t},\mu^t}$ in linear time.

Let $T\in {\rm
LR}(\mu,\nu,\lambda)\cup{\rm
LR}(\mu^{t},\nu^{\bullet t},\lambda^t)$. Because
$\bullet\blacklozenge T$ and ${\blacklozenge\bullet}T$ have  the same   column word $w(T)^{\blacklozenge \bullet}=w(T)^{\bullet \blacklozenge }$, it follows that
\begin{equation}
\begin{array}{ccccc}
\bullet\blacklozenge=\blacklozenge\bullet:{\rm
LR}(\mu,\nu,\lambda)\cup{\rm
LR}(\mu,\nu^\bullet,\lambda)&\longrightarrow &{\rm
LR}(\mu^{t},{\nu^\bullet}^t,\lambda^{t})\cup{\rm
LR}(\mu^t,\nu^t,\lambda^t)\\
w(T)&\mapsto&  w_{col}(T^{\blacklozenge\bullet})=w(T)^{\blacklozenge\bullet}.
\end{array}
\end{equation}
is an involution which transposes the inner shape, outer shape and reverses and transposes  the weight.
\begin{remark} \label{re:crocodile}If $T$ is an LR tableau, $(\bullet\blacklozenge T)^t$ is obtained from $T$ by replacing the horizontal strip $i^{\nu_i}$, from SE to NW, with $\nu_1 \nu_1-1\cdots\nu_1-\nu_i+1$, for all $i$.

 If $T$ is a opposite LR tableau, $(\bullet\blacklozenge T)^t$ is obtained from $T$ by replacing the horizontal strip $(d-i+1)^{\nu_i}$, from  NW to SE, with $12\cdots \nu_i$, for all $i$.
\end{remark}

From the discussion above, it follows the next proposition.

 \begin{proposition} \label{prop:blackbullet}The rotation $\bullet$ and the orthogonal transpose $\blacklozenge$ maps commute on the set of LR and opposite LR tableaux,
$\scriptsize{\blacklozenge\bullet=\bullet\blacklozenge}$,  that is, $(\blacklozenge \bullet)^2=1$.
\end{proposition}

\begin{example} \label{ex:bullblack} Let $n=7$ and $d=3$. Let $\rm
    T=\begin{Young} 1&3&&\cr &1&2&2\cr &&1&1\cr
\end{Young}$ be an LR tableau of shape $\lambda^\vee/\mu$, weight $\nu=421$, where $\lambda=200$, $\lambda^\vee=442$, $\mu=210$, and
   with word $w=1122131$. Then $w^\bullet=3132233$, $w^\blacklozenge=1212314$, column word of $T^\blacklozenge$, and $w^{\blacklozenge\bullet}=1423434=w^{\bullet\blacklozenge}$, $\nu^t=3211$, column word of $T^{\bullet\blacklozenge}$,
    $${\rm
    T}=\begin{Young} 1&3&&\cr &1&2&2\cr &&1&1\cr
\end{Young}\;\ts\overset{\underset{}
    {\longleftrightarrow}}\ts\;\;
    \begin{Young} 4&1&&\cr &3&2&1\cr &&2&1\cr
\end{Young}
\ts\overset{\text{rotate\;\&\;transpose diagram}}{\underset{\text{$\begin{smallmatrix}\end{smallmatrix}$}}{\longleftrightarrow}}\ts
    \begin{Young} 4&&\cr 1&3&\cr &2&2\cr &1&1\cr
\end{Young}={\rm T}^{\blacklozenge}{\underset{\text{$\bullet$}}{\longleftrightarrow}}\ts
    \begin{Young} 4&4&\cr 3&3&\cr &2&4\cr &&1\cr
\end{Young}=\bullet\blacklozenge{\rm T}\;.$$
$$T\longleftrightarrow T^\bullet=\begin{Young} 3&3&&\cr 2&2&3&\cr &&1&3\cr
\end{Young}\;\;{\underset{}{\longleftrightarrow}}\ts\;\;\begin{Young} 4&3&&\cr 4&3&2&\cr &&4&1\cr
\end{Young}\ts\overset{\text{rotate\;\&\;transpose diagram}}{\underset{\text{$\begin{smallmatrix}\end{smallmatrix}$}}{\longleftrightarrow}}\ts\begin{Young} 4&4&\cr 3&3&\cr &2&4\cr &&1\cr
\end{Young}=T^{\bullet\blacklozenge}.$$
$\rm T^{\blacklozenge}$ is an LR tableau with
 shape $(\lambda^\vee/\mu)^{\bullet\,t}$ and column word $w^\blacklozenge=1212314$ of weight $\nu^t$, while
${\rm T}^{\blacklozenge\bullet}={\rm T}^{\bullet\blacklozenge}$ is a opposite LR tableau with
 shape $(\lambda/\mu)^{t}$  and column word $w^{\blacklozenge
 \bullet}=1423434$ of weight $\nu^{t\bullet}$,
  where $U(w)=\begin{Young} 6&&&\cr 3&4&&\cr 1&2&5&7\cr
\end{Young}$,
  $U(w^\blacklozenge)=\begin{Young} 7&&\cr 5&&\cr 2&4&\cr 1&3&6\cr
\end{Young}=U(w)^t$, and $U(w^{\blacklozenge \bullet})=\begin{Young} 2&5&7\cr &4&6\cr &&3\cr &&1\cr
\end{Young}=U(w)^{\bullet t}$.

Let $\rm
    T=\begin{Young} 3&3&&\cr &2&2&3\cr &&1&2\cr
\end{Young}$ be a opposite LR tableau with word $w=2132233$. Following the remark above, $(\blacklozenge \bullet\rm T)^t=\begin{Young} 1&2&&\cr &1&2&3\cr &&1&3\cr
\end{Young}$  and $\blacklozenge \bullet\rm T$ is a LR tableau with column word $w^{\blacklozenge\bullet}=1212313 $.
\end{example}

\subsubsection{Computational complexity of bijection $\blacklozenge$ on LR tableaux}

\label{S:comp}

We now show that the computational complexity of bijection
$\blacklozenge$ is linear on the input where we use skew LR tableaux. Hence, recalling  the definition of a linear cost bijection in Subsection \ref{S:linearcost}, the bijection $\blacklozenge$ is of linear cost.

\begin{alg}\label{alg}[Bijection $\blacklozenge$.]
\qquad \\
\emph{Input:} LR tableau $T$ of skew shape $\lambda / \mu$, with
$\lambda=\(\lambda_{1} \geq \ldots \geq \lambda_{n}\)$, \\
$\mu=\(\mu_{1} \geq \ldots \geq \mu_{n}\)$, and filling
$\nu=\(\nu_{1} \geq \ldots \geq \nu_{n}\)$, having { $A=(a_{i,j})
\in M_{n \times n}(\N)$} \quad ($a_{i,j}=0$ if $j>i$) as (lower
triangular) recording matrix.

Write $\widetilde{A}$, a copy of the matrix $A$.

For $j:=n$ down to $2$ do

\quad For $i:=1$ to $n$ do

\quad \quad Begin

\quad \quad \quad If $i=j$ then
$\widetilde{a}_{i,i}:=\widetilde{a}_{i,i}+\lambda_{1}-\lambda_{i}$

\quad \quad \quad \phantom{If $i=j$} else

\quad \quad \quad \phantom{If $i=j$ else} If $j>i$ then
$\widetilde{a}_{i,j}=0$ else
$\widetilde{a}_{i,j}:=\widetilde{a}_{i,j}+\widetilde{a}_{i,j+1}$.

\quad \quad End

\bigskip

So far the computational cost is $O(n^{2})=O(\langle A
\rangle)$.

\bigskip

Set a matrix $B=(b_{i,j}) \in M_{\lambda_{1} \times
\lambda_{1}}(\N)$ such that $b_{i,j}=0$ for all $i,j$.

For $i:=1$ to $n$ do

\quad Begin

\quad \quad Set $c:=0$.

\quad \quad For $j:=0$ to $n$ do

\quad \quad \quad Begin

\quad \quad \quad \quad $r:=\widetilde{a}_{i+j,i}-a_{i+j,i}$, \quad see Remark~\ref{R:comp}.

\quad \quad \quad \quad For $t:=1$ to $a_{i+j,i}$ do
$b_{r+t,c+t}:=b_{r+t,c+t}+1$.

\quad \quad \quad \quad $c:=c+a_{i+j,i}$.

\quad \quad \quad End

\quad End

This part has total computational cost at most equal to
$$O\(\sum_{1 \leq i.j \leq n}a_{i,j}\)=O\(|\lambda \setminus \mu|\)
=O\(|\lambda|-|\mu|\)=O\(\langle T \rangle\).$$

\emph{Output:} $B$ recording matrix of the output tableau $T^\blacklozenge$.

\end{alg}

\begin{remark} \label{R:comp}
For all $1 \leq i \leq n$ and $0\leq j \leq n-i+1$,
we have
$$\widetilde{a}_{i+j+1,i}-\widetilde{a}_{i+j,i} \geq a_{i+j+1,i}.$$
\end{remark}

\begin{corollary} \label{cor:blackbullet}
 The composition $\scriptsize{\blacklozenge\bullet=\bullet\blacklozenge}$ is a linear cost involution on ${\rm LR}(\mu,\nu^{},\lambda)\cup{\rm LR}(\mu^{t},\nu^{\bullet t},\lambda^t)$. It exhibits the symmetry $c_{\mu,\nu,\lambda}=c_{\mu^t,{\nu^\bullet}^t,\lambda^t}$ in linear time.
\end{corollary}
\subsection{The linear cost involution orthogonal transpose on LR companion  tableaux}\label{lozengecomp}

{If $G$ is the right companion of $T\in \rm LR(\mu,\nu,\lambda)$ and   $\iota(T^\blacklozenge)=:G^\blacklozenge$, where $\iota$ is the linear map bijection \eqref{rightcompanionmap},
  then if $B$ is the recording matrix of $T^\blacklozenge$, Algorithm \ref{alg}, $B^t$ is the recording matrix of $G^\blacklozenge$ and $G\overset{\iota^{-1}}\rightarrow T\overset{\blacklozenge}\rightarrow T^{\blacklozenge}\overset{\iota}\rightarrow G^\blacklozenge$ is a linear cost bijection. It is then  not difficult to define  a linear cost algorithm to directly calculate $G^\blacklozenge$ from $G$} without making recourse of the recording matrix.
We now recall the involution  on LR companion tableaux, also denoted $\blacklozenge$, as described in Steps 2 and 3 in Section 6.1 of \cite{leclen}
\begin{equation}\label{complozenge}\blacklozenge:\rm LR_{\nu,\lambda^\vee/\mu}\rightarrow \rm LR_{\nu^t,(\lambda^{\vee }/\mu)^{t\bullet}}, \;G\mapsto G^\blacklozenge\end{equation} such that $\iota(T^\blacklozenge)=\iota(T)^\blacklozenge$ whenever $T\in \rm LR(\mu,\nu,\lambda)$. We reproduce it below with slightly different notation.

\begin{alg} \label{alg:blacklozenge}[Construction of $G^\blacklozenge$.] Let $G\in \rm LR_{\nu,\lambda^\vee/\mu}$. The construction of $G^{\blacklozenge}$ has the following two steps.

{\em Step 1}. Transpose the tableau $G$ and denote the resulting filling of shape $\nu^t$ by $G^t$.

{\em Step 2}. For each $i=1,\dots, d$, consider in $G^t$ the vertical strip of $i's$ with size $\lambda^\vee_i-\mu_i$, and replace these entries,  from southeast to northwest, with $(\lambda^\bullet)_i+1,(\lambda^\bullet)_i+2,\dots,(\lambda^\bullet)_i+\lambda^\vee_i-\mu_i$ respectively. The resulting tableau is $G^\blacklozenge$ of shape $\nu^t$ and weight $(\lambda^{\vee t }/\mu^t)^{\bullet}$.
\end{alg}
\begin{example} In the previous example, $T\in LR(\mu,\nu,\lambda)$ has recording matrix
 $A=\left(\begin{smallmatrix}
 2&0&0&0\\
 1&2&0&0\\
 1&0&1&0\\
 0&0&0&0\\
  \end{smallmatrix}\right)$ and  $A^t$ is the recording  matrix of  the companion tableau  $G=\begin{Young} 3&&&\cr 2&2&&\cr 1&1&2&3\cr
\end{Young}$  of shape $\nu=421$ and weight $\lambda^\vee/\mu=232 $, where $\lambda=200$, $\lambda^\vee=442$, $\mu=210$. One has $G^t=\begin{Young} 3&&\cr 2&&\cr 1&2&\cr 1&2&3\cr
\end{Young}$,  $\lambda^\bullet=002$ and $(\lambda^\bullet)_1+1=1,(\lambda^\bullet)_1+2=2,$ $(\lambda^\bullet)_2+1=1,(\lambda^\bullet)_2+2=2, (\lambda^\bullet)_2+3=3$,
$(\lambda^\bullet)_3+1=3,(\lambda^\bullet)_3+2=4$. Then $B=A^\blacklozenge=\left(\begin{smallmatrix}
 2&0&0&0\\
 0&2&0&0\\
 1&0&1&0\\
 0&0&0&1\\
  \end{smallmatrix}\right)$ is the recording matrix of $T^\blacklozenge$  and $B^t$ is the recording matrix of $G^\blacklozenge=\begin{Young} 4&&\cr 3&&\cr 2&2&\cr 1&1&3\cr
\end{Young}=\iota(T^\blacklozenge)$ with shape $\nu^t=3211$ and weight $(\lambda^{\vee t}/\mu^t)^\bullet=(3322-2100)_{rev}=2221$ is the companion tableau of $T^\blacklozenge.$

\end{example}

From  the second part of Proposition \ref{prop:burge} and
duality of Burge correspondence, it follows
\begin{proposition}\label{prop:U}
 Let  $w$ and $u$ be two Yamanouchi words such that $w\equiv Y(\nu)$. Then

 $(a)$ $w^{\blacklozenge}\equiv Y(\nu^t)$, and $w\equiv u$ if and only if $ w^\blacklozenge\equiv u^\blacklozenge$.

 $(b)$ $w^{\bullet}\equiv Y(\nu^\bullet)$, and $w^{\bullet\blacklozenge}=w^{\blacklozenge\bullet}\equiv Y(\nu^{t\bullet})$.

 $(c)$ $Q(w^\bullet) ={\rm evac\,}Q(w)={\rm evac\,}U(w)=U(w)^{\bullet \rm
n}= U(w^\bullet)^{\rm n}=U(w)^{\rm a \bullet}$.

$(d)$ ${Q}(w^\blacklozenge)
={Q}(w)^t=U(w)^t$.

$(e)$  $Q(w^{\bullet\blacklozenge})=$ $Q(w^{\blacklozenge\bullet})={\rm evac\,} Q(w)^{\,t}=({\rm evac\,} Q(w))^{t\,}={\rm evac\,} U(w)^{\,t}=({\rm evac\,} U(w))^{t\,}$.
\end{proposition}


\section{ LR transposers  coincidence and  linear equivalence to an LR commuter}\label{sec:commutertransposer}

{In this section we work in $\mathcal{LR}$ as a set of LR tableaux or their companions.}
\subsection{An LR commuter for the symmetry $c_{\mu,\nu,\lambda}=c_{\lambda,\nu,\mu}$}

Let $T\in\rm LR(\mu,\nu,\lambda)$ and $B(T)$ the crystal connected component of $B(\lambda/\mu)$ containing $T$. The highest and lowest weight elements of $B(T)$ are $T$ and $e (T)$ respectively, and henceforth the highest and lowest weight elements of $B(T)^\bullet=B(\bullet T)$ are $\bullet e(T)$ and $\bullet(T)$ respectively. Since highest and lowest weight elements in a connected crystal component are related by reversal $e$, \begin{equation}\label{bullet}e\bullet T=\bullet e T.\end{equation}

 \begin{theorem} \label{fundamentalrev}Let $\rho:=\bullet\,e=e\bullet$. Then the involution $$\rho: {\rm LR}(\mu,\nu,\lambda)\longrightarrow {\rm LR}(\lambda,\nu,\mu), T\mapsto \rho(T)=e\bullet T $$
is an LR commuter that exhibits the symmetry $c_{\mu,\nu,\lambda}=c_{\lambda,\nu,\mu}$.
\end{theorem}

\begin{remark}  We may also use Knuth and dual Knuth equivalence to characterize $\bullet e (T)$ and $e\bullet  (T)$. That is we may use tableau switching as in Algorithm \ref{alg:KD} to calculate $e\bullet  (T)$. In fact ${T}^\bullet\in \rm LR(\lambda,\nu^\bullet,\mu)$, and, from Corollary \ref{cor:reversal},
$T^{\bullet\,e}=[ Y(\nu)]_K\cap [{T^\bullet}]_{dK},$
is the unique  LR tableau in $\rm LR(\lambda,\nu,\mu)$ of the crystal connected component $B({T}^\bullet)$ in $B(\mu^\vee/\lambda)$. 
On the other hand, $T^{e\bullet}=[ Y(\nu^\bullet)^\bullet]_K\cap [{T^\bullet}]_{dK}=[Y(\nu)^{\rm a}]_K\cap [T^\bullet]_{dK}=T^{\bullet e}$.
\end{remark}

\subsection{Coincidence of LR transposers and linear  equivalence to an LR commuter} Let $T$ be an LR tableau.
Recall that $\sigma_0$ coincides with the reversal on   $T$   and on  $T^{\textsf{low}}$ the lowest weight element  of $B(T,d)$. The $\texttt{reversal}(T)=T^{\textsf{low}}$ may be computed by the action of  $\sigma_0$ on $B(T,d)$.  Recall that, for $w$  a Yamanouchi or opposite Yamanouchi word, $\sigma_0 w$ coincides with the reversal $e$ on the LR or opposite LR diagonally-shaped tableau with reading word $w$. Thus the column  word of $\bullet \blacklozenge eT$ is $(\sigma_0 w)^{\blacklozenge\,\bullet}$.

 \begin{proposition}  The following holds on $\mathcal{LR}$ as a set of LR tableaux:

 $(a)$ $\sigma_0 (w^\blacklozenge)=(\sigma_0 w)^\blacklozenge$, where $w$ is a Yamanouchi or opposite Yamanouchi word.

 $(b)$ $e\,\blacklozenge=\blacklozenge\,e$, that is,  $(e\,\blacklozenge)^2=1$.

 $(c)$ the involutions $\blacklozenge$,  $e$, $\bullet$ pairwise commute.

 \end{proposition}
 \begin{proof} $(a)$ Suppose that $w\equiv Y(\nu)$. Then, from Proposition \ref{prop:U}, $w^\blacklozenge\equiv Y(\nu^t)$ and subsection \ref{groupaction}, $\sigma_0 w\equiv w^\bullet\equiv Y(\nu^\bullet)$. Thereby, $(\sigma_0 w)^\diamond\equiv
 w^{\bullet\blacklozenge}=w^{\blacklozenge\,\bullet}\equiv\sigma_0( w^\blacklozenge)\equiv  Y(\nu^{t\bullet})$. On the other hand, $\sigma_0 (w^\blacklozenge)$ and $ (\sigma_0 w)^\blacklozenge$ are dual equivalent because their $Q$-symbol is $Q(w)^t$. In fact $Q(w)^t=Q(w^\blacklozenge)=Q(\sigma_0(w^\blacklozenge))$ and $Q(w)^t=Q(\sigma_0w)^t=Q((\sigma_0w)^\blacklozenge)$. Therefore $\sigma_0 (w^\blacklozenge)=(\sigma_0 w)^\blacklozenge$.

 $(b)$ Let $T$ be an LR tableau 
 where { $T\in B_d(\lambda/\mu)$}.
 The row  and the column reading of the tableaux in  $B(T)$ embeds $B(T)$ into usually different subcrystals of $B^{|\lambda|-|\mu|}$ but  isomorphic. Henceforth, either we consider the row reading or the column reading of  $T$, one always has,  $\sigma_0(w_{col}(T))=w_{col}(T^e)$,
 and $\sigma_0(w(T))=w(T^e)$, with $\omega_0\in \mathfrak{S}_d$ .
 Then the (row reading) word of $T^e$ is $\sigma_0 w(T)$, the column reading word of $T^\blacklozenge\in B_{\nu_1}(\mu^{t\vee}/\lambda^{t\vee})\subseteq B_{n-d}(\mu^{t\vee}/\lambda^{t\vee})$ is $w(T)^\blacklozenge$ on the alphabet $[\nu_1]$, and the column reading word of $T^{\blacklozenge\,e}$ is $\sigma_0(w(T)^\blacklozenge)$ with with $\omega_0\in \mathfrak{S}_{\nu_1}$.  Henceforth, from $(a)$,  $w_{col}(T^{e\,\blacklozenge})=(\sigma_0 w)^\blacklozenge=\sigma_0 (w^\blacklozenge)=\sigma_0(w_{col}(T^{\blacklozenge}))=w_{col}(T^{\blacklozenge\,e})$.

 $(c)$ It follows from $(b)$, Theorem \ref{fundamentalrev} and Proposition \ref{prop:blackbullet}.
 \end{proof}

  The following is an illustration of the previous result.
\begin{example} \label{crystaloperator} Let $w=1111221332$ on the alphabet $[3]$ with weight $\nu=532$, and $w^\blacklozenge=1234125123$ on the alphabet $[\nu_1=5]$ of weight $\nu^t=33211$. Then $\sigma_0
w=$ $\sigma_1\sigma_2\sigma_1w=$ $\sigma_1\sigma_2(2211221332)=$ $\sigma_1(3311221333)=3311222333$ and $(\sigma_0 w)^\blacklozenge =1245345345$. On the other hand,

$\sigma_0(w^\blacklozenge)=\sigma_1.\sigma_2\sigma_1.\sigma_3\sigma_2\sigma_1.\sigma_4\sigma_3\sigma_2\sigma_1 (w^\blacklozenge)=
\sigma_1.\sigma_2\sigma_1.\sigma_3\sigma_2\sigma_1.\sigma_4\sigma_3 1234135123 $

$=\sigma_1.\sigma_2\sigma_1.\sigma_3\sigma_2\sigma_1.\sigma_4 1234145124=
 \sigma_1.\sigma_2\sigma_1.\sigma_3\sigma_2\sigma_1 1235145125$

 $=\sigma_1.\sigma_2\sigma_1.\sigma_3\sigma_2 1235245125=\sigma_1.\sigma_2\sigma_1.\sigma_3 1235345135=\sigma_1.\sigma_2\sigma_1 1245345145$

 $=\sigma_1.\sigma_2 1245345245=\sigma_1 1245345345=1245345345=(\sigma_0w)^\blacklozenge.$

 If one considers the alphabet $[4]$,  one has $\sigma_0w=\sigma_1\sigma_2\sigma_1\sigma_3\sigma_2\sigma_1w=4422333444$ but still $(\sigma_0w)^\blacklozenge=\sigma_0w^\blacklozenge$.
\end{example}
From this discussion,  Section \ref{s:maps}, and  
the computational complexity of the reversal involution $e$ \cite{PV2},    it then follows
 \begin{theorem}\label{plr*}
    Let $\rm T$ be a LR tableau with
shape $\lambda/\mu$ and word $w$.
  Let $\varrho:=\blacklozenge\rho=$ $\blacklozenge\bullet e=$ $\bullet\blacklozenge \,e= \bullet\,e\blacklozenge=\rho\blacklozenge$ where $\rho=\bullet\,e$. Then
$$e: LR(\mu,\nu,\lambda)\rightarrow LR(\mu,\nu^\bullet,\lambda),\, T\mapsto e T,\,w(e T)=\sigma_0 w;$$
$$\rho: LR(\mu,\nu,\lambda)\rightarrow LR(\lambda,\nu,\mu),\, T\mapsto \rho(T)=\bullet e T,\,w(\rho (T))=(\sigma_0 w)^\bullet;$$
and $$\varrho: LR(\mu,\nu,\lambda)\rightarrow LR(\mu^t,\nu^t,\lambda^t),\, T\mapsto \varrho(T)=\blacklozenge\rho(T)=\blacklozenge\bullet e T,\,w_{col}(\varrho(T))=(\sigma_0 w)^{\blacklozenge\,\bullet}$$
are involutions exhibiting the symmetries $c_{\mu\,\nu\,\lambda}=c_{\mu\,\nu^\bullet\,\lambda}$, $c_{\mu\,\nu\,\lambda}=c_{\lambda\,\nu\,\mu}$ and $c_{\mu\,\nu\,\lambda}=c_{\mu^t\,\nu^t\,\lambda^t}$ respectively. The three involutions $e,\rho,$ and $\varrho$ are  linear time equivalent to each other and in particular to the reversal $e$.
\end{theorem}

\begin{remark}\label{remark:evarrho}From the identity $\varrho=\blacklozenge\rho=(\blacklozenge\bullet) e$ we conclude that $\varrho(T)$ can be obtained from $T^e$ by replacing,  for $i=1,\dots, \ell(\nu)$,   from NW to SE, the entries of the $i$-horizontal strip of length $\ell(\nu)-i+1$ in $T^e$,  with $1,\dots,\nu_i-i+1$. This gives $ (\varrho(T))^t=(\blacklozenge\rho(T))^t$ a row semistandard tableau.
\end{remark}
We  may also use Knuth and dual Knuth equivalence to characterize the bijection $\varrho$. This shows that $\varrho (T)$ can also be calculated using tableau switching as in Algorithm \ref{alg:KD}. This is the procedure offered in \cite{bss}.
\begin{corollary} \label{cor:varrho} Let $\rm T$ be a LR tableau with
shape $\lambda/\mu$ and weight $\nu$. Then $\varrho(T)=T^{e\,\blacklozenge\bullet}$ is the unique tableau
Knuth equivalent to $Y(\nu^t)$ and dual Knuth equivalent to $\widehat
T^t$. That is, $\varrho(T)=T^{e\,\blacklozenge\bullet}=[Y(\nu^t)]_K\cap[\widehat{T}^t]_{dK}$.
\end{corollary}
\begin{proof} Let $w$ be the word of $T$ with weight $\nu$. One has $\sigma_0w\equiv w^\bullet$, and $(\sigma_0 w)^{\bullet}\equiv w\equiv Y(\nu)$. Then $(\sigma_0 w)^{\blacklozenge\,\bullet}= ((\sigma_0 w)^\bullet)^{\blacklozenge}=(\sigma_0 (w^\bullet))^{\blacklozenge}\equiv w^\blacklozenge\equiv Y(\nu^t)$. Recall that dual Knuth equivalence between tableaux can be checked either by using row  or column words.

From Haiman's theorem, Theorem \ref{t2}, $[Y(\nu^t)]_K\cap[\widehat{T}^t]_{dK}$ has a sole tableau dual Knuth equivalent to ${\widehat T}^t$ and Knuth equivalent to $Y(\nu^t)$. Since $T^{e\,\blacklozenge\bullet}\equiv Y(\nu^t)$, it is enough to see that the column words of $T^{e\,\blacklozenge\bullet}$ and $\widehat T^t$ have the same $Q$--symbol, that is, $T^{e\,\blacklozenge\bullet}$ is the highest weight element of the connected component  $B(\widehat T^t)$ in  $B((\lambda/\mu)^t)$, on the alphabet $|\lambda|-|\mu|$. Let $\widehat w$ be the word of $\widehat T$. As $rev\,\widehat w$,
the reverse word of $\widehat T$, is the column word of  $\widehat T^t$, we want to show that $Q(rev \widehat w)=Q((\sigma_0w)^{\blacklozenge\, \bullet})$.

We know  that  any word $u$ is dual Knuth equivalent to $\sigma_0 u$, $Q(u)=Q(\sigma_0u)$ and, from Proposition \ref{prop:U}, $Q(u^\bullet)=Q(u)^E$ and $Q(rev\, u)=Q(u)^{Et}$  ~\cite{st}.
Recalling Proposition \ref{prop:U}, $Q(rev\,\widehat w)=Q(\widehat w)^{E\,t}=$ $Q(w)^{E\,t}=Q(\sigma_0 w)^{Et}=Q((\sigma_0w)^{\blacklozenge\, \bullet}).$
\end{proof}

In \cite{bss}, it is observed that the White and the
Hanlon--Sundaram maps ~\cite{white,HS}  produce the same result,  denoted by $\varrho^{WHS}$.
Thus $\varrho^{BSS}(T)$ can be
obtained either by tableau--switching or by the White--Hanlon--Sundaram
transformation $\varrho^{WHS}$ or by $\varrho$.
\begin{theorem}\label{coinctransposer} The LR transposers $\varrho^{BSS}$, $\varrho^{WHS}$ and $\varrho$ are
identical,
and linear time equivalent 
to the reversal involution $e$.

\end{theorem}

\begin{remark} The Sch\"utzenberger's
 \textit{jeu de taquin} formulation of the LR rule says: Fix  $\rm T_\nu \in SYT(\nu)$. The number of $\rm T\in SYT(\lambda^\vee/\mu)$ such that $\rm \textsf{rectification}(T) =T_\nu$
 equals $c_{\mu,\nu,\lambda}$ \cite[Appendix 1]{st}. Equivalently, since $\rm T_\nu$ is Knuth equivalent to $\rm T_\nu^{\rm a}$, the number of $\rm T\in SYT(\lambda^\vee/\mu)$ such that $\rm \textsf{arectification}(T) =T_\nu^{\rm a}$
 equals $c_{\mu,\nu,\lambda}$. The definition only depends on the shapes $\mu,\nu,\lambda$ and not on a particular choice of the filling of the Young diagram of $\nu$.  However, choosing a certain filling of $\rm T_\nu \in SYT(\nu)$ allows to relate this definition with LR tableaux \cite[Appendix 1]{st}. That is, choosing $\rm T_\nu$ to be a standardization of the Yamanouchi tableau $\rm Y(\nu)$ incurs that $\rm \textsf{rectification}(T) =T_\nu$, with $\rm T\in SYT(\lambda^\vee/\mu)$, holds only if the $\nu$-semistandardization of $\rm T$ is an LR tableau of shape $\lambda^\vee/\mu$ and weight $\nu$ \cite[Lemma A1.3.6, Lemma A1.3.7 ]{st}.

  Fix the tableaux $\rm T_\mu\in SYT(\mu)$, $\rm T_\nu\in SYT(\nu)$ and $\rm T_\lambda\in SYT(\lambda)$ to be  the standardizations of the Yamanouchi tableaux $\rm Y(\mu)$, $\rm Y(\nu)$, and $\rm Y(\lambda)$ respectively. Choose $\rm T\in SYT(\lambda^\vee/\mu)$ to be the standardization of an LR tableau that rectifies to $\rm T_\nu$, to initialize a Thomas-Yong \emph{carton filling} \cite{carton},  built upon Fomin's \emph{jeu de taquin} growth-diagrams and  the \emph{infusion involution} \cite{infusion}, a particular case of Benkart-Sottille-Stroomer tableau-switching on pairs of standard tableaux. Let $\texttt{CARTONS}_{\mu,\nu,\lambda}$ be the set of all carton fillings built in this way with initial data $ \rm T_\mu, \rm T_\nu,\rm T_\lambda$.
   The number of carton fillings  is equal to the number of standard tableaux of shape $\lambda^\vee/\nu$ which rectify to $T_\nu$, that is, $c_{\mu,\nu,\lambda}$.
 For this particular choice of $\rm T_\mu$, $\rm T_\nu$ and $\rm T_\lambda$ the carton filling besides to
 showing the $\mathfrak{S}_3$-symmetries of LR coefficients $c_{\mu,\nu,\lambda}=c_{\epsilon,\delta,\gamma}$, where $(\epsilon,\delta,\gamma)$ is any permutation of $(\mu,\nu,\lambda)$, also gives tableau-switching  bijections on LR tableaux exhibiting such symmetries.   More precisely, composing the infusion involution with  the  semistandardization of those standard tableaux in the carton filling , one obtains the tableau-switching on LR tableaux and thus the carton filling gives tableau-switching  bijections on LR tableaux exhibiting such symmetries.

 Any carton filling gives a growth diagram on the face $\emptyset-\nu-\lambda^\vee-\mu$
for which
the edge $\mu-\lambda^\vee$
is a standard tableau $\rm T_{\lambda^\vee/\mu}$ of shape $\lambda^\vee/\mu$
rectifying to
$\rm T_\nu$. By the \textit{jeu de taquin} Littlewood-Richardson rule, fillings of this face count $c_{\mu,\nu,\lambda}$. Any such growth-diagram of this face extends uniquely to  a filling of the entire carton. The carton initialized with $ \rm T_\mu, \rm T_\nu,\rm T_\lambda$ and with
$\rm T_{\lambda^\vee/\mu}$ on the edge $\mu-\lambda^\vee$, by the symmetry of \emph{jeu de taquin}, also contains a standard tableau
 $\rm T_{\lambda^\vee/\nu}$ on the edge $\nu-\lambda^\vee$. Denoting by $\rho_1$ the infusion corresponding to rectification and by $\rho_2$ the infusion corresponding to antirectification, each of the six faces has a pair of skew standard tableaux and a pair of standard tableaux of normal shape or of  antinormal shape:
$\rho_2(\rm T_{\lambda^\vee/\mu},T_\lambda^{\rm a})=(\rm T_{\nu^\vee/\mu},T^{\rm a}_\nu)$,
$\rho_1(\rm T_\mu,\rm T_{\nu^\vee/\mu})=(\rm T_\lambda, T_{\nu^\vee/\lambda})$, $\rho_2( T_{\nu^\vee/\lambda}, \rm T^{\rm a}_\nu)=(T_{\mu^\vee/\lambda},\rm T^{\rm a}_\mu)$
 or
 $\rho_1(\rm T_{\mu},T_{\lambda^\vee/\mu})=(\rm T_{\nu},T_{\lambda^\vee/\nu})$,
$\rho_2(\rm T_{\lambda^\vee/\nu},\rm T^{\rm a}_{\lambda})=( T_{\mu^\vee/\nu},\rm T^{\rm a}_{\mu})$, $\rho_1(\rm T_\nu, T_{\mu^\vee/\nu})=(\rm T_\lambda,T_{\mu^\vee/\lambda})$. That is $\rho_1\rho_2\rho_1(\rm T_{\mu},T_{\lambda^\vee/\mu}, \rm T_{\lambda}=\rho_2\rho_1\rho_2(\rm T_{\mu},T_{\lambda^\vee/\mu}, \rm T_{\lambda})=(\rm T_{\lambda},T_{\mu^\vee/\lambda},\rm T^{\rm a}_\mu)$
$\texttt{CARTONS}_{\mu,\nu,\lambda}\rightarrow \texttt{CARTONS}_{\lambda,\nu,\mu}$
 \end{remark}

\subsection{  LR companion tableaux and Lascoux's  double crystal graph}\label{lascdouble}
{
 Let $T$ be an LR tableau of shape $\lambda/\mu$ with  companion $G$. We recall Lascoux's {\em double crystal graph structure on biwords} \cite{double} where a crystal operator consists of a left and a right operator, and a crystal string of a left and a right string. Let  ${\mathcal B}(T,G)$ be the crystal graph  whose vertices consist of the collection of biwords whose recording tableau is $G$ in the RSK-correspondence, subsection \ref {burgecomp}, with highest weight element the biword $W^{\lambda/\mu}=\left(\begin{smallmatrix}{\bf y}\\
w(T)
\end{smallmatrix}\right)$ on the LHS of \eqref{biword} identified with
  the LR tableau  $T$,   and lowest weight element the biword $\left(\begin{smallmatrix}{\bf y}\\
\sigma_0 w(T)
\end{smallmatrix}\right)$ identified with $\sigma_0T=T^e$. The vertices of $\mathcal{B}(T,G)$ describe simultaneously integer matrices and tableaux. The latter have the former as recording matrices. Instead of biwords we may consider integer matrices which are the recording matrices of the tableaux that they do describe.

 Proposition \ref{prop:burge},
  exhibits the Lascoux's  double crystal graph structure on biwords \cite{double}. It shows that the strings of the crystal graph $\mathcal{B}(T,G)$,
are transformed, by reordering the biwords (or transposing the corresponding matrices) as on the RHS of  \eqref{biword},   into the   strings    of  the {\em cocrystal}  $\mathcal{CB}(T,G)$, the set of biwords whose insertion tableau is $G$ \cite[p.103]{double}, with top     biword $W^{\nu}=\left(\begin{smallmatrix}{w(G)}\\
{\bf x}
\end{smallmatrix}\right)$, and bottom biword    $\left(\begin{smallmatrix}w(G^{\rm a})\\
{\bf x}
\end{smallmatrix}\right)$ with $G^{\rm a}$ the anti--normal form of $G$.   Under this reordering of the biletters, the Kashiwara operators in $\mathcal{B}(T,G)$ are translated to}
the  cocrystal operators or  right operators,  elementary
  {\em jeu de taquin} ({\em reverse  jeu de taquin}) operations on  two--row tableaux (see also \cite{az3}). Again instead of biwords we may consider as vertices of $\mathcal{CB}(T,G)$ the transpose of the matrices as the vertices of the crystal $\mathcal{B}(T,G)$. The matrices as vertices of $\mathcal{CB}(T,G)$ are the recording matrices of the tableaux that they do describe.

 An $i$-string in $\mathcal{CB}(T,G)$
 is an ordered  string of  two--row words in the plactic class of a two--row tableau \cite[p.100]{double}. The two-row words
being ordered according to the length of their bottom row,    such that  the two--row  word on the top is a two-row tableau of  partition shape   $(\nu_i,\nu_{i+1})$, and on the bottom    is the anti-normal  form of the two-row tableau on the top,  
see \cite[Section 2]{double}.
 Let $\theta_i$ be translation of the crystal reflection operators $\sigma_i$ on the cocrystal under the reordering of the billeters in the biwords as explained above. Thus the  symmetric group also acts on the cocrystal through the involutions $\theta_i$ which reflects each $i$-string about the middle, for $i\in{1,\dots,n-1}$. In particular,  $\theta _i$ sends the top of a $i$-string, a two--row tableau, to its anti-normal form in the bottom and vice versa.
Indeed the
entries of the $j$--th row of $\theta_iG$ are precisely the $k$'s
 telling us in which rows of $\sigma_iT$ the $j$'s are filled in.  Indeed
$\iota(\sigma_i T)= \theta_i G$ and put
$\theta_0:=\theta_{i_N}\dots\theta_{i_1}$ where $\sigma_0=\sigma_{i_N}\dots\sigma_{i_1}$. Thus $\iota(\sigma_0 T)=\iota (T^e)=\theta_0G=G^{\rm a }$.
This defines the commutative scheme
$$\begin{matrix}
T&\longleftrightarrow&\sigma_{i_1}T&\longleftrightarrow&\sigma_{i_2}\sigma_{i_1}T&\longleftrightarrow&\cdots&
\longleftrightarrow&\sigma_0T=T^e\\
{\iota}\updownarrow&        &{\iota}\updownarrow&                    &{\iota}\updownarrow           &                    &      &                  &{\iota}\updownarrow\\
G&\longleftrightarrow&\theta_{i_1}G&\longleftrightarrow&\theta_{i_2}\theta_{i_1}G&\longleftrightarrow
&\cdots&\longleftrightarrow&\theta_0G=G^{\rm
a}.\end{matrix}$$


\begin{example} Consider Example \ref{ex:burge}. The following exhibits the action of $\mathfrak{S}_3$ on the key tableaux (straight shape tableaux whose weight is  a reordering of the partition shape) of $B(\nu,G)$ whose weight is  a reordering of the partition shape)

\begin{figure}[!ht]
\centering
\begin{tikzpicture}[line cap=round,line join=round,>=triangle 45,x=1cm,y=1cm, scale=0.88]
\clip(-7,-2) rectangle (7,6);
\draw (0,5) node {$
T\equiv W^{\lambda/\mu}=\left(\begin{smallmatrix}1&1&1&1&2&2&2&3&3&3\cr
1&1&1&1&2&2&1&3&3&2
\end{smallmatrix}\right)$};
\draw [->] (-0.3,4.5)-- (-2,3.5);
\draw [->] (0.3,4.5)-- (2,3.5);
\draw (-0.9,3.8) node {$\sigma_1$};
\draw (0.9,3.8) node {$\sigma_2$};
\draw (-3.2,3) node {$\sigma_1T\equiv \left(\begin{smallmatrix}1&1&1&1&2&2&2&3&3&3\cr
2&2&1&1&2&2&1&3&3&2
\end{smallmatrix}\right)$};
\draw (3.2,3) node {$\sigma_2T\equiv\left(\begin{smallmatrix}1&1&1&1&2&2&2&3&3&3\cr
1&1&1&1&2&2&1&3&3&3
\end{smallmatrix} \right)$};
\draw [->] (-2,2.5)-- (-2,1);
\draw [->] (2,2.5)-- (2,1);
\draw (-1.65,1.7) node {$\sigma_2$};
\draw (1.65,1.7) node {$\sigma_1$};
\draw (-3.2,0.5) node {$\sigma_2\sigma_1 T\equiv\left(\begin{smallmatrix}1&1&1&1&2&2&2&3&3&3\cr
3&3&1&1&2&2&1&3&3&3
\end{smallmatrix}\right)$};
\draw (3.2,0.5) node {$\sigma_1\sigma_2 T\equiv \left(\begin{smallmatrix}1&1&1&1&2&2&2&3&3&3\cr
2&2&1&1&2&2&2&3&3&3
\end{smallmatrix} \right)$};
\draw [->] (-2,-0) -- (-0.3,-1);
\draw [->] (2,-0) -- (0.3,-1);
\draw (-0.9,-0.3) node {$\sigma_1$};
\draw (0.9,-0.3) node {$\sigma_2$};
\draw (0,-1.5) node {$\sigma_0 T\equiv \left(\begin{smallmatrix}1&1&1&1&2&2&2&3&3&3\cr
3&3&1&1&2&2&2&3&3&3
\end{smallmatrix} \right)$};
\end{tikzpicture}
\end{figure}

Reordering the biletters so that the biword bottom row is a row word, equivalently, transposing the  matrices defined by the biwords above, one obtains Figure \ref{fig:doublecrystal}

\begin{figure}[!ht]
\centering
\begin{tikzpicture}[line cap=round,line join=round,>=triangle 45,x=1cm,y=1cm, scale=0.88]
\clip(-7,-2) rectangle (7,6);
\draw (0,5) node {$
G\equiv W^\nu=\left(\begin{smallmatrix}2&1&1&1&1&3&2&2&3&3\cr
1&1&1&1&1&2&2&2&3&3
\end{smallmatrix}\right)$};
\draw [->] (-0.3,4.5)-- (-2,3.5);
\draw [->] (0.3,4.5)-- (2,3.5);
\draw (-0.9,3.8) node {$\theta_1$};
\draw (0.9,3.8) node {$\theta_2$};
\draw (-3.2,3) node {$\theta_1 G\equiv\left(\begin{smallmatrix}2&1&1&3&2&2&1&1&3&3\cr
1&1&1&2&2&2&2&2&3&3
\end{smallmatrix}\right)$};
\draw (3.2,3) node {$\theta_2 G\equiv\left(\begin{smallmatrix}2&1&1&1&1&2&2&3&3&3\cr
1&1&1&1&1&2&2&3&3&3
\end{smallmatrix}\right)$};
\draw [->] (-2,2.5)-- (-2,1);
\draw [->] (2,2.5)-- (2,1);
\draw (-1.65,1.7) node {$\theta_2$};
\draw (1.65,1.7) node {$\theta_1$};
\draw (-3.2,0.5) node {$\theta_2\theta_1 G\equiv\left(\begin{smallmatrix}2&1&1&2&2&3&3&3&1&1\cr
1&1&1&2&2&3&3&3&3&3
\end{smallmatrix}\right)$};
\draw (3.2,0.5) node {$\theta_1\theta_2 G\equiv \left(\begin{smallmatrix}1&1&2&2&2&1&1&3&3&3\cr
1&1&2&2&2&2&2&3&3&3
\end{smallmatrix} \right)$};
\draw [->] (-2,-0) -- (-0.3,-1);
\draw [->] (2,-0) -- (0.3,-1);
\draw (-0.9,-0.3) node {$\theta_1$};
\draw (0.9,-0.3) node {$\theta_2$};
\draw (0,-1.5) node {$\theta_0 G\equiv \left(\begin{smallmatrix}1&1&2&2&2&1&1&3&3&3\cr
1&1&2&2&2&3&3&3&3&3
\end{smallmatrix} \right)$};
\end{tikzpicture}
\caption{}\label{fig:doublecrystal}
\end{figure}

\end{example}

Thus and taking into account Algorithm \ref{alg:blacklozenge}, we have the following result relating the tableaux $T$, $T^e$, $T^{e\bullet}$ and $T^{e\bullet\blacklozenge}$ respectively with their corresponding companion tableaux $G$, $G^{\rm a}$, $G^{E}$ and $G^{E\blacklozenge}=G^{\blacklozenge E}$. See also the construction of $G^{E\blacklozenge}$ in \cite[Section 6.1]{leclen}.

{ \begin{theorem} \label{th:companion} Let $\rm T$ be a LR tableau with
shape $\lambda/\mu$ and right LR companion tableau $G$. Then

 $(a)$ the following diagram is commutative:
$$\begin{matrix}
T&\overset{\text{$e$}}{\longleftrightarrow}&T^e&\overset{\text{\,$\bullet$}}{\longleftrightarrow}
&\rho(T)=T^{e\bullet}&\overset{\text{$\blacklozenge$}}{\longleftrightarrow}&
{\varrho(T)=T^{e\bullet\blacklozenge}}\\
{\iota}\updownarrow&&{\iota}\updownarrow&&{\iota}\updownarrow&&{\iota}\updownarrow\\
G&\underset{\rm a}{\longleftrightarrow}&G^{\rm a}&
\underset{\bullet}{\longleftrightarrow}&G^{{\rm
a}\bullet}=\rm{evac}\,G&\underset{\text{$\blacklozenge$}}\longleftrightarrow &\blacklozenge \rm{evac}\,G.
\end{matrix}
$$
$(b)$ the involution symmetries $\rho$ and $\varrho$ translate to companion LR tableaux as follows $$\begin{array}{ccc}\rho:\rm LR_{\nu,\lambda/\mu}\longrightarrow \rm LR_{\nu,(\lambda/\mu)^\bullet}: G\mapsto {\rm evac}\, G\;\text{such that}\;{\rm evac}\,\iota(T)=\iota(\rho(T));\\
\\
\varrho: \rm LR_{\nu,\lambda/\mu}\longrightarrow \rm LR_{\nu^t,(\lambda/\mu)^t}: G\mapsto \blacklozenge{\rm evac}G={\rm evac}\blacklozenge G \;\text{such that}\;{\rm evac}\,\blacklozenge\iota(T)=\iota(\varrho(T)).
\end{array}$$
\end{theorem}
}


\subsection{Illustration}
Consider the BSS'LR transposer  \cite{bss}
$$\begin{matrix}
\varrho^{BSS}:&\rm LR(\mu,\nu,\lambda)&\rightarrow&\rm LR(\mu^t,\nu^t,\lambda^t)\\
&T&\mapsto&\varrho^{BSS}(T)=[Y(\nu^t)]_K\cap[\widehat{T}^t]_{dK}
\end{matrix}.\qquad 
$$
\noindent  The image of $T$ by the BSS--bijection is the unique
tableau of shape $\lambda^t/\mu^t$ whose rectification is
$Y(\nu^t)$ and the $Q$--symbol of the column reading word is
$Q(w(T))^{Et}$. The idea behind this bijection can be told as follows:
$\widehat T$ constitutes a set of instructions telling where
expanding slides can be applied to $Y(\mu)$. Then $\widehat T^t$ is
a set of instructions telling where expanding slides can be applied
to $Y(\mu)^t$. Tableau--switching provides an algorithm  to give way
to those instructions. In the following, $\textbf{s}$ denotes switching: $$\begin{array}{ccccccccccccccc} {Y(\mu)\cup
T}&\underset{\text{of T}}{\overset{\text{\scriptsize
standardization}}{\longrightarrow}}& {Y(\mu)\cup \widehat
T}&\underset{\text{of $\widehat T$}}{\overset{\text{\scriptsize
transposition}}{\longrightarrow}}&{Y(\mu^t)\cup {\widehat T}^t}&&Y(\mu^t)\cup \varrho^{BSS}(T)&\\
&&&&
 \downarrow\!{\scriptstyle s}&&\uparrow\!{\scriptstyle s}&\\
 &&&&{(\widehat{T}^t)^{\rm n}\cup Z}&\overset{\!{}}\mapsto&Y(\nu^t)\cup Z&
\end{array}.$$
\noindent 
Then
$\varrho^{BSS}(T)\equiv Y(\nu^t)$ and $\varrho^{BSS}(T)\overset{d}{\equiv}
\widehat T^t$.

\begin{example} Let $T$ in  $\rm LR (\mu, \nu,\lambda)$ with $\mu=21$, $\nu=532$ and $\lambda=643$.

$\begin{array}{cccccccccc}
{T=\begin{Young} 2&3&3&&&\cr &1&2&2&&\cr &&1&1&1&1\cr
\end{Young}}& \rightarrow&
{\widehat{T}=\begin{Young} 6&9&10&&&\cr &1&7&8&&\cr &&2&3&4&5\cr
\end{Young}}& \rightarrow&
&{\widehat{T}^t=\begin{Young}5&&\cr 4&&\cr 3&8&\cr 2&7&10\cr &1&9\cr &&6\cr
\end{Young}}&\rightarrow
\end{array}$

${\begin{matrix} \rightarrow&
Y(\mu^t)\cup\widehat{T}^t=\begin{Young}5&&\cr 4&&\cr 3&8&\cr 2&7&10\cr {\color{red}{\bf{2}}}&1&9\cr {\color{red}{\bf{1}}}&{\color{red}{\bf{1}}}&6\cr
\end{Young}
&&
\begin{Young}5&&\cr 4&&\cr 2&3&\cr 1&2&3\cr {\color{red}{\bf{2}}}&1&2\cr {\color{red}{\bf{1}}}&{\color{red}{\bf{1}}}&1\cr
\end{Young}=Y(\nu^t)\cup\varrho^{BSS}(T)\\
&{\scriptstyle s}\!\downarrow&&\uparrow\!{\scriptstyle s}\\
&(\widehat{T}^t)^{\rm n}\cup
Z=\begin{Young}{\color{red}\bf{1}}&&\cr 5&&\cr 4&{\color{red}\bf{2}}&\cr 3&8&{\color{red}\bf{1}}\cr 2&7&10\cr 1&6&9\cr
\end{Young}&\rightarrow&\begin{Young}{\color{red}\bf{1}}&&\cr 5&&\cr 4&{\color{red}\bf{2}}&\cr 3&3&{\color{red}\bf{1}}\cr 2&2&2\cr 1&1&1\cr
\end{Young}=Y(\nu^t)\cup
Z
\end{matrix}}.$
\end{example}

Consider the involution $\varrho=e\bullet\blacklozenge$,

\begin{example} Recall Example \ref{crystaloperator}.
Letting $T$ as in the previous example, we get
\def\Tscale{.78}


$$\begin{array}{ccccccccc}{\rm
T}={\begin{Young}2&3&3&&&\cr &1&2&2&&\cr &&1&1&1&1\cr
\end{Young}}&\overset{\text{e}}\rightarrow&
{\rm
T}^e={\begin{Young}3&3&3&&&\cr &2&2&2&&\cr &&1&1&3&3\cr
\end{Young}}&\overset{\text{$\bullet$}}{\rightarrow}&{\rm
T}^{e\bullet}={\begin{Young}1&1&3&3&&\cr &&2&2&2&\cr &&&1&1&1\cr
\end{Young}}
&\overset{\text{$\blacklozenge$}}{\underset{\text{}}{\rightarrow}}
\end{array}$$
$$\begin{array}{ccccccccc} w=1111221332&\rightarrow&&\sigma_0
w=3311222333&&\overset{\text{$*$}}\rightarrow (\sigma_0
w)^\bullet=1112223311&\overset{\text{$\blacklozenge$}}\rightarrow\end{array}
$$
$\overset{\text{$\blacklozenge$}}{\underset{\text{}}{\rightarrow}}
\varrho(T)={\rm
T}^{e\bullet\blacklozenge}={\begin{Young}5&&\cr 4&&\cr 2&3&\cr 1&2&3\cr &1&2\cr &&1\cr
\end{Young}}.$

\noindent $\overset{\text{$\blacklozenge$}}\rightarrow(\sigma_0
w)^{\bullet\blacklozenge}=1231231245\,\, \text{\scriptsize column
word of }\varrho(T)=\varrho^{BSS}(T).$

On the left hand side one has the biwords associated to $T$ and $T^e$ and on the right hand side the reordered biwords giving the companions $G$, Example \ref{ex:companion}, and $G^{\rm a}$ respectively
$$\begin{array}{cccccccc}W^{\lambda/\mu}&=\left(\begin{smallmatrix}1&1&1&1&2&2&2&3&3&3\cr
1&1&1&1&2&2&1&3&3&2
\end{smallmatrix}\right)&\rightarrow&
W^\nu&=\left(\begin{smallmatrix}2&1&1&1&1&3&2&2&3&3\cr
1&1&1&1&1&2&2&2&3&3
\end{smallmatrix}\right)\\
&\downarrow_{\sigma_0}&&&\downarrow_{\theta_0}
\\
&\left(\begin{smallmatrix}1&1&1&1&2&2&2&3&3&3\cr
3&3&1&1&2&2&2&3&3&3
\end{smallmatrix}\right)&\rightarrow&
&\left(\begin{smallmatrix}1&1&2&2&2&3&3&3&1&1\cr
1&1&2&2&2&3&3&3&3&3
\end{smallmatrix}\right).
\end{array}$$
The companions of $T^e\in \rm LR(\mu,\nu^\bullet,\lambda)$, $T^{e\bullet}\in LR(\lambda,\nu,\mu)$ and $T^{e\bullet\blacklozenge}\in \rm LR(\mu^t,\nu^t,\lambda^t)$ are respectively
$$\begin{array}{cccccc}G^{\rm a}=\begin{Young}&1&1&3&3&3\cr &&&2&2&2\cr &&&&1&1\cr
\end{Young},
&G^{\rm a\bullet}=G^{\bullet\rm n}=G^E=\begin{Young}3&3&&&&\cr 2&2&2&&&\cr 1&1&1&3&3&\cr
\end{Young},&
G^{E\blacklozenge}=
\begin{Young}&&\cr6&&\cr 5&&\cr 3&4&\cr 2&3&4\cr 1&2&3\cr
\end{Young},\,\mu^\bullet=012
\end{array}
$$

\end{example}

\section{The $\mathbb Z_2\times D_3$--symmetries and the subgroup $\mathcal{H}$  of KTW puzzle  dualities and rotations}\label{sec:puzzle}
\subsection{Knutson-Tao-Woodward puzzles and Tao's bijection}
\label{sec:markup}
 A KTW puzzle of size $n$ ~\cite{knutson} is a tiling of an equilateral
triangle of side length $n$ with three kind of puzzle pieces:
$(a)$ unit equilateral triangles with all edges labeled $1$ (here also
represented in blue colour);
$(b)$ unit equilateral triangles with all edges labeled $0$ (here also
represented in pink colour); and
$(c)$ unit rhombi (two equilateral triangles joined together) with the
two edges, clockwise, of acute vertices labeled $0$, and the other two
labeled $1$,

\medskip

\begin{center}
\includegraphics[scale=0.7,trim = 0cm 16.5cm 0cm 0cm,clip]{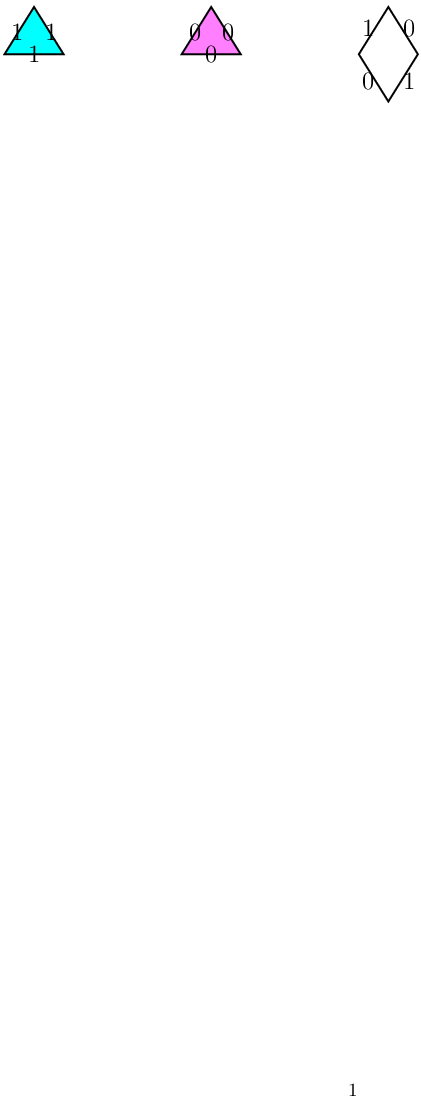}
\end{center}

\noindent such that  whenever two pieces share an edge, the labels
on the edge must agree. Puzzle pieces may be rotated in any
orientation but  rhombi can not be reflected.
The boundary data of the KTW puzzle is the partition triple $(\mu,\nu,\lambda)$ where the
partitions $\mu$, $\nu$ and  $\lambda$ appear clockwise, starting in
the lower-left corner,  as $01$--words. The partitions $\mu$, $ \nu$ and
$\lambda$ as $01$--words have exactly $d$ $1$'s and $n-d$ $0$'s. This
means that the blue unitary triangles constitute a triangle of size
$d$, the {\em $d$-triangle}, and the pink unitary triangles a triangle of size $n-d$, the {\em $(n-d)$-triangle}.
For instance, the following is a puzzle with $n=5$, $d=3$ and  boundary $\mu=01011=(100),\,\nu=01101=(110)$ and
$\lambda=10101=(210)$, read clockwise starting in the lower-left corner.

\begin{equation}\label{puzzle0}
\includegraphics[scale=0.4,trim = 0cm 18cm 0cm 0cm,clip]{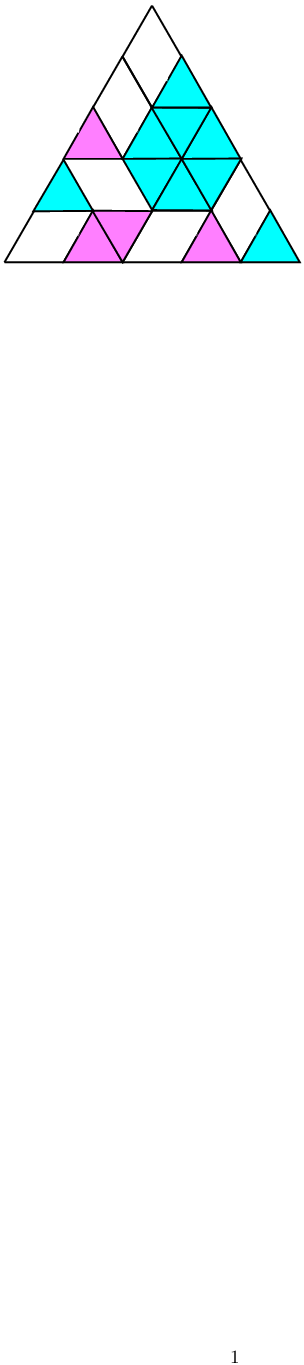}
\end{equation}

The number of puzzles with $\mu$,
$\nu$ and $\lambda$ appearing clockwise as $01$--strings along the
boundary is equal to $c_{\mu\;\nu\;\lambda}$ \cite{knutson}.
KTW puzzles  are in bijection with LR tableaux \cite{knutson}. We use \emph{Tao's bijection "without words" } in~\cite{vakil}  also used in~\cite[Figure 9]{mosaic}.
Tao's bijection defines a one-to-one correspondence between puzzles of size $n$ and boundary $(\mu,\nu,\lambda)$ as $01$ strings with $d$ 1's and $n-d$'s $0$'s, and LR tableaux with boundary $(\lambda,\mu,\nu)$,  inside a rectangle $D$ of size $d\times (n-d)$. From now on we just write puzzle to mean KTW puzzle, as no other puzzles will be considered.
To get the LR tableau filling from a puzzle, as  illustrated in Example \ref{ex:tao}, we  follow Pechenik's wording \cite{Pe16}. We construct disjoint trails of puzzle pieces, one for each $1$ along the bottom side, the $\lambda$-edge. Then these trails will be read to produce the row fillings  of the LR tableau. We think of the puzzle pieces as rooms and the $1$-labeled edges as doorways. We enter through one of the doors on the $\lambda$-edge. Whenever we enter a room, we leave it by a different door edge. We traverse right leaned rhombi from bottom to to top and left leaned rhombi from left to right. When we enter the base of a triangle, we exit through the door on our right, and when we enter the lower left door of a upsidedown triangle, we exit through the door on our left. Thus we will be always moving northeast and eventually we exit through a door on the $\nu$-side. The recording of the rooms together with the filling along this walk gives the track of the initial door on the $\lambda$-side. Reading the filling of the track of each basement door gives the row filling of each row in the LR tableau.
\begin{example}\label{ex:tao}  Tao's bijection on the puzzle below, with $n=20$, $d=4$ and  boundary $(\mu=(10,7,3,2),\nu=8522,\lambda=(11,8,5,1))$, gives, on the right,  the LR tableau $T$ with boundary $(\lambda=(11,8,5,1);$ $\mu=(10,7,3,2);\nu=8522)$,  inside the rectangle $D$ of size $4\times 16$. From bottom to top,  there are exactly $\ell(\mu)=4$ left leaned rhombi corridors, SE to NW, of lengths $\mu_1$, $\mu_2$, $\mu_3$ and $\mu_4$ respectively. Those corridors are filled in with $\mu_i$, $i$'s, respectively.
 In our example, the filling of the track in each  base  door, left to right,  gives the words $1122344$, $112233$, $111222$, and $111$.

\medskip
$$\begin{array}{ccccc}
\colorlet{lightgray}{black!10}\colorlet{lightgray1}{black!5}\colorlet{myred}{magenta!65}
\begin{tikzpicture}[scale=0.7]

\draw [color=myred, fill] (0,0) -- (0.5,0) -- (1.75,2.17) -- (1.49, 2.59)
-- cycle;
\draw [color=myred, fill] (1,0) -- (3,0) -- (3.74,1.3)  -- (3.23,2.17) --
(2.25,2.17) -- cycle;
\draw [color=myred, fill] (3.5,0) -- (5,0) -- (5.5,0.87)  -- (5.25,1.27)
-- (4.25,1.28) -- cycle;
\draw [color=myred, fill] (5.5,0) -- (7,0) -- (7.5,0.86)  -- (6,0.86) --
cycle;
\draw [color=myred, fill] (7.5,0) -- (10,0) -- (9.5,0.86)  -- (8,0.86) --
cycle;
\draw [color=myred, fill] (1.75,3.02) -- (2.01,2.58)  --  (3,2.58) --
(3.24,3.01) -- (2.51,4.33) -- cycle;
\draw [color=myred, fill] (2.75,4.75) -- (3.49,3.49) -- (4.5,3.49) --
(5,4.35) -- (3.75,6.47) -- cycle;
\draw [color=myred, fill] (4,6.92) -- (5.26,4.75) -- (5.75,4.75) --
(4.25,7.35) -- cycle;
\draw [color=myred, fill] (4.5,7.76) -- (6.02,5.18) -- (7,5.18) --
(5,8.65) -- cycle;
\draw [color=myred, fill] (6.25,4.73) -- (7.25,4.73)  -- (7.98,3.48) --
(6.98,3.48) -- cycle;
\draw [color=myred, fill] (8.25,3.02) -- (8.98,1.75) -- (7.5,1.75) --
(7,2.58) -- (7.25,3.02) -- cycle;
\draw [color=myred, fill] (5.75,1.28) -- (7.25,1.28)  -- (6.73,2.16) --
(5.75,2.16) -- (5.5,1.72) -- cycle;
\draw [color=myred, fill] (5,1.73) -- (4,1.73)  -- (3.5,2.58) --
(3.75,3.01) -- (4.75,3.01) -- (5.25,2.16) -- cycle;
\draw [color=myred, fill] (6.5,2.58) -- (5.5,2.58) -- (5,3.48) --
(5.5,4.33) -- (6,4.33) -- (6.74,3.01) -- cycle;

\draw [color=lightgray1] (9.5,0.86)-- (0.5,0.87);
\draw [color=lightgray1] (9,1.73)-- (1.01,1.72);
\draw [color=lightgray1] (8.51,2.58)-- (1.5,2.6);
\draw [color=lightgray1] (8,3.46)-- (1.99,3.48);
\draw [color=lightgray1] (7.5,4.34)-- (2.49,4.35);
\draw [color=lightgray1] (7.01,5.18)-- (3.01,5.18);
\draw [color=lightgray1] (6.51,6.05)-- (3.51,6.04);
\draw [color=lightgray1] (5.51,7.78)-- (4.51,7.77);
\draw [color=lightgray1] (5.51,7.78)-- (1,0);
\draw [color=lightgray1] (3,0)-- (6.51,6.05);

\draw [color=lightgray1] (0.25,0.43) -- (9.75,0.43);
\draw [color=lightgray1] (0.75,1.3) -- (9.25,1.3);
\draw [color=lightgray1] (1.25,2.15) -- (8.75,2.15);
\draw [color=lightgray1] (1.75,3.03) -- (8.25,3.03);
\draw [color=lightgray1] (2.25,3.9) -- (7.75,3.9);
\draw [color=lightgray1] (2.75,4.76) -- (7.25,4.76);
\draw [color=lightgray1] (3.25,5.63) -- (6.75,5.63);
\draw [color=lightgray1] (3.75,6.5) -- (6.25,6.5);
\draw [color=lightgray1] (4.25,7.36) -- (5.75,7.36);
\draw [color=lightgray1] (4.75,8.23) -- (5.25,8.23);

\draw [color=lightgray1] (3.97,6.98)-- (5.97,6.98);
\draw [color=lightgray1] (7.01,5.18)-- (4,0);
\draw [color=lightgray1] (5,0)-- (7.5,4.34);
\draw [color=lightgray1] (6,0)-- (8,3.46);
\draw [color=lightgray1] (8.51,2.58)-- (7,0);
\draw [color=lightgray1] (5.97,6.98)-- (2,0);
\draw [color=lightgray1] (9,1.73)-- (8,0);
\draw [color=lightgray1] (9,0)-- (9.5,0.86);
\draw [color=lightgray1] (0.5,0.87)-- (1,0);
\draw [color=lightgray1] (2,0)-- (1.01,1.72);
\draw [color=lightgray1] (1.5,2.6)-- (3,0);
\draw [color=lightgray1] (4,0)-- (1.99,3.48);
\draw [color=lightgray1] (5,0)-- (2.49,4.35);
\draw [color=lightgray1] (6,0)-- (3.01,5.18);
\draw [color=lightgray1] (7,0)-- (3.51,6.04);
\draw [color=lightgray1] (8,0)-- (3.97,6.98);
\draw [color=lightgray1] (9,0)-- (4.51,7.77);

\draw [color=lightgray1] (0.5,0)-- (0.25,0.43);
\draw [color=lightgray1] (1.5,0)-- (0.75,1.3);
\draw [color=lightgray1] (2.5,0)-- (1.25,2.17);
\draw [color=lightgray1] (3.5,0)-- (1.75,3.03);
\draw [color=lightgray1] (4.5,0)-- (2.25,3.9);
\draw [color=lightgray1] (5.5,0)-- (2.75,4.76);
\draw [color=lightgray1] (6.5,0)-- (3.25,5.63);
\draw [color=lightgray1] (7.5,0)-- (3.75,6.5);
\draw [color=lightgray1] (8.5,0)-- (4.25,7.36);
\draw [color=lightgray1] (9.5,0)-- (4.75,8.23);
\draw [color=lightgray1] (9.5,0)-- (9.75,0.44);
\draw [color=lightgray1] (8.5,0)-- (9.25,1.31);
\draw [color=lightgray1] (7.5,0)-- (8.76,2.14);
\draw [color=lightgray1] (6.5,0)-- (8.25,3.03);
\draw [color=lightgray1] (5.5,0)-- (7.75,3.9);
\draw [color=lightgray1] (4.5,0)-- (7.26,4.74);
\draw [color=lightgray1] (3.5,0)-- (6.75,5.62);
\draw [color=lightgray1] (2.5,0)-- (6.25,6.5);
\draw [color=lightgray1] (1.5,0)-- (5.75,7.37);
\draw [color=lightgray1] (0.5,0)-- (5.25,8.22);
\draw [color=lightgray1] (1,0)-- (2.51,2.61);
\draw [color=lightgray1] (2.51,2.61)-- (1.99,3.48);
\draw [color=lightgray1] (0.5,0)-- (1.98,2.58);
\draw [color=lightgray1] (1.98,2.58)-- (1.75,3.03);

\fill[color=cyan] (7.5,0.87) -- (7.99,0.87) -- (7.49,1.73) -- (7.25,1.3)--
cycle;
\fill[color=cyan] (5.5,0.87) -- (6,0.87) -- (5.75,1.28) -- cycle;
\fill[color=cyan] (5.26,1.28) -- (5.5,1.73) -- (5,1.73) -- cycle;
\fill[color=cyan] (5.75,2.16) -- (5.25,2.17) -- (5.505,2.58) -- cycle;
\fill[color=cyan] (6.5,2.58) -- (6.74,2.18) -- (6.98,2.58) -- cycle;
\fill[color=cyan] (6.73,3.05) -- (7.24,3.05) -- (6.99,3.48) -- cycle;
\fill[color=cyan] (6,4.34) -- (6.25,4.76) -- (6.01,5.18) --
(5.75,4.76)--cycle;
\fill[color=cyan] (4.25,1.3) -- (3.75,1.3) -- (4,1.72) -- cycle;
\fill[color=cyan] (3.24,2.19) -- (3.5,2.59) -- (3,2.59) -- cycle;
\fill[color=cyan] (2.25,2.17) -- (2,2.58) -- (1.75,2.17) -- cycle;
\fill[color=cyan] (3.74,3.03) -- (3.49,3.47) -- (3.24,3.03) -- cycle;
\fill[color=cyan] (4.76,3.02) -- (4.5,3.47) -- (5,3.47) -- cycle;
\fill[color=cyan] (5.49,4.34) -- (5.24,4.77) -- (5,4.34) -- cycle;

\draw[color=black] (0,0) -- (10,0) -- (5,8.66) -- cycle;

\draw (0.5,0)--(1.75,2.16)--(1.5,2.6);
\draw (1,0)--(2.25,2.16)--(1.75,3.03);
\draw (3,0)--(3.75,1.3)--(3,2.58)--(3.24,3.03)--(2.5,4.35);
\draw (3.5,0)--(4.25,1.3)--(3.5,2.6)--(3.75,3.03)--(2.75,4.76);
\draw
(5,0)--(5.5,.87)--(5,1.73)--(5.25,2.17)--(4.5,3.46)--(5,4.34)--(3.75,6.5);
\draw
(5.5,0)--(6,.87)--(5.5,1.73)--(5.75,2.17)--(5,3.46)--(5.5,4.34)--(4,6.93);
\draw (7,0)--(7.5,0.86)--(6.5,2.58)--(6.75,3.03)--(4.25,7.36);
\draw (7.5,0)--(8,0.86)--(7,2.58)--(7.25,3.03)--(4.5,7.77);

\node at (5.25,-.25) {\tiny 1};\node at (7.25,-.25) {\tiny 1};\node at
(.75,-.25) {\tiny 1};\node at (3.25,-.25) {\tiny 1};

\draw
(9.5,.86)--(5.5,.86)--(5.25,1.28)--(3.75,1.3)--(3.25,2.18)--(1.75,2.16);
\draw (9.25,1.28)--(5.75,1.28)--(5.5,1.73)--(3.99,1.73)--(3.5,2.6)--(2,2.58);
\draw (9,1.73)--(7.5,1.73);
\draw (6.75,2.17)--(5.75,2.17);
\draw (6.98,2.58)--(5.5,2.58);
\draw (5,3.48)--(3.5,3.48);
\draw (4.76,3.02)--(3.24,3.02);
\draw (8.25,3.03)--(6.73,3.03);
\draw (8,3.48)--(7,3.48);
\draw (7.26,4.74)--(5.24,4.74);
\draw (7,5.18)--(6,5.18);
\draw (6,4.34)--(5,4.34);
\draw (6.26,3.89) -- (6.51,4.33);
\draw (6.51,3.46) -- (6.75,3.88);

\node at (8.15,1.05) {\tiny 1};\node at (8.65,1.05) {\tiny 1};\node at
(9.15,1.05) {\tiny 1};\node at (9.55,1.1) {\tiny 1};
\node at (8.9,1.5) {\tiny 2};\node at (8.4,1.5) {\tiny 2};\node at
(7.9,1.5) {\tiny 2};\node at (9.3,1.55) {\tiny 1};
\node at (7.15,1.05) {\tiny 1};\node at (6.65,1.05) {\tiny 1};\node at
(6.15,1.05) {\tiny 1};\node at (4.9,1.5) {\tiny 1};
\node at (4.4,1.5) {\tiny 1};\node at (2.9,2.4) {\tiny 1};\node at
(2.4,2.4) {\tiny 1};\node at (1.5,2.9) {\tiny 1};
\node at (6.35,2.4) {\tiny 2};\node at (5.9,2.4) {\tiny 2};\node at
(4.4,3.25) {\tiny 2};\node at (3.9,3.25) {\tiny 2};
\node at (2.5,4.65) {\tiny 1};\node at (5.6,4.55) {\tiny 3};\node at
(6.4,4.95) {\tiny 4};\node at (6.9,4.95) {\tiny 4};
\node at (7.3,5.05) {\tiny 1};\node at (7.4,3.25) {\tiny 3};\node at
(7.9,3.25) {\tiny 3};\node at (8.25,3.35) {\tiny 1};
\node at (3.75,6.8) {\tiny 1};\node at (4.25,7.7) {\tiny 1};

\draw (8.51,1.73)--(8.99,0.86);\draw (8,1.73)--(8.5,0.86);\draw
(6.75,1.28)--(7,0.86);\draw (6.25,1.28)--(6.49,0.86);\draw
(4.75,1.3)--(4.49,1.73);\draw (2.51,2.59)--(2.74,2.17);
\draw (4,3.47)--(4.25,3.03);\draw (6,2.58)--(6.24,2.18);\draw
(7.5,3.47)--(7.75,3.04);
\draw (6.5,5.18)--(6.76,4.75);

\draw (7.75,0.44)--(7.26,0.44);\draw (5.75,0.44)--(5.25,0.44);\draw
(3.75,0.44)--(3.25,0.44);\draw (0.75,0.44)--(1.25,0.44);

\draw (3.5,0.86)--(4,0.86);\draw (1.5,0.86)--(1,0.86);\draw
(1.75,1.29)--(1.25,1.29);\draw (2,1.73)--(1.49,1.73);\draw
(3.5,1.73)--(3.76,2.15);\draw (2.99,4.34)--(2.75,3.9);\draw
(3.26,3.91)--(3,3.48);
\draw (5,2.58)--(5.25,3.03);

\draw (7.25,1.28)--(7.75,1.28)--(7.49,0.87);\draw
(7.25,2.15)--(7,1.75);\draw (4.75,3.9)--(5.25,3.9);
\draw (4.75,7.36)--(4.5,6.94);\draw (4.75,6.48)--(5,6.9);
\draw (4.99,6.06)--(5.25,6.48);\draw (5.47,6.09)--(5.22,5.67);
\draw (5.75,5.61)--(5.49,5.2);\draw (4,6.05)--(4.25,6.5);
\draw (4.5,6.06)--(4.25,5.63);\draw (4.75,5.62)--(4.5,5.18);
\draw (4.73,4.8)--(4.97,5.24);

\end{tikzpicture}&&T=\begin{Young}&1&1&2&2&3&4&4&&&&&&&&\cr &&&&&1&1&2&2&3&3&&&&&\cr &&&&&&&&1&1&1&2&2&2&&\cr
&&&&&&&&&&&1&1&1&&\cr\end{Young}\\
(\mu,\nu,\lambda)&&(\lambda,\mu,\nu).\end{array}$$


\end{example}

\subsection{The  puzzle $H$--subgroup of symmetries }
\label{sec:H}
The reflection of a (upright or left or right leaned) rhombus  in a puzzle swaps 0 and 1 labels  which propagates  in a unique way through the puzzle and eventually swaps  0 and 1 labels in the boundary. To get a puzzle again  we have to make one of three diagonal reflections of the  puzzle  (a diagonal is a  line joining a vertex to the midpoint of the opposite side). Recalling Section \ref{sec:pre}, with respect to operations on 01 words and partitions, this  procedure defines, according to the chosen diagonal,  an involution on puzzles. The vertical reflection while swapping all $01$ labels of a puzzle of boundary $(\mu,\nu,\lambda)$, defines the {\em dual puzzle} of boundary $(\nu^t,\mu^t,\lambda^t)$.  This procedure is an involution, denoted  $\spadesuit$, on puzzles. It 
exhibits the symmetry
 $c_{\mu \nu \lambda}=c_{\nu^t \mu^t  \lambda^t}$ on puzzles.
   Performing the two remaining diagonal left and right reflections,  respectively,  with $0 1$  swapping, these procedures define the involutions $\blacklozenge$ (left) and $\clubsuit$ (right)  respectively. They exhibit the symmetries (puzzle dualities) $c_{\mu \nu \lambda}=c_{\lambda^t \nu^t \mu^t}$ and  $c_{\mu,\nu,\lambda}=c_{\mu^t \lambda^t \nu^t}$ on puzzles, respectively. The duals of the puzzles counted by $c_{\mu \nu \lambda}$ are
 exactly those counted by $c_{\nu^t \mu^t  \lambda^t}$ or $c_{\lambda^t \nu^t \mu^t}$ or $c_{\mu^t \lambda^t \nu^t}$.
     The involutions $\blacklozenge$, $\spadesuit$ and $\clubsuit$  are the unique involutions  which swap
pink (label $0$) and blue (label $1$) colours in a puzzle such that
the resulting tiled triangle is still a puzzle, equivalently, which swap the
blue $d$-triangle  with the pink $(n-d)$-triangle.
As an apart, one observes that the symmetry with respect to complement operation  on partitions, equivalently, reversing the $01$ word, Section \ref{sec:pre}, is indeed not  easy to exhibit due to the  restriction of the reflection operation on rhombi. The same difficulty  has been already observed in \cite[Remarks]{bz} in the case of BZ triangles. The same difficulty applies to the transposition of partitions, reversing and 01 swapping a 01 string.

The action of $H=\langle\tau\varsigma_1,\tau\varsigma_2\rangle=\langle\tau\varsigma_1,\tau\varsigma\rangle=
\langle\tau\varsigma_2,\tau\varsigma\rangle$ on puzzles is defined via any two of the involutions $\blacklozenge, \spadesuit, \clubsuit$.
The compositions $\clubsuit\blacklozenge=\blacklozenge\spadesuit=\spadesuit\clubsuit$  and $\blacklozenge\clubsuit=\spadesuit\blacklozenge=\clubsuit\spadesuit$   rotate  the puzzle  clockwise $2\pi/3$ and $4\pi/3$ radians about the center respectively, and
 realize the corresponding rotational (cycle) symmetries
$c_{\mu,\nu,\lambda}=c_{ \lambda,\mu,\nu}$ and $c_{\mu,\nu,\lambda}=c_{ \nu, \lambda,\mu}$.
The group of symmetries of a puzzle, generated by the action of $H$,
is  realized by the  three diagonal reflections with $01$
swapping  and  clockwise $0$, $2\pi/3$ and $4\pi/3$ radians rotations about the center.

 \begin{example}\label{ex:spade} Illustration of the action of $\tau \varsigma_1\in H$ on puzzles through the involution $\spadesuit$. The involution $\spadesuit$  on the puzzle \eqref{puzzle0}, realizes the action of $\tau\varsigma_1$ on
that puzzle, giving the puzzle

\begin{center}
\includegraphics[scale=0.4,trim = 0cm 18cm 0cm 0cm,clip]{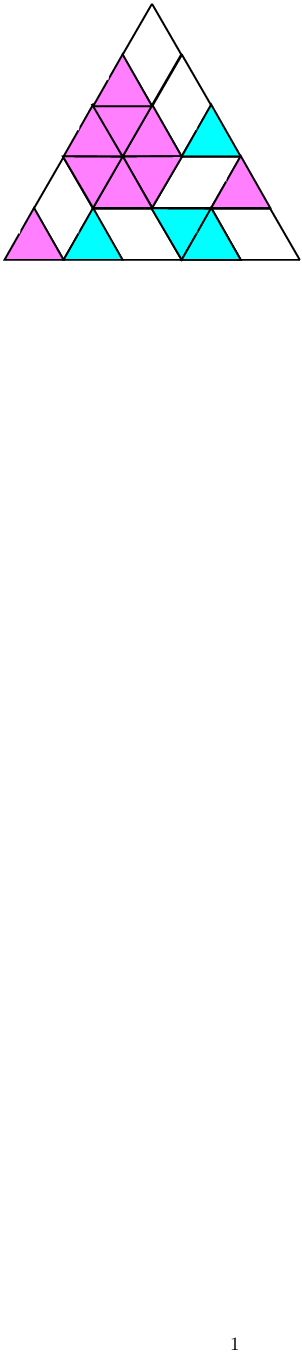}
\end{center}

\noindent whose boundary data is $(\nu^t=01001;\,\mu^t=00101; \lambda^t=01010)$.
\end{example}
  The group action of $H$ on puzzles is faithful and   we  may identify $H$ with its representation in $\mathfrak{S}_{\mathcal{LR}}$, $H\simeq \langle\spadesuit,\blacklozenge\rangle=\langle\clubsuit,\blacklozenge\rangle=\langle\spadesuit,\clubsuit\rangle$, and a possible presentation  is
\begin{equation}\label{hgroup1}
H=\{\spadesuit,\blacklozenge:\spadesuit^2=\blacklozenge^2={\bf 1}=
(\blacklozenge\spadesuit)^3\}\simeq {D}_3\simeq \mathfrak{S}_3,
\end{equation}
 the dihedral group  of an equilateral triangle.
Since through Tao's bijection,  LR tableau and puzzle boundaries  differ  by a cyclic permutation,   first rotate counterclockwise the puzzle $2\pi/3$ radians about the center, follow with the involution $\blacklozenge$, $\spadesuit$ or $\clubsuit$, and then apply Tao's bijection to obtain the involutions $\blacklozenge$, $\spadesuit$ or $\clubsuit$ on LR tableaux respectively. That is, on puzzles, $\clubsuit=\blacklozenge(\blacklozenge\clubsuit)$,  $\blacklozenge=\spadesuit(\spadesuit\blacklozenge)$ and $\spadesuit=\clubsuit(\clubsuit\spadesuit)$ are translated to LR tableaux, by Tao's bijection, to $\blacklozenge$, $\spadesuit$ or $\clubsuit$, respectively.
 In particular, we show that the bijection $\blacklozenge:\rm LR(\mu,\nu,\lambda)\rightarrow LR(\lambda^t,\nu^t,\mu^t)$, defined in Section \ref{blacklozenge}, exhibiting the identity $c_{\mu \nu \lambda}=c_{\lambda^t \nu^t \mu^t}$ on LR tableaux  is translated to puzzles in this fashion by exhibiting the identity $c_{\mu,\nu,\lambda}=c_{\mu^t,\lambda^t,\nu^t}$. We also  define  bijections $\spadesuit: \rm LR(\mu,\nu,\lambda)\rightarrow \rm LR(\nu,\mu,\lambda)$,  $\clubsuit: \rm LR(\mu,\nu,\lambda)\rightarrow LR(\mu,\lambda,\nu)$,  exhibiting the identities $c_{\mu,\nu,\lambda}=c_{\mu^t,\nu^t, \lambda^t }$, $c_{\mu,\nu,\lambda}=c_{\mu^t \lambda^t \nu^t}$  respectively,
 which are also translated to puzzles in the same fashion by exhibiting the identities $c_{\mu,\nu,\lambda}=c_{\lambda^t,\mu^t,\nu^t }$, $c_{\mu,\nu,\lambda}=c_{\nu^t,\mu^t, \lambda^t }$ .

 \begin{example} \label{ex:taoblacklozenge} If we  first rotate clockwise the puzzle $4\pi/3$ radians about the center and then apply $\blacklozenge$, that is, we perform {$\blacklozenge(\blacklozenge\clubsuit)$}, it is not difficult to see that Tao's bijection translates the involution $\blacklozenge$ (the left diagonal reflection with $01$  swapping) on the $4\pi/3$ radians clockwise rotated puzzle  to the involution $\blacklozenge$ on LR tableaux as we have defined in subsection \ref{blacklozenge}.  Consider the puzzle in Example \ref{ex:tao} together with the left leaned rhombi filling. After rotating the puzzle clockwise $\frac{4}{3}\pi$ radians about the center,  to obtain a puzzle of boundary $(\nu,\lambda,\mu)$,   apply the involution $\blacklozenge$ (left diagonal reflection while $01$ label swapping) to give a puzzle of boundary $(\mu^t, \lambda^t,\nu^t)$. Then
 replace East-West: (i) the ten 1's word, inside the left leaned rhombi, with the word $123456789\,(10)$; (ii) the seven 2's word with the word $1234567$; (iii) the three 3's word by $123$; and (iv)  the two $4$'s word with $12$.  The puzzle-LR tableau pair with boundaries { $(\mu,\nu,\lambda)$ and $(\lambda,\mu,\nu)$  } respectively, in  Example \ref{ex:tao}, is transformed into a puzzle-LR tableau pair of boundaries   $(\mu^t,\lambda^t, \nu^t)$ and $(\nu^t, \mu^t,\lambda^t)$, where $\mu^t=4^2 3 2^4 1^3$ and $\ell(\mu^t)=10$, respectively. See figure below

$$\blacklozenge(\blacklozenge\clubsuit)\begin{array}{cccccc}
\colorlet{lightgray}{black!10}\colorlet{lightgray1}{black!5}\colorlet{myred}{magenta!65}
\begin{tikzpicture}[scale=0.8]
\fill[color=cyan] (2.5,4.35) -- (4,4.35) -- (4.25,4.76) -- (3.75,5.63) --
(3.25,5.63) -- cycle;
\fill[color=cyan] (4.75,4.76) -- (5.75,4.76) -- (6.5,6.04) --  (5.5,7.77)
-- (4.25,5.63) -- cycle;
\fill[color=cyan] (3.5,6.04) -- (4,6.04) -- (5.25,8.22) -- (5,8.66);
\fill[color=cyan] (1.25,2.14) -- (3.75,2.14) -- (4.25,3.03) -- (3.75,3.9)
-- (2.25,3.9) -- cycle;
\fill[color=cyan] (4.75,3.03) -- (5.75,3.03) -- (6,3.48) -- (5.5,4.35) --
(4.5,4.35) -- (4.25,3.9) -- cycle;
\fill[color=cyan] (6.5,3.48) -- (7,3.48) -- (7.5,4.35) -- (6.75,5.63) --
(6,4.35);
\fill[color=cyan] (0.75,1.3) -- (3.75,1.3) -- (3.5,1.72) -- (1,1.72) --
cycle;
\fill[color=cyan] (4.25,1.3) -- (5.75,1.3) -- (6,1.72) -- (5.5,2.61) --
(4.5,2.61) -- (4,1.72) -- cycle;
\fill[color=cyan] (6.5,1.72) -- (7.5,1.72) -- (6.75,3.03) -- (6.25,3.03)
-- (6,2.61) -- cycle;
\fill[color=cyan] (8,1.72) -- (8.5,2.61) -- (7.75,3.9) -- (7.25,3.03) --
cycle;
\fill[color=cyan] (0,0) -- (4,0) -- (3.5,0.87) -- (0.5,0.87) -- cycle;
\fill[color=cyan] (4.5,0) -- (6,0) -- (5.5,0.87) -- (4,0.87) -- cycle;
\fill[color=cyan] (6.5,0) -- (8,0) -- (7.25,1.3) -- (6.25,1.3) -- (6,0.87)
-- cycle;
\fill[color=cyan] (9,0) -- (10,0) -- (8.75,2.17) -- (8.25,1.3) -- cycle;
\draw [color=lightgray1] (9.5,0.86)-- (0.5,0.87);
\draw [color=lightgray1] (9,1.73)-- (1.01,1.72);
\draw [color=lightgray1] (8.51,2.58)-- (1.5,2.6);
\draw [color=lightgray1] (8,3.46)-- (1.99,3.48);
\draw [color=lightgray1] (7.5,4.34)-- (2.49,4.35);
\draw [color=lightgray1] (7.01,5.18)-- (3.01,5.18);
\draw [color=lightgray1] (6.51,6.05)-- (3.51,6.04);
\draw [color=lightgray1] (5.51,7.78)-- (4.51,7.77);
\draw [color=lightgray1] (5.51,7.78)-- (1,0);
\draw [color=lightgray1] (3,0)-- (6.51,6.05);

\draw [color=lightgray1] (0.25,0.43) -- (9.75,0.43);
\draw [color=lightgray1] (0.75,1.3) -- (9.25,1.3);
\draw [color=lightgray1] (1.25,2.15) -- (8.75,2.15);
\draw [color=lightgray1] (1.75,3.03) -- (8.25,3.03);
\draw [color=lightgray1] (2.25,3.9) -- (7.75,3.9);
\draw [color=lightgray1] (2.75,4.76) -- (7.25,4.76);
\draw [color=lightgray1] (3.25,5.63) -- (6.75,5.63);
\draw [color=lightgray1] (3.75,6.5) -- (6.25,6.5);
\draw [color=lightgray1] (4.25,7.36) -- (5.75,7.36);
\draw [color=lightgray1] (4.75,8.23) -- (5.25,8.23);

\draw [color=lightgray1] (3.97,6.98)-- (5.97,6.98);
\draw [color=lightgray1] (7.01,5.18)-- (4,0);
\draw [color=lightgray1] (5,0)-- (7.5,4.34);
\draw [color=lightgray1] (6,0)-- (8,3.46);
\draw [color=lightgray1] (8.51,2.58)-- (7,0);
\draw [color=lightgray1] (5.97,6.98)-- (2,0);
\draw [color=lightgray1] (9,1.73)-- (8,0);
\draw [color=lightgray1] (9,0)-- (9.5,0.86);
\draw [color=lightgray1] (0.5,0.87)-- (1,0);
\draw [color=lightgray1] (2,0)-- (1.01,1.72);
\draw [color=lightgray1] (1.5,2.6)-- (3,0);
\draw [color=lightgray1] (4,0)-- (1.99,3.48);
\draw [color=lightgray1] (5,0)-- (2.49,4.35);
\draw [color=lightgray1] (6,0)-- (3.01,5.18);
\draw [color=lightgray1] (7,0)-- (3.51,6.04);
\draw [color=lightgray1] (8,0)-- (3.97,6.98);
\draw [color=lightgray1] (9,0)-- (4.51,7.77);

\draw [color=lightgray1] (0.5,0)-- (0.25,0.43);
\draw [color=lightgray1] (1.5,0)-- (0.75,1.3);
\draw [color=lightgray1] (2.5,0)-- (1.25,2.17);
\draw [color=lightgray1] (3.5,0)-- (1.75,3.03);
\draw [color=lightgray1] (4.5,0)-- (2.25,3.9);
\draw [color=lightgray1] (5.5,0)-- (2.75,4.76);
\draw [color=lightgray1] (6.5,0)-- (3.25,5.63);
\draw [color=lightgray1] (7.5,0)-- (3.75,6.5);
\draw [color=lightgray1] (8.5,0)-- (4.25,7.36);
\draw [color=lightgray1] (9.5,0)-- (4.75,8.23);
\draw [color=lightgray1] (9.5,0)-- (9.75,0.44);
\draw [color=lightgray1] (8.5,0)-- (9.25,1.31);
\draw [color=lightgray1] (7.5,0)-- (8.76,2.14);
\draw [color=lightgray1] (6.5,0)-- (8.25,3.03);
\draw [color=lightgray1] (5.5,0)-- (7.75,3.9);
\draw [color=lightgray1] (4.5,0)-- (7.26,4.74);
\draw [color=lightgray1] (3.5,0)-- (6.75,5.62);
\draw [color=lightgray1] (2.5,0)-- (6.25,6.5);
\draw [color=lightgray1] (1.5,0)-- (5.75,7.37);
\draw [color=lightgray1] (0.5,0)-- (5.25,8.22);
\draw [color=lightgray1] (1,0)-- (2.51,2.61);
\draw [color=lightgray1] (2.51,2.61)-- (1.99,3.48);
\draw [color=lightgray1] (0.5,0)-- (1.98,2.58);
\draw [color=lightgray1] (1.98,2.58)-- (1.75,3.03);

\draw[color=black] (0,0) -- (10,0) -- (5,8.66) -- cycle;


\fill[color=myred] (3.75,5.63) -- (4.25,5.63) -- (4,6.04) -- cycle;
\draw (3.25,5.63) -- (4.25,5.63) -- (5.5,7.77);
\draw (3.5,6.04) -- (4,6.04) -- (5.25,8.22);
\node at (3.2,5.88) {\tiny 0};
\node at (5.57,8.05) {\tiny 0};
\draw (4.75,4.76) -- (4,6.04);
\draw (3.75,5.63) -- (4,6.04);
\draw (4.5,6.04) -- (4.25,6.5);
\draw (4.75,6.5) -- (4.5,6.94);
\draw (5,6.95) -- (4.75,7.36);
\draw (5.25,7.36) -- (5,7.77);
\draw (4,5.18) -- (4.5,5.18);

\node at (4.38,4.95) {\tiny 9};
\node at (4.12,5.38) {\tiny 10};

\fill[color=myred] (3.75,3.9) -- (4.25,3.9) -- (4,4.35) -- cycle;
\fill[color=myred] (4.25,4.76) -- (4.5,4.35) -- (4.75,4.76) -- cycle;
\fill[color=myred] (5.5,4.35) -- (6,4.35) -- (5.75,4.76) -- cycle;
\draw (4.5,4.35) -- (3.75,5.63);
\draw (2.25,3.9) -- (4.25,3.9) -- (4.5,4.35) -- (6,4.35) -- (6.75,5.63);
\draw (2.5,4.35) -- (4,4.35) -- (4.25,4.76) -- (5.75,4.76) -- (6.5,6.04);
\node at (2.22,4.15) {\tiny 0};
\node at (6.78,5.89) {\tiny 0};
\draw (2.75,3.9) -- (3,4.35);
\draw (3.25,3.9) -- (3.5,4.35);
\draw (3.75,3.9) -- (4,4.35);
\draw (4.5,4.35) -- (4.75,4.76);
\draw (5,4.35) -- (5.25,4.76);
\draw (5.5,4.35) -- (5.75,4.76);
\draw (6.5,3.48) -- (5.75,4.76);
\draw (6.25,4.76) -- (6,5.18);
\draw (6.5,5.18) -- (6.25,5.63);
\draw (6,3.48) -- (5.5,4.35);
\draw (5.75,3.9) -- (6.25,3.9);
\draw (4.25,3.03) -- (3.75,3.9);
\draw (4.75,3.03) -- (4.25,3.9);
\draw (4,3.48) -- (4.5,3.48);

\node at (4.35,3.25) {\tiny 6};
\node at (6.13,3.68) {\tiny 7};

\node at (4.12,3.68) {\tiny 7};
\node at (5.88,4.1) {\tiny 8};


\fill[color=myred] (3.5,1.72)-- (4,1.72) -- (3.75,2.14) -- cycle;
\fill[color=myred] (4.5,2.61)-- (4.75,3.03) -- (4.25,3.03) -- cycle;
\fill[color=myred] (5.5,2.61)-- (6,2.61) -- (5.75,3.03) -- cycle;
\fill[color=myred] (6.25,3.03)-- (6.5,3.48) -- (6,3.48) -- cycle;
\fill[color=myred] (6.75,3.03)-- (7.25,3.03) -- (7,3.48) -- cycle;
\draw (1,1.72) -- (4,1.72) -- (4.5,2.61) -- (6,2.61) -- (6.25,3.03) --
(7.25,3.03) -- (7.75,3.9);
\draw (1.25,2.14) -- (3.75,2.14) -- (4.25,3.03) -- (5.75,3.03) -- (6,3.48)
-- (7,3.48) -- (7.5,4.35);
\draw (1.5,1.72) -- (1.75,2.14);
\draw (2,1.72) -- (2.25,2.14);
\draw (2.5,1.72) -- (2.75,2.14);
\draw (3,1.72) -- (3.25,2.14);
\draw (3.5,1.72) -- (3.75,2.14);
\draw (4.5,2.61) -- (4.75,3.03);
\draw (5,2.61) -- (5.25,3.03);
\draw (5.5,2.61) -- (5.75,3.03);
\draw (6.25,3.03) -- (6.5,3.48);
\draw (6.75,3.03) -- (7,3.48);
\draw (4.25,2.14) -- (4,2.61);
\draw (4.5,2.61) -- (4.25,3.03);
\draw (6,2.61) -- (5.75,3.03);
\draw (6.25,3.03) -- (6,3.48);
\draw (7.5,3.48) -- (7.25,3.9);
\draw (6,1.72) -- (5.5,2.61);
\draw (6.5,1.72) -- (6,2.61);
\draw (6.25,2.14) -- (5.75,2.14);
\draw (7.25,2.14) -- (7.75,2.14);
\draw (7,2.61) -- (7.5,2.61);
\node at (0.97,1.98) {\tiny 0};
\node at (7.78,4.18) {\tiny 0};


\fill[color=myred] (3.5,0.87) -- (4,0.87) -- (4.25,1.3) -- (3.75,1.3) --
cycle;
\fill[color=myred] (5.5,0.87) -- (6,0.87) -- (5.75,1.3)  -- cycle;
\fill[color=myred] (6.25,1.3) -- (6.5,1.72) -- (6,1.72) -- cycle;
\fill[color=myred] (7.25,1.3) -- (8.25,1.3) -- (8,1.72) -- (7.5,1.72) --
cycle;
\draw (0.5,0.87) -- (6,0.87) -- (6.25,1.3) -- (8.25,1.3) -- (8.75,2.17);
\draw (0.75,1.3) -- (5.75,1.3) -- (6,1.72) -- (8,1.72) -- (8.5,2.58);
\draw (8.5,1.72) -- (8.25,2.17);
\draw (1,0.87) -- (1.25,1.3);
\draw (1.5,0.87) -- (1.75,1.3);
\draw (2,0.87) -- (2.25,1.3);
\draw (2.5,0.87) -- (2.75,1.3);
\draw (3,0.87) -- (3.25,1.3);
\draw (3.5,0.87) -- (3.75,1.3);
\draw (4,0.87) -- (4.25,1.3);
\draw (4.5,0.87) -- (4.75,1.3);
\draw (5,0.87) -- (5.25,1.3);
\draw (5.5,0.87) -- (5.75,1.3);
\draw (6.25,1.3) -- (6.5,1.72);
\draw (6.75,1.3) -- (7,1.72);
\draw (7.25,1.3) -- (7.5,1.72);
\draw (7.75,1.3) -- (8,1.72);
\node at (0.47,1.12) {\tiny 0};
\node at (8.77,2.43) {\tiny 0};

\draw (4,0) -- (3.5,0.87);
\draw (4.5,0) -- (3.5,1.72);
\draw (4.25,1.3) -- (3.75,2.14);
\node at (4.25,-0.2) {\tiny 0};
\node at (0.25,-0.2) {\tiny 1};
\node at (0.75,-0.2) {\tiny 1};
\node at (1.25,-0.2) {\tiny 1};
\node at (1.75,-0.2) {\tiny 1};
\node at (2.25,-0.2) {\tiny 1};
\node at (2.75,-0.2) {\tiny 1};
\node at (3.25,-0.2) {\tiny 1};
\node at (3.75,-0.2) {\tiny 1};
\node at (4.75,-0.2) {\tiny 1};
\node at (5.25,-0.2) {\tiny 1};
\node at (5.75,-0.2) {\tiny 1};
\node at (6.75,-0.2) {\tiny 1};
\node at (7.25,-0.2) {\tiny 1};
\node at (7.75,-0.2) {\tiny 1};
\node at (9.25,-0.2) {\tiny 1};
\node at (9.75,-0.2) {\tiny 1};

\draw (6,0) -- (5.5,0.87);
\draw (6.5,0) -- (5.75,1.3);
\node at (6.25,-0.2) {\tiny 0};

\draw (8,0) -- (7.25,1.3);
\draw (8.5,0) -- (6.75,3.03);
\draw (9,0) -- (7,3.46);
\node at (8.25,-0.2) {\tiny 0};
\node at (8.75,-0.2) {\tiny 0};

\draw (3.75,0.43) -- (4.25,0.43);
\draw (5.75,0.43) -- (6.25,0.43);
\draw (7.75,0.43) -- (8.75,0.43);
\draw (7.5,0.87) -- (8.5,0.87);

\node at (4.12,0.22) {\tiny 1};
\node at (6.12,0.22) {\tiny 1};
\node at (8.12,0.22) {\tiny 1};
\node at (8.62,0.22) {\tiny 1};
\node at (3.87,0.65) {\tiny 2};
\node at (5.87,0.65) {\tiny 2};
\node at (7.87,0.65) {\tiny 2};
\node at (8.35,0.65) {\tiny 2};

\node at (3.87,1.5) {\tiny 3};
\node at (7.63,1.07) {\tiny 3};
\node at (8.12,1.07) {\tiny 3};

\node at (6.12,1.93) {\tiny 4};
\node at (7.62,1.93) {\tiny 4};

\node at (5.89,2.35) {\tiny 5};
\node at (7.37,2.35) {\tiny 5};
\node at (7.12,2.8) {\tiny 6};

\node at (-0.1,0.25) {\tiny 1};
\node at (0.15,0.66) {\tiny 1};
\node at (0.65,1.52) {\tiny 1};
\node at (1.15,2.36) {\tiny 1};
\node at (1.40,2.82) {\tiny 1};
\node at (1.65,3.28) {\tiny 1};
\node at (1.9,3.72) {\tiny 1};
\node at (2.4,4.57) {\tiny 1};
\node at (2.65,4.98) {\tiny 1};
\node at (2.9,5.39) {\tiny 1};
\node at (3.4,6.27) {\tiny 1};
\node at (3.65,6.72) {\tiny 1};
\node at (3.9,7.18) {\tiny 1};
\node at (4.18,7.59) {\tiny 1};
\node at (4.43,8) {\tiny 1};
\node at (4.68,8.45) {\tiny 1};

\node at (10.08,0.39) {\tiny 1};
\node at (9.83,0.8) {\tiny 1};
\node at (9.58,1.21) {\tiny 1};
\node at (9.33,1.62) {\tiny 1};
\node at (9.11,2) {\tiny 1};
\node at (8.61,2.82) {\tiny 1};
\node at (8.36,3.28) {\tiny 1};
\node at (8.11,3.72) {\tiny 1};
\node at (7.61,4.57) {\tiny 1};
\node at (7.36,4.98) {\tiny 1};
\node at (7.11,5.39) {\tiny 1};
\node at (6.61,6.27) {\tiny 1};
\node at (6.36,6.72) {\tiny 1};
\node at (6.11,7.18) {\tiny 1};
\node at (5.86,7.59) {\tiny 1};
\node at (5.32,8.45) {\tiny 1};
\end{tikzpicture}&&\blacklozenge T=\begin{Young}&&&\cr 10&&&\cr 9&&&\cr 7&&&\cr 6&&&\cr 3&8&&\cr 2&7&&\cr 1&5&&\cr &4&6&\cr &2&5&\cr &1&4&\cr &&3&3\cr &&2&2\cr &&1&1\cr &&&\cr
&&&\cr\end{Young}\\
(\mu^t, \lambda^t,\nu^t) &&\;\;\;\;\;\;\;\;\; (\nu^t, \mu^t,\lambda^t)
\end{array}$$

Performing $\spadesuit(\spadesuit\blacklozenge)=\spadesuit(\blacklozenge\clubsuit)$ on the puzzle in  Example \ref{ex:tao}, or, equivalently, rotate counter clockwise the puzzle just above  by $2\pi/3$ radians about the center $\spadesuit(\spadesuit\blacklozenge)=\spadesuit\blacklozenge(\blacklozenge(\blacklozenge\clubsuit))$ to get a puzzle of boundary $( \nu^t,\lambda^t,\mu^t)$. Tao's bijection translates the action of $\spadesuit$ on puzzles  to LR tableaux.
The puzzle-LR tableau pair with boundaries { $(\mu,\nu,\lambda)$ and $(\lambda,\mu,\nu)$  } respectively, in  Example \ref{ex:tao}, is transformed into a puzzle-LR tableau pair of boundaries  { $(\lambda^t,\nu^t, \mu^t)$} and $(\mu^t, \lambda^t,\nu^t)$ respectively. See figure below, to obtain the LR filling on the right, one has to rotate the puzzle counterclockwise $2\pi/3$ radians about the center, where $\lambda^t=4 3^4 2^3 1^3$, and $\ell(\lambda^t)=11$, and ignore the old filling in gray
$$\begin{array}{cccccc}
\spadesuit(\spadesuit\blacklozenge)\;\;\;\;\;\clubsuit(\clubsuit\spadesuit)&&\\
\colorlet{lightgray}{black!10}\colorlet{lightgray1}{black!5}\colorlet{myred}{magenta!65}
\begin{tikzpicture}[scale=0.8]
\fill[color=cyan] (2.5,4.35) -- (4,4.35) -- (4.25,4.76) -- (3.75,5.63) --
(3.25,5.63) -- cycle;
\fill[color=cyan] (4.75,4.76) -- (5.75,4.76) -- (6.5,6.04) --  (5.5,7.77)
-- (4.25,5.63) -- cycle;
\fill[color=cyan] (3.5,6.04) -- (4,6.04) -- (5.25,8.22) -- (5,8.66);
\fill[color=cyan] (1.25,2.14) -- (3.75,2.14) -- (4.25,3.03) -- (3.75,3.9)
-- (2.25,3.9) -- cycle;
\fill[color=cyan] (4.75,3.03) -- (5.75,3.03) -- (6,3.48) -- (5.5,4.35) --
(4.5,4.35) -- (4.25,3.9) -- cycle;
\fill[color=cyan] (6.5,3.48) -- (7,3.48) -- (7.5,4.35) -- (6.75,5.63) --
(6,4.35);
\fill[color=cyan] (0.75,1.3) -- (3.75,1.3) -- (3.5,1.72) -- (1,1.72) --
cycle;
\fill[color=cyan] (4.25,1.3) -- (5.75,1.3) -- (6,1.72) -- (5.5,2.61) --
(4.5,2.61) -- (4,1.72) -- cycle;
\fill[color=cyan] (6.5,1.72) -- (7.5,1.72) -- (6.75,3.03) -- (6.25,3.03)
-- (6,2.61) -- cycle;
\fill[color=cyan] (8,1.72) -- (8.5,2.61) -- (7.75,3.9) -- (7.25,3.03) --
cycle;
\fill[color=cyan] (0,0) -- (4,0) -- (3.5,0.87) -- (0.5,0.87) -- cycle;
\fill[color=cyan] (4.5,0) -- (6,0) -- (5.5,0.87) -- (4,0.87) -- cycle;
\fill[color=cyan] (6.5,0) -- (8,0) -- (7.25,1.3) -- (6.25,1.3) -- (6,0.87)
-- cycle;
\fill[color=cyan] (9,0) -- (10,0) -- (8.75,2.17) -- (8.25,1.3) -- cycle;
\draw [color=lightgray1] (9.5,0.86)-- (0.5,0.87);
\draw [color=lightgray1] (9,1.73)-- (1.01,1.72);
\draw [color=lightgray1] (8.51,2.58)-- (1.5,2.6);
\draw [color=lightgray1] (8,3.46)-- (1.99,3.48);
\draw [color=lightgray1] (7.5,4.34)-- (2.49,4.35);
\draw [color=lightgray1] (7.01,5.18)-- (3.01,5.18);
\draw [color=lightgray1] (6.51,6.05)-- (3.51,6.04);
\draw [color=lightgray1] (5.51,7.78)-- (4.51,7.77);
\draw [color=lightgray1] (5.51,7.78)-- (1,0);
\draw [color=lightgray1] (3,0)-- (6.51,6.05);

\draw [color=lightgray1] (0.25,0.43) -- (9.75,0.43);
\draw [color=lightgray1] (0.75,1.3) -- (9.25,1.3);
\draw [color=lightgray1] (1.25,2.15) -- (8.75,2.15);
\draw [color=lightgray1] (1.75,3.03) -- (8.25,3.03);
\draw [color=lightgray1] (2.25,3.9) -- (7.75,3.9);
\draw [color=lightgray1] (2.75,4.76) -- (7.25,4.76);
\draw [color=lightgray1] (3.25,5.63) -- (6.75,5.63);
\draw [color=lightgray1] (3.75,6.5) -- (6.25,6.5);
\draw [color=lightgray1] (4.25,7.36) -- (5.75,7.36);
\draw [color=lightgray1] (4.75,8.23) -- (5.25,8.23);

\draw [color=lightgray1] (3.97,6.98)-- (5.97,6.98);
\draw [color=lightgray1] (7.01,5.18)-- (4,0);
\draw [color=lightgray1] (5,0)-- (7.5,4.34);
\draw [color=lightgray1] (6,0)-- (8,3.46);
\draw [color=lightgray1] (8.51,2.58)-- (7,0);
\draw [color=lightgray1] (5.97,6.98)-- (2,0);
\draw [color=lightgray1] (9,1.73)-- (8,0);
\draw [color=lightgray1] (9,0)-- (9.5,0.86);
\draw [color=lightgray1] (0.5,0.87)-- (1,0);
\draw [color=lightgray1] (2,0)-- (1.01,1.72);
\draw [color=lightgray1] (1.5,2.6)-- (3,0);
\draw [color=lightgray1] (4,0)-- (1.99,3.48);
\draw [color=lightgray1] (5,0)-- (2.49,4.35);
\draw [color=lightgray1] (6,0)-- (3.01,5.18);
\draw [color=lightgray1] (7,0)-- (3.51,6.04);
\draw [color=lightgray1] (8,0)-- (3.97,6.98);
\draw [color=lightgray1] (9,0)-- (4.51,7.77);

\draw [color=lightgray1] (0.5,0)-- (0.25,0.43);
\draw [color=lightgray1] (1.5,0)-- (0.75,1.3);
\draw [color=lightgray1] (2.5,0)-- (1.25,2.17);
\draw [color=lightgray1] (3.5,0)-- (1.75,3.03);
\draw [color=lightgray1] (4.5,0)-- (2.25,3.9);
\draw [color=lightgray1] (5.5,0)-- (2.75,4.76);
\draw [color=lightgray1] (6.5,0)-- (3.25,5.63);
\draw [color=lightgray1] (7.5,0)-- (3.75,6.5);
\draw [color=lightgray1] (8.5,0)-- (4.25,7.36);
\draw [color=lightgray1] (9.5,0)-- (4.75,8.23);
\draw [color=lightgray1] (9.5,0)-- (9.75,0.44);
\draw [color=lightgray1] (8.5,0)-- (9.25,1.31);
\draw [color=lightgray1] (7.5,0)-- (8.76,2.14);
\draw [color=lightgray1] (6.5,0)-- (8.25,3.03);
\draw [color=lightgray1] (5.5,0)-- (7.75,3.9);
\draw [color=lightgray1] (4.5,0)-- (7.26,4.74);
\draw [color=lightgray1] (3.5,0)-- (6.75,5.62);
\draw [color=lightgray1] (2.5,0)-- (6.25,6.5);
\draw [color=lightgray1] (1.5,0)-- (5.75,7.37);
\draw [color=lightgray1] (0.5,0)-- (5.25,8.22);
\draw [color=lightgray1] (1,0)-- (2.51,2.61);
\draw [color=lightgray1] (2.51,2.61)-- (1.99,3.48);
\draw [color=lightgray1] (0.5,0)-- (1.98,2.58);
\draw [color=lightgray1] (1.98,2.58)-- (1.75,3.03);

\draw[color=black] (0,0) -- (10,0) -- (5,8.66) -- cycle;


\fill[color=myred] (3.75,5.63) -- (4.25,5.63) -- (4,6.04) -- cycle;
\draw (3.25,5.63) -- (4.25,5.63) -- (5.5,7.77);
\draw (3.5,6.04) -- (4,6.04) -- (5.25,8.22);
\node at (3.2,5.88) {\tiny 0};
\node at (5.57,8.05) {\tiny 0};
\draw (4.75,4.76) -- (4,6.04);
\draw (3.75,5.63) -- (4,6.04);
\draw (4.5,6.04) -- (4.25,6.5);
\draw (4.75,6.5) -- (4.5,6.94);
\draw (5,6.95) -- (4.75,7.36);
\draw (5.25,7.36) -- (5,7.77);
\draw (4,5.18) -- (4.5,5.18);

\node[gray] at (4.38,4.95) {\tiny 9};
\node[gray] at (4.12,5.38) {\tiny 10};

\fill[color=myred] (3.75,3.9) -- (4.25,3.9) -- (4,4.35) -- cycle;
\fill[color=myred] (4.25,4.76) -- (4.5,4.35) -- (4.75,4.76) -- cycle;
\fill[color=myred] (5.5,4.35) -- (6,4.35) -- (5.75,4.76) -- cycle;
\draw (4.5,4.35) -- (3.75,5.63);
\draw (2.25,3.9) -- (4.25,3.9) -- (4.5,4.35) -- (6,4.35) -- (6.75,5.63);
\draw (2.5,4.35) -- (4,4.35) -- (4.25,4.76) -- (5.75,4.76) -- (6.5,6.04);
\node at (2.22,4.15) {\tiny 0};
\node at (6.78,5.89) {\tiny 0};
\draw (2.75,3.9) -- (3,4.35);
\draw (3.25,3.9) -- (3.5,4.35);
\draw (3.75,3.9) -- (4,4.35);
\draw (4.5,4.35) -- (4.75,4.76);
\draw (5,4.35) -- (5.25,4.76);
\draw (5.5,4.35) -- (5.75,4.76);
\draw (6.5,3.48) -- (5.75,4.76);
\draw (6.25,4.76) -- (6,5.18);
\draw (6.5,5.18) -- (6.25,5.63);
\draw (6,3.48) -- (5.5,4.35);
\draw (5.75,3.9) -- (6.25,3.9);
\draw (4.25,3.03) -- (3.75,3.9);
\draw (4.75,3.03) -- (4.25,3.9);
\draw (4,3.48) -- (4.5,3.48);

\node[gray] at (4.35,3.25) {\tiny 6};
\node[gray] at (6.13,3.68) {\tiny 7};

\node[gray] at (4.12,3.68) {\tiny 7};
\node[gray] at (5.88,4.1) {\tiny 8};


\fill[color=myred] (3.5,1.72)-- (4,1.72) -- (3.75,2.14) -- cycle;
\fill[color=myred] (4.5,2.61)-- (4.75,3.03) -- (4.25,3.03) -- cycle;
\fill[color=myred] (5.5,2.61)-- (6,2.61) -- (5.75,3.03) -- cycle;
\fill[color=myred] (6.25,3.03)-- (6.5,3.48) -- (6,3.48) -- cycle;
\fill[color=myred] (6.75,3.03)-- (7.25,3.03) -- (7,3.48) -- cycle;
\draw (1,1.72) -- (4,1.72) -- (4.5,2.61) -- (6,2.61) -- (6.25,3.03) --
(7.25,3.03) -- (7.75,3.9);
\draw (1.25,2.14) -- (3.75,2.14) -- (4.25,3.03) -- (5.75,3.03) -- (6,3.48)
-- (7,3.48) -- (7.5,4.35);
\draw (1.5,1.72) -- (1.75,2.14);
\draw (2,1.72) -- (2.25,2.14);
\draw (2.5,1.72) -- (2.75,2.14);
\draw (3,1.72) -- (3.25,2.14);
\draw (3.5,1.72) -- (3.75,2.14);
\draw (4.5,2.61) -- (4.75,3.03);
\draw (5,2.61) -- (5.25,3.03);
\draw (5.5,2.61) -- (5.75,3.03);
\draw (6.25,3.03) -- (6.5,3.48);
\draw (6.75,3.03) -- (7,3.48);
\draw (4.25,2.14) -- (4,2.61);
\draw (4.5,2.61) -- (4.25,3.03);
\draw (6,2.61) -- (5.75,3.03);
\draw (6.25,3.03) -- (6,3.48);
\draw (7.5,3.48) -- (7.25,3.9);
\draw (6,1.72) -- (5.5,2.61);
\draw (6.5,1.72) -- (6,2.61);
\draw (6.25,2.14) -- (5.75,2.14);
\draw (7.25,2.14) -- (7.75,2.14);
\draw (7,2.61) -- (7.5,2.61);
\node at (0.97,1.98) {\tiny 0};
\node at (7.78,4.18) {\tiny 0};


\fill[color=myred] (3.5,0.87) -- (4,0.87) -- (4.25,1.3) -- (3.75,1.3) --
cycle;
\fill[color=myred] (5.5,0.87) -- (6,0.87) -- (5.75,1.3)  -- cycle;
\fill[color=myred] (6.25,1.3) -- (6.5,1.72) -- (6,1.72) -- cycle;
\fill[color=myred] (7.25,1.3) -- (8.25,1.3) -- (8,1.72) -- (7.5,1.72) --
cycle;
\draw (0.5,0.87) -- (6,0.87) -- (6.25,1.3) -- (8.25,1.3) -- (8.75,2.17);
\draw (0.75,1.3) -- (5.75,1.3) -- (6,1.72) -- (8,1.72) -- (8.5,2.58);
\draw (8.5,1.72) -- (8.25,2.17);
\draw (1,0.87) -- (1.25,1.3);
\draw (1.5,0.87) -- (1.75,1.3);
\draw (2,0.87) -- (2.25,1.3);
\draw (2.5,0.87) -- (2.75,1.3);
\draw (3,0.87) -- (3.25,1.3);
\draw (3.5,0.87) -- (3.75,1.3);
\draw (4,0.87) -- (4.25,1.3);
\draw (4.5,0.87) -- (4.75,1.3);
\draw (5,0.87) -- (5.25,1.3);
\draw (5.5,0.87) -- (5.75,1.3);
\draw (6.25,1.3) -- (6.5,1.72);
\draw (6.75,1.3) -- (7,1.72);
\draw (7.25,1.3) -- (7.5,1.72);
\draw (7.75,1.3) -- (8,1.72);
\node at (0.47,1.12) {\tiny 0};
\node at (8.77,2.43) {\tiny 0};

\draw (4,0) -- (3.5,0.87);
\draw (4.5,0) -- (3.5,1.72);
\draw (4.25,1.3) -- (3.75,2.14);
\node at (4.25,-0.2) {\tiny 0};
\node at (0.25,-0.2) {\tiny 1};
\node at (0.75,-0.2) {\tiny 1};
\node at (1.25,-0.2) {\tiny 1};
\node at (1.75,-0.2) {\tiny 1};
\node at (2.25,-0.2) {\tiny 1};
\node at (2.75,-0.2) {\tiny 1};
\node at (3.25,-0.2) {\tiny 1};
\node at (3.75,-0.2) {\tiny 1};
\node at (4.75,-0.2) {\tiny 1};
\node at (5.25,-0.2) {\tiny 1};
\node at (5.75,-0.2) {\tiny 1};
\node at (6.75,-0.2) {\tiny 1};
\node at (7.25,-0.2) {\tiny 1};
\node at (7.75,-0.2) {\tiny 1};
\node at (9.25,-0.2) {\tiny 1};
\node at (9.75,-0.2) {\tiny 1};

\draw (6,0) -- (5.5,0.87);
\draw (6.5,0) -- (5.75,1.3);
\node at (6.25,-0.2) {\tiny 0};

\draw (8,0) -- (7.25,1.3);
\draw (8.5,0) -- (6.75,3.03);
\draw (9,0) -- (7,3.46);
\node at (8.25,-0.2) {\tiny 0};
\node at (8.75,-0.2) {\tiny 0};

\draw (3.75,0.43) -- (4.25,0.43);
\draw (5.75,0.43) -- (6.25,0.43);
\draw (7.75,0.43) -- (8.75,0.43);
\draw (7.5,0.87) -- (8.5,0.87);

\node[gray] at (4.12,0.22) {\tiny 1};
\node[gray] at (6.12,0.22) {\tiny 1};
\node[gray] at (8.12,0.22) {\tiny 1};
\node[gray] at (8.62,0.22) {\tiny 1};
\node[gray] at (3.87,0.65) {\tiny 2};
\node[gray] at (5.87,0.65) {\tiny 2};
\node[gray] at (7.87,0.65) {\tiny 2};
\node[gray] at (8.35,0.65) {\tiny 2};

\node[gray] at (3.87,1.5) {\tiny 3};
\node[gray] at (7.63,1.07) {\tiny 3};
\node[gray] at (8.12,1.07) {\tiny 3};

\node[gray] at (6.12,1.93) {\tiny 4};
\node[gray] at (7.62,1.93) {\tiny 4};

\node[gray] at (5.89,2.35) {\tiny 5};
\node[gray] at (7.37,2.35) {\tiny 5};
\node[gray] at (7.12,2.8) {\tiny 6};

\node at (-0.1,0.25) {\tiny 1};
\node at (0.15,0.66) {\tiny 1};
\node at (0.65,1.52) {\tiny 1};
\node at (1.15,2.36) {\tiny 1};
\node at (1.40,2.82) {\tiny 1};
\node at (1.65,3.28) {\tiny 1};
\node at (1.9,3.72) {\tiny 1};
\node at (2.4,4.57) {\tiny 1};
\node at (2.65,4.98) {\tiny 1};
\node at (2.9,5.39) {\tiny 1};
\node at (3.4,6.27) {\tiny 1};
\node at (3.65,6.72) {\tiny 1};
\node at (3.9,7.18) {\tiny 1};
\node at (4.18,7.59) {\tiny 1};
\node at (4.43,8) {\tiny 1};
\node at (4.68,8.45) {\tiny 1};

\node at (10.08,0.39) {\tiny 1};
\node at (9.83,0.8) {\tiny 1};
\node at (9.58,1.21) {\tiny 1};
\node at (9.33,1.62) {\tiny 1};
\node at (9.11,2) {\tiny 1};
\node at (8.61,2.82) {\tiny 1};
\node at (8.36,3.28) {\tiny 1};
\node at (8.11,3.72) {\tiny 1};
\node at (7.61,4.57) {\tiny 1};
\node at (7.36,4.98) {\tiny 1};
\node at (7.11,5.39) {\tiny 1};
\node at (6.61,6.27) {\tiny 1};
\node at (6.36,6.72) {\tiny 1};
\node at (6.11,7.18) {\tiny 1};
\node at (5.86,7.59) {\tiny 1};
\node at (5.32,8.45) {\tiny 1};

\node[rotate=-90] at (0.9,1.1) {\tiny 1};
\node[rotate=-90] at (1.4,1.1) {\tiny 2};
\node[rotate=-90] at (1.9,1.1) {\tiny 3};
\node[rotate=-90] at (2.4,1.1) {\tiny 4};
\node[rotate=-90] at (2.9,1.1) {\tiny 5};
\node[rotate=-90] at (3.4,1.1) {\tiny 6};
\node[rotate=-90] at (4.4,1.1) {\tiny 7};
\node[rotate=-90] at (4.9,1.1) {\tiny 8};
\node[rotate=-90] at (5.4,1.1) {\tiny 9};
\node[rotate=-90] at (6.6,1.5) {\tiny 10};
\node[rotate=-90] at (7.1,1.5) {\tiny 11};
\node[rotate=-90] at (1.35,1.92) {\tiny 1};
\node[rotate=-90] at (1.85,1.92) {\tiny 2};
\node[rotate=-90] at (2.35,1.92) {\tiny 3};
\node[rotate=-90] at (2.85,1.92) {\tiny 4};
\node[rotate=-90] at (3.35,1.92) {\tiny 5};
\node[rotate=-90] at (4.85,2.82) {\tiny 6};
\node[rotate=-90] at (5.35,2.82) {\tiny 7};
\node[rotate=-90] at (6.6,3.25) {\tiny 8};
\node[rotate=-90] at (2.6,4.1) {\tiny 1};
\node[rotate=-90] at (3.1,4.1) {\tiny 2};
\node[rotate=-90] at (3.6,4.1) {\tiny 3};
\node[rotate=-90] at (4.85,4.55) {\tiny 4};
\node[rotate=-90] at (5.35,4.55) {\tiny 5};
\node[rotate=-90] at (3.6,5.8) {\tiny 1};

\end{tikzpicture}&
\spadesuit T=\begin{Young}&&&\cr &&&\cr8&11&&\cr 5&10&&\cr 4&7&&\cr 1&6&9&\cr &3&8&\cr &2&7&\cr &1&5&6\cr &&4&5\cr &&3&4\cr &&2&3\cr &&1&2\cr &&&1\cr &&&\cr &&&\cr
\end{Young}\\
&\;\;\;\;\;\;\;\;\;(\mu^t,\lambda^t,\nu^t)
\end{array}$$

\end{example}

Performing  $\clubsuit(\blacklozenge\clubsuit)=\clubsuit(\clubsuit\spadesuit)$ on the puzzle in  Example \ref{ex:tao} (or $\blacklozenge\spadesuit(\blacklozenge(\blacklozenge\clubsuit))$, that is, rotating the previous puzzle clockwise $2\pi/3$ radians about the center), we get a puzzle of boundary $(\nu^t,\mu^t,\lambda^t)$ where $\nu^t=4^22^31^3$. It is not difficult to see that the corridors consisting of upright rhombi, in the previous puzzle, now filled in blue below, give rise to the filling of $\clubsuit T$ of boundary $(\lambda^t,\nu^t,\mu^t)$.

$$\begin{array}{cccc}
\colorlet{lightgray}{black!10}\colorlet{lightgray1}{black!5}\colorlet{myred}{magenta!65}
\begin{tikzpicture}[scale=0.8]
\fill[color=cyan] (2.5,4.35) -- (4,4.35) -- (4.25,4.76) -- (3.75,5.63) --
(3.25,5.63) -- cycle;
\fill[color=cyan] (4.75,4.76) -- (5.75,4.76) -- (6.5,6.04) --  (5.5,7.77)
-- (4.25,5.63) -- cycle;
\fill[color=cyan] (3.5,6.04) -- (4,6.04) -- (5.25,8.22) -- (5,8.66);
\fill[color=cyan] (1.25,2.14) -- (3.75,2.14) -- (4.25,3.03) -- (3.75,3.9)
-- (2.25,3.9) -- cycle;
\fill[color=cyan] (4.75,3.03) -- (5.75,3.03) -- (6,3.48) -- (5.5,4.35) --
(4.5,4.35) -- (4.25,3.9) -- cycle;
\fill[color=cyan] (6.5,3.48) -- (7,3.48) -- (7.5,4.35) -- (6.75,5.63) --
(6,4.35);
\fill[color=cyan] (0.75,1.3) -- (3.75,1.3) -- (3.5,1.72) -- (1,1.72) --
cycle;
\fill[color=cyan] (4.25,1.3) -- (5.75,1.3) -- (6,1.72) -- (5.5,2.61) --
(4.5,2.61) -- (4,1.72) -- cycle;
\fill[color=cyan] (6.5,1.72) -- (7.5,1.72) -- (6.75,3.03) -- (6.25,3.03)
-- (6,2.61) -- cycle;
\fill[color=cyan] (8,1.72) -- (8.5,2.61) -- (7.75,3.9) -- (7.25,3.03) --
cycle;
\fill[color=cyan] (0,0) -- (4,0) -- (3.5,0.87) -- (0.5,0.87) -- cycle;
\fill[color=cyan] (4.5,0) -- (6,0) -- (5.5,0.87) -- (4,0.87) -- cycle;
\fill[color=cyan] (6.5,0) -- (8,0) -- (7.25,1.3) -- (6.25,1.3) -- (6,0.87)
-- cycle;
\fill[color=cyan] (9,0) -- (10,0) -- (8.75,2.17) -- (8.25,1.3) -- cycle;
\draw [color=lightgray1] (9.5,0.86)-- (0.5,0.87);
\draw [color=lightgray1] (9,1.73)-- (1.01,1.72);
\draw [color=lightgray1] (8.51,2.58)-- (1.5,2.6);
\draw [color=lightgray1] (8,3.46)-- (1.99,3.48);
\draw [color=lightgray1] (7.5,4.34)-- (2.49,4.35);
\draw [color=lightgray1] (7.01,5.18)-- (3.01,5.18);
\draw [color=lightgray1] (6.51,6.05)-- (3.51,6.04);
\draw [color=lightgray1] (5.51,7.78)-- (4.51,7.77);
\draw [color=lightgray1] (5.51,7.78)-- (1,0);
\draw [color=lightgray1] (3,0)-- (6.51,6.05);

\draw [color=lightgray1] (0.25,0.43) -- (9.75,0.43);
\draw [color=lightgray1] (0.75,1.3) -- (9.25,1.3);
\draw [color=lightgray1] (1.25,2.15) -- (8.75,2.15);
\draw [color=lightgray1] (1.75,3.03) -- (8.25,3.03);
\draw [color=lightgray1] (2.25,3.9) -- (7.75,3.9);
\draw [color=lightgray1] (2.75,4.76) -- (7.25,4.76);
\draw [color=lightgray1] (3.25,5.63) -- (6.75,5.63);
\draw [color=lightgray1] (3.75,6.5) -- (6.25,6.5);
\draw [color=lightgray1] (4.25,7.36) -- (5.75,7.36);
\draw [color=lightgray1] (4.75,8.23) -- (5.25,8.23);

\draw [color=lightgray1] (3.97,6.98)-- (5.97,6.98);
\draw [color=lightgray1] (7.01,5.18)-- (4,0);
\draw [color=lightgray1] (5,0)-- (7.5,4.34);
\draw [color=lightgray1] (6,0)-- (8,3.46);
\draw [color=lightgray1] (8.51,2.58)-- (7,0);
\draw [color=lightgray1] (5.97,6.98)-- (2,0);
\draw [color=lightgray1] (9,1.73)-- (8,0);
\draw [color=lightgray1] (9,0)-- (9.5,0.86);
\draw [color=lightgray1] (0.5,0.87)-- (1,0);
\draw [color=lightgray1] (2,0)-- (1.01,1.72);
\draw [color=lightgray1] (1.5,2.6)-- (3,0);
\draw [color=lightgray1] (4,0)-- (1.99,3.48);
\draw [color=lightgray1] (5,0)-- (2.49,4.35);
\draw [color=lightgray1] (6,0)-- (3.01,5.18);
\draw [color=lightgray1] (7,0)-- (3.51,6.04);
\draw [color=lightgray1] (8,0)-- (3.97,6.98);
\draw [color=lightgray1] (9,0)-- (4.51,7.77);

\draw [color=lightgray1] (0.5,0)-- (0.25,0.43);
\draw [color=lightgray1] (1.5,0)-- (0.75,1.3);
\draw [color=lightgray1] (2.5,0)-- (1.25,2.17);
\draw [color=lightgray1] (3.5,0)-- (1.75,3.03);
\draw [color=lightgray1] (4.5,0)-- (2.25,3.9);
\draw [color=lightgray1] (5.5,0)-- (2.75,4.76);
\draw [color=lightgray1] (6.5,0)-- (3.25,5.63);
\draw [color=lightgray1] (7.5,0)-- (3.75,6.5);
\draw [color=lightgray1] (8.5,0)-- (4.25,7.36);
\draw [color=lightgray1] (9.5,0)-- (4.75,8.23);
\draw [color=lightgray1] (9.5,0)-- (9.75,0.44);
\draw [color=lightgray1] (8.5,0)-- (9.25,1.31);
\draw [color=lightgray1] (7.5,0)-- (8.76,2.14);
\draw [color=lightgray1] (6.5,0)-- (8.25,3.03);
\draw [color=lightgray1] (5.5,0)-- (7.75,3.9);
\draw [color=lightgray1] (4.5,0)-- (7.26,4.74);
\draw [color=lightgray1] (3.5,0)-- (6.75,5.62);
\draw [color=lightgray1] (2.5,0)-- (6.25,6.5);
\draw [color=lightgray1] (1.5,0)-- (5.75,7.37);
\draw [color=lightgray1] (0.5,0)-- (5.25,8.22);
\draw [color=lightgray1] (1,0)-- (2.51,2.61);
\draw [color=lightgray1] (2.51,2.61)-- (1.99,3.48);
\draw [color=lightgray1] (0.5,0)-- (1.98,2.58);
\draw [color=lightgray1] (1.98,2.58)-- (1.75,3.03);

\draw[color=black] (0,0) -- (10,0) -- (5,8.66) -- cycle;


\fill[color=myred] (3.75,5.63) -- (4.25,5.63) -- (4,6.04) -- cycle;
\draw (3.25,5.63) -- (4.25,5.63) -- (5.5,7.77);
\draw (3.5,6.04) -- (4,6.04) -- (5.25,8.22);
\node at (3.2,5.88) {\tiny 0};
\node at (5.57,8.05) {\tiny 0};
\draw (4.75,4.76) -- (4,6.04);
\draw (3.75,5.63) -- (4,6.04);
\draw (4.5,6.04) -- (4.25,6.5);
\draw (4.75,6.5) -- (4.5,6.94);
\draw (5,6.95) -- (4.75,7.36);
\draw (5.25,7.36) -- (5,7.77);
\draw (4,5.18) -- (4.5,5.18);

\node[gray] at (4.38,4.95) {\tiny 9};
\node[gray] at (4.12,5.38) {\tiny 10};

\fill[color=myred] (3.75,3.9) -- (4.25,3.9) -- (4,4.35) -- cycle;
\fill[color=myred] (4.25,4.76) -- (4.5,4.35) -- (4.75,4.76) -- cycle;
\fill[color=myred] (5.5,4.35) -- (6,4.35) -- (5.75,4.76) -- cycle;
\draw (4.5,4.35) -- (3.75,5.63);
\draw (2.25,3.9) -- (4.25,3.9) -- (4.5,4.35) -- (6,4.35) -- (6.75,5.63);
\draw (2.5,4.35) -- (4,4.35) -- (4.25,4.76) -- (5.75,4.76) -- (6.5,6.04);
\node at (2.22,4.15) {\tiny 0};
\node at (6.78,5.89) {\tiny 0};
\draw (2.75,3.9) -- (3,4.35);
\draw (3.25,3.9) -- (3.5,4.35);
\draw (3.75,3.9) -- (4,4.35);
\draw (4.5,4.35) -- (4.75,4.76);
\draw (5,4.35) -- (5.25,4.76);
\draw (5.5,4.35) -- (5.75,4.76);
\draw (6.5,3.48) -- (5.75,4.76);
\draw (6.25,4.76) -- (6,5.18);
\draw (6.5,5.18) -- (6.25,5.63);
\draw (6,3.48) -- (5.5,4.35);
\draw (5.75,3.9) -- (6.25,3.9);
\draw (4.25,3.03) -- (3.75,3.9);
\draw (4.75,3.03) -- (4.25,3.9);
\draw (4,3.48) -- (4.5,3.48);

\node[gray] at (4.35,3.25) {\tiny 6};
\node[gray] at (6.13,3.68) {\tiny 7};

\node[gray] at (4.12,3.68) {\tiny 7};
\node[gray] at (5.88,4.1) {\tiny 8};


\fill[color=myred] (3.5,1.72)-- (4,1.72) -- (3.75,2.14) -- cycle;
\fill[color=myred] (4.5,2.61)-- (4.75,3.03) -- (4.25,3.03) -- cycle;
\fill[color=myred] (5.5,2.61)-- (6,2.61) -- (5.75,3.03) -- cycle;
\fill[color=myred] (6.25,3.03)-- (6.5,3.48) -- (6,3.48) -- cycle;
\fill[color=myred] (6.75,3.03)-- (7.25,3.03) -- (7,3.48) -- cycle;
\draw (1,1.72) -- (4,1.72) -- (4.5,2.61) -- (6,2.61) -- (6.25,3.03) --
(7.25,3.03) -- (7.75,3.9);
\draw (1.25,2.14) -- (3.75,2.14) -- (4.25,3.03) -- (5.75,3.03) -- (6,3.48)
-- (7,3.48) -- (7.5,4.35);
\draw (1.5,1.72) -- (1.75,2.14);
\draw (2,1.72) -- (2.25,2.14);
\draw (2.5,1.72) -- (2.75,2.14);
\draw (3,1.72) -- (3.25,2.14);
\draw (3.5,1.72) -- (3.75,2.14);
\draw (4.5,2.61) -- (4.75,3.03);
\draw (5,2.61) -- (5.25,3.03);
\draw (5.5,2.61) -- (5.75,3.03);
\draw (6.25,3.03) -- (6.5,3.48);
\draw (6.75,3.03) -- (7,3.48);
\draw (4.25,2.14) -- (4,2.61);
\draw (4.5,2.61) -- (4.25,3.03);
\draw (6,2.61) -- (5.75,3.03);
\draw (6.25,3.03) -- (6,3.48);
\draw (7.5,3.48) -- (7.25,3.9);
\draw (6,1.72) -- (5.5,2.61);
\draw (6.5,1.72) -- (6,2.61);
\draw (6.25,2.14) -- (5.75,2.14);
\draw (7.25,2.14) -- (7.75,2.14);
\draw (7,2.61) -- (7.5,2.61);
\node at (0.97,1.98) {\tiny 0};
\node at (7.78,4.18) {\tiny 0};


\fill[color=myred] (3.5,0.87) -- (4,0.87) -- (4.25,1.3) -- (3.75,1.3) --
cycle;
\fill[color=myred] (5.5,0.87) -- (6,0.87) -- (5.75,1.3)  -- cycle;
\fill[color=myred] (6.25,1.3) -- (6.5,1.72) -- (6,1.72) -- cycle;
\fill[color=myred] (7.25,1.3) -- (8.25,1.3) -- (8,1.72) -- (7.5,1.72) --
cycle;
\draw (0.5,0.87) -- (6,0.87) -- (6.25,1.3) -- (8.25,1.3) -- (8.75,2.17);
\draw (0.75,1.3) -- (5.75,1.3) -- (6,1.72) -- (8,1.72) -- (8.5,2.58);
\draw (8.5,1.72) -- (8.25,2.17);
\draw (1,0.87) -- (1.25,1.3);
\draw (1.5,0.87) -- (1.75,1.3);
\draw (2,0.87) -- (2.25,1.3);
\draw (2.5,0.87) -- (2.75,1.3);
\draw (3,0.87) -- (3.25,1.3);
\draw (3.5,0.87) -- (3.75,1.3);
\draw (4,0.87) -- (4.25,1.3);
\draw (4.5,0.87) -- (4.75,1.3);
\draw (5,0.87) -- (5.25,1.3);
\draw (5.5,0.87) -- (5.75,1.3);
\draw (6.25,1.3) -- (6.5,1.72);
\draw (6.75,1.3) -- (7,1.72);
\draw (7.25,1.3) -- (7.5,1.72);
\draw (7.75,1.3) -- (8,1.72);
\node at (0.47,1.12) {\tiny 0};
\node at (8.77,2.43) {\tiny 0};

\draw (4,0) -- (3.5,0.87);
\draw (4.5,0) -- (3.5,1.72);
\draw (4.25,1.3) -- (3.75,2.14);
\node at (4.25,-0.2) {\tiny 0};
\node at (0.25,-0.2) {\tiny 1};
\node at (0.75,-0.2) {\tiny 1};
\node at (1.25,-0.2) {\tiny 1};
\node at (1.75,-0.2) {\tiny 1};
\node at (2.25,-0.2) {\tiny 1};
\node at (2.75,-0.2) {\tiny 1};
\node at (3.25,-0.2) {\tiny 1};
\node at (3.75,-0.2) {\tiny 1};
\node at (4.75,-0.2) {\tiny 1};
\node at (5.25,-0.2) {\tiny 1};
\node at (5.75,-0.2) {\tiny 1};
\node at (6.75,-0.2) {\tiny 1};
\node at (7.25,-0.2) {\tiny 1};
\node at (7.75,-0.2) {\tiny 1};
\node at (9.25,-0.2) {\tiny 1};
\node at (9.75,-0.2) {\tiny 1};

\draw (6,0) -- (5.5,0.87);
\draw (6.5,0) -- (5.75,1.3);
\node at (6.25,-0.2) {\tiny 0};

\draw (8,0) -- (7.25,1.3);
\draw (8.5,0) -- (6.75,3.03);
\draw (9,0) -- (7,3.46);
\node at (8.25,-0.2) {\tiny 0};
\node at (8.75,-0.2) {\tiny 0};

\draw (3.75,0.43) -- (4.25,0.43);
\draw (5.75,0.43) -- (6.25,0.43);
\draw (7.75,0.43) -- (8.75,0.43);
\draw (7.5,0.87) -- (8.5,0.87);

\node[gray] at (4.12,0.22) {\tiny 1};
\node[gray] at (6.12,0.22) {\tiny 1};
\node[gray] at (8.12,0.22) {\tiny 1};
\node[gray] at (8.62,0.22) {\tiny 1};
\node[gray] at (3.87,0.65) {\tiny 2};
\node[gray] at (5.87,0.65) {\tiny 2};
\node[gray] at (7.87,0.65) {\tiny 2};
\node[gray] at (8.35,0.65) {\tiny 2};

\node[gray] at (3.87,1.5) {\tiny 3};
\node[gray] at (7.63,1.07) {\tiny 3};
\node[gray] at (8.12,1.07) {\tiny 3};

\node[gray] at (6.12,1.93) {\tiny 4};
\node[gray] at (7.62,1.93) {\tiny 4};

\node[gray] at (5.89,2.35) {\tiny 5};
\node[gray] at (7.37,2.35) {\tiny 5};
\node[gray] at (7.12,2.8) {\tiny 6};

\node at (-0.1,0.25) {\tiny 1};
\node at (0.15,0.66) {\tiny 1};
\node at (0.65,1.52) {\tiny 1};
\node at (1.15,2.36) {\tiny 1};
\node at (1.40,2.82) {\tiny 1};
\node at (1.65,3.28) {\tiny 1};
\node at (1.9,3.72) {\tiny 1};
\node at (2.4,4.57) {\tiny 1};
\node at (2.65,4.98) {\tiny 1};
\node at (2.9,5.39) {\tiny 1};
\node at (3.4,6.27) {\tiny 1};
\node at (3.65,6.72) {\tiny 1};
\node at (3.9,7.18) {\tiny 1};
\node at (4.18,7.59) {\tiny 1};
\node at (4.43,8) {\tiny 1};
\node at (4.68,8.45) {\tiny 1};

\node at (10.08,0.39) {\tiny 1};
\node at (9.83,0.8) {\tiny 1};
\node at (9.58,1.21) {\tiny 1};
\node at (9.33,1.62) {\tiny 1};
\node at (9.11,2) {\tiny 1};
\node at (8.61,2.82) {\tiny 1};
\node at (8.36,3.28) {\tiny 1};
\node at (8.11,3.72) {\tiny 1};
\node at (7.61,4.57) {\tiny 1};
\node at (7.36,4.98) {\tiny 1};
\node at (7.11,5.39) {\tiny 1};
\node at (6.61,6.27) {\tiny 1};
\node at (6.36,6.72) {\tiny 1};
\node at (6.11,7.18) {\tiny 1};
\node at (5.86,7.59) {\tiny 1};
\node at (5.32,8.45) {\tiny 1};

\node[rotate=-90] at (0.9,1.1) {\tiny 1};
\node[rotate=-90] at (1.4,1.1) {\tiny 2};
\node[rotate=-90] at (1.9,1.1) {\tiny 3};
\node[rotate=-90] at (2.4,1.1) {\tiny 4};
\node[rotate=-90] at (2.9,1.1) {\tiny 5};
\node[rotate=-90] at (3.4,1.1) {\tiny 6};
\node[rotate=-90] at (4.4,1.1) {\tiny 7};
\node[rotate=-90] at (4.9,1.1) {\tiny 8};
\node[rotate=-90] at (5.4,1.1) {\tiny 9};
\node[rotate=-90] at (6.6,1.5) {\tiny 10};
\node[rotate=-90] at (7.1,1.5) {\tiny 11};
\node[rotate=-90] at (1.35,1.92) {\tiny 1};
\node[rotate=-90] at (1.85,1.92) {\tiny 2};
\node[rotate=-90] at (2.35,1.92) {\tiny 3};
\node[rotate=-90] at (2.85,1.92) {\tiny 4};
\node[rotate=-90] at (3.35,1.92) {\tiny 5};
\node[rotate=-90] at (4.85,2.82) {\tiny 6};
\node[rotate=-90] at (5.35,2.82) {\tiny 7};
\node[rotate=-90] at (6.6,3.25) {\tiny 8};
\node[rotate=-90] at (2.6,4.1) {\tiny 1};
\node[rotate=-90] at (3.1,4.1) {\tiny 2};
\node[rotate=-90] at (3.6,4.1) {\tiny 3};
\node[rotate=-90] at (4.85,4.55) {\tiny 4};
\node[rotate=-90] at (5.35,4.55) {\tiny 5};
\node[rotate=-90] at (3.6,5.8) {\tiny 1};
\node[rotate=135,blue] at (5.25,7.75) {\tiny 1};
\node[rotate=135,blue] at (5,7.35) {\tiny 2};
\node[rotate=135,blue] at (4.75,6.95) {\tiny 3};
\node[rotate=135,blue] at (4.5,6.5) {\tiny 4};
\node[rotate=135,blue] at (4.25,6.1) {\tiny 5};
\node[rotate=135,blue] at (4.25,4.35) {\tiny 6};
\node[rotate=135,blue] at (4.25,2.6) {\tiny 7};
\node[rotate=135,blue] at (4,2.15) {\tiny 8};
\node[rotate=135,blue] at (6.5,5.65) {\tiny 1};
\node[rotate=135,blue] at (6.25,5.2) {\tiny 2};
\node[rotate=135,blue] at (6,4.75) {\tiny 3};
\node[rotate=135,blue] at (6,3.05) {\tiny 4};
\node[rotate=135,blue] at (6,1.3) {\tiny 5};
\node[rotate=135,blue] at (7.5,3.9) {\tiny 1};
\node[rotate=135,blue] at (7.25,3.45) {\tiny 2};

\node[rotate=135,blue] at (8.5,2.15) {\tiny 1};
\node[rotate=135,blue] at (8.25,1.75) {\tiny 2};

\end{tikzpicture}&
\clubsuit T=\begin{Young}&&&\cr &&&\cr 5&&&\cr 2&8&&\cr 1&7&&\cr &4&&\cr &2&&\cr &1&6&\cr &&3&\cr &&2&\cr &&1&5\cr &&&4\cr &&&3\cr &&&2\cr &&&1\cr &&&\cr
\end{Young}\\
&\;\;\;\;\;\;\;\;(\lambda^t,\nu^t,\mu^t)
\end{array}$$

Tao's bijection translates the involutions $\spadesuit$,   $\clubsuit$  and $\blacklozenge$
   on puzzles to LR tableaux which in turn can also be explained using Puhrboo's mosaics ~\cite{mosaic} naturally in bijection with puzzles. In ~\cite[Section 5.1]{mosaic}  it is discussed how the operation \emph{migration} of a single
rhombus in a mosaic is related with {\em jeu de taquin} slides on
tableaux.
\emph{Migration} is
an invertible operation on tableau-like  structures on the  rhombi of a mosaic, called flocks, that allows    to identify a mosaic (equivalently, a puzzle) with an LR tableau. It gives a bijection between mosaics (equivalently puzzles) and LR tableaux, and, with appropriate orientation of the flocks in the mosaic, it coincides with Tao's bijection. More importantly, \emph{migration} allows to relate operations on puzzles with {\em jeu de taquin} operations, like tableau-switching \cite{bss},  on LR tableaux.
This explains the correspondence between  the action of
$H$ on puzzles and on LR tableaux. Technical details and illustrations on mosaics are deferred to Appendix \ref{a:mosaic}.
Our concern next is to show that although  involutions $\spadesuit$,   $\clubsuit$  and $\blacklozenge$  may be executed using {\em jeu de taquin slides} on LR tableaux,   those slides do not need to be performed upon a scan of the neighbours. The slides are independent of the relative size of the neighbours and are reduced to simple procedures defining linear cost involutions.

\subsection{The LR tableau $H$--symmetries }
\label{sec:names}
 Migration on Puhrboo's mosaics \cite{mosaic} is related with {\em jeu de taquin} and  translates H-symmetries on puzzles to LR tableaux through tableau switching and {\em standardization}. This translation coincides with Tao's bijection as illustrated in Example \ref{ex:taoblacklozenge}. We want to avoid the computational complexity of  standardization and {\em jeu de taquin} procedure which we succeed  through hybrid tableau pairs.

Given a decomposition of the   rectangle $D$ into shapes
$\mu$, $\lambda^\vee/\mu$, and $D/\lambda^\vee$, a triple of tableaux
$(U_1,U_2,U_3)$ is said to be a three--fold multitableau of shape $(\mu,\lambda^\vee/\mu,D/\lambda^\vee)$
if $U_1$ is  a filling of
the shape $\mu$,  $U_2$ is a filling of $\lambda^\vee/\mu$, and
$U_3$ is a filling of $D/\lambda^\vee$.
A three--fold LR multitableau of  boundary data $(\mu,\nu,\lambda)$ is a three--fold multitableau where the inner tableau is the Yamanouchi tableau $Y(\mu)$, the middle
one is the LR tableau of shape $\lambda^\vee/\mu$ and content $\nu$, and the outer tableau is  $Y(\lambda)^{\rm a}$ the {{\em anti-normal  form} of $Y(\lambda)$, defined to be the filling of the anti-normal shape of $\lambda$ such that each column, right to left, is filled with consecutive integers bottom to top  starting with $1$}.
 For instance, see \eqref{anti}.
{ Tableau switching can be adjusted to move a pair of tableaux through each other where one is a column strict tableau and the other is a row strict tableau see~\cite[Section 2, pp. 22 ]{bss}. Our next definitions are just a particular case where the left or right tableau is the transpose of a (or antinormal) Yamanouchi tableau row strict) and therefore switch moves can a priori be prescribed.}

We define the bijection
$\spadesuit: \rm LR(\mu,{\nu},\lambda)\longrightarrow
LR({\nu}^t,\mu^t,\lambda^t)$ as a five step procedure.

\begin{definition}[Map $\spadesuit$] \label{D:spade} Let $T\in \rm LR(\mu,{\nu},\lambda)$ inside the rectangle $D$ of size $d\times (n-d)$.

\begin{enumerate}
\item Fill the inner shape $\mu$, using a completely ordered alphabet different from the ``numerical'' filling of $T$, so that its transpose is the Yamanouchi tableau $Y(\mu^t)$.
\item For $i=1,\dots, d$, slide down vertically the $i$'s in the filling of $T$ to the $i$th row.
\item For $i=1,\dots, d$, slide horizontally all the numbers $i$'s to the left so that we get the Yamanouchi tableau $Y(\nu)$. Erase $Y(\nu)$.
\item { Transpose  the resulting filling to obtain $T^\spadesuit\in \rm LR({\nu}^t,\mu^t,\lambda^t)$.}
\end{enumerate}
\end{definition}
Clearly, the last step can   also be the first step with obvious adaptations in the next steps.
An illustration of this procedure follows.
\begin{example} \label{E:spade} Let $T\in LR(\mu,\nu,\lambda)$ with $d=4$, $n=11$, $\mu=4210$, $\nu=5420$ and  $\lambda=5320$. Then considering the twofold hybrid tableau
$$\begin{array}{ccccccccccccc}
T=\begin{Young}
1&3&&&&&\cr &2&2&3&&&\cr &&1&2&2&&\cr &&&&1&1&1\cr
\end{Young}
\rightarrow
([Y(\mu)]^t,T)=\begin{Young}
1&3&&&&&\cr { a}&2&2&3&&&\cr { a}&{ b}&1&2&2&&\cr { a}&{ b}&{ c}&{ d}&1&1&1\cr
\end{Young}
\rightarrow
\begin{Young}
1&&&\cr
1&&&\cr
1&2&&\cr
{\color{blue}d}&2&{\small 3}&\cr
{\color{blue}c}&1&2&\cr
 {\color{blue}b}&{\color{blue}b}&2&{\small3}\cr
{\color{blue}a}&{\color{blue}a}&{\color{blue}a}&1\cr
\end{Young}=(Y(\mu^t),  T^t)\\
\\
\rightarrow
\begin{Young}
1&&&\cr
1&&&\cr
1&2&&\cr
{\color{blue}d}&2&{\small 3}&\cr
1&2&{\color{blue}c}&\cr
{\color{blue}b}&2&{\small3}&{\color{blue}b}\cr
1&{\color{blue}a}&{\color{blue}a}&{\color{blue}a}\cr
\end{Young}
\rightarrow
\begin{Young}
{\color{blue}d}&     &&\cr
{\color{blue}b}&     &&\cr
1&{\color{blue}a}&               &\cr
1              &2&{\color{blue}c}&\cr
1              &2 &{\color{blue}a}&\cr
1 &2&{\small 3} &{\color{blue}b}\cr
1              &2 &{\small 3}&{\color{blue}a}\cr
\end{Young}=( [Y(\nu)]^t,T^\spadesuit)
\rightarrow
\begin{Young}
{\color{blue}d}&     &&\cr {\color{blue}b}&     &&\cr
              &{\color{blue}a}&               &\cr
              &              &{\color{blue}c}&\cr
              &              &{\color{blue}a}&\cr
              &              &             &{\color{blue}b}\cr
              &              &             &{\color{blue}a}\cr
\end{Young}={\rm T}^\spadesuit.

\end{array}$$

\end{example}
The procedure is clearly reversible and  an involution on $\mathcal{LR}$. Next we check that it yields
the desired tableau.
Let ${\bf s}_ i$, $i=1,2$, denote the tableau--switching
operation on the LR--multitableau of boundary data
$(\mu,\nu,\lambda)$ (recall subsection \ref{subsec:switching}) which switches the first two LR tableaux and the
last two  respectively, see~\cite{bss}. Compare the procedure with the explanation given by the migration for the operation $\spadesuit$ on mosaics in Appendix.
\begin{proposition} \label{P:spade} The map $\spadesuit$ is such that
$$ T
{\longrightarrow}(Y(\mu^t),  T^t, Y(\lambda^t)^{\rm
a})\longrightarrow {\bf s}_1(Y(\mu^t),  T^t, Y(\lambda^t)^{\rm
a})=( [Y(\nu)]^t, {T}^\spadesuit, Y(\lambda^t)^{\rm
a}).$$

\end{proposition}

\begin{proof} Since $T$ is column strict then $T^t$ is row strict
The second and third steps of the definition of the map $\spadesuit$
coincides  with the action of the switching operation ${\bf s}_1$ on the hybrid twofold tableau $(Y(\mu^t),
 T^t)$, that is, ${\bf s}_1(Y(\mu^t),
 T^t)=( [Y(\nu)]^t, {T}^\spadesuit)$, see \cite[Section 2]{bss}.
\end{proof}


The  bijection
$\rm LR(\mu,{\nu},\lambda)\overset{\text{$\clubsuit$}}\longrightarrow
LR(\mu^t,{\lambda^t},\nu^t)$ is defined similarly as a five step
procedure.

\begin{definition}[Map $\clubsuit$] \label{D:club} 
Let $T\in LR(\mu,{\nu},\lambda)$  inside the rectangle $D$ of size $d\times (n-d)$ and $\nu^t=(\nu_1^t,\nu_2^t,\dots,\nu_{\nu_1}^t)$.

\begin{enumerate}

\item Fill the  outer shape $\lambda$, using a completely ordered alphabet different from the ``numerical'' filling of $T$, so
that its transpose is  $Y(\lambda^t)^{\rm a}$.

\item { For $i=1,\dots,\nu_1^t$, slide horizontally the rightmost $i$ of $T$  to the $(n-d)$th column of $D$; for $i=1,\dots,\nu_2^t$, slide horizontally the rightmost $i$ in the first $(n-d-1)$ columns of $T$ to the $(n-d-1)$th column of $D$;  $\dots$, lastly slide horizontally the remaining $\nu_1$th $1$ to the $(n-d-\nu_1+1)$th column of $D$.}

\item Slide up vertically the numbers along each column so that we get $Y(\nu)^{\rm a}$. Erase $Y(\nu)^{\rm a}$.

\item { Transpose the resulting filling to obtain $T^\clubsuit \in LR(\mu^t,{\lambda^t},\nu^t)$.}

\end{enumerate}
\end{definition}
The example illustrates the procedure.
\begin{example}\label{E:club}Consider again the  LR tableau $T$ of Example \ref{E:spade} with $d=4$ and $n-d=7$. One has $\nu=542$, $\nu^t=33221$, $\lambda=532$, $\lambda^t=33211$,
and
 $$\begin{array}{ccccccccccccc}
{\rm T}=\begin{Young} 1&3&&&&&\cr &2&2&3&&&\cr &&1&2&2&&\cr
&&&&1&1&1\cr
\end{Young}
\rightarrow
\begin{Young}
1&3&{\color{blue}a}&{\color{blue}b}&{\color{blue}c}&{\color{blue}d}&{\color{blue}e}\cr
&2&2&3&{\color{blue}a}&{\color{blue}b}&{\color{blue}c}\cr
&&1&2&2&{\color{blue}a}&{\color{blue}b}\cr &&&&1&1&1\cr
\end{Young}
\rightarrow \cr
\begin{Young}
{\color{blue}a}&{\color{blue}b}&1&{\color{blue}c}&{\color{blue}d}&3&{\color{blue}e}\cr
&{\color{blue}a}&{\color{blue}b}&2&2&{\color{blue}c}&3\cr
&&{\color{blue}a}&1&{\color{blue}b}&2&2\cr &&&&1&1&1\cr
\end{Young}\rightarrow
\begin{Young}
{\color{blue}a}&{\color{blue}b}&1&2&2&3&3\cr
&{\color{blue}a}&{\color{blue}b}&1&1&2&2\cr
&&{\color{blue}a}&{\color{blue}c}&{\color{blue}d}&1&1\cr
&&&&{\color{blue}b}&{\color{blue}c}&{\color{blue}e}\cr
\end{Young}\rightarrow
\begin{Young}
{\color{blue}e}&&&\cr {\color{blue}c}&&&\cr
{\color{blue}b}&{\color{blue}d}&&\cr
              &{\color{blue}c}&&\cr
              &{\color{blue}a}&{\color{blue}b}&\cr
              &               &{\color{blue}a}&{\color{blue}b}\cr
               &                 &                &{\color{blue}a}\cr
\end{Young}={\rm T}^{\clubsuit}.
\end{array}$$
\end{example}
The procedure is clearly an involution on $\mathcal{LR}$ and as before we check that it
yields the desired tableau. Compare the procedure with the explanation given by the migration on mosaics for the operation $\clubsuit$  in Appendix. Again we have avoided the standardization of $T$.

\begin{proposition}\label{P:club} The map
$\clubsuit$  is such that
$$ T
{\longrightarrow}
(Y(\mu^t),  T^t,
Y(\lambda^t)^{\rm a}){\longrightarrow}{\bf s}_2(Y(\mu^t),   T^t,
Y(\lambda^t)^{\rm a})
=(Y(\mu^t), {T}^\clubsuit,
[Y(\nu)^{\rm a}]^t ).$$

\end{proposition}
\begin{proof}
The second and third steps of the definition of the map $\clubsuit$
correspond exactly to the action of the switching operation ${\bf s}_2$ on the hybrid two-fold tableau
$({ T}^t, Y(\lambda^t)^{\rm a})$.
\end{proof}

\begin{theorem} \label{th:lcost}
The involutions $\spadesuit$ and $\clubsuit$ on $\mathcal{LR}$ have linear cost.
\end{theorem}
\begin{proof}
The use of hybrid tableaux show clearly  that the maps
are performed using {\em jeu de taquin} without   need to scan  the neighbours.
\end{proof}

The  rotation symmetries on puzzles { are explicit} and exhibit the cyclic symmetries of LR coefficients $c_{\mu\,\nu\,\lambda}=c_{\lambda\,\mu\,\nu}=c_{\nu\,\lambda\,\mu}$. On LR tableaux  although less explicitly  they are easily performed,
noting that $\clubsuit\blacklozenge=\blacklozenge\spadesuit=\spadesuit\clubsuit$ and
$\blacklozenge\clubsuit=\spadesuit\blacklozenge=\clubsuit\spadesuit$.  Migration on mosaics explains with \textit{jeu de taquin} the rotation symmetries of LR tableaux using standardisation. For instance, the rotation below  $\blacklozenge\clubsuit=\spadesuit\blacklozenge=\clubsuit\spadesuit$ on puzzles is translated to LR tableaux in \cite[Corollary 5.3]{mosaic}.  Again we may avoid standardisation operation  with hybrid tableau pairs as shown next.

\begin{definition}[$2\pi/3$ radians counterclockwise rotation symmetry map $\blacklozenge\clubsuit=\spadesuit\blacklozenge=\clubsuit\spadesuit$] \label{D:rotation}
Let $T\in \rm LR(\mu,\nu,\lambda)$ inside the rectangle $D$ of size $d\times (n-d)$.
\begin{enumerate}

\item Rotate $T$ by $\pi$ radians.

\item Fill the  inner shape $\lambda$  using a completely ordered alphabet different
from the ``numerical'' filling of $T$ so that we get the  Yamanouchi tableau $Y(\lambda)$.

\item For each $i$, replace the $\nu_i$ $i$'s with $1,2,\dots,\nu_i$ according to the standard order on the boxes.

\item Slide horizontally each $1$ to the first column, each $2$ to the second column, each $3$ to the third column and so on.

\item Slide down vertically along each column the numbers.

\item Erase the numerical filling to obtain $\spadesuit\blacklozenge{\rm T}\in \rm LR(\nu,\lambda,\mu)$.

\end{enumerate}
\end{definition}

The example illustrates the procedure from which we see that it  is
reversible.
\begin{example}

$$\begin{array}{ccccccccccccc}
{\rm T}=\begin{Young} 1&3&&&&&\cr &2&2&3&&&\cr &&1&2&2&&\cr
&&&&1&1&1\cr
\end{Young}
\overset{\pi-rotation}\longrightarrow
\begin{Young}
1&1&1&&&&\cr \color{blue}c&\color{blue}c&2&2&1&&\cr
\color{blue}b&\color{blue}b&\color{blue}b&3&2&2&\cr
\color{blue}a&\color{blue}a&\color{blue}a&\color{blue}a&\color{blue}a&3&1\cr
\end{Young}
\rightarrow
\begin{Young}
1&2&3&&&&\cr \color{blue}c&\color{blue}c&1&2&4&&\cr
\color{blue}b&\color{blue}b&\color{blue}b&1&3&4&\cr
\color{blue}a&\color{blue}a&\color{blue}a&\color{blue}a&\color{blue}a&2&5\cr
\end{Young}=Y(\lambda)\cup (\blacklozenge T)^t
\cr
\rightarrow\begin{Young}
1&2&3&&&&\cr 1&2&\color{blue}c&4&\color{blue}c&&\cr
1&\color{blue}b&3&4&\color{blue}b&\color{blue}b&\cr
\color{blue}a&2&\color{blue}a&\color{blue}a&5&\color{blue}a&\color{blue}a\cr
\end{Young}
\rightarrow
\begin{Young}
\color{blue}a&\color{blue}b&\color{blue}c&&&&\cr
1&2&\color{blue}a&\color{blue}a&\color{blue}c&&\cr
1&2&3&4&\color{blue}b&\color{blue}b&\cr
1&2&3&4&5&\color{blue}a&\color{blue}a\cr
\end{Young}
\rightarrow
\begin{Young}
\color{blue}a&\color{blue}b&\color{blue}c&&&&\cr
&&\color{blue}a&\color{blue}a&\color{blue}c&&\cr
&&&&\color{blue}b&\color{blue}b&\cr
&&&&&\color{blue}a&\color{blue}a\cr
\end{Young}
={\spadesuit\blacklozenge{\rm T}}
\end{array}$$
\end{example}
The $4\pi/3$ radians counterclockwise rotation symmetry $\clubsuit\blacklozenge=\blacklozenge\spadesuit$ can be performed in a
similar way, considering  $Y(\mu)^{\rm a}$ as the filling of the
outer shape $\mu$ after rotating $T$ by $\pi$ radians.
We may now observe that all  these symmetries can be exhibited using  the {\em hybrid} tableau--switching
involution  where the   slides are  executed without scanning the neighbours.

The
$H$ action on puzzles (equivalently mosaics) or LR--tableaux is defined by the group with presentation
$\eqref{hgroup1}$ or $\begin{matrix}
\langle\clubsuit,\spadesuit:\clubsuit^2 =\spadesuit^2=1=(\spadesuit\clubsuit)^3\rangle\end{matrix}$
 $=\{{\bf 1},\spadesuit,\clubsuit,
\spadesuit\clubsuit,$ $ \clubsuit\spadesuit,$ $\clubsuit\spadesuit\clubsuit=\spadesuit\clubsuit\spadesuit\}\simeq {D}_3,$
 where $\blacklozenge$, $\spadesuit$ and $\clubsuit$  are given for LR tableaux
 in \eqref{orthogonal} and definitions \ref{D:spade} and \ref{D:club} respectively. From
 Theorem \ref{th:lcost} (which agrees with the computational cost of Algorithm \ref{alg}), we may say that $H$ is  a linear time subgroup  of
index $2$ of $\mathbb Z_2\times \mathfrak{S}_3$. In particular, the $
\mathfrak{S}_3$ index two subgroup of symmetries  given by the cyclic group
$R\simeq\{\spadesuit\clubsuit:(\spadesuit\clubsuit)^3=1\}=\{{\bf 1},\spadesuit\clubsuit, (\spadesuit\clubsuit)^2=\clubsuit\spadesuit\}$  is of linear cost as already shown in \cite{PV1}.
Since $\clubsuit\spadesuit\clubsuit=\spadesuit\clubsuit\spadesuit$ is equivalent to $(\spadesuit\clubsuit)^3=1$, from
propositions~\ref{P:spade} and~\ref{P:club}, we may conclude the following identity  in the hybrid threefold
 $$\textbf{s}_1\textbf{s}_2\textbf{s}_1( Y(\mu^t), { T}^t, Y(\lambda^t)^{\rm
a})=\textbf{s}_2\textbf{s}_1\textbf{s}_2( Y(\mu^t), { T}^t, Y(\lambda^t)^{\rm
a})=( Y(\lambda^t), T^{\clubsuit\spadesuit\clubsuit}=T^{\spadesuit\clubsuit\spadesuit},  Y(\mu^t)^{\rm
a}).$$   Therefore, switching in a three fold multi-tableau consisting of Yamanouchi tableaux on the left and right and a standard tableau in the middle satisfy braid relations, $\textbf{s}_1\textbf{s}_2\textbf{s}_1( Y(\alpha),$ $ U, Y(\gamma)^{\rm
a})$ $=\textbf{s}_2 \textbf{s}_1 \textbf{s}_2( Y(\alpha), U, Y(\gamma)^{\rm
a})$ where $U$ is a standard tableau. \begin{remark}In general, braid relations are not satisfied and the result depends on the factorization of the permutation, see \cite[Lemma 3.2, Section 3]{bss}.
\end{remark}

\subsection{The  LR companion tableau $H$--symmetries} \label{compsym} As in the case of the map $\blacklozenge$,  subsection \ref{blacklozenge}, the linear cost of the maps $\spadesuit$ and $\clubsuit$ allows to give  bijections between the companion tableaux of $T$ and $T^\spadesuit$  or $T^\clubsuit$ whenever $T\in \rm LR(\mu,\nu,\lambda^\vee)$. We describe the procedure for the map $\spadesuit$ on companion tableaux.  We construct a bijection between the sets of companion LR tableaux  of shape $\nu$ and weight $\lambda/\mu$ and those   of shape $\mu^t$ and content $\lambda^t/\nu^t$,

\begin{equation}\label{compspade}\spadesuit: \rm LR_{\nu,\lambda/\mu}\rightarrow \rm LR_{\mu^t,\lambda^t/\nu^t}, G_\nu\mapsto G_{\mu^t}^{\spadesuit}\end{equation}
  where $G$ is the  right LR companion of some $T\in \rm LR(\mu,\nu,\lambda^\vee)$  and $\iota(T^\spadesuit)=(\iota(T))^\spadesuit$, that is, the recording matrix of $\spadesuit G$ is the transpose of the recording matrix of $\spadesuit T$, equivalently, $G$ is the LR companion of $T$ if and only if $\spadesuit G$ is the LR companion of $\spadesuit T$. The construction has the following three steps.

\begin{alg}\label{alg:spade}[Construction of $G^\spadesuit$.] Let $G\in \rm LR_{\nu,\lambda/\mu}$ and $T\in \rm LR(\mu,\nu,\lambda^\vee)$ with right LR companion $G$. The construction of $G^{\spadesuit}$ has the following three steps on the track of  Definition  \ref{D:spade} to construct $T^\spadesuit$.

{\em Step 1}: For $i=1,\dots, \ell(\lambda)$, consider the $i$-horizontal strip (the horizontal strip consisting of all boxes filled with $i$) of $G$ of size $\lambda_i-\mu_i$ and replace the entries with
$ \mu_i+1,\mu_i+2,\dots,\lambda_i$, scanned from SE to NW. In the resulting filling of shape $\nu$ and content $\lambda^t/\mu^t$, sort  by decreasing order the entries of the rows to obtain the companion plane partition $C$ of $T^t$.

The entries in row $r$ of $C$ tell  which columns of $T$ contain $r$ as an entry.

{\em Step 2}:  For $k=1,\dots, \ell(\lambda^t)$, let $R_k$ consists
 of the row indices of $C$ containing the entry $k$,
$$R_k:=\{r\in\{1,\dots,\ell(\nu)\}: \text{$k$ is an entry in row $r$ of $C$}\},$$
with $\#R_k=\lambda^t_k-\mu_k^t$. $R_k$ consists of the entries in the column $k$ of $T$.
{ Put
$$F_k:=[ \lambda_k^t]\setminus R_k,$$
with $\#F_k=\mu_k^t$. Note that $F_k=\emptyset$, for $k>\ell(\mu^t)$.}

{\em Step 3}:  Define the filling of shape $\mu^t=(|F_1|,\dots,|F_{\mu_k^t}|,0,\dots,0)$ whose $i$-th row consists of the elements in $F_i$ by decreasing order, for $i=1,\dots,\ell(\mu^t)$. { Its content is
$(\epsilon_j)_{j=1}^{\lambda_1^t}$  the multiplicity vector} of the multiset union $\bigcup_{i=1}^{\ell(\mu^t)} F_i$.
 For each $j=1,\dots,\lambda_1^t$, replace the entries of the $j$-vertical strip of size $\epsilon_j$ with $\nu_j+1,\dots,\nu_j+\epsilon_j$ scanned from bottom to top. The resulting tableau is $G^\spadesuit$ of shape $\mu^t$ and content $\lambda^t/\nu^t$.
 \end{alg}

 It is easy to see that if $G$ is the companion of $T\in LR(\mu,\nu,\lambda)$, then our construction of $G^\spadesuit$ gives the companion of $T^\spadesuit\in LR(\nu^t,\mu^t,\lambda^t)$.

 \begin{example}\label{G:spade}Consider the LR tableau $T$ of shape $\lambda^\vee/\mu$  and content $\nu$  in
  Example \ref{E:spade},  where $d=4$,  $\mu=4210$, $\nu=5420$, $\lambda=5320$, $\lambda^\vee=7542$; $\mu^t=3211000$, $\nu^t=3322100$, ${\lambda^\vee}^t=4433211$, $\ell(\lambda^{\vee t})=7$, and  its  right companion tableau $G \in LR_{\nu,\lambda^\vee/\mu} $ of shape $\nu$ and content $\lambda^\vee/\mu$,
 $$\begin{array} {cccc}T=\begin{Young}
1&3&&&&&\cr &2&2&3&&&\cr &&1&2&2&&\cr &&&&1&1&1\cr
\end{Young}&&G=\begin{Young}&&&&&&\cr 3&4&&&&&\cr 2&2&3&3&&&\cr 1&1&1&2&4&&\cr
\end{Young} \end{array}.$$
 We explain Algorithm \ref{alg:spade} as a track of Example \ref{E:spade}.

 Step $1$ produces the companion plane partition $C$ of $T^t$ of shape $\nu$ and content $\alpha:=$ $\lambda^{\vee t}/\mu^t=1222211$,
$$\begin{array}{cccccc}T^t=\begin{Young}
1&&&\cr
1&&&\cr
1&2&&\cr
&2&{\small 3}&\cr
&1&2&\cr
 &&2&{\small3}\cr
&&&1\cr
\end{Young}&& C=\begin{Young} &&&&&&\cr 4&2&&&&&\cr 5&4&3&2&&&\cr 7&6&5&3&1&&\cr
\end{Young}\end{array}$$
The entries in row $r$ of $C$ tell  which rows of $T^t $ contain $r$ as an entry.

Step 2 produces the sets $R_1=\{1\}, R_2=\{2,3\}$, $R_3=\{1,2\}$, $R_4=\{2,3\}$, $R_5=\{1,2\}$, $R_6=\{1\}$, $R_7=\{1\}$ where $R_k$ consists of the row indices of $C$ containing the entry $k$, for $k=1,\dots,7=\ell(\lambda^{\vee t})$.
$R_k$ also consists of the entries in  column $k$ of $T$ (that is, in row $k$ of $T^t$).

Let
$F_1=[4]\setminus R_1=\{2,3,4\}, F_2=[4]\setminus R_2=\{1,4\}$, $F_3=[3]\setminus \{1,2\}=\{3\}$, $F_4=[3]\setminus\{2,3\}=\{1\}$, $\ell(\mu^t)=4$.
That is, $\mu^t=(|F_1|, |F_2|, |F_3|,|F_4|,0,0,0)$.

Step 3 constructs  the filling of shape $\mu^t$, below on the right, whose row $i$ entries,  bottom to top,  consist of the elements in $F_i$, $i=1,\dots,\ell(\mu^t)=4$, with content $(1^2,2,3^2,4^2,0^3)$,
$$\begin{array}{ccccc}\begin{Young}
1&&&\cr
1&&&\cr
1&2&&\cr
{\color{blue}d}&2&{\small 3}&\cr
1&2&{\color{blue}c}&\cr
{\color{blue}b}&2&{\small3}&{\color{blue}b}\cr
1&{\color{blue}a}&{\color{blue}a}&{\color{blue}a}\cr
\end{Young}&&&\begin{Young}&&&\cr &&&\cr  &&&\cr  1&&&\cr 3&&&\cr 4&1&&\cr 4&3&2&\cr
\end{Young}
\end{array}$$
The number of $j$'s in the $j$th column of the tableau on the left is $\nu_j$, for $j=1,\dots,\ell(\nu)$. The second part of Step 3, equivalent to push down the $\nu_j$'s $j$'s in column  $j$ of the tableau above on the left
$$\begin{array}{ccccccc}\begin{Young}
{\color{blue}d}&     &&\cr
{\color{blue}b}&     &&\cr
1&{\color{blue}a}&               &\cr
1              &2&{\color{blue}c}&\cr
1              &2 &{\color{blue}a}&\cr
1 &2&{\small 3} &{\color{blue}b}\cr
1              &2 &{\small 3}&{\color{blue}a}\cr
\end{Young}, &&\text{constructs}&& G^\spadesuit=\begin{Young}&&&\cr &&&\cr &&&\cr 7&&&\cr 4&&&\cr 2&6&&\cr 1&3&5&\cr
\end{Young} \end{array}$$ of shape $\mu^t=3211$ and content $\lambda^{\vee t}/\nu^t=1^7$,
 the companion tableau of $T^\spadesuit.$ Note that $A=\left(\begin{smallmatrix}
 3&0&0&0\\
 1&2&0&0\\
 0&2&1&0\\
 1&0&1&0\\
  \end{smallmatrix}\right)$ is the recording matrix of $T$, $A^t$ is the recording matrix of $G$, $B=\left(\begin{smallmatrix}
 1&0&0&0&0&0&0\\
 0&1&0&0&0&0&0\\
 0&0&1&0&0&0&0\\
 1&0&0&0&0&0&0\\
 0&0&0&1&0&0&0
  \end{smallmatrix}\right)$ is the recording matrix of $\spadesuit T$ and $B^t$ is the recording matrix of $\spadesuit G$, that is $\spadesuit\iota(T)=\iota (\spadesuit T)$.
 \end{example}

\subsection{The symmetries outside of $H$}
\label{sec:check}  We have  discussed the $H$-symmetries of  puzzles, LR tableaux, and LR companion  tableaux and how do they do translate to each other. We
  now discuss the symmetries under the action of the other coset in $\mathbb{Z}_2\times\mathfrak{S}_3/H$, that is, the coset $\mathbb{Z}_2\times\mathfrak{S}_3-H=\zeta H=H\zeta\neq H$, for $\zeta=\varsigma_1,$ $\varsigma_2,\tau$, $\varsigma
=\varsigma_1\varsigma_2\varsigma_1$, $\tau\varsigma_2\varsigma_1$, $\tau\varsigma_1\varsigma_2$.  On KTW puzzles the symmetries outside of $H$ are also explained by the \emph{migration} on Purbhoo mosaics which correspond to \textit{jeu de taquin} slides and tableau switching on LR tableaux.
 Recall theorems \ref{fundamentalrev}, \ref{plr*}, and the LR commuter involution $ \rho=\bullet e$,  $\rho: \rm LR(\mu,\nu,\lambda)\rightarrow \rm LR(\lambda,\nu,\mu)$, $\;\rho(T) =\bullet e\,T$, and the LR transposer (involution) $\varrho=\blacklozenge\rho$, $\varrho:\rm LR(\mu,\nu,\lambda)\rightarrow \rm LR(\mu^t,\nu^t,\lambda^t)$, $\;\varrho(T) =\blacklozenge\bullet eT$.
Let us now denote by $\rho_1$ and $\rho_2$ the  LR commutativity bijections
$\rho_1: \rm LR(\mu,\nu,\lambda) \rightarrow \rm LR(\nu,\mu,\lambda)$ and $\rho_2:\rm LR(\mu,\nu,\lambda)\rightarrow \rm LR(\mu,\lambda,\nu)$ defined by the tableau switching involutions $s_1$ and $s_2$ respectively.
In \cite{az17} it has been shown that all known   LR commuters, known in the literature, exhibiting the identity $c_{\mu\,\nu\,\lambda}=c_{\nu\,\mu\,\lambda}$ coincide with the tableau switching involution, that is, with $\rho_1$.
{ In \cite[Corollary 3.4]{mosaic}, the proof  of the commutativity for LR tableaux, using \emph{migration} on mosaics, is equivalent to \emph{tableau switching}.}
We have  seen in theorems \ref{plr*} and \ref{coinctransposer} that the LR commuter $\rho$ and the LR transposer $\varrho$ on LR tableaux are related through the linear cost involution $\blacklozenge$, $\varrho=\blacklozenge\rho=\rho\blacklozenge$. We  next show that the same holds for the LR commuters $\rho_1$, $\rho_2$ and $\rho$ and the LR transposer $\varrho$ via the linear cost involutions $\spadesuit$, $\clubsuit$ and $\blacklozenge$ in $H$,
\begin{equation}\label{fundamental}
\spadesuit\,\rho_{1}=\rho_{1}\,\spadesuit=\blacklozenge\,\rho=\rho\,\blacklozenge=\varrho=\clubsuit\,\rho_2=\rho_2\clubsuit.
\end{equation}

\begin{theorem}\label{th:fundamental} Consider the set up as above. Then

$(a)$ $\rho_{1}=\spadesuit\blacklozenge\,\rho=\clubsuit\spadesuit\,\rho=$ $\blacklozenge\clubsuit\rho=\spadesuit\varrho$.

$(b)$ $\rho_{2}=\blacklozenge\spadesuit\,\rho=\spadesuit\clubsuit\,\rho=\clubsuit\blacklozenge\,\rho=\clubsuit\varrho=
\blacklozenge\spadesuit\blacklozenge\,\varrho$.

All known LR commuters and LR transposers are linear time reducible to each other and to the tableau switching involution, in particular, to the   reversal involution, equivalently, Sch\"utzenberger evacuation involution.
\end{theorem}
\begin{proof} We use the LR tableau model.

 $(a)$ We prove $\spadesuit\,\rho_{1}=\blacklozenge\,\rho=\varrho$. The key ingredient is to show that, in Algorithm \ref{alg:KD}, the second stage of the calculation of the reversal  of an LR tableau, performed by the switching involution { $s_1$} on the right hand side of \eqref{linearcoststepnew} below,  can be performed in a special way, without scanning the neighbours, and thus is a linear cost involution. That is, when $W=Y(\mu)$ is a Yamanouchi tableau,  $Q\equiv Y(\nu)$ is an LR tableau, and   $V=Q^{\rm n E}=Y(\nu^\bullet)$, Algorithm \ref{alg:KD} calculates $Q^e=[Q^{\rm n E}]_K\cap [Q]_{dK}=[Y(\nu^\bullet)]_K\cap [Q]_{dK}$,

\begin{equation}\begin{matrix}\label{linearcoststepnew}
Y(\mu)\cup Q &&Y(\mu)\cup Q^e\\
{\scriptstyle \bf s_1}\!\downarrow&&\uparrow\!{\scriptstyle \bf s_1}\\
Y(\nu)\cup \rho_1(Q)&\rightarrow&Y(\nu^\bullet)\cup \rho_1(Q).
\end{matrix}\end{equation}

First, we observe that the last step of \eqref{linearcoststepnew}, $Y(\nu^\bullet)\cup \rho_1(Q)\overset{\text{$s_1$}}\longleftrightarrow Y(\mu)\cup Q^e$  can be performed by a linear cost map as follows: first, for $i=1,\dots,\ell(\mu)$, the $i$-horizontal strip of length $\mu_i$ of $\rho_1(Q)$, an LR tableau of weight $\mu$, slides down to the $i$th row (see Example \ref{ex:fish}); then, for $i=1,\dots,\ell(\mu)$, sliding horizontally, justify to the left the $\mu_i$, $i$'s, to get $Y(\mu)$. Simultaneously it produces $Q^e$. { This is possible thanks to the filling of inner tableau, the reverse Yamanouchi  $Y(\nu^\bullet)$.

Let $\nu^t=(n_1,n_2,\dots,n_{\nu_1})$. The column $i$ of $Y(\nu^\bullet)$ is $C_i:=n_1>n_1-1>\cdots>n_1-n_i+1$ for $i=1,\dots,\nu_1$. When the $1$ in column $i$ of $\rho_1(Q)$ slides down to the first row of $Y(\nu^\bullet)$, column $i$ of $Y(\nu^\bullet)$ is shifted one box up. If there is an entry in column $i-1$ of $\rho_1(Q)$, an LR tableau, next to the left of $1$, then this entry is also $1$. In this case, $C_{i-1}=C_i$ and $C_{i-1}$ is also shifted one box up. If there is no entry in column $i-1$ of $\rho_1(Q)$, next to the left of $1$, then $\ell(C_{i-1}>\ell(C_i)$ and the entry of $C_{i-1}$ next to the left of $1$ is $\le n_1$. Thus when $C_i$ is shifted one box up the semistandardness along rows is preserved. Since $\ell(C_i)\ge \ell(C_{i+1})$, when column $C_i$ is shifted one box up the semistandardness along rows is preserved. Therefore when the $1$'s of $\rho_1(Q)$ slide down to the first row of $Y(\nu^\bullet)$, we get a perforated tableau pair. This perforated tableau pair has the following property: ignoring the first row the result is a tableau pair consisting of an opposite Yamanouchi tableau and the LR tableau obtained from $\rho_1(Q)$ removing the border strip consisting of $1$'s. By induction on $\ell(\mu)$ we conclude the validity of the procedure above. This is illustrated in \eqref{lozenge1} and \eqref{elozenge}.}

Second, thanks to remarks \ref{re:crocodile} and \ref{remark:evarrho}, for $i=1,\dots, \ell(\nu)$, replacing,  from NW to SE, the entries of the $i$-horizontal strip of length $\ell(\nu)-i+1$ in $Q^e$,  with $1,\dots,\ell(\nu)-i+1$, gives $ (\blacklozenge\rho(Q))^t$. We then get the sequence
\begin{equation}\label{lozenge0}Y(\nu^\bullet)\cup \rho_1(Q)\overset{\text{$s_1$}}\longleftrightarrow Y(\mu)\cup Q^e\longleftrightarrow Y(\mu)\cup (\blacklozenge\rho(Q))^t.\end{equation}

In fact, to obtain $(\blacklozenge\rho(Q))^t$,  we do not need the last step $Q^e \longleftrightarrow (\blacklozenge\rho(Q))^t$ in the sequence \eqref{lozenge0}.
Replace on the left hand side of \eqref{lozenge0}
 $Y(\nu^\bullet)$ with $[Y(\nu^t)]^t$ (they have  the same shape $\nu$) and  apply the map involution  $\spadesuit$ without transposing (see Definition \ref{D:spade}) to $[Y(\nu^t)]^t\cup\rho_1(Q)$. The  sequence of instructions   dictated by  $\spadesuit$, without transposing, on $[Y(\nu^t)]^t\cup \rho_1(Q)$, is exactly the same as the ones described above to realize  $Y(\nu^\bullet)\cup\rho_1(Q)\overset{\text{$s_1$}}\longleftrightarrow Y(\mu)\cup Q^e$ in \eqref{lozenge0}. At the end $\spadesuit$ produces  $Y(\mu)\cup (\blacklozenge\rho(Q))^t$ instead of $Y(\mu)\cup Q^e$, that is,
$$[Y(\nu^t)]^t\cup\rho_1(Q) \overset{\spadesuit}\rightarrow Y(\mu)\cup (\spadesuit\rho_1(Q))^t= Y(\mu)\cup (\blacklozenge\rho(Q))^t.$$
Henceforth, $\spadesuit\rho_1 =\blacklozenge\rho=\varrho$.

$(b)$ We prove $\blacklozenge\,\rho=\clubsuit\,\rho_2$.  Thanks to Remark \ref{re:KD}, Algorithm \ref{alg:KD} also calculates $Q^e=[Y(\nu^\bullet)^{\rm a}]_K\cap [Q]_{dK}$,

\begin{equation}\begin{matrix}\label{linearcoststep2b}
Q\cup Y(\lambda)^{\rm a} && Q^e\cup Y(\lambda)^{\rm a}\\
{\scriptstyle \bf s_2}\!\downarrow&&\uparrow\!{\scriptstyle \bf s_2}\\
\rho_2(Q)\cup Y(\nu)^{\rm a}&\rightarrow& \rho_2(Q)\cup Y(\nu^\bullet)^{\rm a}.
\end{matrix}\end{equation}

First, we observe that the last step of \eqref{linearcoststep2b}, $ \rho_2(Q)\cup Y(\nu^\bullet)^{\rm a}\overset{\text{$s_2$}}\longleftrightarrow  Q^e\cup Y(\lambda)^{\rm a}$  can be performed by a linear cost map as follows. Let $\lambda^t=(m_1,m_2,\dots,m_{\lambda_1})$.
For $i=1,\dots,m_1$, slide horizontally the rightmost $i$  { of $\rho_2(Q)$ } to the $(n-d)$th column { of $D$}, and put $C_1:=12\cdots m_1$; for $i=1,\dots,m_2$, slide horizontally the rightmost $i$ { of $\rho_2(Q)$ }, in the first $(n-d-1)$ columns,  to the $(n-d-1)$th column { of $D$}, and put $C_2:=12\cdots m_2$;  $\dots$; lastly, slide horizontally the remaining $1,2,\dots,m_{\lambda_1}$ { of $\rho_2(Q)$ } to the $(n-d-\lambda_1+1)$th column { of $D$, and put $C_{\lambda_1}:=12\cdots m_{\lambda_1}$}. Then, for each for $j=1,2,\dots,\lambda_1$,  slide up vertically the column word $C_j:=12\cdots m_j$ in the $n-d-j+1$ of $D$, to justify on the Northeast corner of $D$ the $\lambda_i$, $i$'s { of $\rho_2(Q)$ } so that we get $Y(\lambda)^{\rm a}$. Simultaneously one produces $Q^e$. The procedure is illustrated in Example \ref{ex:dolphin}.


Second, for $i=1,\dots, \ell(\nu)$, replacing,  from NW to SE, the entries of the $i$-horizontal strip of length $\nu_{\ell(\nu)-i+1}$ in $Q^e$,  with $1,\dots,\nu_{\ell(\nu)-i+1}$, gives $ (\blacklozenge\rho(Q))^t$. We then get the sequence
\begin{equation}\label{lozengeb} \rho_2(Q)\cup Y(\nu^\bullet)^{\rm a}\overset{\text{$s_2$}}\longleftrightarrow
Q^e\cup Y(\lambda)^{\rm a} \longleftrightarrow  (\blacklozenge\rho(Q))^t\cup Y(\lambda)^{\rm a}.\end{equation}

In fact, to obtain $(\blacklozenge\rho(Q))^t$,  we do not need the last step $Q^e \longleftrightarrow (\blacklozenge\rho(Q))^t$ in the sequence \eqref{lozengeb}.
Replace on the left hand side of \eqref{lozengeb}
 $Y(\nu^\bullet)^{\rm a}$ with $[Y(\nu^t)^{\rm a}]^t$ (they have  the same anti normal shape $rev \,\nu$) and  apply the map involution $\clubsuit$ without transposing (see  Definition \ref{D:club}) to $\rho_2(Q)$. The  sequence of instructions   dictated by  $\clubsuit$, without transposing, on $ \rho_2(Q)\cup [Y(\nu^t)^{\rm a}]^t$, is exactly the same as the ones described above to produce  $\rho_2(Q)\cup Y(\nu^\bullet)^{\rm a}\overset{\text{$s_2$}}\longleftrightarrow
Q^e\cup Y(\lambda)^{\rm a}$ in \eqref{lozengeb}. At the end $\clubsuit$ returns $(\blacklozenge\rho(Q))^t\cup Y(\lambda)^{\rm a}$ instead of $Q^e\cup Y(\lambda)^{\rm a}$, that is,
$$\rho_2(Q)\cup [Y(\nu^t)^{\rm a}]^t \overset{\clubsuit}\rightarrow  (\clubsuit\,\rho_2(Q))^t\cup Y(\lambda)^{\rm a}=  (\blacklozenge\rho_2(Q))^t\cup Y(\lambda)^{\rm a}.$$

\end{proof}

\begin{example}\label{ex:fish} Let $\nu=552$ and  $\mu=5321$, $Q\equiv Y(\nu)$ and
$Y(\mu)\cup Q\overset{\text{$s_1$}}\longleftrightarrow Y(\nu)\cup \rho_1(Q)$, where $\rho_1(Q)\equiv Y(\mu)$. Then

\begin{equation}Y(\nu^\bullet)\cup\rho_1(Q)=
\begin{Young}
 \olga{1}&\olga{2}&\olga{3}&\olga{4}&&&\cr
 3&3&\olga{1}&\olga{1}&\olga{2}&\olga{3}&\cr
 2&2&3&3&3&\olga{2}&\cr
  1&1&2&2&2&\olga{1}&\olga{1}\cr
\end{Young}\leftrightarrow \begin{Young}
 3&\olga{2}&\olga{3}&\olga{4}&&&\cr
 2&3&3&3&\olga{2}&\olga{3}&\cr
 1&2&2&2&3&\olga{2}&\cr
  \olga{1}&1&\olga{1}&\olga{1}&2&\olga{1}&\olga{1}\cr
\end{Young}\leftrightarrow
\begin{Young}
 3&3&\olga{3}&\olga{4}&&&\cr
 2&2&3&3&3&\olga{3}&\cr
 1&\olga{2}&2&2&\olga{2}&\olga{2}&\cr
  \olga{1}&1&\olga{1}&\olga{1}&2&\olga{1}&\olga{1}\cr
\end{Young}\label{lozenge1}\end{equation}
\begin{equation}\leftrightarrow\begin{Young}
 3&3&3&\olga{4}&&&\cr
 2&2&\olga{3}&3&3&\olga{3}&\cr
 1&\olga{2}&2&2&\olga{2}&\olga{2}&\cr
  \olga{1}&1&\olga{1}&\olga{1}&2&\olga{1}&\olga{1}\cr
\end{Young}
\leftrightarrow Y(\mu)\cup Q^e=\begin{Young}
 \olga{4}&3&3&3&&&\cr
 \olga{3}&\olga{3}&2&2&3&3&\cr
 \olga{2}&\olga{2}&\olga{2}&1&2&2&\cr
  \olga{1}&\olga{1}&\olga{1}&\olga{1}&\olga{1}&1&2\cr
\end{Young}\label{elozenge}\end{equation}
On the righthand of \eqref{elozenge}, replacing, from NW to SE, the entries of the $(\ell(\nu)-i+1)$-horizontal strip of $Q^e$, with $1,\dots,\nu_i$, for $i=1,\dots, \ell(\nu)$, gives $$\leftrightarrow \begin{Young}
 \olga{4}&1&2&3&&&\cr
 \olga{3}&\olga{3}&1&2&4&5&\cr
 \olga{2}&\olga{2}&\olga{2}&1&3&4&\cr
  \olga{1}&\olga{1}&\olga{1}&\olga{1}&\olga{1}&2&5\cr
\end{Young}=Y(\mu)\cup (\blacklozenge\rho(Q))^t=Y(\mu)\cup (\varrho(Q))^t.$$

This is equivalent to replace, on the left hand side of \eqref{lozenge1},
 $Y(\nu^\bullet)$ with $[Y(\nu^t)]^t$ (they have  the same shape $\nu$) and  apply the map involution  $\spadesuit$, Definition \ref{D:spade}, to $[Y(\nu^t)]^t\cup\rho_1(Q)$. It produces the same sequence of steps as in \eqref{lozenge1} and \eqref{elozenge} with the appropriate replacement of $Y(\nu^\bullet)$ by $[Y(\nu^t)]^t$,
$$[Y(\nu^t)]^t\cup\rho_1(Q)=
\begin{Young}
 \olga{1}&\olga{2}&\olga{3}&\olga{4}&&&\cr
 1&2&\olga{1}&\olga{1}&\olga{2}&\olga{3}&\cr
 1&2&3&4&5&\olga{2}&\cr
 1&2&3&4&5&\olga{1}&\olga{1}\cr
\end{Young}
\leftrightarrow \begin{Young}
 1&\olga{2}&\olga{3}&\olga{4}&&&\cr
 1&2&3&4&\olga{2}&\olga{3}&\cr
 1&2&3&4&5&\olga{2}&\cr
  \olga{1}&2&\olga{1}&\olga{1}&5&\olga{1}&\olga{1}\cr
\end{Young}\leftrightarrow
\begin{Young}
 1&2&\olga{3}&\olga{4}&&&\cr
 1&2&3&3&5&\olga{3}&\cr
 1&\olga{2}&3&4&\olga{2}&\olga{2}&\cr
  \olga{1}&2&\olga{1}&\olga{1}&5&\olga{1}&\olga{1}\cr
\end{Young}$$
$$\leftrightarrow\begin{Young}
 1&2&3&\olga{4}&&&\cr
 1&2&\olga{3}&4&5&\olga{3}&\cr
 1&\olga{2}&3&4&\olga{2}&\olga{2}&\cr
  \olga{1}&2&\olga{1}&\olga{1}&5&\olga{1}&\olga{1}\cr
\end{Young}
\leftrightarrow \begin{Young}
 \olga{4}&1&2&3&&&\cr
 \olga{3}&\olga{3}&1&2&4&5&\cr
 \olga{2}&\olga{2}&\olga{2}&1&3&4&\cr
  \olga{1}&\olga{1}&\olga{1}&\olga{1}&\olga{1}&2&5\cr
\end{Young}=Y(\mu)\cup (\spadesuit\rho_1(Q) )^t =Y(\mu)\cup (\blacklozenge\rho(Q))^t.
$$
Henceforth, $\spadesuit\rho_1(Q) =\blacklozenge\rho(Q)=\varrho(Q)$.

\end{example}

\begin{example}\label{ex:dolphin} Let $\lambda=542$, $\lambda^t=33221$ and  $\nu^\bullet=235$, $Q\equiv Y(\nu)$ and
$ \rho_2(Q)\cup Y(\nu^\bullet)^{\rm a}\overset{\text{$s_2$}}\longleftrightarrow  Q^e\cup Y(\lambda)^{\rm a}$, where $\rho_2(Q)\equiv  Y(\lambda)$. Then

\begin{equation}\rho_2(Q)\cup Y(\nu^\bullet)^{\rm a}=
\begin{Young}
 \olga{1}&\olga{3}&3&3&3&3&3\cr
 &\olga{2}&\olga{2}&\olga{3}&2&2&2\cr
  &&\olga{1}&\olga{2}&\olga{2}&1&1\cr
  &&&&\olga{1}&\olga{1}&\olga{1}\cr
\end{Young}\leftrightarrow \begin{Young}
 \olga{1}&\olga{3}&3&3&3&3&3\cr
 &\olga{2}&\olga{2}&2&2&2&\olga{3}\cr
  &&\olga{1}&\olga{2}&1&1&\olga{2}\cr
  &&&&\olga{1}&\olga{1}&\olga{1}\cr
\end{Young}\leftrightarrow
\begin{Young}
 \olga{1}&3&3&3&3&\olga{3}&3\cr
 &\olga{2}&\olga{2}&2&2&2&\olga{3}\cr
  &&\olga{1}&1&1&\olga{2}&\olga{2}\cr
  &&&&\olga{1}&\olga{1}&\olga{1}\cr
\end{Young}\label{dolphin1}\end{equation}
\begin{equation}\leftrightarrow\begin{Young}
 \olga{1}&3&3&3&3&\olga{3}&3\cr
 &\olga{2}&2&2&\olga{2}&2&\olga{3}\cr
  &&\olga{1}&1&1&\olga{2}&\olga{2}\cr
  &&&&\olga{1}&\olga{1}&\olga{1}\cr
\end{Young}
\leftrightarrow \begin{Young}
 \olga{1}&3&3&3&3&\olga{3}&3\cr
 &2&2&\olga{2}&\olga{2}&2&\olga{3}\cr
  &&1&\olga{1}&1&\olga{2}&\olga{2}\cr
  &&&&\olga{1}&\olga{1}&\olga{1}\cr
\end{Young}\leftrightarrow\begin{Young}
 3&3&\olga{1}&3&3&\olga{3}&3\cr
 &2&2&\olga{2}&\olga{2}&2&\olga{3}\cr
  &&1&\olga{1}&1&\olga{2}&\olga{2}\cr
  &&&&\olga{1}&\olga{1}&\olga{1}\cr
\end{Young}\nonumber \end{equation}
\begin{equation}\leftrightarrow Q^e\cup Y(\lambda)^{\rm a}=\begin{Young}
 3&3&\olga{1}&\olga{2}&\olga{2}&\olga{3}&\olga{3}\cr
 &2&2&\olga{1}&\olga{1}&\olga{2}&\olga{2}\cr
  &&1&3&3&\olga{1}&\olga{1}\cr
  &&&&1&2&3\cr
\end{Young}
 \label{dolphin2}\end{equation}
On the righthand of \eqref{dolphin2}, replacing, from NW to SE, the entries of the $i$-horizontal strip of $Q^e$, with $1,\dots,\nu_{\ell(\nu)-i+1}$, for $i=1,\dots, \ell(\nu)$, gives $$\leftrightarrow \begin{Young}
 1&2&\olga{1}&\olga{2}&\olga{2}&\olga{3}&\olga{3}\cr
 &1&2&\olga{1}&\olga{1}&\olga{2}&\olga{2}\cr
  &&1&3&4&\olga{1}&\olga{1}\cr
  &&&&2&3&5\cr
\end{Young}= (\blacklozenge\rho(Q))^t \cup Y(\lambda)^{\rm a}= (\varrho(Q))^t  \cup Y(\lambda)^{\rm a}.$$

This is equivalent to replace, on the left hand side of \eqref{dolphin1},  $Y(\nu^\bullet)^{\rm a}$ with $[Y(\nu^t)^{\rm a}]^t$
 (they have  the same antinormal shape $rev\,\nu$) and  apply the map involution  $\clubsuit$, Definition \ref{D:club}, to $\rho_2(Q)\cup[Y(\nu^t)]^t$. It produces the same sequence of steps as in \eqref{dolphin1} and \eqref{dolphin2} with the appropriate replacement of $Y(\nu^\bullet)^{\rm a}$ by $[Y(\nu^t)^{\rm a}]^t$,
$$\rho_2(Q)\cup [Y(\nu^t)^{\rm a}]^t=
\begin{Young}
 \olga{1}&\olga{3}&1&2&3&4&5\cr
 &\olga{2}&\olga{2}&\olga{3}&1&2&3\cr
  &&\olga{1}&\olga{2}&\olga{2}&1&2\cr
  &&&&\olga{1}&\olga{1}&\olga{1}\cr
\end{Young}
\leftrightarrow \begin{Young}
 \olga{1}&\olga{3}&1&2&3&4&5\cr
 &\olga{2}&\olga{2}&1&2&3&\olga{3}\cr
  &&\olga{1}&\olga{2}&1&2&\olga{2}\cr
  &&&&\olga{1}&\olga{1}&\olga{1}\cr
\end{Young}\leftrightarrow
\begin{Young}
 \olga{1}&1&2&3&4&\olga{3}&5\cr
 &\olga{2}&\olga{2}&1&2&3&\olga{3}\cr
  &&\olga{1}&1&2&\olga{2}&\olga{2}\cr
  &&&&\olga{1}&\olga{1}&\olga{1}\cr
\end{Young}
$$
$$\leftrightarrow
 \begin{Young}
 \olga{1}&1&2&3&4&\olga{3}&5\cr
 &\olga{2}&1&2&\olga{2}&3&\olga{3}\cr
  &&\olga{1}&1&2&\olga{2}&\olga{2}\cr
  &&&&\olga{1}&\olga{1}&\olga{1}\cr
\end{Young}\leftrightarrow\begin{Young}
 \olga{1}&1&2&3&4&\olga{3}&5\cr
&1&2&\olga{2}&\olga{2}&3&\olga{3}\cr
  &&1&\olga{1}&2&\olga{2}&\olga{2}\cr
  &&&&\olga{1}&\olga{1}&\olga{1}\cr
\end{Young} \leftrightarrow\begin{Young}
 \olga{1}&1&2&3&4&\olga{3}&5\cr
&1&2&\olga{2}&\olga{2}&3&\olga{3}\cr
  &&1&\olga{1}&2&\olga{2}&\olga{2}\cr
  &&&&\olga{1}&\olga{1}&\olga{1}\cr
\end{Young}\leftrightarrow\begin{Young}
 1&2&\olga{1}&3&4&\olga{3}&5\cr
&1&2&\olga{2}&\olga{2}&3&\olga{3}\cr
  &&1&\olga{1}&2&\olga{2}&\olga{2}\cr
  &&&&\olga{1}&\olga{1}&\olga{1}\cr
\end{Young}$$
$$
\leftrightarrow\begin{Young}
 1&2&\olga{1}&\olga{2}&\olga{2}&\olga{3}&\olga{3}\cr
 &1&2&\olga{1}&\olga{1}&\olga{2}&\olga{2}\cr
  &&1&3&4&\olga{1}&\olga{1}\cr
  &&&&2&3&5\cr
\end{Young}=(\clubsuit\,\rho_2(Q))^t\cup Y(\lambda)^{\rm a}=  (\blacklozenge\rho_2(Q))^t\cup Y(\lambda)^{\rm a}.
$$
Henceforth, $\clubsuit\rho_2(Q) =\blacklozenge\rho_2(Q)=\varrho(Q)$.
\end{example}


\begin{corollary}\label{cor:fund} 
Consider the symmetries outside of $H$.
 The following  holds  in $\mathcal{LR}$:

$(a)$ $\varrho=\blacklozenge\rho=\rho\blacklozenge$, and
$\rho=\varrho\blacklozenge=\blacklozenge\varrho$.

$(b)$ $\rho_1=\spadesuit\blacklozenge\rho=\rho\blacklozenge\spadesuit=\clubsuit\spadesuit\rho=\rho\spadesuit\clubsuit=\blacklozenge\clubsuit\rho$.
$=\rho\clubsuit\blacklozenge=\spadesuit\varrho=\varrho\spadesuit$.

$(c)$ $\rho_2=\rho\spadesuit\blacklozenge=\blacklozenge\spadesuit\rho=\rho\clubsuit\spadesuit=\spadesuit\clubsuit\rho=
\rho\blacklozenge\clubsuit=\clubsuit\blacklozenge\rho=$ $\varrho\clubsuit=\clubsuit\varrho$.

$(d)$ $\rho\spadesuit=\spadesuit\blacklozenge\spadesuit\rho=\clubsuit\rho=\blacklozenge\spadesuit\varrho=\varrho\blacklozenge\spadesuit$.

$(e)$ $\spadesuit\rho=\rho\clubsuit=\spadesuit\blacklozenge\varrho=\varrho\spadesuit\blacklozenge$.

$(f)$ $\rho_1\rho_2\rho_1=\rho_2\rho_1\rho_2=(\spadesuit\blacklozenge)^3\rho=\rho=\bullet\,e$.

$(g)$   $(\rho_1\rho_2)^3=1$.

$(h)$  $\rho_1$ and $\rho_2$ generate a  representation of $\mathfrak{S}_3$ in $\mathfrak{S}_{\mathcal{LR}}$.

\end{corollary}

\begin{proof} $(a)$  It  follows from Theorem \ref{plr*}.  $(b)$ and $(c)$   follow  from Theorem \ref{th:fundamental}, $(a)$ and  $(b)$, respectively: $\rho_1^2=$ $1=\spadesuit\blacklozenge\rho\spadesuit\blacklozenge\rho=\clubsuit\spadesuit\rho\clubsuit\spadesuit\rho=\blacklozenge\clubsuit
\rho\blacklozenge\clubsuit\rho=\varrho\spadesuit\varrho\spadesuit$ $\Leftrightarrow $ $ \rho\blacklozenge\spadesuit=\spadesuit\blacklozenge\rho=\clubsuit\spadesuit\rho=\spadesuit\clubsuit\rho=\blacklozenge\clubsuit\rho=
\rho\clubsuit\blacklozenge=\varrho\spadesuit=\spadesuit\varrho$;
 $\rho_2^2=$ $1=\blacklozenge\spadesuit\rho\blacklozenge\spadesuit\rho$ $\Leftrightarrow $ $ \blacklozenge\spadesuit\rho=\rho\spadesuit\blacklozenge$ and $\varrho\clubsuit=\rho\spadesuit\blacklozenge=\blacklozenge\spadesuit\rho$ $=\clubsuit\varrho$.
$(d)$From $(c)$,  $\spadesuit\rho\spadesuit=\spadesuit\blacklozenge\spadesuit\rho\blacklozenge=\clubsuit \rho \blacklozenge=\clubsuit  \blacklozenge\rho=\blacklozenge\spadesuit\rho.$
$(e)$ From $(d)$ and $(c)$, $\spadesuit\rho=\spadesuit\clubsuit\rho\spadesuit=\rho\blacklozenge\clubsuit\spadesuit=\rho\clubsuit\spadesuit\spadesuit=\rho\clubsuit$.
\end{proof}

\begin{remark}  From the previous corollary, we conclude that the tableau-switching composition $\textbf{s}_1 \textbf{s}_2 \textbf{s}_1=\textbf{s}_2 \textbf{s}_1 \textbf{s}_2$ on (three fold tableau model) $\mathcal{LR}$ satisfy braid relations and coincides with reversal composed with rotation. That is, if $T\in LR(\mu,\nu,\lambda)$ then \begin{eqnarray}\label{eq:switchbraid}&&\textbf{s}_1 \textbf{s}_2 \textbf{s}_1(Y(\mu)\cup T\cup Y(\lambda)^{\rm a})= \textbf{s}_2 \textbf{s}_1 \textbf{s}_2(Y(\mu)\cup T\cup Y(\lambda)^{\rm a})\nonumber\\
&=& Y(\lambda)\cup \rho_1\rho_2\rho_1(T)\cup Y(\mu)^{\rm a}=Y(\lambda)\cup \rho_2\rho_1\rho_2(T)\cup Y(\mu)^{\rm a}\nonumber\\
&=&Y(\lambda)\cup \rho(T)\cup Y(\mu)^{\rm a}
=Y(\lambda)\cup \bullet\,e(T)\cup Y(\mu)^{\rm a}.\end{eqnarray}
The same happens for the tableau switching on pairs of Yamanouchi  and antinormal Yamanouchi tableaux inside $D$. From \eqref{eq:switchbraid}, $\textbf{s}_1 \textbf{s}_2 \textbf{s}_1(\emptyset\cup Y(\mu)\cup Y(\mu^\vee)^{\rm a})=Y(\mu^\vee)\cup \bullet{\rm evac} Y(\mu)\cup \emptyset$
$=\textbf{s}_2 \textbf{s}_1 \textbf{s}_2(\ Y(\mu)\cup Y(\mu^\vee)^{\rm a}\cup \emptyset)=\emptyset\cup {\rm evac}\bullet Y(\mu^\vee)^{\rm a}\cup Y(\mu)^{\rm a}$. Therefore,
$Y(\mu^\vee)\cup \bullet{\rm evac} Y(\mu)={\rm evac}\bullet Y(\mu^\vee)^{\rm a}\cup Y(\mu)^{\rm a}$, and $Y(\mu)^{\rm a}=\bullet {\rm evac} Y(\mu)\Leftrightarrow Y(\mu) ={\rm evac}\bullet  Y(\mu)^{\rm a}$.
On the other hand, $\textbf{s}(Y(\mu)\cup Y(\mu^\vee)^{\rm a})=Y(\mu^\vee)\cup Y(\mu)^{\rm a}=Y(\mu^\vee)\cup \bullet {\rm evac} Y(\mu)$.
\end{remark}

Taking into account Theorem \ref{th:companion} and algorithms \ref{alg:blacklozenge} and  \ref{alg:spade}, replacing $\rho$ with ${\rm evac}$ in the previous corollary,  in particular, one has
\begin{corollary}
The following holds on $LR$ companions

$(a)$ $\blacklozenge{\rm evac}={\rm evac}\blacklozenge$.

$(b)$ $\spadesuit\blacklozenge{\rm evac}={\rm evac}\blacklozenge\spadesuit=\clubsuit\spadesuit{\rm evac}={\rm evac}\spadesuit\clubsuit=\blacklozenge\clubsuit{\rm evac}$
$={\rm evac}\clubsuit\blacklozenge$.

$(c)$ ${\rm evac}\spadesuit\blacklozenge=\blacklozenge\spadesuit{\rm evac}={\rm evac}\clubsuit\spadesuit=\spadesuit\clubsuit{\rm evac}=
{\rm evac}\blacklozenge\clubsuit=\clubsuit\blacklozenge{\rm evac}$.

$(d)$ ${\rm evac}\spadesuit=\spadesuit\blacklozenge\spadesuit{\rm evac}=\clubsuit{\rm evac}$.

$(e)$ $\spadesuit{\rm evac}={\rm evac}\clubsuit$.

\end{corollary}
\begin{remark}

From  Theorem \ref{th:companion} and algorithms \ref{alg:blacklozenge} and  \ref{alg:spade}, if $G$ is the companion of $T$, the companion of $\rho_1(T)=\spadesuit\blacklozenge \rho(T)=\rho\blacklozenge\spadesuit(T)$ is ${\rm evac}\blacklozenge\spadesuit G $ $=$ $\spadesuit\blacklozenge{\rm evac}G$. Similarly, for $\rho_2$ and $\varrho$. This identity does not depend on the straight shape tableau $G$ because every tableau of straight shape is the right companion of some LR tableau with  appropriate inner shape and outer shape. Henceforth, the identities above are valid for any tableaux of straight shape.
 \end{remark}

The cyclic group
$R\simeq\{\spadesuit\clubsuit:(\spadesuit\clubsuit)^3=1\}=\{{\bf 1},\spadesuit\clubsuit, (\spadesuit\clubsuit)^2=\clubsuit\spadesuit\}$ of rotation involutions, where $\clubsuit\blacklozenge=\blacklozenge\spadesuit=\spadesuit\clubsuit$ is the $4\pi/3$ radians counterclockwise rotation symmetry involution, and
$\blacklozenge\clubsuit=\spadesuit\blacklozenge=\clubsuit\spadesuit$ is the $2\pi/3$ radians counterclockwise rotation symmetry involution, is an index two subgroup of $\mathfrak{S}_3\simeq<\rho_1,\rho_2:\rho_1^2=\rho_2^2=(\rho_1\rho_2)^3=1>=\{1,\rho_1,\rho_2,\rho_1\rho_2,\rho_2\rho_1, \rho_1\rho_2\rho_1=\rho_2\rho_1\rho_2=\rho\}$ $=R\cup \rho_1 R$ with $\rho_1R=R\rho_1=\rho_2R=R\rho_2=$ $\rho R=R\rho=\{\rho_1,\rho_2,\rho\}$ $=\{\spadesuit\blacklozenge\rho,\blacklozenge\spadesuit\rho,\rho\}$ $=\{\spadesuit\blacklozenge\bullet e,\blacklozenge\spadesuit\bullet e,\bullet e\}$. That is,
\begin{equation}\label{s3}
\mathfrak{S}_3\simeq\{1,\spadesuit\blacklozenge, \blacklozenge\spadesuit,\spadesuit\blacklozenge\bullet e,\blacklozenge\spadesuit\bullet e,\bullet e\},
\end{equation}
and, from Corollary \ref{cor:fund}, $\blacklozenge\mathfrak{S}_3=\mathfrak{S}_3\blacklozenge=\blacklozenge\mathfrak{S}_3=\mathfrak{S}_3\spadesuit=
\spadesuit\mathfrak{S}_3=\mathfrak{S}_3\clubsuit=\clubsuit\mathfrak{S}_3$, where
\begin{equation}\label{outsides3}
\mathbb{ Z}_2\times \mathfrak{S}_3-\mathfrak{S}_3=\blacklozenge\mathfrak{S}_3=\mathfrak{S}_3\blacklozenge=\{ \blacklozenge,\clubsuit,\spadesuit,  \blacklozenge\spadesuit\varrho,\spadesuit\blacklozenge\varrho,\varrho\}.
\end{equation}
The faithful  permutation representation of $\mathbb{ Z}_2\times \mathfrak{S}_3$ affords the following presentation
$$\mathbb Z_2\times \mathfrak{S}_3=<\blacklozenge,\rho_1,\rho_2:\blacklozenge^2=\rho_1^2=\rho_2^2=(\rho_1\rho_2)^3=(\blacklozenge\rho_2)^2=(\blacklozenge\rho_1)^2=1>.$$

 Therefore, $\varsigma_1R=R\varsigma_1=\varsigma_2R$ $=R\varsigma_2=\varsigma R=R\varsigma=\mathfrak{S}_3-R$.

The symmetries outside of $H$,  commutativity and conjugation symmetries,  \eqref{outH}, are given by the involutions $\rho$, $\rho_1$, $\rho_2$ and $\varrho$, and
 the two  others combining rotation and conjugation,
 \eqref{outH1},
 $c_{\mu\,\nu\,\lambda}=c_{\lambda^t\,\mu^t\,\nu^t},$   $c_{\mu\,\nu\,\lambda}=c_{\nu^t\,\lambda^t\,\mu^t}$,
and are given by the bijections $\clubsuit\spadesuit\varrho=\blacklozenge\clubsuit\varrho=\spadesuit\blacklozenge\varrho=\spadesuit\rho$ and
$\spadesuit\clubsuit\varrho=\blacklozenge\clubsuit\varrho=\blacklozenge\spadesuit\varrho=\clubsuit \rho$ respectively. That is,
\begin{eqnarray}\label{outideH0}\mathbb Z_2\times \mathfrak{S}_3-H=\rho_1H=H\rho_1=\rho_2H=H\rho_2=\rho H=H\rho=\varrho H\nonumber\\
=H\varrho= \spadesuit\rho H=H\spadesuit\rho =\clubsuit\rho H=H\clubsuit\rho\nonumber\\ =\{\rho_1,\rho_2,\rho,\varrho,\clubsuit\rho,\spadesuit\rho\}=\{\rho_1,\rho_2,\rho,\blacklozenge\rho,\clubsuit\rho,\spadesuit\rho\}
\end{eqnarray}
Equivalently,
\begin{equation}\label{outsideH}\mathbb Z_2\times \mathfrak{S}_3-H=\{\spadesuit\blacklozenge\bullet e,\blacklozenge\spadesuit\bullet e,\bullet e,\blacklozenge\bullet e,\clubsuit\bullet e,\spadesuit\bullet e\}.\end{equation}
The
symmetries outside of $H$ are therefore linearly reducible to each
other, and, in particular, to the reversal involution $e$. On its
turn, as shown in Section 3.2, the reversal map $e$ is linearly
reducible to  the Sch\"utzenberger involution $E$.
The linear time maps
$\iota$ and $\blacklozenge\spadesuit\iota$ define bijections between tableaux of normal (straight) shape and LR tableaux, see~\cite{lee2,ouch,PV2}.
We can thus state the following result.

\begin{theorem}
LR transposers and commuters are linearly reducible to each other, in particular, to the   reversal involution $e$, equivalently, Sch\"utzenberger involution
$\rm evac$.
\end{theorem}
Hence, we may use two involutions of $H$ and an LR commuter or an LR transposer to realize the action of $\mathbb Z_2\times \mathfrak{S}_3$ on $\mathcal{LR}$. For instance, $\mathbb Z_2\times \mathfrak{S}_3$ has the following presentations
$$\mathbb Z_2\times \mathfrak{S}_3\simeq \langle\spadesuit,\blacklozenge,\varrho\rangle= \langle\spadesuit,\blacklozenge,{\varrho}:{\varrho}^2=\spadesuit^2=\blacklozenge^2=
(\spadesuit\blacklozenge)^3=(\spadesuit{\varrho})^2=(\blacklozenge{\varrho})^2=1\rangle$$
$$\;\qquad\qquad=\langle\blacklozenge,\spadesuit,\rho:\rho^2=\spadesuit^2=\blacklozenge^2=
(\spadesuit\blacklozenge)^3=(\spadesuit\blacklozenge\rho)^2=(\blacklozenge\rho)^2=1\rangle,$$
$$\;\qquad\qquad=\langle\blacklozenge,\spadesuit,\rho_i:\rho_i^2=\spadesuit^2=\blacklozenge^2=
(\spadesuit\blacklozenge)^3=(\spadesuit\blacklozenge\rho_i)^2=(\blacklozenge\rho_i)^2=1\rangle,\;i=1,2$$
$$\;\qquad\qquad=\langle\blacklozenge,\spadesuit,\varrho\spadesuit\blacklozenge:\spadesuit^2=\blacklozenge^2=
(\spadesuit\blacklozenge)^3=(\spadesuit\varrho\spadesuit\blacklozenge)^2=(\blacklozenge\varrho\spadesuit\blacklozenge)^2=1\rangle,$$
\noindent where $\spadesuit\varrho,\blacklozenge\,\varrho$ determine  an action of
$\mathfrak{S}_3$  on puzzles and LR tableaux.

\section{The action of $\mathbb{Z}_2\times\mathfrak{S}_3$  on LR companion pairs and hives}\label{sec:hives}
 Theorems \ref{fundamentalrev} and \ref{th:fundamental} show that  LR commuters and LR transposers are linear time reducible to each other, in particular, to the Lusztig-Sch\"utzenberger involution. As for the LR companion tableaux, the Henriques-Kamnitzer LR commuter version of $\rho_1$ shows that if $(L,G)$ is the companion pair (left and right) of $T\in LR(\mu,\nu,\lambda)$ then $({\rm evac}\, G,{\rm evac} \,L)$ is the companion pair (left and right) of $\rho_1(T)\in LR(\nu,\mu,\lambda)$. On the other hand, Theorem \ref{th:companion} shows that ${\rm evac}\, G$ is also the right companion of $\rho(T)$. We now translate the action  of $\mathbb{Z}_2\times\mathfrak{S}_3$ on LR tableaux to their companion pairs.

\begin{proposition} \label{th:rho1} Let $(L_\mu, G_\nu)$ be the companion  pair  of  $T$ in $LR(\mu,\nu,\lambda)$. The following holds

$(a)$  the   companion pair of $\rho_1(T)\in  LR(\nu,\mu,\lambda)$ is $$(\blacklozenge \spadesuit{\rm evac}\,L_\mu,{\rm evac}\,L_\mu)=({\rm evac}\,G_\nu,\spadesuit\blacklozenge {\rm evac}\,G_\nu).$$

$(b)$ $G_\nu=\spadesuit\blacklozenge L_\mu$, equivalently, $L_\mu=\blacklozenge\spadesuit G_\nu$.

 \end{proposition}
 \begin{proof} From Henriques-Kamnitzer LR commuter version of $\rho_1$ and Theorem \ref{th:companion},
if $G$ is the right companion of $T$ in $LR(\mu,\nu,\lambda)$ then ${\rm evac}G$ is simultaneously the left companion of $\rho_1(T)\in LR(\nu,\mu,\lambda)$ and the right companion of $\rho(T)\in LR(\lambda,\nu,\mu)$. Being ${\rm evac}G$ the right companion of $\rho(T)\in LR(\lambda,\nu,\mu)$, Algorithm \ref{alg:blacklozenge} says that $\blacklozenge {\rm evac}G$ is the right companion of $\blacklozenge\rho(T)\in  LR(\mu^t,\nu^t,\lambda^t)$. Then, from Algorithm \ref{alg:spade} and Corollary \ref{cor:fund}, $\spadesuit\blacklozenge {\rm evac}G={\rm evac}\blacklozenge \spadesuit G$ is the right companion of $\rho_1(T)=\spadesuit\blacklozenge\rho(T)=\rho\blacklozenge\spadesuit(T)\in  LR(\nu,\mu,\lambda)$. Finally from Henriques-Kamnitzer LR commuter version of $\rho_1$, $({\rm evac}G_\nu,{\rm evac}L_\mu)$ is the  companion pair of $\rho_1(T)$. Henceforth,
${\rm evac}L_\mu=\spadesuit\blacklozenge {\rm evac}G_\nu={\rm evac}\blacklozenge\spadesuit G_\nu$ and thereby $L_\mu=\blacklozenge\spadesuit G_\nu$ and $({\rm evac}G_\nu,\spadesuit\blacklozenge {\rm evac}G_\nu)$ is the companion pair of $\rho_1(T)$ in $LR(\nu,\mu,\lambda)$. Similarly, since $G_\nu=\spadesuit\blacklozenge L_\mu$,  from Corollary \ref{cor:fund} one has  ${\rm evac } G_\nu=\blacklozenge \spadesuit{\rm evac}L_\mu$, and then   $(\blacklozenge \spadesuit{\rm evac}L_\mu,{\rm evac}L_\mu)=({\rm evac}G_\nu,\spadesuit\blacklozenge {\rm evac}G_\nu)$ is the companion pair of $\rho_1(T)$.
 \end{proof}

 \begin{corollary} \label{cor:rho1} Let $T\in LR(\mu,\nu,\lambda)$  and $(L,G)$ be a pair of straight semistandard
 tableaux of shape $\mu$, weight ${\rm rev}(\lambda/\nu)$, and shape $\nu$, weight $\lambda^\vee/\mu$ respectively. The pair $(L,G)$ is  a  companion pair of  $T$  if and only if $L=\blacklozenge\spadesuit G$ and $L$ or $G$ is the left or right companion of $T$ respectively.
 \end{corollary}

 \begin{proof} The "only if " part is the content of  previous proposition.
We now prove the "if" part. We  assume that $L=\blacklozenge\spadesuit G$ and $G$ is the right companion of $T$. 
Then from previous proposition $({\rm evac} G,\spadesuit\blacklozenge {\rm evac}G)$ is the companion pair of $\rho_1(T)$. In particular, $\spadesuit\blacklozenge {\rm evac}G={\rm evac}\blacklozenge\spadesuit G={\rm evac}L$ is the right companion of $\rho_1(T)$. From the Henriques-Kamnitzer commuter, $L$ is the left companion of $T$.
\end{proof}

 \begin{example} Considering  Example \ref{E:spade}, one has $\blacklozenge\spadesuit G=\begin{Young} 3&&&\cr  2&4&&\cr 1&1&2&4\cr
\end{Young}=L_\mu$ of shape $\mu$ and content $(\lambda^\vee/\nu)^\bullet=(7542-5420)_{rev}=2212$. The left companion tableau of $T$ in Example \ref{E:spade}.
\end{example}
 \begin{equation}\label{comprho1}\rho_1:LR_{\nu,\lambda/\mu}\rightarrow LR_{\mu,\lambda/\nu}: G\mapsto\spadesuit\blacklozenge\,\rm{evac}\, G
 \end{equation}
such that $\iota(\rho_1(T))=\rho_1(\iota(T))=\spadesuit\blacklozenge \rm{evac}\, G={\rm evac}\, L$.

\begin{equation}\label{comprho2}\rho_2:LR_{\nu,\lambda/\mu}\rightarrow LR_{\lambda^\vee,\nu^\vee/\mu}: G\mapsto\blacklozenge\spadesuit\,\rm{evac}\, G
 \end{equation}
 such that $\iota(\rho_2(T))=\rho_2(\iota(T))=\blacklozenge \spadesuit\rm{evac}\, G$.

\subsection{The action of $\mathbb{Z}_2\times\mathfrak{S}_3$ on  hives}
\label{hives}

{\it Hives\/} were
first introduced by Knutson and Tao \cite{knutson0} with properties described in more detail by Buch \cite{Bu}.
Hives have also a {\it edge representation\/} as introduced by~\cite{KTT1} and used in~\cite{KTT2}. In this representation, a hive is specified by superposing three Gelfand-Tsetlin patterns where two of them constitute the companion pair of an LR tableau and the third is a consequence of the triangle condition on the edge labels of the hive. On the other hand, superposing a companion pair of an LR tableau always specifies a unique LR hive (for details see \cite{akt}). Let
${\mathcal H}(\mu,\nu,\lambda)$
be the set of LR hives whose  left, right and lower boundary edge
labels are specified by the parts of the partitions  $\mu$, $\nu$ and $\lambda$  fitting a $d\times n-d$ rectangle.
Hives are thereby naturally in bijection with LR tableaux:  $T\mapsto (\blacklozenge\spadesuit G,G)$ where $G$ is the right LR companion of $T\in LR(\mu,\nu,\lambda)$. We define the action of $\mathbb{Z}_2\times\mathfrak{S}_3$ on $\mathcal{H}(\mu,\nu,\lambda)$ via its  action on LR companion pairs. Henceforth, the following theorem exhibits the twelve symmetries of LR coefficients in the  hive model via LR companion pairs.

\begin{theorem}\label{actionLRcompanion} Let $(L=\blacklozenge \spadesuit G,G)$ be the LR companion pair of $T\in LR(\mu,\nu,\lambda)$. Then we have the following LR companion pairs under the action of $\mathbb{Z}_2\times\mathfrak{S}_3$:

$\begin{array}{cccccccc}
(1)& \spadesuit T&\mapsto&\spadesuit G\mapsto(\blacklozenge  G,\spadesuit G)\\
(2)&  \blacklozenge T& \mapsto&\blacklozenge G\mapsto(\spadesuit L=\blacklozenge \spadesuit\blacklozenge   G, \blacklozenge  G)\\
 (3)&\blacklozenge \spadesuit\blacklozenge T&\mapsto &\blacklozenge \spadesuit\blacklozenge G\mapsto( \spadesuit   G, \blacklozenge \spadesuit\blacklozenge  G=\spadesuit L)\\
 (4)&\blacklozenge \spadesuit T&\mapsto &\blacklozenge \spadesuit  G\mapsto(\blacklozenge\spadesuit L=\spadesuit\blacklozenge G,\blacklozenge \spadesuit G=L)\\
(5) &\spadesuit\blacklozenge T&\mapsto &\spadesuit\blacklozenge G\mapsto(G, \spadesuit\blacklozenge G)\\
(6)&\rho(T)&\mapsto &\rho G\mapsto(\blacklozenge \spadesuit {\rm evac}\,G, {\rm evac}\,G)\\
(7)& \rho_1(T)&\mapsto &\rho_1 G \mapsto({\rm evac}\,G, \spadesuit\blacklozenge {\rm evac}\,G={\rm evac} L)\\
(8)&\rho_2(T)&\mapsto &\rho_2 G\mapsto(\spadesuit\blacklozenge {\rm evac}\,G={\rm evac} L, \blacklozenge\spadesuit  {\rm evac}\,G)\\
(9)&\varrho(T)&\mapsto &\varrho G \mapsto(\blacklozenge{\rm evac}\,L, \blacklozenge {\rm evac}\,G)\\
(10)&\spadesuit \rho(T)&\mapsto &\spadesuit \rho G \mapsto(\blacklozenge  {\rm evac}\,G, \spadesuit {\rm evac}\,G)\\
(11)&\blacklozenge\spadesuit\blacklozenge \rho(T)&\mapsto &\blacklozenge\spadesuit\blacklozenge \rho G \mapsto
(\spadesuit  {\rm evac}\,G, \blacklozenge\spadesuit\blacklozenge {\rm evac}\,G)=(\spadesuit  {\rm evac}\,G, \blacklozenge {\rm evac}\,L).
\end{array}$

\end{theorem}
\begin{proof} The proof is a direct consequence of Theorem \ref{th:companion}, Algorithm \ref{alg:spade}, Theorem \ref{th:fundamental} and Theorem \ref{th:rho1}. We leave the calculations for the interested  reader.
\end{proof}

\appendix
\section{Purbhoo mosaics, a rhombus--square--triangle model, and migration}
\label{a:mosaic}
 In this section, we follow closely~\cite{mosaic} (to which we refer the reader for more details) { except  that,  an LR tableau of boundary $(\lambda,\mu,\nu)$ here  it is written there $(\mu,\lambda,\nu^\vee)$.} Consider a
puzzle of side length $n$ and { boundary $(\mu,\nu,\lambda)$}, and replace the rhombi by unitary squares.
The puzzle will be distorted and a  convex  diagram can be recovered
 by adding thin rhombi with angles of $\frac{5}{6}\pi$ and $\frac{\pi}{6}$ radians to the three distorted edges of
the puzzle. If one ignores the  labels on the puzzle pieces, the resulting diagram, called {\em  mosaic}, is a tiled
hexagon by  three shapes:
  equilateral triangles with side length $1$;  squares with side length $1$; and   rhombi with side length $1$ and angles $\frac{\pi}{6}$ and $\frac{5}{6}\pi$
radians, in a way that all rhombi are packed into three nests $ A$, $ B$ and $ C$ of the hexagon. See the picture
 below, where { the mosaic was built on  a puzzle with boundary $\mu=2111,$ $\nu=2110$ and $\lambda=2100$}, the colours should be looked at this point as  decoration. The mosaic has side lengths $A'A=B'B=C'C=n-d$
 and $AB'=BC'=CA'=d$. The collections of rhombi in the  nests $A$, $B$ and $C$,
 denoted respectively by $\alpha, \beta$, and $\gamma$, define the boundary data $(\alpha, \beta, \gamma)$ of the mosaic.
  {  The collections of extra  rhombi $\alpha$, $\beta$ and $\gamma$, with the standard orientation, given by the edges of the nests  clockwise, can be
regarded as the three Young diagrams $\mu$, $\nu$, and $\lambda$, respectively, clockwise encoded by the $01$--words on
the boundary of the puzzle.}

\medskip
\hspace{4cm}\includegraphics[scale=0.7,trim = 0cm 13.5cm 0cm
0cm,clip]{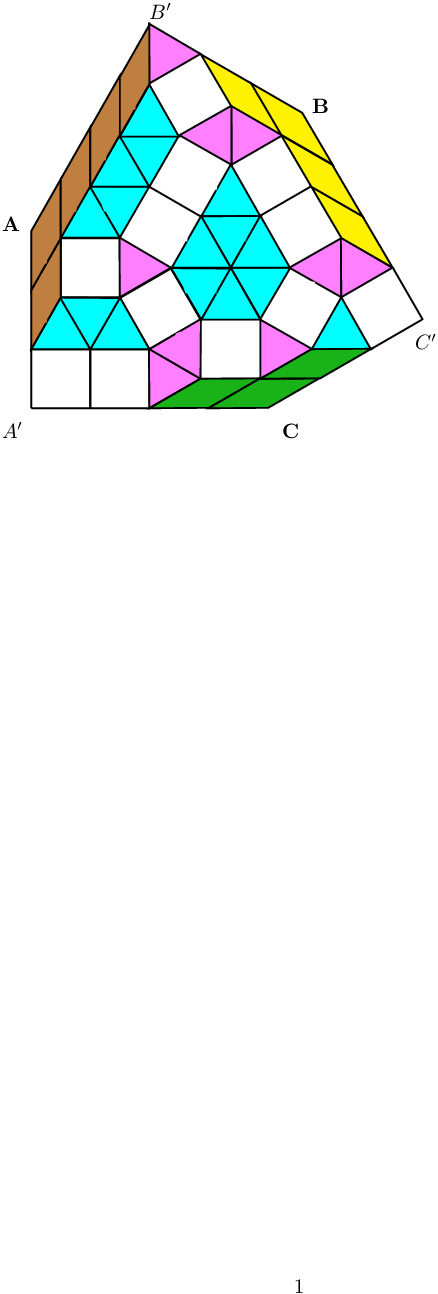}

This construction exhibits a natural bijection between mosaics
and puzzles.
 Removing the extra rhombi and
straightening the resulting shape, we can go from a mosaic to a puzzle. Walking from
$A'$ to $B'$, the shape that is left, by removing $\alpha$, turns into the
string of $0$'s and $1$'s, $0$ for each unit step west, and a $1$
for each unit step north; and  straightening  the squares they will become $\frac{\pi}{6}$
/ $\frac{\pi}{3}$ radians rhombi. This will determine the remain labels of the puzzle pieces.
(Similarly, walking anticlockwise from
$C'$ to $B'$, we get  the dual puzzle, that is, the one obtained by
mirror reflection along the vertical axis and label swapping.) In the standard orientation,
that is, read clockwise, a mosaic of boundary $(\alpha, \beta,
\gamma)$  can be identified with the corresponding puzzle of boundary data $(\mu,
\nu, \lambda)$,  where $\alpha$ is identified with the Young diagram of $\mu$,
$\beta$ with $\nu$, and $\gamma$ with $\lambda$.  The number of mosaics with boundary data $(\alpha,
\beta, \gamma)$ is equal to the number of puzzles of boundary
$(\mu,\nu,\lambda)$.

One of the
advantages of mosaics over puzzles is that we can give different
orientations to the nests $A$, $B$, and $C$. This
allows us to relate  the  symmetry bijections on puzzles and  LR tableaux.
Define unit vectors $E_A$, $N_A$, $E_B$, $N_B$, $E_C$, $N_C$ in
the directions of $AA'$,
$AB'$,
$BB'$,
$BC'$,
$CC'$,
$CA'$
respectively, and fix orientations
$(E_A,N_A)$, $(E_B,N_B)$, $(E_C,N_C)$ on the nests at $A$, $B$, and $C$ respectively.
The letters $E$, $ N$, $-E$, $-N$ are thought as compass directions east, north, west and
south, respectively. Thus the orientations $(E,N)$ and $(N, E)$ in  a nest means  the standard or clockwise
orientation, and the counterclockwise orientation respectively. Flocks are (skew) tableau--like structures, defined on the thin rhombi in a
mosaic, packed into one of the nests $A$, $B$ or $C$. Four orientations can be given to a nest. Each orientation uniquely
determines the flock as a  LR tableau. Fix  a nest, say $A$, the rhombi in $\alpha$ under the orientation $(E, N)$ define
the Yamanouchi tableau $Y(\mu)$; under $(N,E)$, $Y(\mu^t)$; under $(-E, -N)$, $Y(\mu)^{\rm a}$; and
under $(-N, -E)$, $Y(\mu^t)^{\rm a}$. {(The second compass direction indicates in which direction the entries of the LR tableau strictly increase. In the standard orientation $(E,N)$, the entries of an LR tableau weakly increase eastward and strictly increase northward. This is consistent with the representation of  partitions and their linear transformations as  $01$ words in \eqref{rectangle} and \eqref{01word},  Section \ref{sec:pre}.)}
\emph{ Migration} is
an invertible operation on mosaics that takes a flock from a nest to a new nest whose shape  can be interpreted as  an LR tableau. In this case the initial mosaic is identified with that LR tableau. This operation gives a bijection between mosaics (equivalently puzzles) and LR tableaux. More precisely, with appropriate orientation of the flocks in the mosaic, \emph{migration} coincides with Tao's bijection and allows to relate operations on puzzles with {\em jeu de taquin} operations (or switching)  on LR tableaux.
The rhombi
must move in the standard order of  a tableau (recall the definition in Section \ref{subsec:standard}).
 Choose the target nest. Rhombi move in the chosen direction of migration, inside a smallest hexagon in which the thin
 rhombus $\diamondsuit$ is
contained:

\vskip0.2cm
\includegraphics[scale=0.7,trim = 0cm 20cm 0cm
0cm,clip]{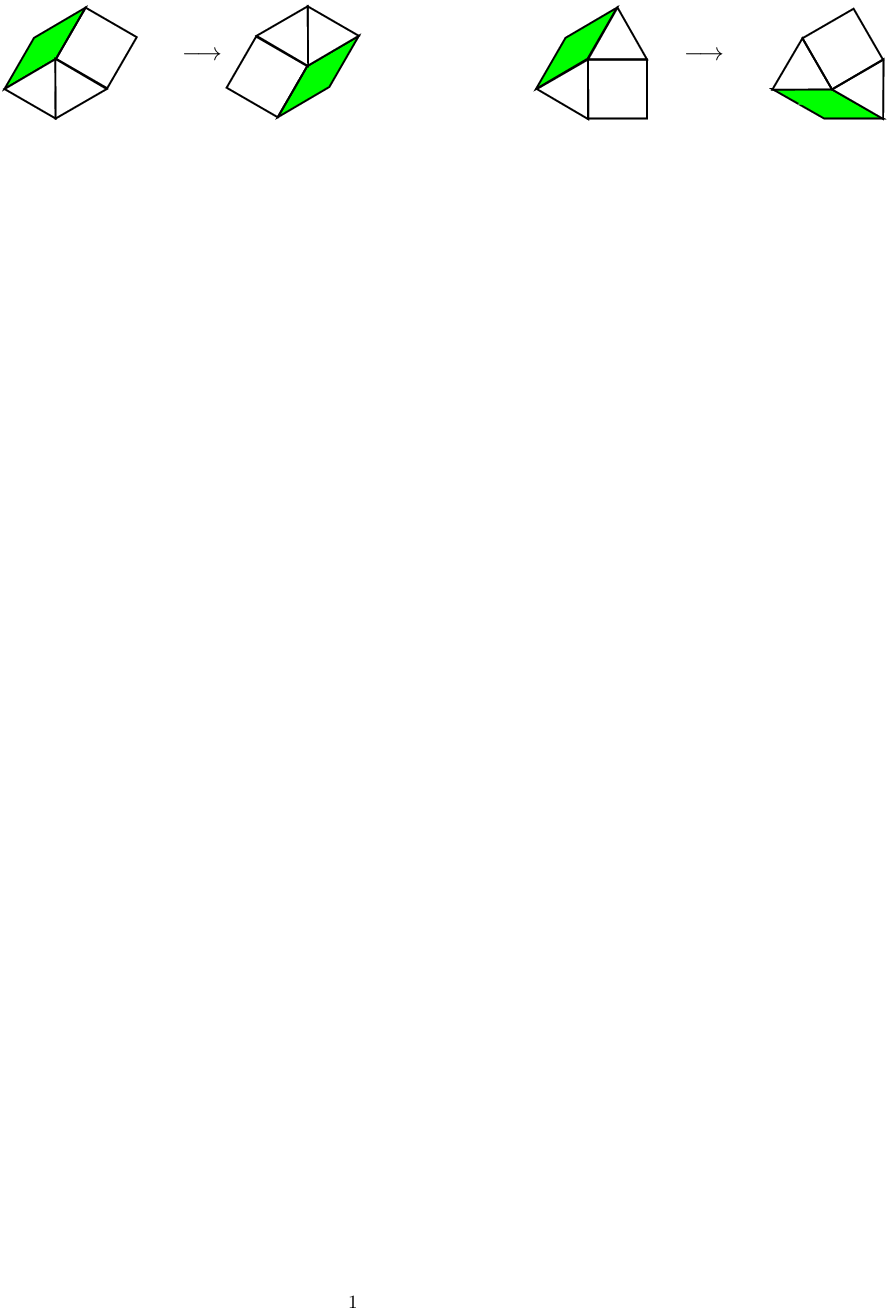}

\vskip0.5cm
\emph{Migration} from the nest $A$ to the  nest $B$ of the mosaic, of the flock $\alpha$ with standard orientation $(E,N)$,
 gives a bijection
between   mosaics of boundary $(\alpha,\beta,\gamma)$ and LR
tableaux of boundary $(\lambda,\mu,\nu)$ where $\nu=\beta$,
$\lambda=\gamma$ and $ \mu=\alpha$. This bijection coincides
with the Tao's bijection, "without words"  in~\cite{vakil} (see also  \cite[Fig. 9]{mosaic}), between   puzzles of boundary
$(\mu,\nu,\lambda)$, or the corresponding  mosaics, and LR tableaux of
boundary data $(\lambda,\mu,\nu)$. On the other hand, migration from the left nest $A$  to the bottom nest C of the mosaic, of the flock $\alpha$ with the orientation $(-E,-N)$,
 gives the same bijection \cite[Proposition 5.1]{mosaic}. This coincidence  with respect to the two orientations given to  the flock $\alpha$ is consistent with the definition of LR tableau. An LR tableau of boundary $(\lambda,\mu,\nu)$ rectifies to the Yamanouchi tableau $Y(\mu)$, and by reverse {\em jeu de taquin}, to its antinormal form $Y(\mu)^{\rm a}$.
An illustration  of  Tao's bijection  which coincides with the
migration  of the flock $Y(\mu)$  to the  nest $B$ of the mosaic, and with the migration of the flock $Y(\mu)^{\rm a}$ to the nest $C$:

\hspace{3cm}\includegraphics[scale=0.8,trim = 0cm 17cm 0cm
0cm,clip]{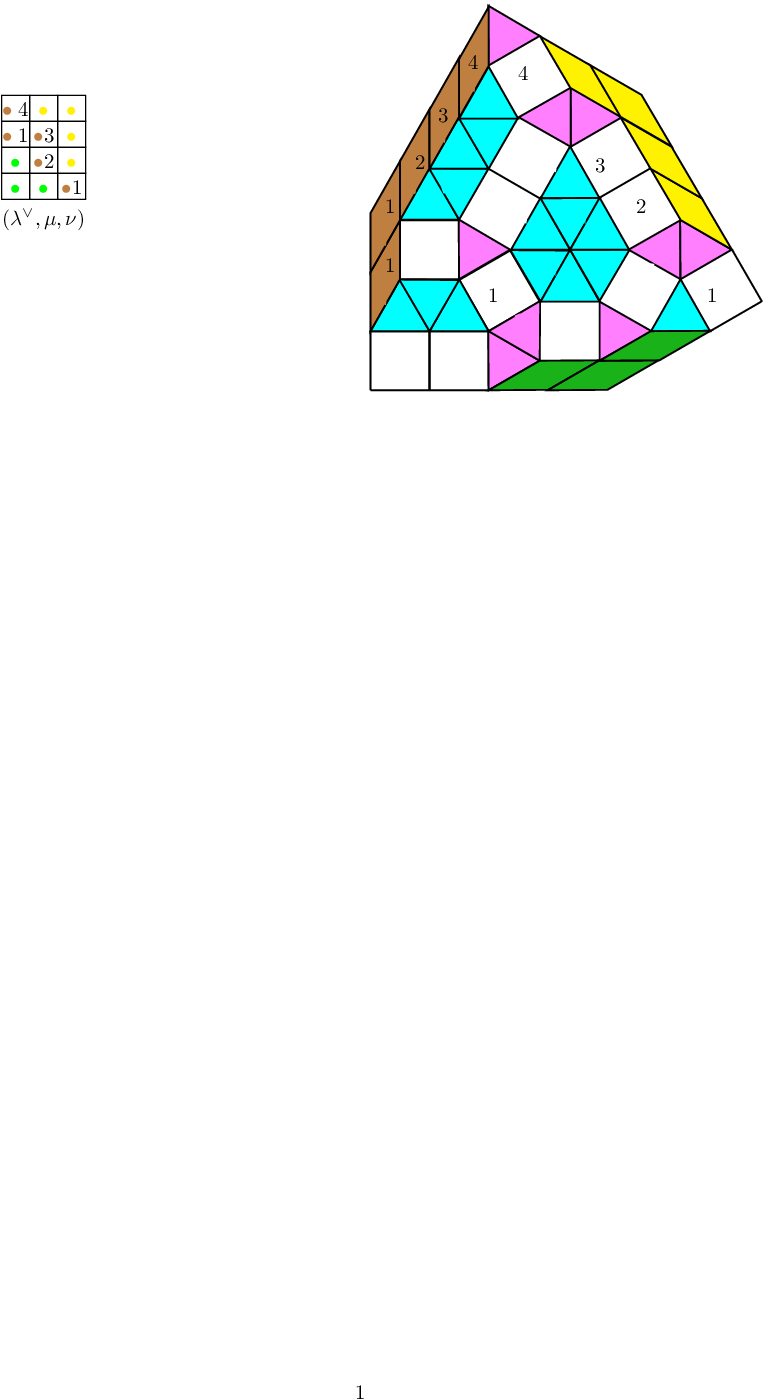}

Following \cite[Table 1]{mosaic}, we next  see that the involutions $\clubsuit,\spadesuit,\blacklozenge$ defined by the action of the group $H$ on puzzles
and on LR tableaux are exactly what we get with migration on mosaics. In this discussion it is better to consider the
presentation $H=<\clubsuit,\clubsuit\blacklozenge>$.
Migration from the nest $B$  to the nest $A$ of the mosaic, of the flock $\beta$ with orientation $(N,E)$
(read counterclockwise),  thus identified with $Y(\nu^t)$ coincides with  Tao's bijection on the back side of the mosaic, that is, on the mosaic of boundary $(\beta,\alpha,\gamma)$  defined  by the reflection of our current mosaic along the vertical axis (also  corresponding to the puzzle obtained by reflecting along the vertical axis and swapping the colours). This  defines
the $\spadesuit$ involution on puzzles which  translates to the $\clubsuit$ involution  on LR tableaux, as illustrated below.
\vskip 0.5cm

\hspace{4cm}\includegraphics[scale=0.8,trim = 0cm 14cm 0cm
0cm,clip]{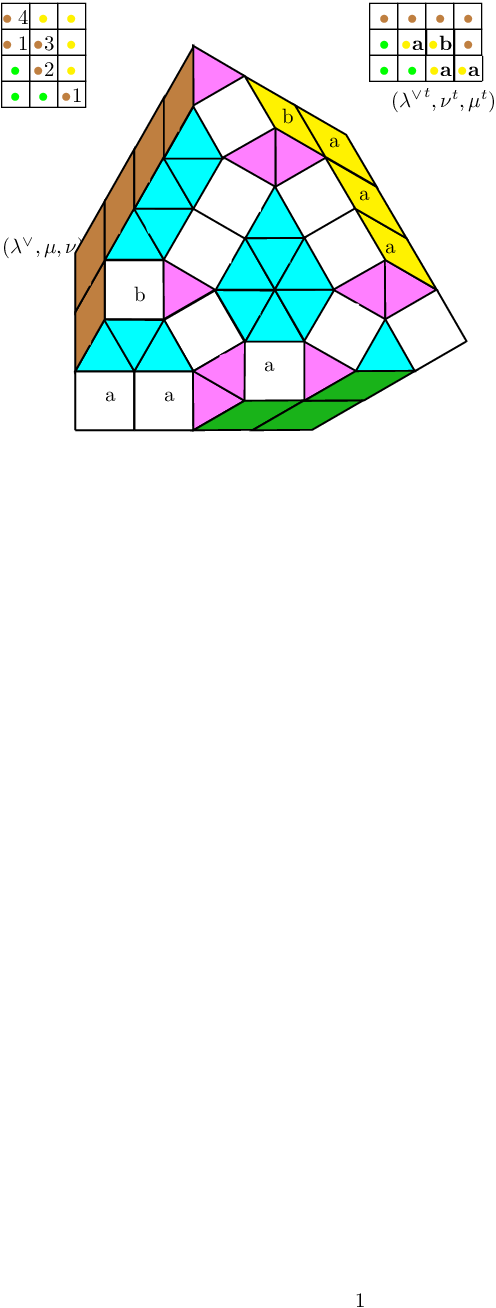}

Migration of the flock $\alpha$ with orientation $(N,E)$ (read counterclockwise), from  the nest $A$ to the nest $C$ of the
mosaic,  coincides with what Tao's bijection gives when applied
on the back side of the mosaic after rotating it $\frac{2}{3}\pi$ radians clockwise (reflection of the corresponding puzzle,  after rotating it $\frac{2}{3}\pi$ radians clockwise, along the vertical axis and swapping the colours).
This defines the involution $\clubsuit$ on puzzles which translates to the $\blacklozenge$ involution on LR
tableaux as  illustrated below with the LR tableau on the right.

Migration from the nest $C$ to the nest $A$  of the mosaic, of the flock $\gamma$
with standard orientation $(E,N)$,  coincides with what Tao's bijection gives when applied to the mosaic  after a rotation of $\frac{2}{3}\pi$ radians clockwise. This   defines the $\frac{2}{3}\pi$ clockwise  rotation
$\clubsuit\blacklozenge$ on puzzles which is translated to $\blacklozenge\spadesuit$ on LR tableaux.  Illustrated below with the LR tableau on the left.
The bijections are illustrated on the mosaic at the same time.

\hspace{4cm}\includegraphics[scale=0.8,trim = 0cm 13cm 0cm
0cm,clip]{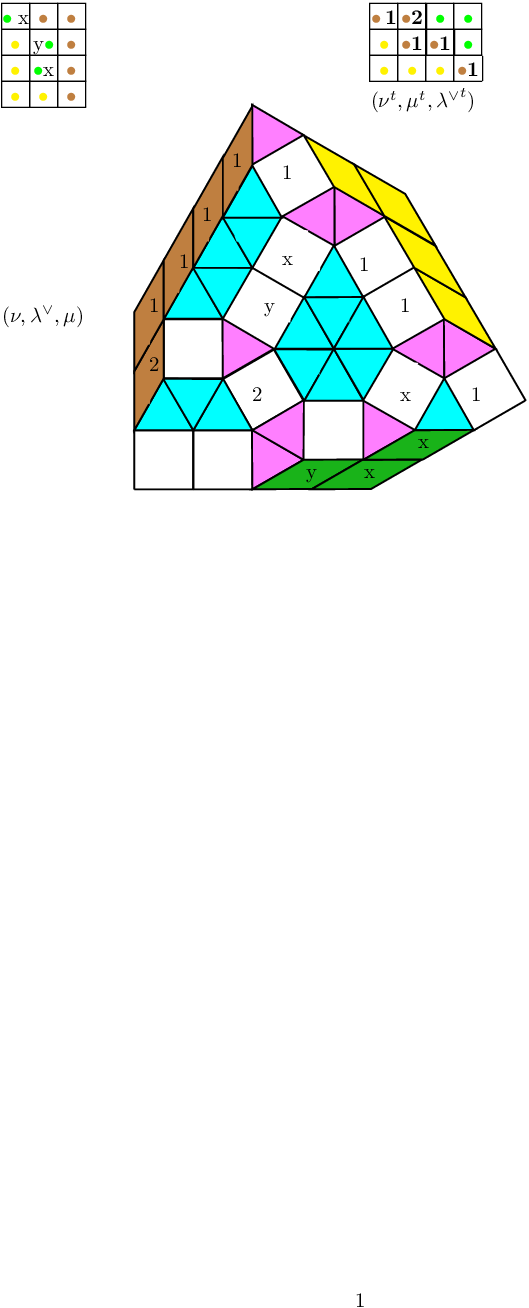}

Similarly, migration of the flock $\gamma$ with orientation $(N,E)$ (read counterclockwise), from the nest $C$ to the nest $B$ of the
mosaic,  coincides  with Tao's bijection
on the back side of the mosaic after rotating it $\frac{2}{3}\pi$ radians counterclockwise (or $\frac{4}{3}\pi$ radians clockwise). It defines the involution $\blacklozenge$ on puzzles which translates to the involution $\spadesuit$ on LR tableaux.

In~\cite{mosaic} it is discussed how the \emph{migration} of a single
rhombus in a mosaic is related with {\em jeu de taquin} slides on
tableaux. This explains the correspondence between  the action of
$\mathbb Z_2\times \mathfrak{S}_3$ on puzzles and on LR tableaux.
We have seen that $\clubsuit$ { on  LR tableaux} can be described as the migration of
the flock $\beta$ with orientation $(N,E)$, thus identified with $Y(\nu^t)$, from the nest $B$ to the nest $A$ of the mosaic.
However $\clubsuit$  can  also be described on mosaics as the migration of the flock $\beta$ with
orientation $(-N,-E)$, thus identified with $Y(\nu^t)^{\rm a}$, from the nest $B$ to the
nest $C$ of the mosaic.  On the back  mosaic, it is the migration
of the  flock $\nu^t$ with orientation $(-E,-N)$, thus
identified with $Y(\nu^t)^{\rm a}$, to the nest $C$. Combining the Tao's bijection on the the mosaic of
boundary $(\alpha, \beta, \gamma)$, in  standard orientation, giving a LR tableau of
boundary $(\lambda, \mu,\nu)$, with  the migration of the flock $\beta$ with orientation $(-N,-E)$, identified with $Y(\nu^t)^{\rm a}$, to the nest $C$, we get
Proposition~\ref{P:club}. This is illustrated below.
 Consider our current mosaic  with the LR tableau of boundary $(\lambda,\mu,\nu)$  produced  by the migration of the flock $Y(\mu)$ to the nest $B$ or Tao's bijection.

\begin{center}
\parbox{6.5cm}{\includegraphics[scale=0.5,trim = 0cm 17cm 0cm 0cm,clip]{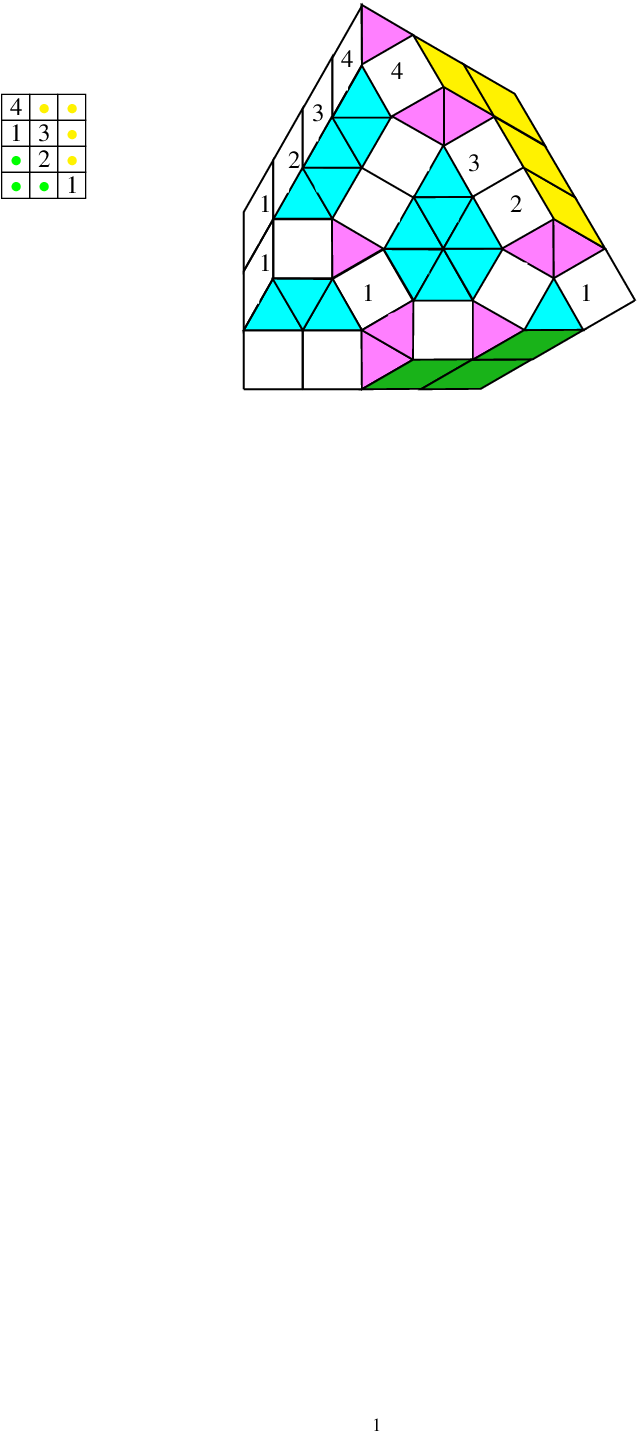}}
\end{center}

 Reflecting our mosaic vertically we get the "back" mosaic with boundary $(\beta,\alpha,\gamma)$  naturally in bijection with the  dual puzzle of boundary $(\nu^t,\mu^t,\lambda^t)$ of the current mosaic. Tao's bijection on the back mosaic gives the LR tableau
 \begin{Young}
 &&&\cr
 &a&b&\cr
 &&a&a\cr
 \end{Young} of boundary $(\lambda^t,\nu^t,\mu^t)$ which coincides with the migration on the back mosaic of $Y\nu^t)^{\rm a}$ to $C$.

  Transposing the standardized  LR tableau  produced by Tao's bijection on the mosaic of boundary $(\alpha,\beta,\gamma)$ just above and filling the outer  partition $\nu^t$ in the  antinormal form to obtain $Y(\nu^t)^{\rm a}$,
 \begin{equation}\label{fish} \begin{Young}
 2&a&a&b\cr
 &3&4&a\cr
 &&1&5\cr
 \end{Young},
  \end{equation} we verify that the migration    of the flock $Y(\nu^t)^{\rm a}$ (flock $\beta$ with orientation $(-E,-N)$) to the nest $C$,  on the back mosaic, is also explained by the tableau switching on \eqref{fish}. Begin with the minimal rhombus in the standard order of $Y(\nu^t)^{\rm a}$ and proceed in standard order  with the remain rhombi in $Y(\nu^t)^{\rm a}$. (Migration preserves the standard order of the flock.) The migrated flock packed in $C$ is identified with an LR tableau of boundary $(\lambda^t,\nu^t,\mu^t)$ which is  the same LR tableau of boundary $(\lambda^t,\nu^t,\mu^t)$ as Tao's bijection on the dual puzzle. (Note that some of the moves in the migration are
omitted.)

\parbox{6.5cm}{\includegraphics[scale=0.5,trim = 0cm 17cm 0cm 0cm,clip]{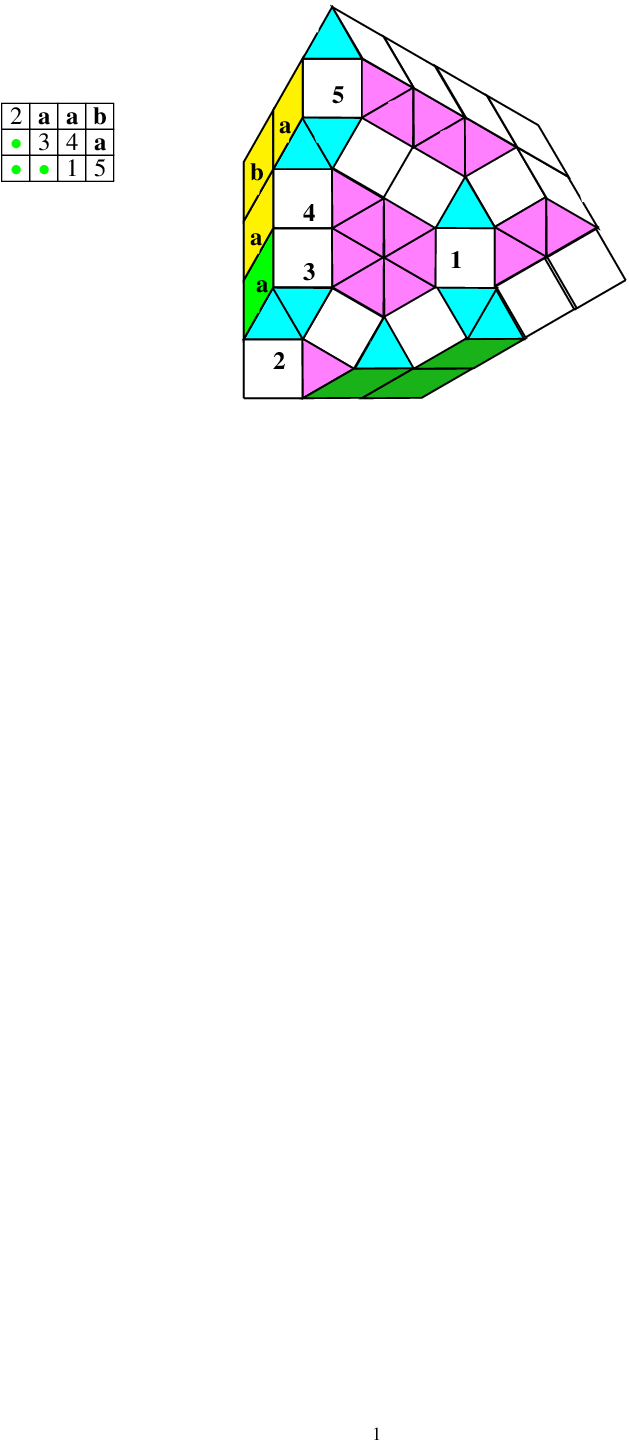}}
\parbox{6.5cm}{\includegraphics[scale=0.5,trim = 0cm 17cm 0cm 0cm,clip]{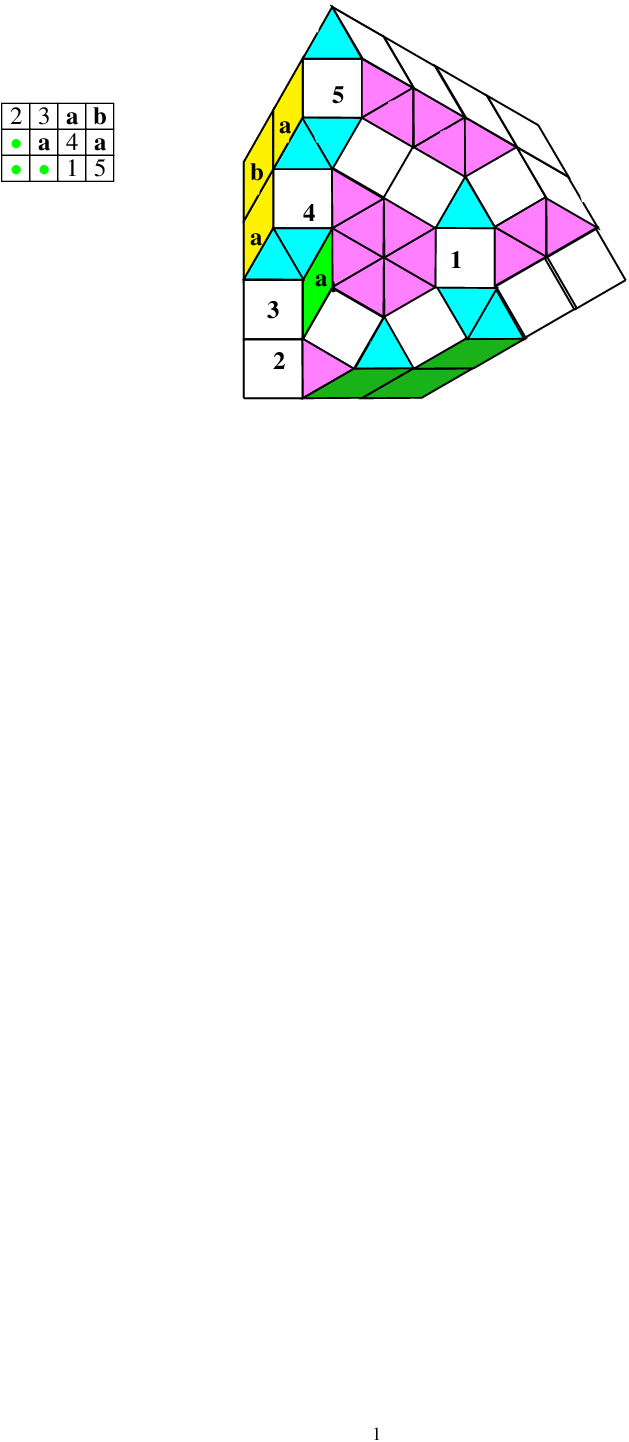}}

\parbox{6.5cm}{\includegraphics[scale=0.5,trim = 0cm 17cm 0cm 0cm,clip]{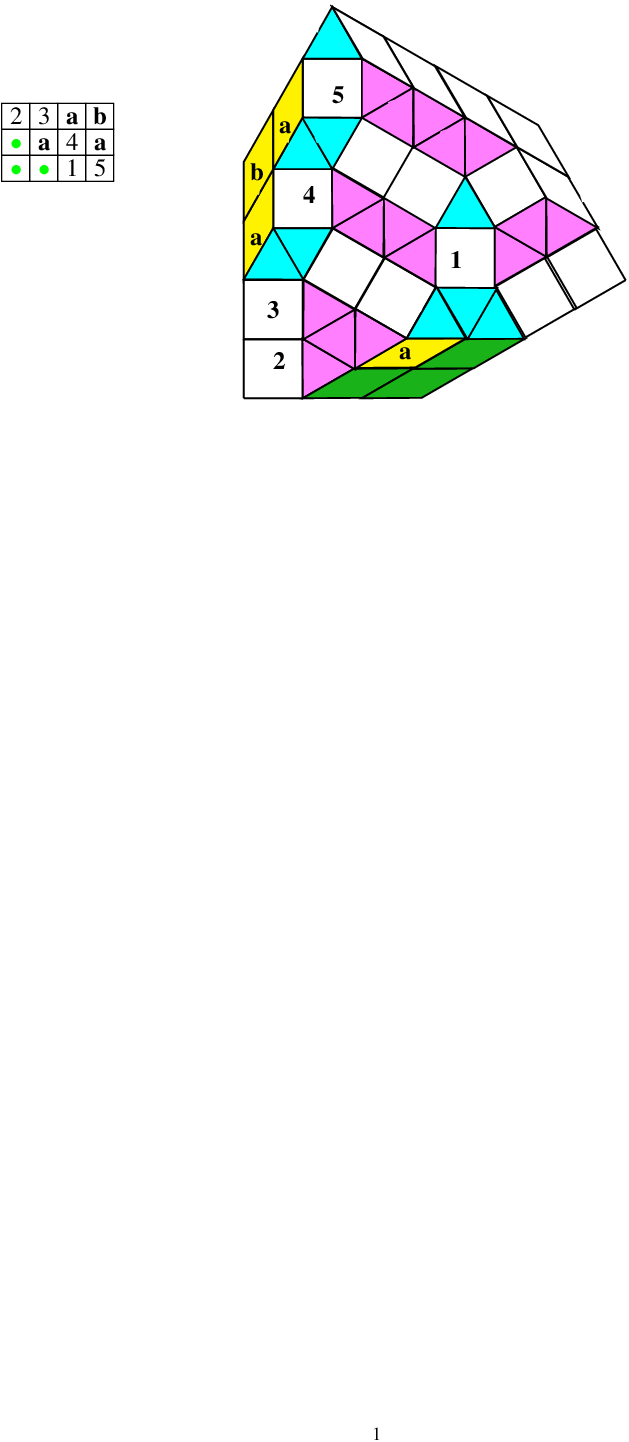}}
\parbox{6.5cm}{\includegraphics[scale=0.5,trim = 0cm 17cm 0cm 0cm,clip]{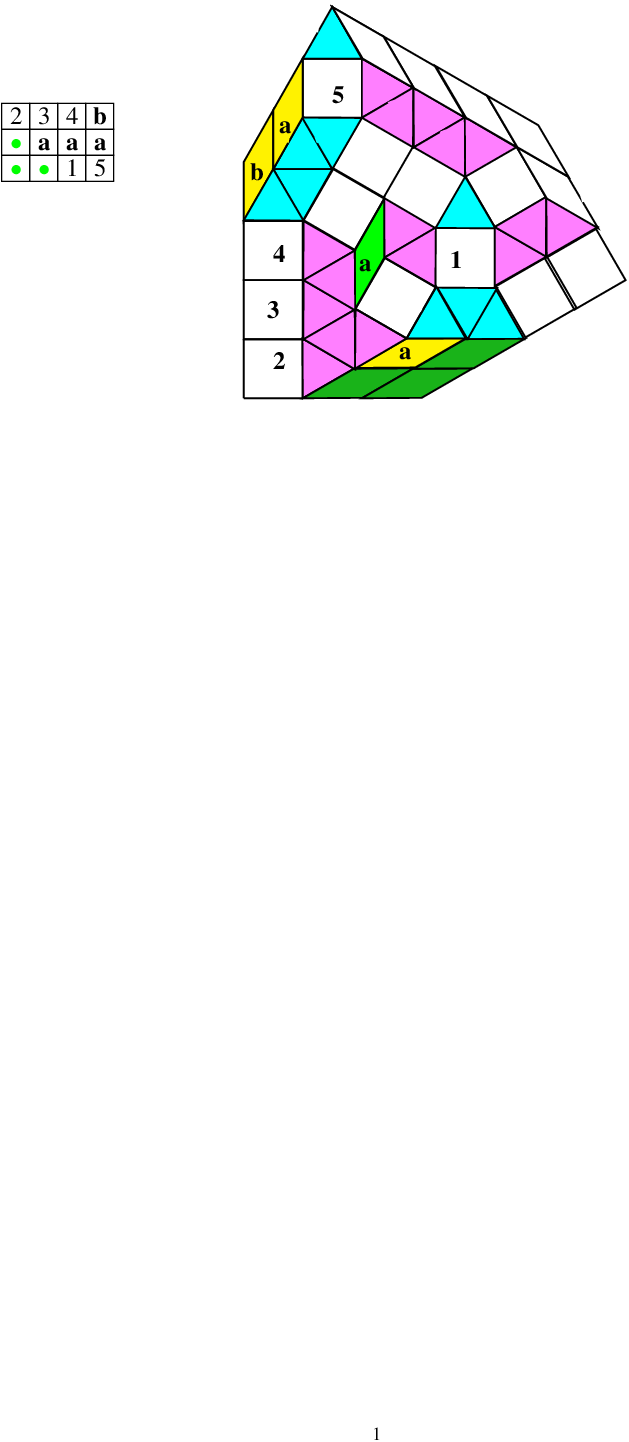}}

\parbox{6.5cm}{\includegraphics[scale=0.5,trim = 0cm 17cm 0cm 0cm,clip]{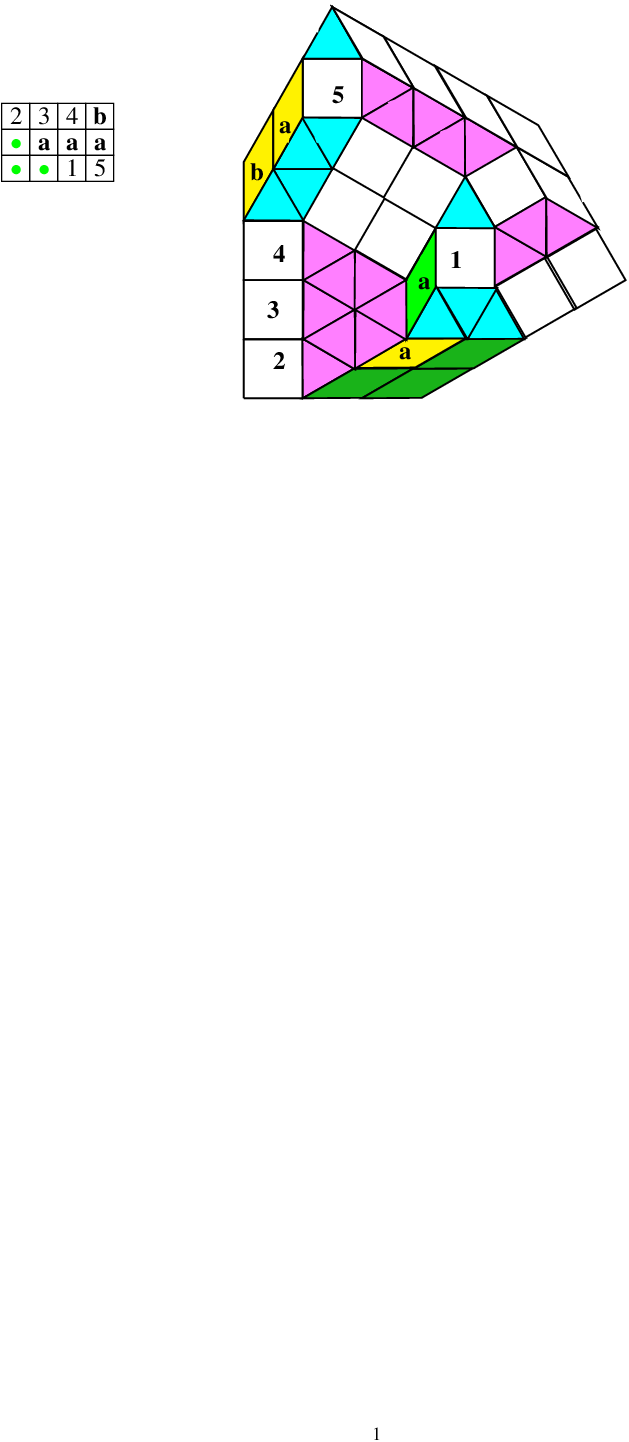}}
\parbox{6.5cm}{\includegraphics[scale=0.5,trim = 0cm 17cm 0cm 0cm,clip]{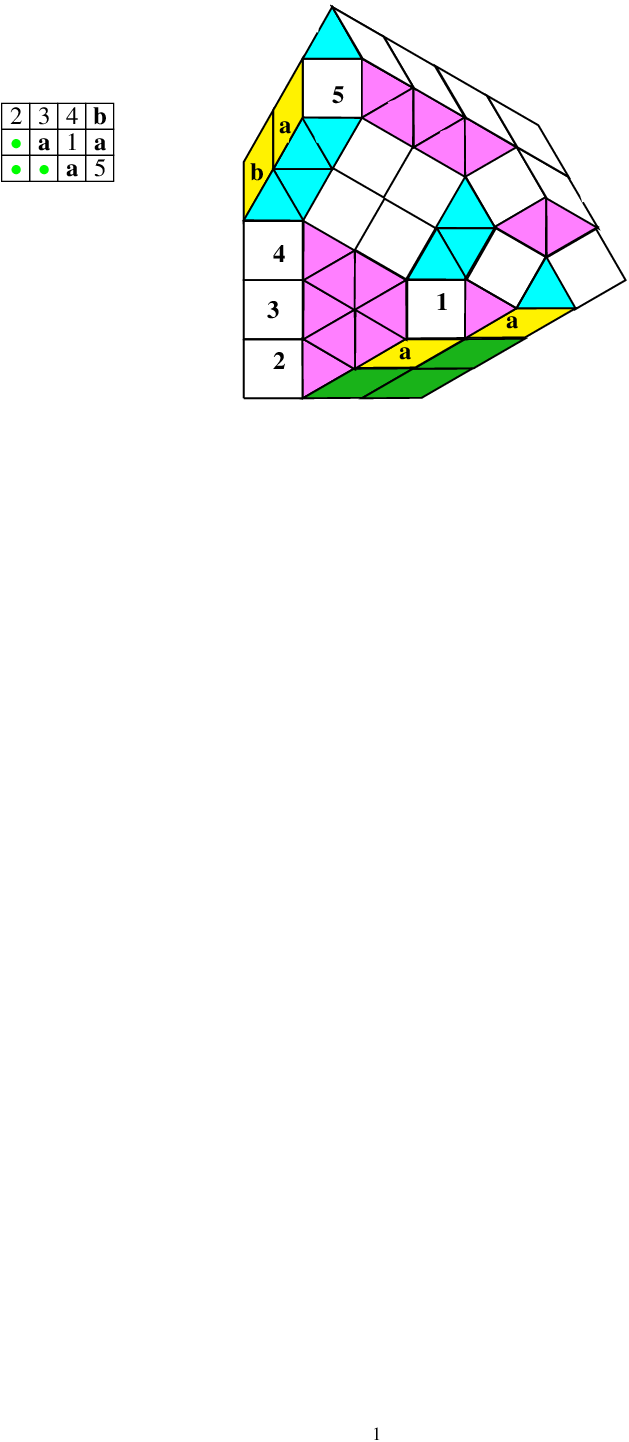}}

\parbox{6.5cm}{\includegraphics[scale=0.5,trim = 0cm 17cm 0cm 0cm,clip]{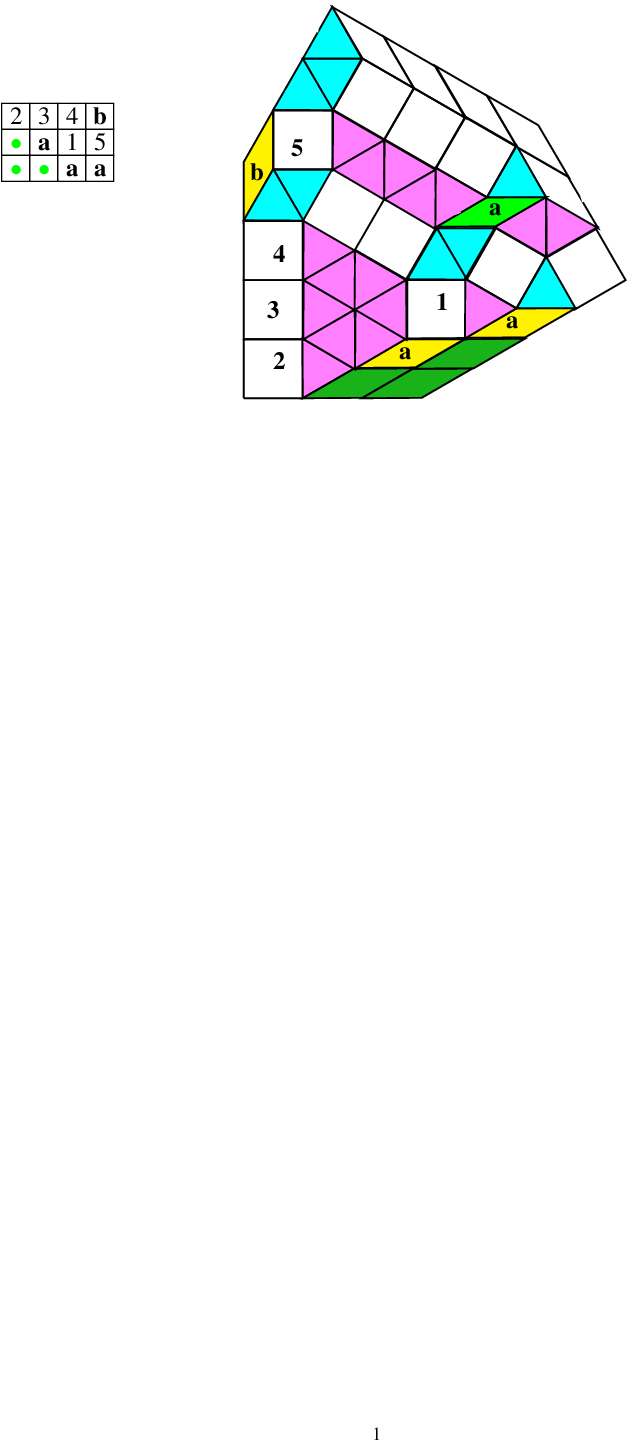}}
\parbox{6.5cm}{\includegraphics[scale=0.5,trim = 0cm 17cm 0cm 0cm,clip]{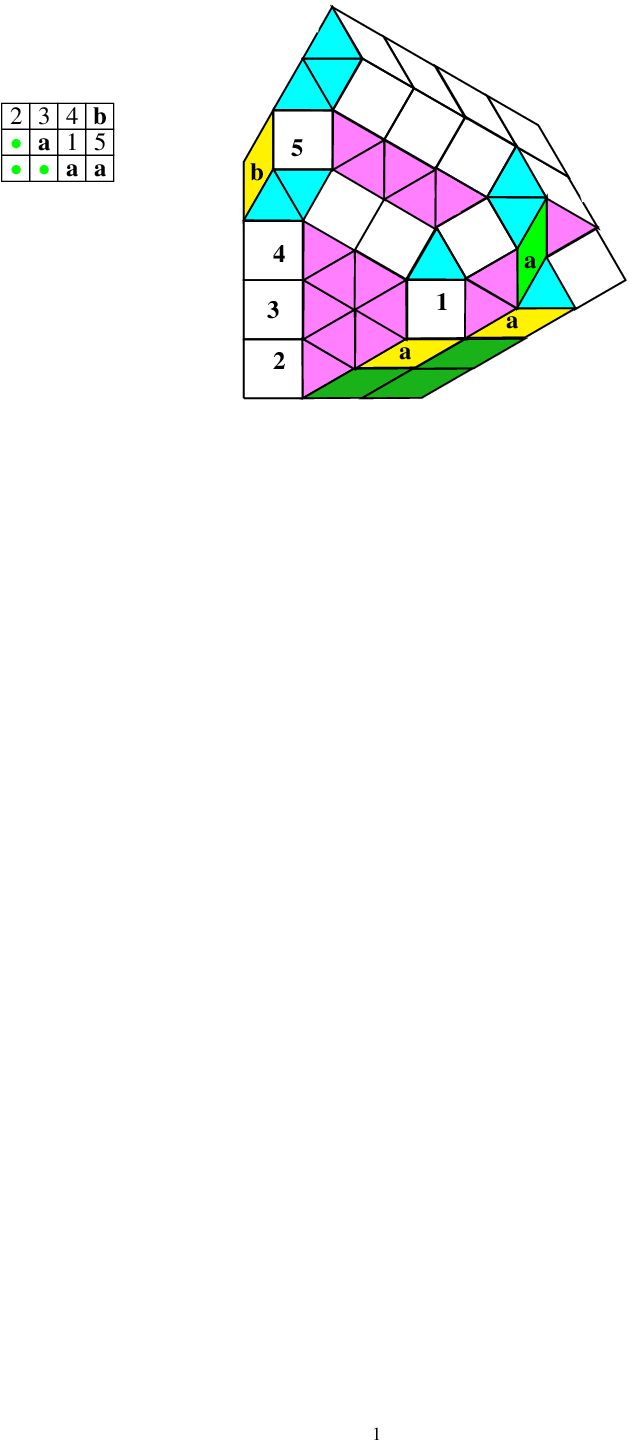}}

\parbox{6.5cm}{\includegraphics[scale=0.5,trim = 0cm 17cm 0cm 0cm,clip]{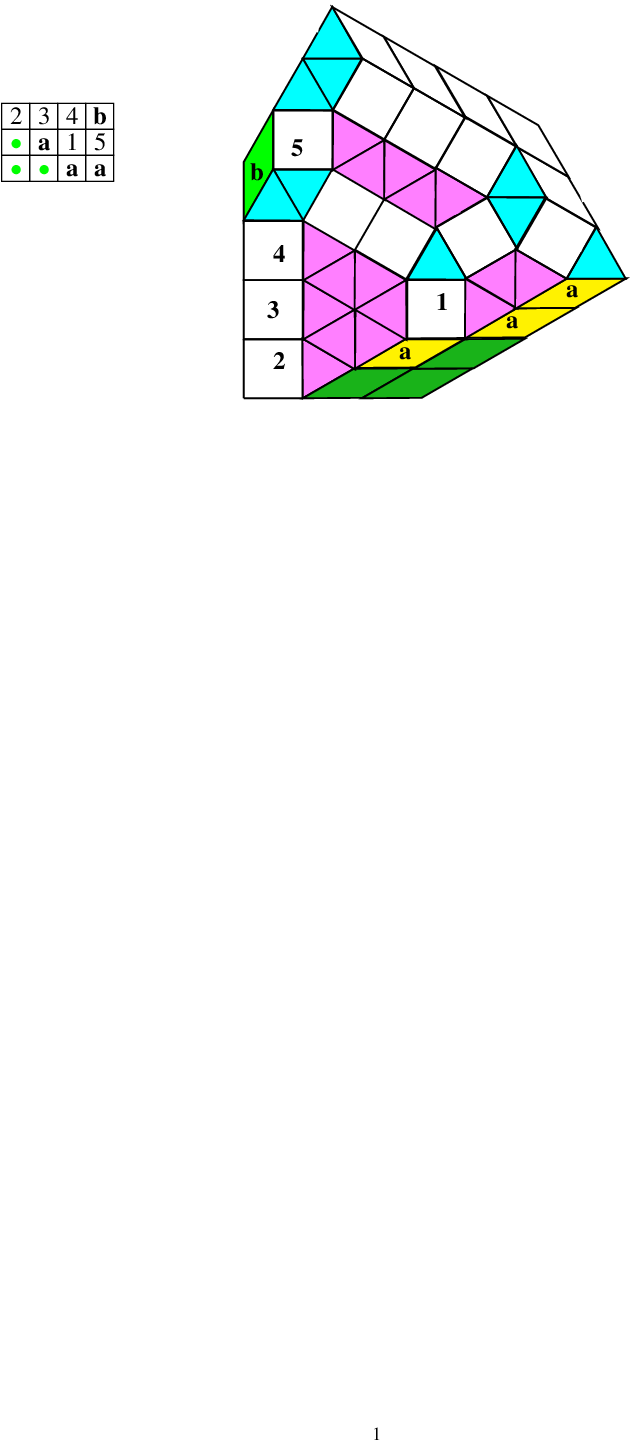}}
\parbox{6.5cm}{\includegraphics[scale=0.5,trim = 0cm 17cm 0cm 0cm,clip]{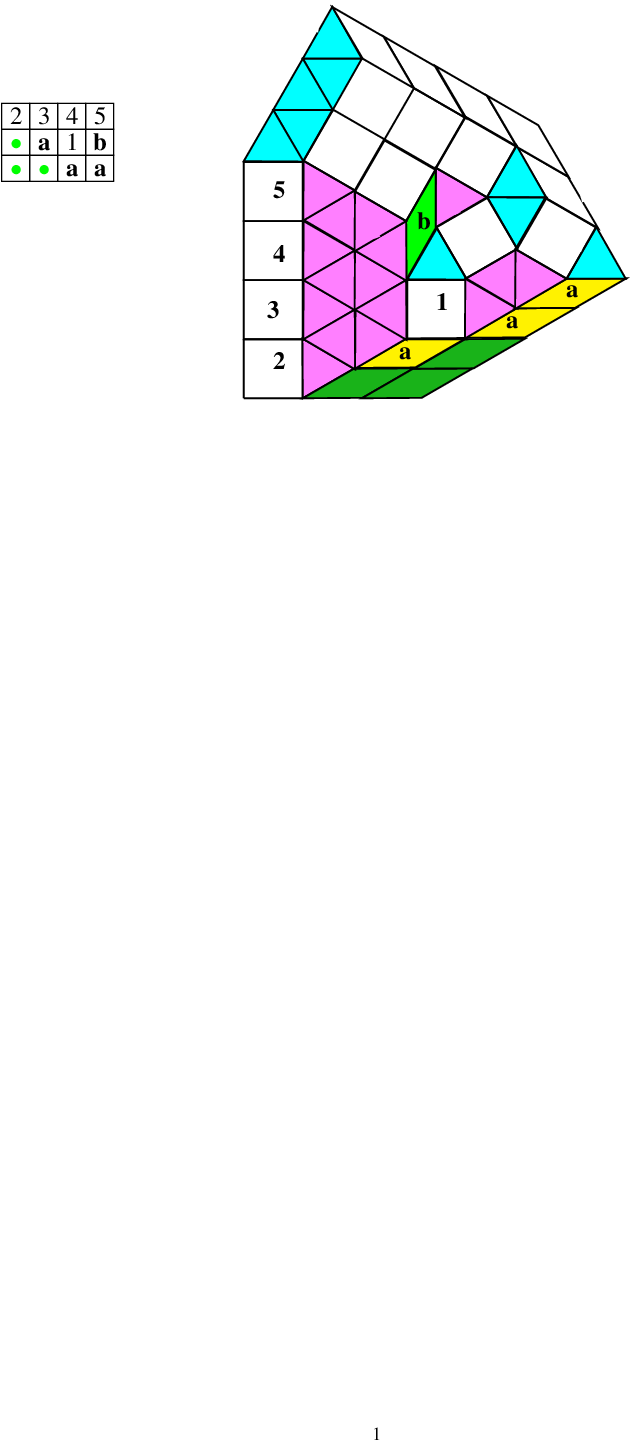}}

\parbox{6.5cm}{\includegraphics[scale=0.5,trim = 0cm 17cm 0cm 0cm,clip]{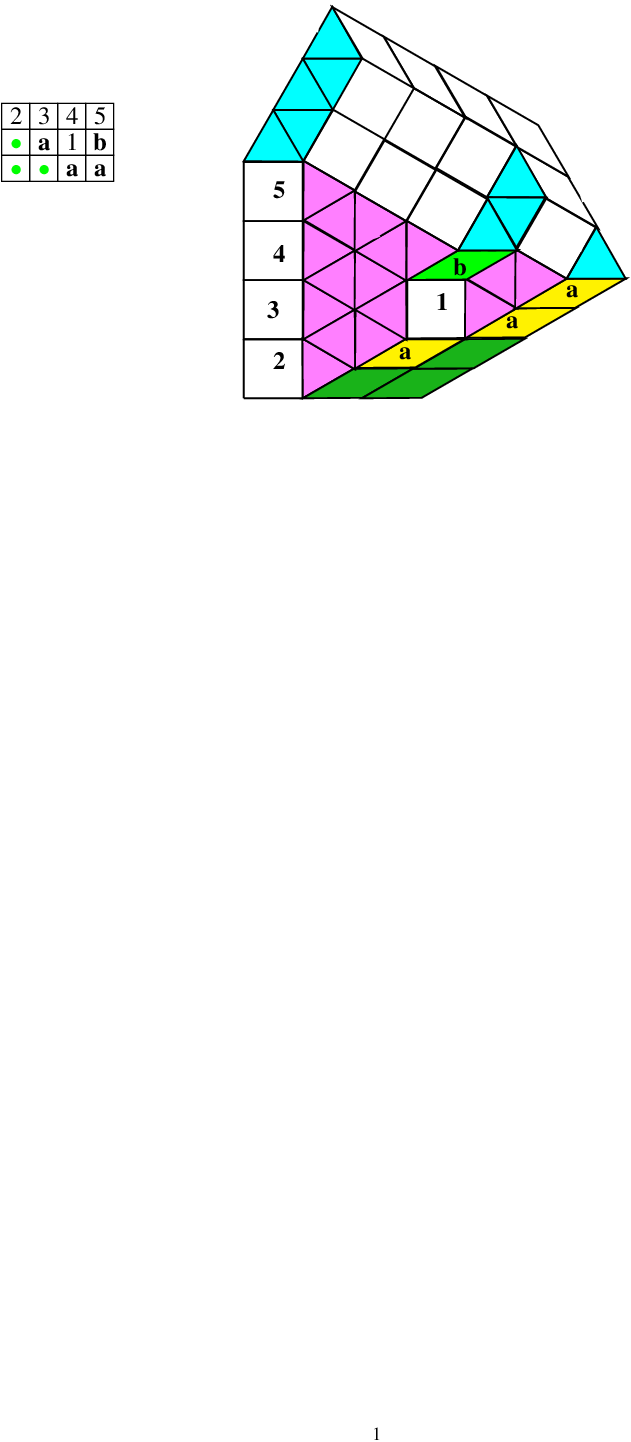}}
\parbox{6.5cm}{\includegraphics[scale=0.5,trim = 0cm 17cm 0cm 0cm,clip]{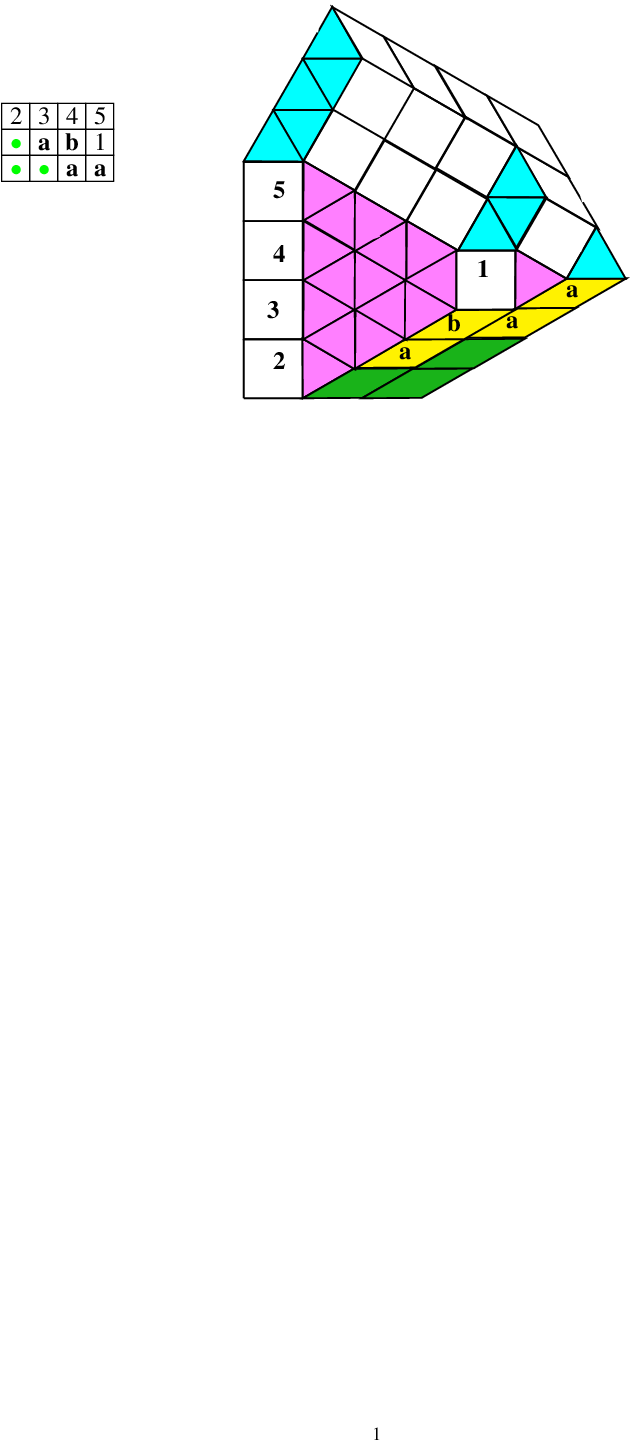}}

The last mosaic, on the right, has boundary $(\emptyset,\alpha,\gamma')$ where, in standard orientation, $\alpha$ is identified with $\mu^t$, and $\gamma'$ identified with $(\mu^t)^\vee $ gives  the  LR tableau of boundary $(\lambda^t,\nu^t,\mu^t)$ already obtained with Tao's bijection above.

\bigskip

\noindent{\bf Acknowledgement:}
The first author wishes to thank  R. C. King, C\'{e}dric Lecouvey,  Cristian Lenart, Igor Pak, Itaru Terada,  and Ernesto Vallejo for helpful discussions at different times along the years.

\end{document}